\documentclass[10pt,reqno]{article}

\usepackage{dsfont}
\usepackage{amsthm,amsfonts,latexsym,amsmath, amssymb, mathrsfs, verbatim}

\usepackage{tocloft}

\usepackage{mathrsfs}
\usepackage{bm}
\usepackage{bbm}
\usepackage{color}

\usepackage{geometry}
\usepackage{graphicx}
\usepackage{hyperref}
\hypersetup{
    linktoc=page
}
\usepackage{aligned-overset}
\usepackage{yhmath}
\usepackage{appendix}
\usepackage{comment}
\usepackage{xcolor}
\usepackage{tikz}
\usepackage{scalerel}
\usepackage{tkz-euclide}
\usetikzlibrary{patterns,arrows.meta}

\usepackage{pgfplots}
\usepackage{float}
\usepackage{listings}
\usepackage{indentfirst}
\usepackage[cmtip,all]{xy}

\usepackage{cite}

\numberwithin{equation}{section}
\numberwithin{figure}{section}

\graphicspath{{./figure/}}
\DeclareGraphicsExtensions{.png}

\theoremstyle{definition}
\newtheorem{definition}{Definition}[section]

\theoremstyle{plain}
\newtheorem{proposition}[definition]{Proposition}

\newtheorem{lemma}[definition]{Lemma}
\newtheorem{theorem}[definition]{Theorem}
\newtheorem{corollary}[definition]{Corollary}
\newtheorem{remark}[definition]{Remark}

\newcommand{\ovl}{\overline}
\newcommand{\dif}{\mathrm{d}}
\newcommand{\jpnlambda}{\langle\lambda\rangle}
\newcommand{\abs}[1]{\left|#1\right|}
\newcommand{\norm}[1]{\left\|#1\right\|}
\newcommand{\longsquiggly}{\xymatrix{{}\ar@{~>}[r]&{}}}

\newcommand{\pt}{\partial_{\tilde t}}

\newcommand{\wtt}{\tilde{\partial}_{\tilde t}}
\newcommand{\mb}{\mathbb}

\DeclareMathOperator{\Supp}{Supp}

\allowdisplaybreaks

\pgfplotsset{compat=1.18}

\title{Shock Formation for Compressible Euler Equations on $\mathbb{S}^2$}

\author{Xinliang An\thanks{National University of Singapore, Singapore. matax@nus.edu.sg}, Haoyang Chen\thanks{National University of Singapore, Singapore. hychen@nus.edu.sg}, Fulin Qi\thanks{National University of Singapore, Singapore. fulin\_qi@u.nus.edu}, and Wenze Su\thanks{Tsinghua University, China. wenzesu@mail.tsinghua.edu.cn}}

\begin{document}
\maketitle
\begin{abstract}
In this paper, we prove the finite-time shock formation for the compressible Euler equations on the two-dimensional sphere $\mathbb{S}^2$. In contrast to the flat Euclidean case $\mathbb{R}^2$, the geometry of $\mathbb S^2$ imposes new difficulties, and the fluid dynamics are affected by the curved background. To overcome these challenges, we modify the existing modulation method and employ a set of carefully constructed, time-dependent coordinates that precisely track the shock formation on $\mathbb{S}^2$. In particular, we first perform a time-dependent rotation of $\mathbb S^2$, then apply the stereographic projection to the sphere, straighten the steepening shock front, and finally construct shock-adapted coordinates. In the shock-adapted coordinates, the compressible Euler equations on $\mathbb{S}^2$ can be recast into a form suitable for self-similar analysis. Within this framework, we implement a detailed bootstrap argument and establish global well-posedness for the self-similar system. After transferring these results back to the original physical system, we thereby demonstrate the finite-time shock formation on $\mathbb{S}^2$.
\end{abstract}

\tableofcontents

\section{Introduction}
\indent One central question in the mathematical theory of compressible fluids is the shock-formation problem for the compressible Euler equations. While shock formation for the Euler equations has been extensively studied in flat Euclidean space $\mathbb{R}^n$, much less is known about the dynamics of compressible flows on curved manifolds. At the same time, fluid equations posed on curved surfaces, such as the rotating sphere, are relevant in modeling atmospheric and oceanic phenomena on planetary scales, where curvature plays an important role in shaping large-scale flow behavior. One notable example is the tropical cyclone, whose horizontal dynamics can be effectively approximated by compressible flows on $\mathbb{S}^2$. To advance our understanding of the Euler equations on non-flat manifolds, it is important to develop a systematic framework that accurately captures the curvature effects and clarifies how the underlying geometry influences the fluid dynamics at each stage of the analysis. As such, the problem addressed in this paper is both mathematically challenging and physically rich.  

In this paper, we prove the \textit{first} shock formation result for the compressible Euler equations on $\mathbb{S}^2$. As pointed out by Arnold and Khesin in the book \cite{Arnold-Khesin} (Chapter VI, §2), the isentropic compressible Euler equations on a Riemannian manifold $(M, g)$ take the form:
\begin{equation}
\label{the Euler equations}
\left\{\begin{array}{l}
\partial_{\mathring t} \rho+\nabla_{v} \rho+\rho \operatorname{div} v=0, \\
\partial_{\mathring t} v+\nabla_{v} v+\frac{\nabla p}{\rho}=0, \\
p=\frac{1}{\gamma} \rho^{\gamma}.
\end{array}\right.
\end{equation}
Here, $\mathring t \in \mathbb{R}$ denotes the physical time, $\nabla_v$ stands for the covariant derivative on $M$ along the direction of the vector field $v$, and $\nabla p$ is the gradient of the pressure $p$. The unknowns of the system are the velocity field $v \in \Gamma(TM)$, the density function $\rho \in C^\infty(M)$, and the pressure function $p \in C^\infty(M)$. The parameter $\gamma > 1$ is a physical constant known as the adiabatic index. We refer interested readers to \cite{Arnold-Khesin} for a derivation of the above system.

Throughout the paper, for a fixed $r_0>0$, we set $M=\mathbb{S}^2=\{\mathring x\in\mathbb{R}^3: |\mathring x|=r_0\}$, and we use the standard spherical metric $g$ on $\mathbb S^2$. Our result below establishes shock formation for system \eqref{the Euler equations} on manifold $(\mathbb{S}^2,g)$. It marks the \emph{first} rigorous result of shock formation for compressible Euler equations on a curved manifold. In addition, accurate descriptions of the blow-up time, the blow-up location, and the regularity of the resulting shock solution are provided. Owing to the physical interest and other potential applications, we also independently investigate shock formation for equivariant Euler equations on $\mathbb{S}^2$ using different coordinates and modulation variables. The detailed proofs are provided in Appendix \ref{append}\footnote{The proofs presented in this appendix are based on Fulin Qi's \emph{Undergraduate Research Opportunities Programme in Science} (UROPS) research project.}.

\subsection{Main Results}
The rough version of our main theorem is stated below. For the detailed version, readers are referred to Theorem \ref{thm: main theorem}.

\begin{theorem}[Shock Formation for Compressible Euler Equations on $\mathbb{S}^2$]
\label{thm: main theorem rough version}
Consider the isentropic compressible Euler equations \eqref{the Euler equations} on the 2-dimensional sphere $\mathbb{S}^2$
\begin{equation*}
\left\{\begin{array}{l}
\partial_{\mathring t} \rho+\nabla_{v} \rho+\rho \operatorname{div} v=0, \\
\partial_{\mathring t} v+\nabla_{v} v+\frac{\nabla p}{\rho}=0, \\
p=\frac{1}{\gamma} \rho^{\gamma},
\end{array}\right.
\end{equation*}
where the adiabatic index $\gamma >1$. There exists a family of compactly supported smooth initial data, given as perturbations around the self-similar Burgers profile, such that the corresponding solutions develop a point shock in finite time $T_*$. 

At the shock point, the first derivatives of $\rho$ and $v$ blow up at the rate $\frac{1}{T_*-t}$, while the fluid variables themselves remain bounded. In particular, the solution exhibits $C^{1/3}$ regularity at the shock. Moreover, the vorticity is bounded up to the shock and is allowed to be non-trivial.
\end{theorem}

\begin{remark}
We obtain detailed information about the shock singularity formed, including precise descriptions of the blow-up location and time, the blow-up rate, and the regularity of the shock solution. Compared to the Euclidean case, the isentropic Euler equations on the sphere contain additional terms arising from the curved geometry, and the lack of a global Cartesian structure precludes the use of standard translations to construct co-moving coordinates. Moreover, the natural choice of spherical coordinates can lead to expressions that are prohibitively complicated for our analysis. To overcome these difficulties, we employ an approach that combines time-dependent rotations with stereographic projection, thereby establishing a suitable coordinate system well adapted to the shock.
\end{remark}

\begin{remark}
Our theorem also fits into the broader perspective that background geometry can strongly influence shock formation in compressible fluids. In expanding cosmological spacetimes, sufficiently fast expansion induces a damping effect that suppresses shock formation and yields global regularity. This stabilization was established for relativistic Euler equations with the equation of state $p=K\rho$, $0<K<1/3$, on linearly expanding backgrounds $a(t)=t$ by Fajman–Ofner–Oliynyk–Wyatt \cite{non-acc-stability} and was further extended in \cite{linear-expanding-stability}. Numerical experiments in \cite{phase-transition} suggest a sharp transition between stability and shock formation for general power-law expansions $a(t)=t^\alpha$ with $\alpha\in(0,1)$, depending on the competition between the expansion rate and the sound speed: stability is observed when $K<K_{\mathrm{crit}}(\alpha)=1-\tfrac{2}{3\alpha}$ and shock formation occurs when $K>K_{\mathrm{crit}}(\alpha)$. The endpoint instability for radiation fluids $(K=1/3)$ on $a(t)=t$ was rigorously proved earlier by Speck \cite{speck2013stabilizing}. In contrast, the static geometry of $\mathbb{S}^2$ provides no expansion-induced damping. Our result shows that without cosmological expansion, curvature alone cannot prevent shock formation.
\end{remark}

Owing to the physical interest and other potential applications, we also provide an independent and different proof of shock formation for equivariant flows on $\mathbb{S}^2$. In this case, the fluid variables depend only on the latitude (i.e., the angle from the symmetry axis) and are invariant under longitude rotations. In this proof, we employ a different coordinate system on $\mathbb{S}^2$ and design a different set of modulation variables. The statement of this result is as follows. 

\begin{theorem}[Equivariant Shock Formation for Compressible Euler Equations on $\mathbb{S}^2$]
\label{thm: equiv theorem rough version}
For the isentropic compressible Euler equations \eqref{the Euler equations}, there exists a family of equivariant initial data that evolve to solutions ceasing to be smooth in a finite time $T_*$. At $t=T_*$, they develop $C^{1/3}$ shock singularities equipped with $C^1$ blow-up at speed $O(\frac{1}{T_*-t})$.
\end{theorem}

\begin{remark}\label{remark: importance of the appendix}
    Physically speaking, this class of equivariant solutions captures the structure of axisymmetric flows with no swirl, i.e., flows exhibiting no rotation about the axis of symmetry. Such flows serve as an idealized physical model for the zonally symmetric motion of atmosphere. 
\end{remark}

\subsection{Related Literature}
In this section, we review relevant works on both the mathematical theory of shock formation for the Euler equations and the study of fluid mechanics on curved manifolds.

\subsubsection{Shock Formation for the Euler Equations}

In what follows, we provide a brief survey of key results and developments concerning shock formation for the compressible Euler equations.

In one-dimensional flows, shock formation rests on a well-established theoretical foundation and is extensively explored via the method of characteristics. The foundational work is due to Riemann \cite{riemann1860fortpflanzung}, who introduced the powerful method of Riemann invariants for the 1D isentropic Euler equations (a $2 \times 2$ system); these invariants remain constant along characteristic curves. Using Riemann invariants, Lax \cite{Lax1964} proved finite-time shock formation for general $2 \times 2$ strictly-hyperbolic genuinely-nonlinear systems, which in particular covers the small data scenarios for the 1D isentropic Euler equations. For a comprehensive bibliography and review of the early developments in the 1D setting, we refer readers to the monograph by Dafermos \cite{Dafermos_hyperbolic_book}.

The situation becomes substantially more complex for multi-dimensional flows. In \cite{sideris1984formation, sideris1985formation, sideris1997delayed}, Sideris proved finite-time blowup of smooth solutions to the compressible Euler equations for certain classes of small initial data. By using virial-type integral inequalities, Sideris derived upper bounds of the lifespan of smooth solutions, which qualitatively indicates that breakdown must occur in finite time. In \cite{alinhac1999blowup, alinhac1999}, Alinhac developed a more constructive approach to show shock formation for irrotational flows, for which the compressible Euler equations can be reduced to a single scalar quasilinear wave equation for the velocity potential. Without any symmetry assumption, Alinhac also proved that an isolated shock singularity can form from smooth initial data satisfying a non-degeneracy condition via Nash-Moser iteration.

A significant breakthrough was later achieved by Christodoulou in \cite{christodoulou2007formation} regarding the 3D \textit{irrotational} ($\nabla \times v = 0$) relativistic Euler equations, which introduced robust energy methods applicable to a broader class of shock singularities. For the corresponding non-relativistic result, see Christodoulou-Miao \cite{christodoulou2014}. The key innovation in \cite{christodoulou2007formation} was the introduction of a fully geometric framework, which precisely characterizes shock formation via the vanishing inverse foliation density. For the compressible Euler equations with non-trivial vorticity, Luk and Speck \cite{luk2018shock} proved shock formation for the 2D isentropic case with spacetime decompositions. This result was later extended to three dimensions, allowing for entropy variations, in \cite{luk2021stability}. Their analysis builds on a reformulation of the compressible Euler equations developed in \cite{luk2020hidden,speck2019formulation}, which features favorable null structures relative to the acoustic geometry and yields refined regularity properties for vorticity and entropy. For recent developments in this direction, we refer the reader to Abbrescia-Speck \cite{leospeck22,leospeckCQG,leospeckARMA}, Holzegel-Klainerman-Speck-Wong \cite{specksurvey}, Luo-Yu \cite{luoyu1,luoyu2},  Miao-Yu \cite{miao2017formation}, Speck \cite{speckbk,speck2018}, and Speck-Holzegel-Luk-Wong \cite{speck2016stable}. For the elastic waves and the ideal MHD system, An-Chen-Yin \cite{an2020low,an2021low,an2022elasticwaves,an2022MHD,AnJMP} further employed this geometric shock formation mechanism to establish low-regularity ill-posedness. We also refer to the ill-posedness works by Bourgain-Li \cite{illposedness-2d-incompressible,illposedness-3d-incompressible,illposedness-incompressible-Cm} for the incompressible Euler equations.

Meanwhile, an independent and inspiring self-similar framework for proving shock formation was developed by Buckmaster, Shkoller, and Vicol, which focuses on constructing solutions arising from well-prepared initial data that are close to the self-similar Burgers blow-up profile. In self-similar variables, their formulation can be viewed as a renormalized evolution around the Burgers profile, where several modulation variables are used to adjust shifts and scalings dynamically to track the self-similar dynamics. With this approach, a detailed description of the shock-type singularity can be obtained: the (first) shock is shown to be asymptotically self-similar. Furthermore, associated stability properties can also be established (e.g., the construction of stable shocks with $C^{1/3}$ Hölder continuity). In \cite{Buckmaster2022}, they implemented this self-similar approach and proved shock formation for the 2D isentropic Euler equations under azimuthal symmetry. They subsequently extended the self-similar method to prove shock formation for open sets of initial data in 3D without symmetry assumptions in \cite{Buckmaster2023}, and for the full (non-isentropic) Euler equations in \cite{buckmaster2023shock}. This self-similar framework has also motivated further developments, including the unstable shock formation result by Buckmaster and Iyer \cite{buckmaster2022formation}, the works on damped Euler equations by Chen \cite{chen2025shifted}, on the fractal Burgers equation by Chickering-Moreno-Vasquez-Pandya \cite{chickering2023asymptotically} , on the general dispersive or dissipative perturbations of the Burgers equation by Oh-Pasqualotto \cite{gradient-blowup}, on the Euler-Poisson equations by Qiu-Zhao \cite{qiu2021shock}, on the compressible 2D Euler equations by Su \cite{su2023shock}, and on the Burgers-Hilbert equation by Yang \cite{yang2021shock}.

\subsubsection{Fluid Mechanics on Manifolds}
We also review relevant works studying fluid dynamics on curved manifolds. 

The validity of the Euler and Navier-Stokes equations, which govern fluid motions, can be extended to the setting of a curved background. This generalization is motivated by applications in atmospheric and oceanic sciences, where fluids exhibit spherical dynamics. For example, stationary solutions of the Euler equations on a rotating sphere provide building blocks for traveling-wave models that describe stratospheric flows on Earth and giant planets such as Jupiter and Saturn. Constantin-Germain \cite{Constantin-Germain} studied such stationary solutions and proved rigidity and stability results. In addition, Constantin-Germain-Lin-Zhu \cite{constantin2025onset} analyzed the instability of 3-jet zonal flows on a rotating sphere.

To formulate fluid equations on a manifold, one needs to interpret each term in the fluid equations by its corresponding analogue on manifolds. In this paper, we adopt the formulation of the compressible Euler equations for barotropic fluids in \cite{Arnold-Khesin}:
$$\left\{\begin{aligned}
    &\partial_{\mathring t}\rho+\operatorname{div}(\rho v)=0,\\
    &\rho\partial_{\mathring t}v+\rho\nabla_vv+\nabla h(\rho)=0,
\end{aligned}\right.$$
where the pressure function $h(\rho)$ depends on the physical properties of the fluid. For further insights and discussions, see Schutz \cite{Schutz_1980}, Arnold-Khesin \cite{Arnold-Khesin}, and Taylor \cite{Taylor-2023}.  Extending this approach to the Navier-Stokes equations requires choosing an appropriate version of the Laplace operator on the manifold, for which several alternatives exist. We do not delve into these technical details here, and we refer interested readers to Chan-Czubak-Disconzi \cite{Chan-Czubak-Disconzi} and the references therein. 

In a landmark paper by Arnold \cite{Arnold-1966}, he observed that the incompressible Euler equations on a manifold $M$ can be interpreted as geodesic equations. See also the work of Ebin and Marsden \cite{Ebin-Marsden-1970}. An analogous geometric formulation for the compressible Euler equations was developed in Guillemin-Sternberg \cite{Guillemin-Sternberg} and Marsden-Ratiu-Weinstein \cite{Marsden-Ratiu-Weistein-1984}. This geometric point of view was also applied by Lin-Zeng \cite{Lin-Zeng} in the study of unstable manifolds of the Euler equations in a fixed bounded domain.

The study of fluid equations on manifolds has also led to intriguing and impactful mathematical discoveries. For instance, Tao \cite{Tao-universality-I,Tao-universality-II} proved the universality of the incompressible Euler equation on compact manifolds, showing that any finite-dimensional dynamical system satisfying a certain cancellation property can be embedded within the incompressible Euler equation. 

\subsection{Proof Strategy for the Main Theorem} \label{strategy}
In this section, we present the main ideas behind the proof of Theorem~\ref{thm: main theorem rough version}. We employ the framework by Buckmaster-Shkoller-Vicol \cite{Buckmaster2023} in Euclidean space and extend it to the setting of $\mathbb{S}^2$, which requires new geometric insights and analytical adaptations. In the following, we outline the key steps of the proof. 

\bigskip

\noindent\underline{\textbf{Step 1. Construction of geometrically adapted moving coordinates:}}

\indent To analyze how a point shock forms on the sphere $\mathbb{S}^2$, we use the modulation technique. This involves constructing a suitable moving coordinate system that tracks the evolution toward the singularity. This framework enables us to capture the dynamics near the point of steepest gradient. 
\vspace{\baselineskip}

\underline{\emph{Modulation variables.}} Before giving the coordinates, we first introduce the modulation variables, which trace various properties of the shock. In the subsequent construction of coordinates, we will appeal to these modulation variables. 

Let $w(\tilde t, p)$ be a time-dependent scalar function on the sphere $\mathbb{S}^2$ (in this paper, it will be one of the Riemann variables defined in Section \ref{subsection: Riemann Variables}). Assume that $w$ develops an increasingly steep gradient, indicating the formation of the first shock. To analyze the behavior of this emerging singularity, it is essential to introduce coordinates that are centered at the point of maximal gradient $\xi(\tilde t)$ and adapted to the level set of $w$ passing through $\xi(\tilde t)$. This requires tracking several dynamical geometric quantities:
\begin{itemize}
    \item The \emph{predicted shock location} $\xi(\tilde t) \in \mathbb{S}^2$: the point where the gradient magnitude $|\nabla w|$ attains its maximum.
    \item The \emph{predicted shock amplitude} $\kappa(\tilde t)$: a time-dependent function $\kappa(\tilde t)=w(\tilde t,\xi(\tilde t))$ obtained by evaluating $w$ at $\xi(\tilde t)$. 
    \item The \emph{predicted shock normal} $n(\tilde t) \in T_{\xi(\tilde t)}\mathbb{S}^2$: a unit vector at $\xi(\tilde t)$ pointing in the direction of the gradient $\nabla w(\tilde t,\xi(\tilde t))$, perpendicular to the level set $\Gamma(\tilde t)$ of $w$ which passes through $\xi(\tilde t)$.
    \item The \emph{level set curvature} $\psi(\tilde t)$: a scalar function quantifying the curvature of the level set $\Gamma(\tilde t)$ at the point $\xi(\tilde t)$.
    \item The \emph{predicted blow-up time} $\tau(\tilde t)$:
    \[ \tau(\tilde t) = \underbrace{\tilde t}_{\text{current time}} + \underbrace{|\nabla w(\tilde t, \xi(\tilde t))|^{-1}}_{\text{estimated remaining time}}. \]
    The quantity $\tau(\tilde t)$ plays a crucial role in the self-similar rescaling.
\end{itemize}

\underline{\emph{Construction of shock-adapted coordinates.}} Our approach integrates the aforementioned time-dependent modulation variables into the construction of a dynamically adapted coordinate system:
\begin{enumerate}
    \item \textbf{Co-moving coordinates:} We begin by constructing a time-dependent coordinate system that evolves with the shock formation and is centered at the point where the Riemann variable $w$ has the steepest gradient.
    \begin{itemize}
        \item \textbf{Time-dependent rotation of the sphere:} We apply a time-dependent rotation $O(\tilde t) \in SO(3)$ to the sphere, mapping the original Cartesian coordinates $\mathring{x}\in\mathbb{R}^3$ (restricted to $\mathbb{S}^2$) to new coordinates $x = O(\tilde t)^T \mathring x$. The rotation is chosen so that the critical point $\xi(\tilde t)$ maps to a fixed reference location (e.g., the south pole $(0,0,-r_0)$), and the normal vector $n(\tilde t)$ aligns with a fixed direction (e.g., the $(1,0,0)$ direction on the tangent plane at the south pole).
        \item \textbf{Stereographic projection of the rotated sphere:} From the rotated frame $x$, we then employ stereographic projection (from the north pole $(0,0,r_0)$ onto the tangent plane at the south pole $\{x_3 = -r_0\}$) to define a planar coordinate system $u = (u_1, u_2)$. Direct calculation shows
        $$u_{1} = \frac{2 r_{0} x_{1}}{r_0 - x_{3}},\quad u_{2} = \frac{2 r_{0} x_{2}}{r_0 - x_{3}}.$$

        By construction, the origin $u=0$ now corresponds to the predicted shock location $\xi(\tilde t)$, and the $u_1$-axis aligns with the predicted shock normal $n(\tilde t)$. This setup simplifies the spatial derivatives of the Riemann variable $w$ at the origin:
        $$\left.w\right|_{u=0}=\ast,\quad\left.\Big(\frac{\partial w}{\partial u_1},\,\frac{\partial w}{\partial u_2}\Big)\right|_{u=0}=(\ast,0),\quad\left. \begin{pmatrix}
        \displaystyle \frac{\partial^2 w}{\partial u_1^2} & \displaystyle \frac{\partial^2 w}{\partial u_1 \partial u_2} \\
        \displaystyle \frac{\partial^2 w}{\partial u_1 \partial u_2} & \displaystyle \frac{\partial^2 w}{\partial u_2^2}
        \end{pmatrix} \right|_{u = 0}
        =\begin{pmatrix}
        0&0\\
        0&\ast
        \end{pmatrix}.$$
        Here, the symbol $\ast$ denotes that the corresponding quantity (specifically, $w$, $\partial_{u_1} w$, and $\partial_{u_2}^2 w$ evaluated at the origin) is not fixed at this stage. 
    \end{itemize}
    \item \textbf{Level set straightening:} To further adapt the coordinates to the geometry of the shock, we perform a coordinate transformation $(u_1,u_2) \leadsto (\tilde u_1,\tilde u_2)$, realized by a quadratic adjustment incorporating the curvature $\psi(\tilde{t})$ of the level set: 
    $$\tilde u_1=u_1-\frac12\psi(\tilde t)u_2^2,\quad\tilde u_2=u_2.$$
    This transformation flattens the level set geometry near the origin in the $\tilde u$ coordinates, achieving full normalization of the second derivatives:
    $$\left.w\right|_{\tilde{u}=0}=\ast,\quad\left.\Big(\frac{\partial w}{\partial \tilde u_1},\,\frac{\partial w}{\partial \tilde u_2}\Big)\right|_{\tilde {u}=0}=(\ast,0),\quad\left. \begin{pmatrix}
    \displaystyle \frac{\partial^2 w}{\partial \tilde u_1^2} & \displaystyle \frac{\partial^2 w}{\partial \tilde u_1 \partial \tilde u_2} \\
    \displaystyle \frac{\partial^2 w}{\partial \tilde u_1 \partial \tilde u_2} & \displaystyle \frac{\partial^2 w}{\partial \tilde u_2^2}
\end{pmatrix} \right|_{\tilde{u}= 0}=\begin{pmatrix}
    0&0\\
    0&0
\end{pmatrix}.$$

\end{enumerate}

The sequence of coordinate transformations described above can be illustrated in the following figure:

\begin{figure}[H]
    \centering
    \input{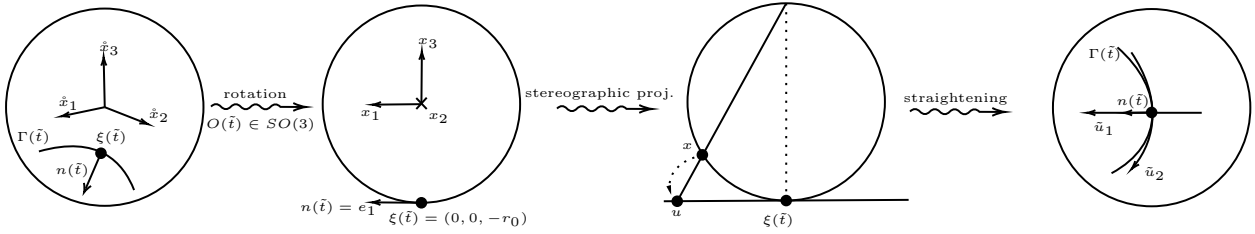} 
    \caption{Coordinate transformations}
    \label{fig:coordinate transformations}
\end{figure}

For clarity, the transformation chain can be schematically expressed as follows:
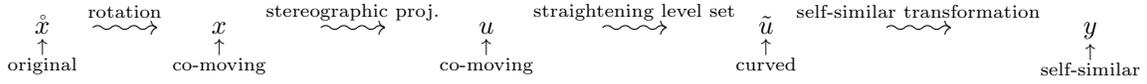
\begin{align*}
\underset{\substack{\uparrow\\\text{original}}}{\mathring{x}}
\overset{\text{rotation}}{\longsquiggly}
\underset{\substack{\uparrow\\\text{co-moving}}}{x}
\overset{\text{stereographic proj.}}{\longsquiggly}
\underset{\substack{\uparrow\\\text{co-moving}}}{u}
\overset{\text{straightening level set}}{\longsquiggly}
\underset{\substack{\uparrow\\\text{curved}}}{\tilde u}
\overset{\text{self-similar transformation}}{\longsquiggly}
\underset{\substack{\uparrow\\\text{self-similar}}}{y}    
\end{align*}

For technical convenience, we take the rotation $O(\tilde t)$ itself as the modulation variable, rather than separately tracking the shock location $\xi(\tilde t)$ and normal direction $n(\tilde t)$. This choice streamlines the expression of coordinate transformations, which would otherwise lead to more complicated formulations in terms of $\xi(\tilde{t})$ and $n(\tilde{t})$.

\bigskip

\noindent\underline{\textbf{Step 2. Symmetrization and diagonalization of the system:}}

\indent To match the transport speed with that of the standard Burgers equation, we set $\alpha=\frac{\gamma-1}{2}$, and rescale the original time $\mathring{t}$ by defining $\tilde{t}=\frac{1+\alpha}{2}\mathring t$.

Expressing the Euler equations \eqref{the Euler equations} directly in the co-moving stereographic $u$-coordinates yields a system \eqref{equations in u_i coordinates} that is not symmetric, risking derivative loss in $H^k$ energy estimates. To address this, we introduce an orthonormal frame ${E_i}$ (see Section \ref{subsection: Reformulation as a Symmetric Hyperbolic System}), and decompose the velocity field as $v = \sum_i V_i E_i$. Denoting $P = (V_1, V_2, \sigma)$, where $\sigma=\frac{1}{\alpha}\rho^\alpha$ is the rescaled sound speed, we obtain a symmetric hyperbolic system:
$$\partial_{\tilde t}P+D_P P+A_{P,u_1}\partial_{u_1}P+A_{P,u_2}\partial_{u_2}P=F_P.$$
Here, the coefficient of the damping term $D_P$ is given by \eqref{def of DP}, the coefficients of the transport terms $A_{P,u_1}$ and $A_{P,u_2}$ are given by \eqref{def of APu1} and \eqref{def of APu2}, and the forcing term $F_P$ is defined in \eqref{def of FP}. We see that $A_{P,u_1}$ and $A_{P,u_2}$ are symmetric matrices. This symmetric hyperbolic structure is preserved under the coordinate transformation $u \leadsto \tilde{u}$.

To effectively analyze shock formation, we further diagonalize the system by introducing Riemann variables $R = (w, z, a)$ (defined in Section \ref{subsection: Riemann Variables}). In these variables, the system becomes
$$\tilde\partial_{\tilde t}R+D_R R+A_{R,\tilde u_1}\partial_{\tilde u_1}R+A_{R,\tilde u_2}\partial_{\tilde u_2}R=F_R,$$
where the coefficient matrix in the $\tilde u_1$-direction is diagonal:
$$A_{R,\tilde{u}_1}=J\mathrm{diag}(w+\beta_2z,\beta_2w+z,\beta_1w+\beta_1z)+\left(g_1-\lambda g_2-\frac{1}{2}\dot\psi \tilde u_2^2\right)I_3. $$
The matrix $A_{R,\tilde u_2}$ and the forcing term $F_R$ are given by
$$A_{R,\tilde u_2}=(2\beta_1\varphi^{-1}V_2+g_2)I_3+2\beta_3\varphi^{-1}\jpnlambda^{-1}S\begin{pmatrix}
-\lambda&0&1\\
0&\lambda&-1\\
\frac{1}{2}&-\frac{1}{2}&0
\end{pmatrix}$$
and
$$\begin{aligned}
F_R=&\begin{pmatrix}
-\frac{\beta_1}{r_0^2\jpnlambda}|V|^2(u_1-\lambda u_2)-\frac{\varphi(V_1-\lambda V_2)}{2r_0^2\jpnlambda}u_j(g_j-2\beta_1\varphi^{-1}V_j)+\frac{\beta_3}{r_0^2}u_jV_jS\\
-\frac{\beta_1}{r_0^2\jpnlambda}|V|^2(u_1-\lambda u_2)-\frac{\varphi(V_1-\lambda V_2)}{2r_0^2\jpnlambda}u_j(g_j-2\beta_1\varphi^{-1}V_j)-\frac{\beta_3}{r_0^2}u_jV_jS\\
-\frac{\beta_1}{r_0^2\jpnlambda}|V|^2(\lambda u_1+u_2)-\frac{\varphi(\lambda V_1+V_2)}{2r_0^2\jpnlambda}u_j(g_j-2\beta_1\varphi^{-1}V_j)
\end{pmatrix}\\
&+\begin{pmatrix}
-\left[2\beta_1\varphi^{-1}(V_2-\alpha\lambda\jpnlambda^{-1}S)+g_2\right]a+\beta_3\varphi^{-1}\jpnlambda^{-1}S(w+z)\\
-\left[2\beta_1\varphi^{-1}(V_2+\alpha\lambda\jpnlambda^{-1}S)-g_2\right]a-\beta_3\varphi^{-1}\jpnlambda^{-1}S(w+z)\\
(2\beta_1\varphi^{-1}V_2+g_2)\frac{w+z}{2}
\end{pmatrix}\jpnlambda^{-2}\partial_{\tilde u_2}\lambda\\
&+\jpnlambda^{-2} \begin{pmatrix}
-a\\
-a\\
\frac{w+z}{2}
\end{pmatrix}\tilde u_2\pt\psi.
\end{aligned}$$

Therefore, The scalar equation for $w$ is
$$\boxed{\tilde\partial_{\tilde t}w}
+[J(\boxed{w}+\beta_2 z)+g_1-\lambda g_2-\frac12\dot\psi\tilde u_2^2]\boxed{\partial_{\tilde u_1}w}
+\cdots=F_R^1,$$
where $F_R^1$ denotes the first component of the vector $F_R$. The boxed terms coincide with those appearing in the one-dimensional inviscid Burgers equation and will dominate the dynamics near the singularity. Moreover, we will see that $J= 1+o(1)$ in the bootstrap analysis, and therefore the factor does not interfere with the dominant Burgers-type dynamics.

\bigskip

\noindent\underline{\textbf{Step 3. Self-similar transformation and modulation ODEs:}}

\indent
To analyze the singularity, we employ a self-similar transformation, which dynamically rescales the solution and effectively ``zooms in" on the singularity. The solution is expected to remain regular in the self-similar coordinates. 

We define the self-similar coordinates and the self-similar variables as follows:
$$s=-\ln(\tau(\tilde t)-\tilde t),\quad y_1=e^{\frac{3}{2}s}\tilde u_1,\quad y_2=e^{\frac s2}\tilde u_2.$$
$$W(s,y)=e^{\frac{s}{2}}(w(\tilde t,\tilde u)-\kappa(\tilde t)),\quad Z(s,y)=z(\tilde t,\tilde u),\quad A(s,y)=a(\tilde t,\tilde u).$$
This transformation yields a system of three transport-type equations for the variables $W$, $Z$, and $A$ in the self-similar coordinate system:
\begin{equation}
\left\{\begin{aligned}
&\left(\partial_{s}-\frac{1}{2}\right)W+\left(\frac{3}{2} y_{1}+\beta_\tau JW+G_{W}\right) \partial_{y_1} W+\left(\frac{1}{2} y_{2}+h_{W}\right) \partial_{y_2}W =F_{W}-\beta_{\tau} e^{-\frac{s}{2}} \partial_{\tilde t} \kappa, \\
&\partial_{s} Z+\left(\frac{3}{2} y_{1}+\beta_2\beta_\tau JW+G_{Z}\right) \partial_{y_1} Z+\left(\frac{1}{2} y_{2}+h_{Z}\right) \partial_{y_2} Z =F_{Z}, \\
&\partial_{s} A+\left(\frac{3}{2} y_{1}+\beta_1\beta_\tau JW+G_{A}\right) \partial_{y_1} A+\left(\frac{1}{2} y_{2}+h_{A}\right) \partial_{y_2} A =F_{A}.
\end{aligned}\right. 
\end{equation}
For any given variable $\mathcal{R}\in\{W,Z,A\}$, the coefficients of $\partial_{y_1}\mathcal{R}$ and $\partial_{y_2}\mathcal{R}$ in the transport equations are given by 
$$\frac{3}{2} y_{1} + \beta_{\mathcal{R}}\beta_{\tau} J W + G_{\mathcal{R}} \quad \text{and} \quad \frac{1}{2} y_{2} + h_{\mathcal{R}},$$
respectively, where $\beta_{W} = 1$, $\beta_{Z} = \beta_2$, and $\beta_{A} = \beta_1$. The term $F_{\mathcal{R}}$ denotes the corresponding forcing term.

In our analysis, the shock formation roots in the self-similar type blow-up for Burgers equation of $W$. To capture this mechanism, we employ the following self-similar profile $\ovl W$, which solves the two-dimensional self-similar Burgers equation:
$$-\frac12\ovl W+\left(\frac32y_1+\ovl W\right)\partial_{y_1}\ovl W+\frac12y_2\partial_{y_2}\ovl W=0.$$
In particular, we select the following $\ovl W$ as the blow-up profile:
$$\ovl W(y_1,y_2)=\langle y_2\rangle \ovl W_{1d}\left(\langle y_2\rangle^{-3}y_1\right),$$
where $\langle y_2\rangle=\sqrt{1+y_2^2}$, and $\ovl W_{1d}(y)$ is implicitly defined by the relation $y=-\ovl W_{1d}(y)-\ovl W_{1d}(y)^3$ for $y\in\mathbb{R}$.

In our construction, we normalize $W$ at the origin $y=0$ up to its second derivatives such that
$$\left.W\right|_{y=0}=0,\quad\left.\Big(\frac{\partial W}{\partial  y_1},\,\frac{\partial W}{\partial  y_2}\Big)\right|_{y=0}=(-1,0),\quad\left. \begin{pmatrix}
    \displaystyle \frac{\partial^2 W}{\partial y_1^2} & \displaystyle \frac{\partial^2 W}{\partial y_1 \partial y_2} \\
    \displaystyle \frac{\partial^2 W}{\partial y_1 \partial y_2} & \displaystyle \frac{\partial^2 W}{\partial y_2^2}
\end{pmatrix} \right|_{y= 0}=\begin{pmatrix}
    0&0\\
    0&0
\end{pmatrix}.$$
These values are identical to those of the profile $\ovl W$:
$$W(s,0)=\ovl W(0),\quad\nabla_yW(s,0)=\nabla_y\ovl W(0),\quad\nabla_y^2 W(s,0)=\nabla_y^2\ovl W(0).$$
This normalization is a direct consequence of our selection of modulation variables. These constraints, in turn, enable the derivation of the evolution ODEs for the parameters $\kappa$, $\tau$, $\psi$, and $O$. These equations are enumerated in \eqref{eqn of kappa}\eqref{eqn of tau}\eqref{eqn of Q12}\eqref{eqn of Q13,Q23}\eqref{eqn of psi}.

We now arrive at a system of three PDEs for $W$, $Z$, and $A$, coupled with six ODEs for the modulation variables $\kappa$, $\tau$, $\psi$, and $O$. Provided that the self-similar transformation remains non-degenerate, this system is equivalent to the compressible Euler's equations on the sphere $\mathbb{S}^2$.

\bigskip

\noindent\underline{\textbf{Step 4. Global well-posedness of the self-similar system:}}

\indent Our aim is to prove the global well-posedness of the self-similar system. Since the self-similar coordinates become singular in finite time, establishing global existence in this frame implies finite-time shock formation when translated back to the original physical coordinates. 

To achieve this, we employ a bootstrap argument. The bootstrap assumptions are organized into four categories, asserting that
\begin{enumerate}
    \item The modulation variables behave similarly to those of the classical Burgers equation.
    \item The dominant Riemann variable $W$ remains close to a known blow-up profile $\ovl W$.
    \item The secondary Riemann variable $Z$ and the tangential velocity component $A$ stay sufficiently small, allowing them to be treated perturbatively.
    \item The energy of $W,Z$, and $A$ are appropriately controlled.
\end{enumerate}

In what follows, we outline the main ideas on how to close these bootstrap assumptions.

\noindent\underline{\emph{Step 4a. Estimates for the modulation variables:}}

\indent We use the bootstrap assumptions to estimate each term on the right-hand side of the modulation ODEs \eqref{eqn of kappa}\eqref{eqn of tau}\eqref{eqn of Q12}\eqref{eqn of Q13,Q23}\eqref{eqn of psi}, thereby improving the assumptions for the modulation variables:
$$|\pt\kappa|\lesssim(\ln M)^{\frac{1}{10}}\tau_0^{\frac16},\quad |\pt\tau|\lesssim M^{\frac14}e^{-s},\quad\left|\pt \psi-\left(2\psi^2+\frac{3}{r_0^2}\right)\beta_3\kappa_0\right|\lesssim\tau_0^{\frac16}\ln M,$$
$$|Q_{12}|\lesssim \tau_0^{\frac13}\quad,|Q_{23}|\lesssim M^{\frac34}\tau_0^{\frac12},\quad\left|Q_{13}-\frac{2\beta_3\kappa_0}{r_0}\right|\lesssim\tau_0^{\frac16}.$$
See Section \ref{section: Estimates for Modulation Variables} for the details.

\noindent\underline{\emph{Step 4b. Estimates for $W$ and $\widetilde W$:}}

This step forms the centerpiece of the analysis, as the behavior of $W$ governs the shock formation process. We focus on improving the bootstrap assumptions for the dominant Riemann variable $W$, and its deviation from the explicit blow-up profile $\ovl W$. The primary goal is to control the perturbation $\widetilde W = W - \ovl W$, such that the solution remains close to the known self-similar profile.

Our method centers on a detailed analysis of the transport equations for $W$, $\widetilde W$, and their derivatives. Indeed, for $\mathcal{R}\in\{W,\widetilde W\}$ and $\gamma\ge0$, we can systematically write that
$$\partial_s\partial_y^\gamma\mathcal{R}+D_{\mathcal{R}}^{(\gamma)}\partial_y^\gamma\mathcal{R}+\mathcal{V}_W\cdot\nabla_y\partial_y^\gamma\mathcal{R}=F_{\mathcal{R}}^{(\gamma)}.$$
Since these transport equations share the same transport coefficients, it is convenient to define the Lagrangian trajectories (or characteristics) $\Phi_W$:
$$\left\{\begin{aligned}
    &\partial_s \Phi_{W}(s;s_1,y_0) = \mathcal{V}_{W}(s,\Phi_{W}(s;s_1,y_0)), \\
    &\Phi_{W}(s_1;s_1,y_0) = y_0.
\end{aligned}\right.$$
By integrating along these characteristics, we establish pointwise bounds on $\partial^\gamma W$ or $\partial^\gamma\widetilde W$. An important property of $\Phi_W$ is that if the starting point $(s_1,y_0)$ is slightly away from the origin, i.e. $|y_0|\ge l$ for a small constant $l=(\ln M)^{-5}$, then $|\Phi_W(s;s_1,y_0)|$ grows exponentially:
$$|\Phi_W(s;s_1,y_0)|\ge|y_0|e^{\frac{s-s_1}{5}}.$$
This enables us to carry out analysis separately in two distinct spatial regions, distinguishing the dominant dynamics near the origin ($y = 0$) from those farther away:
\begin{itemize}
    \item Estimates in $\{|y|\le l\}$: In this neighborhood of $y=0$, the choice of modulation variables ensures that $W$ remains close to $\ovl W$ by enforcing the constraints 
    $$W(s,0)=\ovl W(0),\quad\nabla_y W(s,0)=\nabla_y\ovl W(0),\quad \nabla_y^2 W(s,0)=\nabla_y^2\ovl W(0).$$
    Furthermore, the damping term for $D^4\widetilde W$ admits a positive lower bound:
    $$D_{\widetilde W}^{(\gamma)}(s,y)\ge\frac13\quad\text{for }|\gamma|=4\text{ and }|y|\le l,$$
    which contributes to maintaining this closeness. We show that the derivatives $\partial^\gamma W$ remain close to $\partial^\gamma\ovl W$ for all multi-indices $|\gamma|\le4$, enabling us to improve the bootstrap assumptions in this region.
    \item Estimates in $\{|y|> l\}$: In the far-field, we perform weighted $L^\infty$ estimates for $\mathcal{R}\in\{W,\widetilde W\}$ by employing a weighted function of the form $q = \eta^{\mu(\gamma)}\partial^\gamma_y\mathcal{R}$, which satisfies a transport equation
    $$\partial_sq+D_qq+\mathcal{V}_W\cdot\nabla_yq=F_q.$$
    Solving $q$ along the Lagrangian trajectories $\Phi_W$, we find
    \begin{align*}
        \eta^{(\gamma)}(y)|\partial_y^\gamma \mathcal{R}(s,y)|\le\ &\underbrace{\eta^{(\gamma)}(y_0)|\partial_y^\gamma \mathcal{R}(s_1,y_0)|\exp\left(\int_{s_1}^sD_q(s',\Phi_W(s';s_1,y_0))\dif s'\right)}_{\text{contributions from the initial data or the boundary valus at $\{|y|=l\}$}}\\
        &+\underbrace{\int_{s_1}^s|F_q(s',\Phi_W(s';s_1,y_0))|\exp\left(\int_{s_1}^sD_q(s'',\Phi_W(s'';s_1,y_0))\dif s''\right)\dif s'}_{\text{contributions from the forcing term}}.
    \end{align*}
    Here, $(s_1,y_0)$ is the starting point of the trajectory $\Phi_W(s';s_1,y_0)$, and $\Phi(s;s_1,y_0)=y$. Due to the property of $\Phi_W$, we can require either $s_1=s_0$ or $|y_0|=l$. The critical insight is that the smallness of $\widetilde W$ in the inner region $\{|y|\le l\}$ is propagated outward along these trajectories. The exponential growth of $|\Phi_W|$ then converts the spatial decay of the forcing terms $|F_W^{(\gamma)}|$ into rapid time decay, ensuring their contributions remain dominated by boundary values at $|y| = l$. This behavior contrasts with the estimates for $Z$ and $A$, where the bootstrap assumptions are primarily determined by the magnitude of the forcing terms themselves. Consequently, the main analytical task reduces to establishing sharp, pointwise estimates for these forcing terms $F_W^{(\gamma)}$. As detailed in Section \ref{section: Estimates for Transport Equations}, these forcing terms are complicated functions of the primary variables and the geometric quantities, and controlling their derivatives is a crucial part of the process. 
\end{itemize}

By systematically combining the estimates from both the inner and outer regions, we successfully improve the bootstrap assumptions \eqref{asmp of W} and \eqref{asmp of W tilde} for $W$ and $\widetilde W$, respectively. The rigorous derivations for these bounds are presented in detail in Section \ref{section: Estimates for W}.

\noindent\underline{\emph{Step 4c. Estimates for $Z$ and $A$:}}

In contrast to the dominant variable $W$, the secondary Riemann variable $Z$ and the tangential velocity component $A$ are treated as perturbations that are expected to remain small throughout the evolution. The primary strategy for controlling these variables and their derivatives is to analyze their respective transport equations, as given in system \eqref{eqn of W,Z,A}.

The main idea of this step is to integrate the equations for $\partial_y^\gamma Z$ and $\partial_y^\gamma A$ along their Lagrangian trajectories. As in Step 4b, the main task is to estimate the forcing terms, $F_Z^{(\gamma)}$ and $F_A^{(\gamma)}$. This is detailed in Section \ref{section: Estimates for Transport Equations}. This approach successfully establishes bounds on $Z$, $A$, and all of their derivatives required by the bootstrap assumptions, with the exception for $\partial_{y_1}A$.

However, an exception occurs for the normal derivative of the tangential velocity $\partial_{y_1}A$. For this term, direct analysis via using the transport equation  fails to close the bootstrap argument. Additionally, one of the bootstrap assumptions--specifically the estimate of the tangential velocity derivative $\partial_1 A$ (where $A$ denotes the tangential velocity component, and $\partial_1$ indicates differentiation in the normal direction $\tilde u_1$)--cannot be closed via directly using transport estimates and the main energy estimates. To overcome this, we analyze the evolution of the \emph{vorticity}. Following \cite{Taylor-2023}, we define the vorticity scalar $\omega$ on the sphere $\mathbb{S}^2$ as:
\begin{equation}
    \omega=\ast\dif v^\flat=\epsilon^{ab}\nabla_a v_b.
\end{equation}
Here, $\ast$ is the Hodge dual operator, $\dif$ denotes the exterior derivative, and $\flat:TM\rightarrow T^*M$ is the musical isomorphism induced by the metric $g$. The tensor $\epsilon_{ab}$ represents the volume form associated with $g$. Then we derive that $\omega$ satisfies the evolution equation:
\begin{equation}
    \partial_{\mathring t}\omega+\nabla_v\omega+\omega\operatorname{div}v=0.
\end{equation}
This equation indicates that vorticity is not merely transported by the flow, but is also affected by compressibility through the divergence term. In contrast, the \emph{specific vorticity}, defined as $\zeta = \omega / \rho$ (where $\rho$ is the density), satisfies a pure transport equation:
\begin{equation}
\label{eqn of specific vorticity}
    \partial_{\mathring t}\zeta+\nabla_v\zeta=0.
\end{equation}

To effectively utilize the vorticity equation, we need to express the vorticity $\omega$ in terms of the variables we are working with. We derive the following key identity:
$$\omega=\varphi^{-2}\left[\varphi\jpnlambda\partial_{\tilde u_1}a-\partial_{\tilde u_2}(\varphi V_1)\right],$$
where $a$ is the tangential velocity component, $V_1$ is the velocity component in the $E_1$ direction (aligned with $\tilde u_1$), and $\varphi, \jpnlambda$ are auxiliary quantities determined by the coordinate system and metric, which are explicitly defined in Section \ref{section: Derivation of Equations}. 

This identity connects the vorticity $\omega$ (whose associated specific vorticity $\zeta$ satisfies a favorable transport equation) to the term $\partial_{\tilde u_1} a$, which corresponds to $\partial_1 A$ under self-similar scaling. By exploiting the transport structure of equation \eqref{eqn of specific vorticity}, we gain control on $\|\zeta\|_{L^\infty}$, and hence $\|\omega\|_{L^\infty}$ via using bounds on $\rho$). Together with the above identity and estimates on additional terms like $\partial_{\tilde u_2} (\varphi V_1)$, we can finally close the problematic bootstrap assumption for $\partial_1 A$. 

\noindent\underline{\emph{Step 4d. Energy estimates:}}

As mentioned above, most bootstrap assumptions impose $L^\infty$ -type bounds on the solution and its derivatives, which are typically closed via transport equations. However, since
$$A_{R,\tilde u_2}=(2\beta_1\varphi^{-1}V_2+g_2)I_3+2\beta_3\varphi^{-1}\jpnlambda^{-1}S\begin{pmatrix}
-\lambda&0&1\\
0&\lambda&-1\\
\frac{1}{2}&-\frac{1}{2}&0
\end{pmatrix},$$
which is not diagonal, the equation of $w$ includes the term $2\beta_3\varphi^{-1}\jpnlambda^{-1}\sigma\partial_{\tilde u_2}a$. This term involves the derivative of $a$, resulting in a loss of derivatives in the $L^\infty$ estimates. To address this difficulty, we supplement the $L^\infty$ analysis with $H^k$ energy estimates:
$$e^{-s}\|W\|_{\dot H^k_y}^2+\|Z\|_{\dot H^k_y}^2+\|A\|_{\dot H^k_y}^2\lesssim M^{3k-5}e^{-s}.$$

\bigskip

\noindent\underline{\textbf{Step 5. Recovering the physical variables and proving the main theorem:}}

\noindent
With the bootstrap argument successfully completed in self-similar coordinates, the final step is to translate the resulting estimates in the original physical coordinates, thereby proving the main theorem. Specifically, we have:
\begin{align*}
    \partial_{u_1}w|_{u=0}&=\partial_{\tilde u_1}w|_{\tilde u_=0}=e^s\partial_{y_1}W|_{y=0}=-e^s=-\frac{1}{\tau(\tilde t)-\tilde t}\approx-\frac{1}{\widetilde T_*-\tilde t},
\end{align*}
where $\widetilde T_*$ is the blow-up time. This computation highlights that the shock formation roots in the singularity of the self-similar coordinate transformation. The full procedure for this translation, together with the final argument, is detailed in Section \ref{subsec: Bootstrap Assumptions and Framework of Bootstrap Argument}.

\subsection{Organization of the Paper}

The rest of the paper is structured as follows:

\begin{itemize}
    \item In Section \ref{section: Derivation of Equations}, we perform a sequence of coordinate transformations, deriving the governing equations at each step. Then, we introduce the Riemann variables, and define the self-similar transformation. Finally, we present the blow-up profile and derive the modulation ODEs that govern the evolution of its parameters.

    \item In Section \ref{section: Main Theorem and Reduction to a Bootstrap Argument}, we specify the initial data, including the initial values of the modulation variables and the initial data for self-similar Riemann variables. We then state the main theorem, formulate the bootstrap assumptions, and outline the strategy for closing the argument.

    \item In Section \ref{section: Estimates for Transport Equations}, we employ the bootstrap assumptions to establish estimates on various terms appearing in the transport equations. These estimates lay the groundwork for subsequent analysis.

    \item Sections \ref{section: Estimates for Modulation Variables}, \ref{section: Estimates for Z and A}, and \ref{section: Estimates for W} are devoted to improving the $L^\infty$-bounds in the bootstrap argument. Specifically, Section \ref{section: Estimates for Modulation Variables} focuses on the estimates of modulation variables. Sections \ref{section: Estimates for Z and A} and \ref{section: Estimates for W} address the estimates for self-similar variables $Z$, $A$, and $W$, respectively.

    \item Section \ref{section: energy estimates} completes the argument by establishing energy estimates that allow us to prevent derivative loss and complete the analysis.
\end{itemize}

\subsection{Notation}
For the reader's convenience, we record below the key notation that will be employed repeatedly in the remainder of this work:

\begin{itemize}
    \item The Greek indices $\mu, \nu, \dots$ range over $1,2,3$, while the Latin indices $a,b,c,i,j,\dots$ range over $1,2$.
    
    \item Physical constants: $\alpha=\frac{\gamma-1}{2}>0$, $\beta_{1}=\frac{1}{1+\alpha}$, $\beta_{2}=\frac{1-\alpha}{1+\alpha}$, $\beta_{3}=\frac{\alpha}{1+\alpha}$.

    \item Rescaled sound speed: $\sigma=\frac{1}{\alpha}\rho^\alpha$; steady-state value outside its support: $\sigma_\infty$.

    \item Original time-space coordinates: $(\mathring t,\mathring x)\in\mathbb{R}\times\mathbb{R}^3$; rescaled time: $\tilde{t}=\frac{1+\alpha}{2}\mathring t$.

    \item Stereographic projection coordinates: $u_i$; Shock-adapted coordinates: $\tilde u_i$.

    \item Metric factor: $\varphi(u)=\frac{4 r_{0}^{2}}{u_1^2 + u_2^2 + 4 r_{0}^{2}}$; sphere metric: $g_{\mathbb{S}^2}= \varphi(u)^2 \sum_i\dif u_i\otimes\dif u_i$.

    \item Modulation variables: $O(\tilde t)$, $\psi(\tilde t)$, $\tau(\tilde t)$, $\kappa(\tilde t)$.

    \item Time derivatives: $\mathring\partial_{\tilde t}=(\partial_{\tilde t})_{\mathring x}$ (for fixed $\mathring{x}$), $\partial_{\tilde t}=(\partial_{\tilde t})_{x,u}$ (for fixed $x,u$), $\tilde \partial_{\tilde t}=(\partial_{\tilde t})_{\tilde u}$ (for fixed $\tilde u$).

    \item Auxiliary quantities: $\lambda=\psi\tilde{u}_2$, $\langle\lambda\rangle=\sqrt{1+\lambda^2}$, $J=|\nabla \tilde{u}_1|=\varphi^{-1}\langle\lambda\rangle$.
    \item Vector fields:
    $N=\frac{E_1-\lambda E_2}{\langle\lambda\rangle}$, $T=\frac{\lambda E_1+E_2}{\langle\lambda\rangle}$ with $E_i=\varphi^{-1}\partial_{u_i}$.

    \item Multi-index notations: for $\gamma=(\gamma_1,\gamma_2)\in\mathbb{Z}_{\ge0}^2$, we denote $|\gamma|=\gamma_1+\gamma_2$ and define $\partial_{\tilde u}^\gamma=\partial_{\tilde u_1}^{\gamma_1}\partial_{\tilde u_2}^{\gamma_2}$. In self-similar coordinates ($y$-coordinates), we set $\partial^\gamma=\partial_{y_1}^{\gamma_1}\partial_{y_2}^{\gamma_2}=e^{-\frac{3\gamma_1+\gamma_2}{2}s}\partial_{\tilde u}^\gamma$. For brevity, we also adopt the notation $\partial_{12}=\partial_{y_1}\partial_{y_2}$, $\partial_{11}=\partial_{y_1}^2$, etc.

    \item Velocity fields: $v=v_1\partial_{u_1}+v_2\partial_{u_2}=V_1E_1+V_2E_2$.

    \item Riemann variables: $w=v\cdot N+\sigma$, $z=v\cdot N-\sigma$, $a=v\cdot T$.

    \item Self-similar coordinates: $s=-\ln(\tau(\tilde t)-\tilde t)$, $y_1=e^{\frac{3}{2}s}\tilde u_1$, $y_2=e^{\frac{s}{2}}\tilde u_2$.

    \item Self-similar variables: $w=e^{-\frac{s}{2}}W+\kappa(\tilde t)$, $z=Z$, $a=A$; and $\sigma=S$.

    \item Vector notation: $P=(V_1,V_2,\sigma)=(V_1,V_2,S)$, $\widetilde P=(V_1,V_2,\sigma-\sigma_\infty)=(V_1,V_2,S-\sigma_\infty)$, $R=(w,z,a)$.

    \item We use calligraphic letters (e.g., $\mathcal{R}$) to represent a generic self-similar variable among $W,Z$ and $A$.

    \item Superscript ``$0$" denotes evaluation at $y=0$, for example: $Z^0:=Z(s,0)$, $G_W^0:=G_W(s,0)$, $\partial_1F_W^0:=\partial_1F_W(s,0)$, etc.
\end{itemize}

\subsection{Acknowledgements}
Xinliang An is supported by MOE Tier 1 grants A-0004287-00-00, A-0008492-00-00, A-8002933-00-00, and MOE Tier 2 grant A-8000977-00-00. Haoyang Chen is supported by MOE Tier 1 grant A-0008492-00-00, Tier 2 grant A-8000977-00-00, and NSFC (Grant No. 12171097).

\section{Preliminaries}
\label{section: Derivation of Equations}
In this section, we introduce the equations, coordinate systems, and transformations used throughout the paper. We begin with the Euler equations in their original form and introduce a sequence of coordinate transformations--as illustrated in Figure \ref{fig:coordinate transformations}--to track shock formation. These include a time rescaling aligned with the Burgers equation, a rotation of the sphere, a stereographic projection, an adaptation to shock geometry. We then define the Riemann variables, apply a self-similar ansatz to capture the blow-up structure, and identify the modulation parameters to formulate the transformed system.

\subsection{Rescaled Time and Sound Speed}
First, we introduce the rescaled time and the sound speed variable, and rewrite the system with respect to them. Setting $\alpha=\frac{\gamma-1}{2}$, we define the rescaled time $\tilde t$ by $\tilde{t}=\frac{1+\alpha}{2}\mathring t$, which yields $\partial_{\mathring t}=\frac{\dif \tilde{t}}{\dif \mathring t} \partial_{\tilde t}=\frac{1+\alpha}{2} \partial_{\tilde t}=\frac{1}{2\beta_1}\partial_{\tilde t}$. To simplify the notation, we also define the following three constants
\begin{equation} \label{beta param}
\beta_{1}=\frac{1}{1+\alpha},\quad \beta_{2}=\frac{1-\alpha}{1+\alpha},\quad \beta_{3}=\frac{\alpha}{1+\alpha}.
\end{equation}
After rescaling time, the Euler equations take the form
\begin{equation}
\left\{\begin{aligned}
&\frac{1}{2 \beta_{1}}\partial_{\tilde{t}} \rho+\nabla_{v} \rho+\rho \operatorname{div} v=0, \\
&\frac{1}{2 \beta_{1}}\partial_{\tilde{t}} v+\nabla_{v} v+\frac{\nabla p}{\rho}=0, \\
&p=\frac{1}{\gamma} \rho^{\gamma}.
\end{aligned}\right.
\end{equation}
Similarly, we introduce the rescaled sound speed $\sigma=\frac{1}{\alpha}\rho^\alpha$ to rewrite the above equations as follows:
\begin{equation}
\label{rescaled Euler equations}
\left\{\begin{aligned}
&\frac{1}{2 \beta_{1}}\partial_{\tilde{t}} \sigma+\nabla_{v} \sigma+\alpha\sigma \operatorname{div} v=0, \\
&\frac{1}{2 \beta_{1}}\partial_{\tilde{t}} v+\nabla_{v} v+\alpha\sigma\nabla \sigma=0.
\end{aligned}\right.
\end{equation}

\subsection{Rotation and Stereographic Projection}
\label{subsect: Rotation and Stereographic Projection}
In this section, we first perform an appropriate rotation of the sphere and then apply stereographic projection to construct the time-dependent coordinates. Suppose the original spatial variables of $\mathbb{R}^3$ are $\mathring{x}$, and consider the sphere $\mathbb{S}^2=\{\mathring x\in\mathbb{R}^3: |\mathring x|=r_0\}\subset\mathbb{R}^3$. We introduce a time-dependent rotation matrix $O(\tilde t) = (O_{\mu\nu}) \in SO(3)$, which transforms the original coordinates $\mathring x$ into $x = O(\tilde t)^T \mathring x$ (in components, $x_\mu = O_{\nu\mu} \mathring x_\nu$). Next, we define the  $3\times3$ skew-symmetric matrix $Q = \dot O(\tilde t)^T O(\tilde t)$ (equivalently, $Q_{\mu\nu} = \dot O_{\delta\mu} O_{\delta\nu}$), which represents the angular velocity of the rotation. This yields the relations $\mathring\partial_{\tilde t} x_\mu = Q_{\mu\nu} x_\nu$ and $\mathring\partial_{\tilde t} \frac{\partial}{\partial x_\mu} = Q_{\mu\nu} \frac{\partial}{\partial x_\nu}$. Here, we distinguish $\mathring\partial_{\tilde t} := (\partial_{\tilde t})_{\mathring x}$, the time derivative holding $\mathring x$ fixed, from $\partial_{\tilde t} := (\partial_{\tilde t})_{x}$, the time derivative holding $x$ fixed.

We then map $\mathbb{S}^2$ onto the plane $\{x \in \mathbb{R}^3 : x_3 = -r_0\}$ via stereographic projection from the north pole $(0, 0, r_0)$. Specifically, for each $p \in \mathbb{S}^2$, we define $q(p)$ as the intersection of the line connecting $p$ to the north pole with this plane. By setting $u_i(p) = x_i(q)$, we obtain a coordinate chart on $\mathbb{S}^2$. Straightforward calculations show that
\begin{align}
\label{def of u coordinates}
    u_{i} = \frac{2 r_{0} x_{i}}{r_0 - x_{3}} = \frac{x_{i}}{\varphi(u)},\quad \text{for}\,\,i=1,2,
\end{align}
where
\begin{align}
\label{def of varphi}
    \varphi(u) := \frac{4 r_{0}^{2}}{|u|^{2} + 4 r_{0}^{2}} = \frac{4 r_{0}^{2}}{u_1^2 + u_2^2 + 4 r_{0}^{2}}.
\end{align}
The metric on $\mathbb{S}^2$ is then given by
$$ g_{\mathbb{S}^2} = \sum_\mu \dif x_\mu \otimes \dif x_\mu |_{T\mathbb{S}^2} = \varphi(u)^2 \sum_i \dif u_i \otimes \dif u_i. $$
Moreover, we deduce
\begin{align}
\label{def of gi}
g_i := \mathring\partial_{\tilde t} u_i = -\frac{1}{2 r_{0}} (Q_{13} u_{1} + Q_{23} u_{2}) u_{i} + Q_{i j} u_{j} + Q_{i 3} r_{0} (\varphi^{-1} - 2),
\end{align}
and $$\mathring\partial_{\tilde t} \frac{\partial}{\partial u_i} = h_{ij} \frac{\partial}{\partial u_j},$$
with
\begin{equation}
\label{def of hij}
\begin{array}{c}
h_{11} = h_{22} = \frac{1}{2 r_{0}} (Q_{13} u_{1} + Q_{23} u_{2}), \\
h_{12} = -h_{21} = Q_{12} + \frac{1}{2 r_{0}} (Q_{13} u_{2} - Q_{23} u_{1}).
\end{array}    
\end{equation}
The operator $\mathring\partial_{\tilde t}$ relates to $\partial_{\tilde t}$ through
\[ \mathring\partial_{\tilde t} = \partial_{\tilde t} + \mathring\partial_{\tilde t} u_j \frac{\partial}{\partial u_j} = \partial_{\tilde t} + g_j \frac{\partial}{\partial u_j}. \]

Now, we can express the Euler equations in the $u_i$ coordinates as follows:
\begin{equation}
\label{equations in u_i coordinates}
\left\{
\begin{aligned}
&\partial_{\tilde t} v_{i} + h_{j i} v_{j} + (2\beta_1 v_{j} + g_j) \frac{\partial v_{i}}{\partial u_{j}} + 2\beta_1 \varphi^{-1} \left( 2 v_{i} v_{j} \frac{\partial \varphi}{\partial u_{j}} - \sum_{j} v_j^2 \frac{\partial \varphi}{\partial u_{i}} \right) + 2\beta_3  \varphi^{-2}\sigma \frac{\partial \sigma}{\partial u_{i}} = 0, \\
&\partial_{\tilde t} \sigma + (2\beta_1 v_{j} + g_j) \frac{\partial \sigma}{\partial u_{j}} + 2\beta_3 \sigma \left( \frac{\partial v_{j}}{\partial u_{j}} + 2 \varphi^{-1} v_{j} \frac{\partial \varphi}{\partial u_{j}} \right) = 0.
\end{aligned}
\right.
\end{equation}

\subsection{Reformulation as a Symmetric Hyperbolic System}
\label{subsection: Reformulation as a Symmetric Hyperbolic System}
To avoid loss of derivatives in the energy estimates, we reformulate the system into symmetric hyperbolic form.

At $p\in\mathbb S^2$, we define an orthonormal frame $E_i=\varphi^{-1}\partial_{u_i}$ for the tangent space $T_p\,\mathbb S^2$. A velocity vector $v\in T_p\,\mathbb S^2$ can then be expanded as $v=V_iE_i=v_i\partial_{u_i}$, so that
$$V_i=\varphi v_i.$$
We also introduce the notation
\begin{equation*}
S=\sigma,\quad P=(V_1,V_2,S)^T.
\end{equation*}
The system \eqref{equations in u_i coordinates} can then be rewritten as
\begin{equation}
\partial_{\tilde t}P+D_P P+A_{P,u_1}\partial_{u_1}P+A_{P,u_2}\partial_{u_2}P=F_P,
\label{system in u coordinates}
\end{equation}
where
\begin{equation}
\label{def of DP}
D_P=\begin{pmatrix}
    h_{11}&h_{21}&0\\
    h_{12}&h_{22}&0\\
    0&0&0
\end{pmatrix},
\end{equation}
\begin{equation}
\label{def of APu1}
A_{P,u_1}=(2\beta_1\varphi^{-1}V_1+g_1)I_3+2\beta_3\varphi^{-1}S\begin{pmatrix}
    0&0&1\\
    0&0&0\\
    1&0&0
\end{pmatrix},
\end{equation}
\begin{equation}
\label{def of APu2}
A_{P,u_2}=(2\beta_1\varphi^{-1}V_2+g_2)I_3+2\beta_3\varphi^{-1}S\begin{pmatrix}
    0&0&0\\
    0&0&1\\
    0&1&0
\end{pmatrix},
\end{equation}
and
\begin{equation}
\label{def of FP}
\begin{aligned}
F_P=&\begin{pmatrix}
\varphi^{-1}V_1(g_j-2\beta_1\varphi^{-1}V_j)\frac{\partial\varphi}{\partial u_j}+2\beta_1\varphi^{-2}|V|^2\frac{\partial\varphi}{\partial u_1}\\
\varphi^{-1}V_2(g_j-2\beta_1\varphi^{-1}V_j)\frac{\partial\varphi}{\partial u_j}+2\beta_1\varphi^{-2}|V|^2\frac{\partial\varphi}{\partial u_2}\\
-2\beta_3\varphi^{-2}SV_j\frac{\partial\varphi}{\partial u_j}
\end{pmatrix}\\
=&\begin{pmatrix}
-\frac{\beta_1}{r_0^2}u_1|V|^2-\frac{\varphi}{2r_0^2}V_1u_j(g_j-2\beta_1\varphi^{-1}V_j)\\
-\frac{\beta_1}{r_0^2}u_2|V|^2-\frac{\varphi}{2r_0^2}V_2u_j(g_j-2\beta_1\varphi^{-1}V_j)\\
\frac{\beta_3}{r_0^2}u_jV_jS
\end{pmatrix}.
\end{aligned}
\end{equation}
Here, $I_3$ denotes the $3\times3$ identity matrix, and $|V|^2=V_1^2+V_2^2$. 

\subsection{Coordinates Adapted to the Shock}
To track the evolution of the shock more effectively in the $u_i$ coordinate system, we introduce a new set of coordinates $\tilde{u}_i$ adapted to the shock formation. Let $\psi(\tilde t) \in \mathbb{R}$ be a time-dependent parameter controlling the curvature of the shock front. The new coordinates are defined by:
\begin{equation}
\label{def of tilde u coordinates}
\left\{\begin{aligned}
    &\tilde{u}_1 = u_1 - \frac{1}{2} \psi(\tilde t) u_2^2, \\
    &\tilde{u}_2 = u_2.
\end{aligned}\right.
\end{equation}
Here, $\psi(\tilde t)$ adjusts the transformation dynamically. The corresponding differential operators become
\begin{equation}
\begin{cases}
\partial_{\tilde t} = \tilde\partial_{\tilde t} - \frac{1}{2} \dot{\psi} \tilde{u}_2^2 \partial_{\tilde{u}_1}, \\
\partial_{u_1} = \partial_{\tilde{u}_1}, \\
\partial_{u_2} = \partial_{\tilde{u}_2} - \psi \tilde{u}_2 \partial_{\tilde{u}_1},
\end{cases}\nonumber
\end{equation}
where $\tilde\partial_{\tilde t} := (\partial_{\tilde t})_{\tilde{u}}$ denotes the time derivative with $\tilde{u}_i$ fixed. 

Next, we introduce auxiliary quantities to describe the geometry of the shock. Define $\langle a \rangle = \sqrt{1 + a^2}$ for $a \in \mathbb{R}$, and set
\begin{equation}
\lambda = \psi \tilde{u}_2, \quad J = |\nabla \tilde{u}_1| = \varphi^{-1} \langle \lambda \rangle, \quad N = J^{-1} \nabla \tilde{u}_1 = \langle \lambda \rangle^{-1} (E_1 - \lambda E_2), \quad T = \langle \lambda \rangle^{-1} (\lambda E_1 + E_2).\nonumber
\end{equation}
Note that $\{N, T\}$ forms an orthonormal frame, with $N$ aligned to the gradient of $\tilde{u}_1$ and $T$ orthogonal to it. Moreover, by defining the 1-form $\theta(X)=\langle\nabla_XN,T\rangle_g$, one can verify that
\begin{equation}
    \label{deriv of N,T}\nabla_NN=\theta(N)T,\quad\nabla_NT=-\theta(N)N,\quad\nabla_TN=\theta(T)T,\quad\nabla_TT=-\theta(T)N,
\end{equation}
and
\begin{equation}
\label{def of 1-form theta}
    \theta(N)=-\frac{u_2+\lambda u_1}{2r_0^2\langle\lambda\rangle}+\frac{\lambda\psi}{\varphi\langle\lambda\rangle^3},\quad\theta(T)=\frac{u_1-\lambda u_2}{2r_0^2\langle\lambda\rangle}+\frac{\psi}{\varphi\langle\lambda\rangle^3}.
\end{equation}

Applying these transformations, we rewrite system \eqref{system in u coordinates} in the $\tilde{u}_i$ coordinates as
\begin{equation}
\label{system in tilde u coordinates}
\tilde\partial_{\tilde t} P + D_P P + A_{P,\tilde{u}_1} \partial_{\tilde{u}_1} P + A_{P,\tilde{u}_2} \partial_{\tilde{u}_2} P = F_P,
\end{equation}
where the coefficient matrices are given by
\begin{equation}
\label{def of AP tilde u1}
\begin{aligned}
A_{P,\tilde u_1}&=A_{P,u_1}-\lambda A_{P,u_2}-\frac{1}{2}\dot\psi\tilde u_2^2 I_3\\
&=\left[2\beta_1\varphi^{-1}(V_1-\lambda V_2)+(g_1-\lambda g_2)-\frac{1}{2}\dot\psi \tilde u_2^2\right]I_3+2\beta_3\varphi^{-1}S\begin{pmatrix}
    0&0&1\\
    0&0&-\lambda\\
    1&-\lambda&0
\end{pmatrix},
\end{aligned}
\end{equation}
and
\begin{equation}
\label{def of AP tilde u2}
A_{P,\tilde u_2}=A_{P,u_2}=(2\beta_1\varphi^{-1}V_2+g_2)I_3+2\beta_3\varphi^{-1}S\begin{pmatrix}
    0&0&0\\
    0&0&1\\
    0&1&0
\end{pmatrix}.
\end{equation}

\subsection{Riemann Variables}
\label{subsection: Riemann Variables}
Building on the reformulated system in $\tilde{u}_i$ coordinates, we now introduce Riemann variables to diagonalize the equations. We define the Riemann variables as below:
\begin{equation}
\label{def of w,z,a}
\left\{
\begin{aligned}
&w=v\cdot N+\sigma=\jpnlambda^{-1}(V_1-\lambda V_2)+S,\\
&z=v\cdot N-\sigma=\jpnlambda^{-1}(V_1-\lambda V_2)-S,\\
&a=v\cdot T=\jpnlambda^{-1}(\lambda V_1+V_2),
\end{aligned}
\right.
\end{equation}
or equivalently
\begin{equation}
    R:=(w,z,a)^T=BP,\nonumber
\end{equation}
where
\begin{equation}
B=\begin{pmatrix}
\jpnlambda^{-1}&-\lambda\jpnlambda^{-1}&1\\
\jpnlambda^{-1}&-\lambda\jpnlambda^{-1}&-1\\
\lambda\jpnlambda^{-1}&\jpnlambda^{-1}&0
\end{pmatrix}.\nonumber
\end{equation}
The Riemann variables satisfy the following system:
\begin{equation}
\label{eqn of wza, u tilde}
\tilde\partial_{\tilde t}R+D_R R+A_{R,\tilde u_1}\partial_{\tilde u_1}R+A_{R,\tilde u_2}\partial_{\tilde u_2}R=F_R.
\end{equation}
Here, the damping term is
\begin{equation}
D_R=BD_PB^{-1}=\begin{pmatrix}
h_{11}&0&0\\
0&h_{22}&0\\
0&0&0
\end{pmatrix}+h_{12}\begin{pmatrix}
0&0&-1\\
0&0&-1\\
\frac{1}{2}&\frac{1}{2}&0
\end{pmatrix}.\nonumber
\end{equation}
The transport terms are given by
\begin{equation}
\label{def of AR tilde u1}
\begin{aligned}
A_{R,\tilde u_1}&=BA_{P,\tilde u_1}B^{-1} \\
&=\left[2\beta_1\varphi^{-1}(V_1-\lambda V_2)+(g_1-\lambda g_2)-\frac{1}{2}\dot\psi \tilde u_2^2\right]I_3+2\beta_3\varphi^{-1}\jpnlambda S\begin{pmatrix}
    1&0&0\\
    0&1&0\\
    0&0&0
\end{pmatrix}\\
&=J\mathrm{diag}(w+\beta_2z,\beta_2w+z,\beta_1w+\beta_1z)+\left(g_1-\lambda g_2-\frac{1}{2}\dot\psi \tilde u_2^2\right)I_3,
\end{aligned}\nonumber
\end{equation}
\begin{equation}
\label{def of AR tilde u2}
A_{R,\tilde u_2}=BA_{P,\tilde u_2}B^{-1} 
=(2\beta_1\varphi^{-1}V_2+g_2)I_3+2\beta_3\varphi^{-1}\jpnlambda^{-1}S\begin{pmatrix}
-\lambda&0&1\\
0&\lambda&-1\\
\frac{1}{2}&-\frac{1}{2}&0
\end{pmatrix}.\nonumber
\end{equation}
And the forcing term is
\begin{equation}
\label{def of FR, vector form}
\begin{aligned}
F_R=&BF_P+(\wtt B+A_{R,\tilde u_1}\underbrace{\partial_{\tilde u_1}B}_{=0}+A_{R,\tilde u_2}\partial_{\tilde u_2}B)P\\
=&BF_P+\jpnlambda^{-2} \begin{pmatrix}
-a\\
-a\\
\frac{w+z}{2}
\end{pmatrix}\tilde u_2\pt\psi+\jpnlambda^{-2} A_{R,\tilde u_2}\begin{pmatrix}
-a\\
-a\\
\frac{w+z}{2}
\end{pmatrix}\partial_{\tilde u_2}\lambda\\
=&\begin{pmatrix}
-\frac{\beta_1}{r_0^2\jpnlambda}|V|^2(u_1-\lambda u_2)-\frac{\varphi(V_1-\lambda V_2)}{2r_0^2\jpnlambda}u_j(g_j-2\beta_1\varphi^{-1}V_j)+\frac{\beta_3}{r_0^2}u_jV_jS\\
-\frac{\beta_1}{r_0^2\jpnlambda}|V|^2(u_1-\lambda u_2)-\frac{\varphi(V_1-\lambda V_2)}{2r_0^2\jpnlambda}u_j(g_j-2\beta_1\varphi^{-1}V_j)-\frac{\beta_3}{r_0^2}u_jV_jS\\
-\frac{\beta_1}{r_0^2\jpnlambda}|V|^2(\lambda u_1+u_2)-\frac{\varphi(\lambda V_1+V_2)}{2r_0^2\jpnlambda}u_j(g_j-2\beta_1\varphi^{-1}V_j)
\end{pmatrix}\\
&+\begin{pmatrix}
-\left[2\beta_1\varphi^{-1}(V_2-\alpha\lambda\jpnlambda^{-1}S)+g_2\right]a+\beta_3\varphi^{-1}\jpnlambda^{-1}S(w+z)\\
-\left[2\beta_1\varphi^{-1}(V_2+\alpha\lambda\jpnlambda^{-1}S)-g_2\right]a-\beta_3\varphi^{-1}\jpnlambda^{-1}S(w+z)\\
(2\beta_1\varphi^{-1}V_2+g_2)\frac{w+z}{2}
\end{pmatrix}\jpnlambda^{-2}\partial_{\tilde u_2}\lambda\\
&+\jpnlambda^{-2} \begin{pmatrix}
-a\\
-a\\
\frac{w+z}{2}
\end{pmatrix}\tilde u_2\pt\psi.
\end{aligned}\nonumber
\end{equation}

\subsection{Self-Similar Transformation}
We now introduce the self-similar coordinates to resolve the shock formation. Let $\tau(\tilde t)>0$ be a modulation variable representing the expected blow-up time. The self-similar coordinates are defined by
\begin{equation}
s=-\ln(\tau(\tilde t)-\tilde t),\quad y_1=e^{\frac{3}{2}s}\tilde u_1,\quad y_2=e^{\frac{s}{2}}\tilde u_2.\nonumber
\end{equation}
In the self-similar coordinates, the system \eqref{eqn of wza, u tilde} becomes
\begin{equation}
\partial_{s} R+\beta_{\tau} e^{-s} D_R R+\left(\frac{3}{2} y_{1}+\beta_{\tau} e^{\frac{s}{2}} A_{R,\tilde u_1}\right) \partial_{1} R+\left(\frac{1}{2} y_{2}+\beta_{\tau} e^{-\frac{s}{2}} A_{R, \tilde u_2}\right) \partial_{2} R=\beta_{\tau} e^{-s} F_{R},\nonumber
\end{equation}
Next, we define the self-similar Riemann variables as
\begin{equation}
\left\{
\begin{aligned}
&w=e^{-\frac{s}{2}}W+\kappa(\tilde t),\\
&z=Z,\\
&a=A.
\end{aligned}\right.\nonumber
\end{equation}
Then, we derive the transport equations for $W,Z$ and $A$:
\begin{equation}
\label{eqn of W,Z,A}
\left\{\begin{aligned}
&\left(\partial_{s}-\frac{1}{2}\right)W+\left(\frac{3}{2} y_{1}+g_{W}\right) \partial_{1} W+\left(\frac{1}{2} y_{2}+h_{W}\right) \partial_{2}W =F_{W}-\beta_{\tau} e^{-\frac{s}{2}} \partial_{\tilde t} \kappa, \\
&\partial_{s} Z+\left(\frac{3}{2} y_{1}+g_{Z}\right) \partial_{1} Z+\left(\frac{1}{2} y_{2}+h_{Z}\right) \partial_{2} Z =F_{Z}, \\
&\partial_{s} A+\left(\frac{3}{2} y_{1}+g_{A}\right) \partial_{1} A+\left(\frac{1}{2} y_{2}+h_{A}\right) \partial_{2} A =F_{A}.
\end{aligned}\right. 
\end{equation}
Here $\beta_\tau:=\frac{1}{1-\partial_{\tilde t}\tau}$, and $\partial_i:=\partial_{y_i}$. The transport terms are
\begin{equation}
\label{def of gR,GR,G}
\left\{
\begin{aligned}
&g_{W}=\beta_{\tau} J W+G_{W}=\beta_{\tau} J W+ \beta _ { \tau } e ^ { \frac { s } { 2 } } \left[ J(\beta_2Z + \kappa) + G \right], \\
&g_{Z}=\beta_2 \beta_{\tau} J W+G_{Z}=\beta_2 \beta_{\tau} J W+\beta _ { \tau } e ^ { \frac { s } { 2 } } \left[ J (Z + \beta_2 \kappa) + G\right], \\
&g_{A}=\beta_1 \beta \tau J W+G_{A}=\beta_1 \beta \tau J W+\beta _ { \tau } e ^ { \frac { s } { 2 } } \left[\beta_1J( Z + \kappa) + G\right],
\end{aligned}\right.
\end{equation}
where
\begin{equation}
\label{def of G}
G:=g_1-\lambda g_2- \frac { 1 } { 2 } \dot { \psi } \tilde { u } _ { 2 } ^ { 2 } ,
\end{equation}
and
\begin{equation}
\label{def of hR}
\left\{
\begin{aligned}
&h_W=\beta_\tau e^{-\frac{s}{2}}\left[2\beta_1\varphi^{-1}\left(V_2-\alpha\varphi^{-1}\lambda\jpnlambda^{-1}S\right)+g_2\right],\\
&h_Z=\beta_\tau e^{-\frac{s}{2}}\left[2\beta_1\varphi^{-1}\left(V_2+\alpha\varphi^{-1}\lambda\jpnlambda^{-1}S\right)+g_2\right],\\
&h_A=\beta_\tau e^{-\frac{s}{2}}\left(2\beta_1\varphi^{-1}V_2+g_2\right).
\end{aligned}\right.
\end{equation}
The forcing terms take the form
\begin{equation}
\label{def of FR}
\left\{\begin{aligned}
F_W=&-\beta_\tau e^{-\frac{s}{2}}\left[\frac{\beta_1 |V|^2}{r_0^2\jpnlambda}(u_1 - \lambda u_2) + \frac{V_1 - \lambda V_2}{2r_0^2J} u_j \left(g_j - 2\beta_1 \varphi^{-1} V_j \right) - \frac{\beta_3}{r_0^2} u_jV_jS\right]\\
&-e^{\frac{s}{2}} \jpnlambda^{-2} \left[\beta_\tau e^{-s}\pt\psi\tilde u_2 A + h_W A\partial_2 \lambda - e^{-\frac{s}{2}} h_{W,A} v\cdot N \partial_2\lambda \right]\\
&-\beta_\tau e^{-\frac{s}{2}}\left[h_{11}(e^{-\frac s2}W+\kappa)-h_{12} A\right] - h_{W, A} \partial_2 A,\\
F_Z=&-\beta_\tau e^{-s}\left[\frac{\beta_1 |V|^2}{r_0^2\jpnlambda}(u_1 - \lambda u_2) + \frac{V_1 - \lambda V_2}{2r_0^2J} u_j \left(g_j - 2\beta_1 \varphi^{-1} V_j \right) - \frac{\beta_3}{r_0^2} u_jV_jS\right]\\
&-\jpnlambda^{-2} \left[\beta_\tau e^{-s}\pt\psi\tilde u_2A + h_Z A \partial_2 \lambda- h_{Z,A} v\cdot N \partial_2\lambda \right]\\
&-\beta_\tau e^{-s}\left(h_{22}Z-h_{12} A\right) - h_{Z, A} \partial_2 A,\\
F_A=&-\beta_\tau e^{-s}\left[\frac{\beta_1|V|^2}{r_0^2\jpnlambda}(\lambda u_1+u_2)+\frac{\lambda V_1+V_2}{2r_0^2J}u_j(g_j-2\beta_1\varphi^{-1}V_j)\right]\\
&-\beta_\tau e^{-s}h_{12} v\cdot N- h_{A, W} \partial_2 W-h_{A, Z} \partial_2 Z + \beta_\tau\jpnlambda^{-2}e^{-s}\pt\psi\tilde u_2v\cdot N,
\end{aligned}\right.
\end{equation}
with
\begin{equation}
\label{def of hRR}
\left\{
\begin{aligned}
&h_{W, A} = 2\beta_3 \beta_\tau \varphi^{-1} \jpnlambda^{-1} S, \quad h_{Z, A} = - e^{\frac{s}{2}} h_{W, A},\\
&h_{A, W} = \frac{1}{2} e^{-s} h_{W, A}, \quad h_{A, Z} = -e^{\frac{s}{2}} h_{A, W} = h_{Z, A}.
\end{aligned}\right.
\end{equation}

For each $\mathcal{R}\in\{W,Z,A\}$, by introducing the notation $\mathcal{V}_{\mathcal{R}}=(\mathcal{V}_{\mathcal{R}}^1,\mathcal{V}_{\mathcal{R}}^2)=(\frac{3}{2}y_1+g_{\mathcal{R}},\frac{1}{2}y_2+h_{\mathcal{R}})$, we can write the equation of $\mathcal{R}$ in a compact form:
\begin{equation}
\label{eqn of R, cpt form}
\partial_s\mathcal{R}+D_{\mathcal{R}}\mathcal{R}+\mathcal{V}_{\mathcal{R}}\cdot\nabla\mathcal{R}=F_{\mathcal{R}}-\mathbbm{1}_{\mathcal{R}=W}\beta_\tau e^{-\frac{s}{2}}\partial_{\tilde t}\kappa,
\end{equation}
where $D_{\mathcal{R}}=-\frac{1}{2}\mathbbm{1}_{\mathcal{R}=W}$ denotes the damping term. Here, the indicator function is defined by
$$\mathbbm{1}_{\mathcal{R}=W}=\begin{cases}
    1,&\mathcal{R=W},\\
    0,&\mathcal{R}=Z\text{ or }A.
\end{cases}$$

For any multi-index $\gamma>0$, by applying $\partial^\gamma=\partial_1^{\gamma_1}\partial_2^{\gamma_2}$ to both side of (\ref{eqn of W,Z,A}), we get
\begin{equation}
	\label{eqns of deriv of W,Z,A}
	\left\{\begin{aligned}
		&\left(\partial_s+\frac{3\gamma_1+\gamma_2-1}{2}+\beta_\tau(1+\gamma_1\mathbbm{1}_{|\gamma|\ge2})J\partial_1W\right)\partial^\gamma W+\mathcal{V}_W\cdot\nabla\partial^\gamma W=F_W^{(\gamma)}\\
		&\left(\partial_s+\frac{3\gamma_1+\gamma_2}{2}+\beta_2\beta_\tau\gamma_1J\partial_1W\right)\partial^\gamma Z+\mathcal{V}_Z\cdot\nabla\partial^\gamma Z=F_Z^{(\gamma)}\\
		&\left(\partial_s+\frac{3\gamma_1+\gamma_2}{2}+\beta_1\beta_\tau\gamma_1J\partial_1W\right)\partial^\gamma A+\mathcal{V}_A\cdot\nabla\partial^\gamma A=F_A^{(\gamma)},
	\end{aligned}
	\right.
\end{equation}
with the forcing terms being
\begin{equation}
	\begin{aligned}
		F_W^{(\gamma)}=&\partial^\gamma F_W-\beta_\tau\partial_1W[\partial^\gamma,J]W-\gamma_2\beta_\tau\mathbbm{1}_{|\gamma|\ge2}\partial_2(JW)\partial_1^{\gamma_1+1}\partial_2^{\gamma_2-1}W\\
		&-\beta_\tau\mathbbm{1}_{|\gamma|\ge3}\sum_{\substack{1\le|\beta|\le|\gamma|-2\\\beta\le\gamma}}\binom{\gamma}{\beta}\partial^{\gamma-\beta}(JW)\partial_1\partial^\beta W-\sum_{0\le\beta<\gamma}\binom{\gamma}{\beta}\left(\partial^{\gamma-\beta}G_W\partial_1\partial^\beta W+\partial^{\gamma-\beta}h_W\partial_2\partial^\beta W\right),
	\end{aligned}
	\label{def of FW^gamma}
\end{equation}
\begin{equation}
	\begin{aligned}
		F_Z^{(\gamma)}=&\partial^\gamma F_Z-\gamma_2\beta_2\beta_\tau\mathbbm{1}_{|\gamma|\ge2}\partial_2(JW)\partial_1^{\gamma_1+1}\partial_2^{\gamma_2-1}Z\\
		&-\beta_2\beta_\tau\mathbbm{1}_{|\gamma|\ge2}\sum_{\substack{0\le|\beta|\le|\gamma|-2\\\beta\le\gamma}}\binom{\gamma}{\beta}\partial^{\gamma-\beta}(JW)\partial_1\partial^\beta Z-\sum_{0\le\beta<\gamma}\binom{\gamma}{\beta}\left(\partial^{\gamma-\beta}G_Z\partial_1\partial^\beta Z+\partial^{\gamma-\beta}h_Z\partial_2\partial^\beta Z\right),
	\end{aligned}
	\label{def of FZ^gamma}
\end{equation}
\begin{equation}
	\begin{aligned}
		F_A^{(\gamma)}=&\partial^\gamma F_A-\gamma_2\beta_1\beta_\tau\mathbbm{1}_{|\gamma|\ge2}\partial_2(JW)\partial_1^{\gamma_1+1}\partial_2^{\gamma_2-1}A\\
		&-\beta_2\beta_\tau\mathbbm{1}_{|\gamma|\ge2}\sum_{\substack{0\le|\beta|\le|\gamma|-2\\\beta\le\gamma}}\binom{\gamma}{\beta}\partial^{\gamma-\beta}(JW)\partial_1\partial^\beta A-\sum_{0\le\beta<\gamma}\binom{\gamma}{\beta}\left(\partial^{\gamma-\beta}G_A\partial_1\partial^\beta A+\partial^{\gamma-\beta}h_A\partial_2\partial^\beta A\right).
	\end{aligned}
	\label{def of FA^gamma}
\end{equation}

Here and hereafter, we use several different abbreviations for the indicator function. For example, we write
$$\mathbbm{1}_{|\gamma|\ge2}:=\mathbbm{1}_{\{\beta\in\mathbb{N}_0^2:|\beta|\ge2\}}(\gamma)=\begin{cases}
    1,&|\gamma|\ge2,\\
    0,&|\gamma|<2.
\end{cases}$$

For $\gamma \ge 0$, the equations for higher-order derivatives of $\mathcal{R} \in \{W, Z, A\}$ can be systematically written as follows:  
\begin{equation}
\partial_s \partial^\gamma \mathcal{R} + D_{\mathcal{R}}^{(\gamma)} \partial^\gamma \mathcal{R} + \mathcal{V}_\mathcal{R} \cdot \nabla \partial^\gamma \mathcal{R} = F_{\mathcal{R}}^{(\gamma)}.\nonumber
\end{equation}
Here, we have for any multi-index $\gamma\ge0$ that
\begin{align}
    \label{def of DRgamma}
    D_{\mathcal{R}}^{(\gamma)}=\begin{cases}
        \frac{3\gamma_1+\gamma_2-1}{2}+\beta_\tau(\mathbbm{1}_{|\gamma|>0}+\gamma_1\mathbbm{1}_{|\gamma|\ge2})J\partial_1W,&\mathcal{R}=W,\\
        \frac{3\gamma_1+\gamma_2}{2}+\beta_2\beta_\tau\gamma_1J\partial_1W,&\mathcal{R}=Z,\\
        \frac{3\gamma_1+\gamma_2}{2}+\beta_1\beta_\tau\gamma_1J\partial_1W,&\mathcal{R}=A.
    \end{cases}
\end{align}
For $\gamma > 0$, the forcing terms $F_W^{(\gamma)}, F_Z^{(\gamma)}$, and $F_A^{(\gamma)}$ are defined in \eqref{def of FW^gamma}, \eqref{def of FZ^gamma}, and \eqref{def of FA^gamma}, respectively. For the case $\gamma = 0$, we set 
$$F_W^{(0,0)} := F_W - \beta_\tau e^{-s/2} \partial_{\tilde{t}} \kappa, \quad \text{and} \quad F_{\mathcal{R}}^{(0,0)} := F_{\mathcal{R}} \quad \text{for } \mathcal{R} \in \{Z, A\}.$$

We also record the following relations:
\begin{equation}
\label{P expressed by R}
\left\{\begin{aligned}
    &V_1=\frac{1}{2}\jpnlambda^{-1}(e^{-\frac{s}{2}}W+\kappa+Z)+\lambda \jpnlambda^{-1}A,\\
    &V_2=-\frac{1}{2}\lambda\jpnlambda^{-1}(e^{-\frac{s}{2}}W+\kappa+Z)+\jpnlambda^{-1}A,\\
    &S=\frac{1}{2}(e^{-\frac{s}{2}}W+\kappa-Z),\\
    &v\cdot N=\jpnlambda^{-1}(V_1-\lambda V_2)=\frac{1}{2}(e^{-\frac{s}{2}}W+\kappa+Z),
\end{aligned}\right.
\end{equation}
\begin{equation}
\label{R expressed by P}
\left\{\begin{aligned}
    &W=e^{\frac{s}{2}}\left[\jpnlambda^{-1}(V_1-\lambda V_2)+S-\kappa\right],\\
    &Z=\jpnlambda^{-1}(V_1-\lambda V_2)-S,\\
    &A=\jpnlambda^{-1}(\lambda V_1+V_2),
\end{aligned}\right.
\end{equation}
along with the equation of $P=(V_1,V_2,S)^T$ in self-similar coordinates:
\begin{equation}
\label{eqn of P, self-similar}
    \partial_{s} P+\beta_{\tau} e^{-s} D_P P+\left(\frac{3}{2} y_{1}+\beta_{\tau} e^{\frac{s}{2}} A_{P,\tilde u_1}\right) \partial_{1} P+\left(\frac{1}{2} y_{2}+\beta_{\tau} e^{-\frac{s}{2}} A_{P, \tilde u_1}\right) \partial_{2} P=\beta_{\tau} e^{-s} F_{P}.
\end{equation}

\subsection{Blow-up Profile}\label{subsection: blow-up profile}
In this section, we analyze the blow-up profile by closely following the construction introduced in \cite{collot2018singularity}. While our profile largely builds upon this prior work, we adapt it to suit the self-similar framework developed in the previous sections and emphasize its relevance to our shock formation scenario.

We first introduce the 1D self-similar Burgers profile
\begin{equation}
	\ovl W_{1d}(y_1)=\left(-\frac{y_1}{2}+\left(\frac{1}{27}+\frac{y_1^2}{4}\right)^{\frac{1}{2}}\right)^{\frac{1}{3}}-\left(\frac{y_1}{2}+\left(\frac{1}{27}+\frac{y_1^2}{4}\right)^{\frac{1}{2}}\right)^{\frac{1}{3}},\nonumber
\end{equation}
which solves the 1D self-similar Burgers equation: 
\begin{equation}
	-\frac{1}{2}\ovl W_{1d}+\left(\frac{3}{2}y_1+\ovl W_{1d}\right)\partial_{y_1}\ovl W_{1d}=0.\nonumber
\end{equation}
Then, define $\ovl W(y_1,y_2)=\langle y_2\rangle \ovl W_{1d}\left(\langle y_2\rangle^{-3}y_1\right)$, where $\langle y_2\rangle=\sqrt{1+y_2^2}$. One can verify that $\overline{W}$ is a solution to the 2D self-similar Burgers equation:
\begin{equation}
	\label{self-similar Burgers eqn}
	-\frac{1}{2}\overline{W}+\left(\frac{3}{2}y_1+\overline{W}\right)\partial_{y_1}\overline{W}+\frac{1}{2}y_2\partial_{y_2}\overline{W}=0,
\end{equation}
and it serves as the shock profile in this paper. From the definition of $\ovl W$, one can verify that
\begin{equation}
	\ovl W(0)=0,\ \ \nabla\ovl W(0)=(-1,0)^T,\ \ \nabla^2\ovl W(0)=\begin{pmatrix}
		0&0\\
		0&0
	\end{pmatrix},\ \ \nabla^2\partial_1\ovl W(0)=\begin{pmatrix}
	6&0\\
	0&2
	\end{pmatrix}.
	\label{evaluation of Wbar and derivs at 0}
\end{equation}
Furthermore, the profile $\ovl W$ satisfies the following estimates:
\begin{proposition}
Let $\eta(y)=1+y_1^2+y_2^6$, then we have
\begin{equation}
	\left|\overline{W}\right|\le\eta^{\frac{1}{6}},\ 
	\left|\partial_1\ovl{W}\right|\le\eta^{-\frac{1}{3}}, \ \left|\partial_2\ovl{W}\right|\le\frac{\sqrt{3}}{3}.
	\label{est of ovl W and D ovl W}
\end{equation}
For any multi-index $\gamma\ge0$, it holds that
\begin{equation}
	\label{est of Dk ovl W}
    \left|\partial^\gamma\ovl{W}\right|\lesssim_\gamma\eta^{\frac{1}{6}-\frac{\gamma_1}{2}-\frac{\gamma_2}{6}}.
\end{equation}
\end{proposition}
\begin{proof}
    Note that $t=-\ovl W_{1d}(t)-\ovl W^3_{1d}(t)$. Using this identity, we deduce that $|\ovl W_{1d}(t)|\le\min(|t|,|t|^{\frac{1}{3}})$ and $\ovl W_{1d}'=-\frac{1}{1+3\ovl W_{1d}^2}$. It then follows that
    $$\ovl W\le\langle y_2\rangle\min(\langle y_2\rangle^{-3}|y_1|,\langle y_2\rangle^{-1}|y_1|^{\frac{1}{3}})\le|y_1|^{\frac{1}{3}}\le\eta^{\frac{1}{6}},$$
    and
    $$|\partial_2\ovl W|=2|y_2|\langle y_2\rangle^{-1}\frac{|\ovl W_{1d}|}{1+3\ovl W_{1d}^2}\le\frac{\sqrt{3}}{3}.$$
    To estimate $\partial_1\ovl W$, we employ the identity $\partial_1\ovl W=-\frac{1}{\langle y_2\rangle^2+3\ovl W^2}$. Hence, it suffices to show that $\eta^{\frac{1}{3}}\le\langle y_2\rangle^2+3\ovl W^2$, which is equivalent to $1+(\ovl W_{1d}+\ovl W_{1d}^3)^2+y_2^6\le(1+y_2^2)^3(1+3\ovl W_{1d}^2)^3$. Expanding the left-hand side and moving all terms to the right-hand side, we obtain a polynomial in $y_2^2$ and $\overline{W}_{1d}^2$ with positive coefficients, thus confirming the inequality.

    Using the formula $\partial_2\ovl W=2y_2\ovl W\partial_1\ovl W$, one can then show the estimate for higher derivatives $|\partial^\gamma\ovl{W}|\lesssim_\gamma\eta^{\frac{1}{6}-\frac{\gamma_1}{2}-\frac{\gamma_2}{6}}$, proceeding by first induction on $|\gamma|$ and for each fixed $|\gamma|$, then by induction on $\gamma_2$.
\end{proof}

\subsection{Perturbation Equation}
To control the deviation of $W$ from the profile $\ovl W$, we derive the perturbation equation. Let $\widetilde W=W-\ovl W$ denote the error between $W$ and the blow-up profile $\ovl W$. We can thus deduce the equation of $\widetilde{W}$:
\begin{equation}
	\left(\partial_s-\frac{1}{2}+\beta_\tau J\partial_1\ovl{W}\right)\widetilde{W}+\mathcal{V}_W\cdot\nabla\widetilde{W}=F_{\widetilde{W}},
	\label{eqn of tilde W}
\end{equation}
where
\begin{equation}
F_{\widetilde{W}}=F_W+\left[(1-\beta_\tau J)\ovl{W}-G_W\right]\partial_1\ovl{W}-h_W\partial_2\ovl{W}-\beta_\tau e^{-\frac{s}{2}}\partial_{\tilde t}\kappa.
	\label{def of F_Wtilde}
\end{equation}
The higher derivatives of $\widetilde{W}$ satisfy
\begin{equation}
	\left[\partial_s+\frac{3\gamma_1+\gamma_2-1}{2}+\beta_\tau J(\partial_1\ovl{W}+\gamma_1\partial_1W)\right]\partial^\gamma\widetilde{W}+\mathcal{V}_W\cdot\nabla\partial^\gamma\widetilde{W}=F_{\widetilde{W}}^{(\gamma)},
	\label{eqn of derivatives of tilde W}
\end{equation}
where
\begin{equation}
	\begin{aligned}
		F_{\widetilde{W}}^{(\gamma)}=&\partial^\gamma F_{\widetilde{W}}-\sum_{0\le\beta<\gamma}\binom{\gamma}{\beta}\left[\partial^{\gamma-\beta}G_W\partial_1\partial^\beta \widetilde{W}+\partial^{\gamma-\beta}h_W\partial_2\partial^\beta \widetilde{W}+\beta_\tau\partial^{\gamma-\beta}(J\partial_1\ovl{W})\partial^\beta \widetilde{W}\right]\\
		&-\gamma_2\beta_\tau\partial_2(JW)\partial_1^{\gamma_1+1}\partial_2^{\gamma_2-1}\widetilde{W}-\gamma_1\beta_\tau\partial_1JW\partial^\gamma\widetilde{W}-\beta_\tau\mathbbm{1}_{|\gamma|\ge2}\sum_{\substack{0\le|\beta|\le|\gamma|-2\\\beta\le\gamma}}\binom{\gamma}{\beta}\partial^{\gamma-\beta}(JW)\partial_1\partial^\beta\widetilde{W}.
	\end{aligned}
	\label{def of F_tildeW^gamma}
\end{equation}

\subsection{Modulation ODEs}
In this section, we provide the ordinary differential equations that the modulation variables satisfy. Using modulation variables, we can match the values of $\partial^\gamma W$ and $\partial^\gamma\ovl W$ ($|\gamma|\le2$) at the origin. Specifically, we impose
\begin{equation}
\label{match W and ovl W}
	W(0,s)=\ovl W(0)=0,\ \ 
	\nabla W(0,s)=\nabla\ovl W(0)=(-1,0)^T,\ \ 
	\nabla^2 W(0)=\nabla^2\ovl W(0)=\begin{pmatrix}
		0&0\\
		0&0
	\end{pmatrix}.
\end{equation}
The values of modulation variables are determined based on these conditions. 

To derive the evolution equations of the modulation variables, we evaluate the equation of $\partial^\gamma W$, as well as the auxiliary quantities, at the origin. For convenience, throughout the paper, we adopt the convention that superscript 0 denotes the values of the corresponding quantities at $y=0$:
$$Z^0:=Z(0,s),\ \ G_W^0:=G_W(0,s),\ \ \partial_1F_W^0:=\partial_1F_W(0,s)\ldots.$$
In particular, there holds
\begin{proposition}
\label{evaluation at origin}
At $y=0$, we have that
$$
\partial_{\tilde{u}_{1}} u_{1}^{0}=1, \quad \partial_{\tilde{u}_{2}} u_{1}^{0}=0, \quad \partial_{\tilde{u}_{2}}^{2} u_{1}^{0}=\psi, \quad\partial_{\tilde{u}_{2}} u_{2}^{0}=1,\quad \partial_{\tilde{u}_{2}} \lambda^{0}=\psi;
$$
$$
\partial_{\tilde{u}_{1}}\left(u_{1}^{2}\right)^{0}=0, \quad \partial_{\tilde{u}_{1}}^{2}\left(u_{1}^{2}\right)^{0}=2,\quad \partial_{\tilde{u}_{2}}\left(u_{1}^{2}\right)^{0}=0, \quad \partial_{\tilde{u}_{1} \tilde{u}_{2}}^{2}\left(u_{1}^{2}\right)^{0}=0, \quad \partial_{\tilde{u}_{2}}^{2}\left(u_{1}^{2}\right)^{0}=0 ;
$$
$$
\partial_{\tilde{u}_{1}}\left(|u|^{2}\right)^{0}=\partial_{\widetilde{u}_{2}}\left(|u|^{2}\right)^{0}=0,\quad \partial_{\tilde{u}_{1}}^{2}\left(|u|^{2}\right)^{0}=2 ,\quad \partial_{u_{1} \tilde{u}_{2}}^2\left(| u|^{2}\right)^{0}=0 ,\quad \partial_{\tilde{u}_{2}}^{2}\left(|u|^{2}\right)^{0}=2;
$$
$$
\left(\varphi^{-1}\right)^{0}=1,\quad \partial_{\tilde u_1}\left(\varphi^{-1}\right)^{0}=\partial_{\tilde u_1}\left(\varphi^{-1}\right)^{0}=0,\quad \partial_{\tilde u_1}^2\left(\varphi^{-1}\right)^{0}=\partial_{\tilde u_2}^2\left(\varphi^{-1}\right)^{0}=\frac{1}{2r_0^2},\quad \partial_{\tilde u_1\tilde u_2}^2\left(\varphi^{-1}\right)^{0}=0;$$
$$
\jpnlambda^0=1,\quad \partial_{\tilde u_2}\jpnlambda^0=0, \quad\partial_{\tilde u_2}^2\jpnlambda^0=\psi^2;$$
$$
J^{0}=1, \quad \partial_{\tilde{u}_{1}} J^{0}=\partial_{\tilde{u}_{2}} J^{0}=0,\quad \partial_{\tilde{u}_{1}}^{2} J^{0}=\frac{1}{2 r_{0}^{2}}, \quad \partial_{\tilde{u}_{1} \bar{u}_{2}}^{2} J^{0}=0, \quad \partial_{\tilde{u}_{2}}^{2} J^{0}=\frac{1}{2 r_{0}^{2}}+\psi^{2};
$$
$$
g_{1}^{0}=-r_{0} Q_{13}, \quad \partial_{\tilde{u}_{1}} g_{1}^{0}=0, \quad \partial_{\tilde{u}_{2}} g_{1}^{0}=Q_{12},\quad \partial_{\tilde{u}_{1}}^{2} g_{1}^{0}=-\frac{Q_{13}}{2 r_{0}} , \quad \partial_{\tilde{u}_{1} \tilde{u}_{2}}^{2} g_{1}^{0}=-\frac{Q_{23}}{2 r_{0}} ,\quad \partial_{\tilde{u}_{2}}^{2} g_{1}^{0}=\frac{Q_{13}}{2 r_{0}} ;
$$
$$
g_{2}^{0}=-r_{0} Q_{23}, \quad \partial_{\tilde{u}_{1}} g_{2}^{0}=-Q_{12}, \quad \partial_{\tilde{u}_{2}} g_{2}^{0}=0,\quad \partial_{\tilde{u}_{1}}^{2} g_{2}^{0}=\frac{Q_{23}}{2 r_{0}} , \quad \partial_{\tilde{u}_{1} \tilde{u}_{2}}^{2} g_{2}^{0}=-\frac{Q_{13}}{2 r_{0}} , \quad \partial_{\tilde{u}_{2}}^{2} g_{2}^{0}=-\frac{Q_{23}}{2 r_{0}};
$$
$$
(g_1-\lambda g_2)^{0}=-r_{0} Q_{13}, \quad \partial_{\tilde{u}_{1}}(g_1-\lambda g_2)^{0}=0, \quad \partial_{\tilde{u}_{2}}(g_1-\lambda g_2)^{0}=Q_{12}+\psi r_{0} Q_{23},
$$
$$
\partial_{\tilde{u}_{1}}^{2}(g_1-\lambda g_2)^{0}=-\frac{Q_{23}}{2 r_{0}} ,\quad \partial_{\tilde{u}_{1} \tilde{u}_{2}}^{2}(g_1-\lambda g_2)^{0}=\psi Q_{12}-\frac{Q_{23}}{2 r_{0}} ,\quad \partial_{\tilde{u}_{2}}^{2}(g_1-\lambda g_2)^{0}=\frac{Q_{13}}{2 r_{0}} ;
$$
$$
h_{11}^0=0,\quad h_{12}^0=Q_{12}.
$$
\end{proposition}

We are now ready to derive the evolution equations for modulation variables.

By setting $y=0$ in the equation of $W$, we have that
\begin{equation}
	\boxed{\partial_{\tilde t}\kappa=\frac{e^{\frac{s}{2}}}{\beta_\tau}(G_W^0+F_W^0).}
	\label{eqn of kappa}
\end{equation}
Evaluating the equation of $\partial_1W$ at $y=0$, we deduce that
\begin{equation}
	\boxed{\partial_{\tilde t}\tau=\frac{1}{\beta_\tau}\left(\partial_1F_W^0+\partial_1G_W^0\right). }
	\label{eqn of tau}
\end{equation}
In a similar fashion, for $\partial_2W$, we obtain
\begin{equation}
\label{relation btwn D2FW and D2GW}
\partial_2F_W^0+\partial_2G_W^0=0. 
\end{equation}
Note that by Proposition \ref{evaluation at origin}, we have
$$
\begin{aligned}
    \partial_2G_W^0&=\beta_\tau e^{\frac{s}{2}}\left(\beta_2\partial_2Z^0+\partial_2(g_1-\lambda g_2)^0\right)\\
    &=\beta_2\beta_{\tau} e^{\frac{s}{2}} \partial_{2} Z^{0}+ \beta_{\tau}\left(Q_{12}+\psi r_{0} Q_{23}\right). 
\end{aligned}
$$
Consequently, the evolution equation of $Q_{12}$ is
\begin{equation}
	\boxed{Q_{12}=-\frac{1}{\beta_{\tau}} \partial_{2} F_{W}^{0}-\beta_2 e^{\frac{s}{2}} \partial_{2} Z^{0}-\psi r_{0} Q_{23}}.
	\label{eqn of Q12}
\end{equation}

Setting $y=0$ in the equation of $\partial_1\nabla W$, we get
$$
	\left\{\begin{aligned}
		&G_W^0\partial_{111}W^0+h_W^0\partial_{112}W^0=F_W^{(2,0),0}\\
		&G_W^0\partial_{112}W^0+h_W^0\partial_{122}W^0=F_W^{(1,1),0}. 
	\end{aligned}\right.
$$
It then follows that
\begin{equation}
	\begin{pmatrix}
		G_W^0\\ h_W^0
	\end{pmatrix}=\left(\nabla^2\partial_1W^0\right)^{-1}
	\begin{pmatrix}
		F_W^{(2,0),0}\\ F_W^{(1,1),0}
	\end{pmatrix},
	\label{relation of G_W h_W and F_W}
\end{equation}
provided that the matrix $\nabla^2\partial_1W^0$ is invertible. Observing that 
$$
\begin{aligned}
    G_{W}^{0}&=\beta_{\tau} e^{\frac{s}{2}}\left(\beta_2 Z^{0}+\kappa-r_{0} Q_{13}\right),\\
    h_W^0&=\beta_{\tau} e^{-\frac{s}{2}}\left(2\beta_1V_{2}^{0}+g_{2}^{0}\right)=\beta_{\tau} e^{-\frac{s}{2}}\left(2\beta_1A^{0}-r_{0} Q_{23}\right),
\end{aligned}$$
we get the equation for $Q_{13}$ and $Q_{23}$:
\begin{equation}
\label{eqn of Q13,Q23}
\boxed{\begin{aligned}
    Q_{13}-\frac{2\beta_3\kappa}{r_0}&=-\frac{e^{-\frac{s}{2}}}{ \beta_{\tau}r_{0}} G_W^0+\frac{\beta_2}{ r_{0}} (Z^{0}+\kappa), \\
    Q_{23}&=-\frac{e^{\frac{s}{2}}}{ \beta_{\tau}r_{0}} h_W^0+\frac{2\beta_1}{r_{0}} A^{0}.
\end{aligned}}
\end{equation}
Finally, for $\partial_{22}W$, we evaluate its equation at $y=0$ as below:
$$G_W^0\partial_{122}W^0+h_W^0\partial_{222}W^0=F_W^{(0,2),0}=\partial_{22} F_{W}^{0}+\partial_{22} G_{W}^{0}.$$
Direct computation yields that
$$\begin{aligned}
    \partial_{22} G_{W}^{0}&=\beta_{\tau} e^{\frac{s}{2}}\left[\beta_2\left(\partial_{22}Z^{0}+\left(\frac{1}{2 r_{0}^{2}}+\psi^{2}\right) e^{-s} Z^{0}\right)+\kappa\left(\frac{1}{2 r_{0}^{2}}+\psi^{2}\right) e^{-s}+\frac{Q_{13} e^{-s}}{2 r_{0}} -e^{-s}\partial_{\tilde{t}} \psi\right].
\end{aligned}$$
Therefore, we derive the equation of $\psi$:
\begin{equation}
\label{eqn of psi}
\boxed{\begin{aligned}
    \pt \psi-2\beta_3\kappa\left(\psi^2+\frac{1}{r_0^2}\right)=&\frac{e^{\frac{s}{2}}}{\beta_{\tau}} \left(\partial_{22} F_{W}^{0}-G_{W}^{0} \partial_{122} W^{0}-h_{W}^0 \partial_{222} W^{0}\right)+\beta_2e^{s} \partial_{22} Z^{0}\\
    &+\left(\frac{1}{2 r_{0}^{2}}+\psi^{2}\right)\beta_2(Z^0+\kappa)+\frac{1}{r_{0}}\left(Q_{13}-\frac{2\beta_3\kappa}{r_0}\right).
\end{aligned}}
\end{equation}

\section{Initial Setup and Main Theorem}
\label{section: Main Theorem and Reduction to a Bootstrap Argument}
In this section, we specify the initial data and formulate the main theorem in detail. The proof of the main theorem builds on a bootstrap argument, whose resolution is carried out across the subsequent sections.

We introduce two parameters $\tau_0 \in (0,1)$ and $M > 1$. The parameter $M$ is chosen sufficiently large to dominate all constants arising from the initial data and the subsequent analysis, while the parameter $\tau_0$ is then taken sufficiently small so that $\tau_0^{-1}$ dominates all quantities depending solely on $M$. This choice yields the following hierarchy:
\begin{equation}
\label{hierarchy of constants}
1 \ll \ln M \ll M \ll e^M\ll e^{e^M} \ll \tau_0^{-1}.
\end{equation}
Based on these, we further define two constants:
\begin{equation}
\label{def of l and L}
l = (\ln M)^{-5} \ll 1, \quad L = \tau_0^{-\frac{1}{10}} \gg 1.
\end{equation}

\subsection{Initial Data}
We now prescribe the initial data that lead to shock formation. These conditions are given by assigning initial values for the modulation variables, which determine the coordinate transformations, and for the unknowns $W, Z, A$ in the self-similar coordinates. This equivalence holds since the transformation from the original coordinates to the self-similar coordinates is entirely governed by the modulation variables.

\paragraph{Initial Data of the Modulation Variables} Let us first prescribe initial data for the modulation variables. At $\tilde t=0$, we set
\begin{equation}
    \label{init values of mod variables}
    O(0)=I_3=\begin{pmatrix}
        1&0&0\\
        0&1&0\\
        0&0&1
    \end{pmatrix},\quad\tau(0)=\tau_0,\quad\kappa(0)=\kappa_0,\quad\psi(0)=\psi_0.
\end{equation}
Assume that the constants $\tau_0,\kappa_0,\psi_0$ satisfy
\begin{equation}
    \label{init conditions for mod variables}
    0<\tau_0\ll1,\quad|\kappa_0-\sigma_\infty|\le\tau_0^{\frac{1}{5}},\quad|\psi_0|\le10,
\end{equation}
where $\sigma=\sigma_\infty$ is the background steady state with
\begin{equation}
    \sigma_\infty\ge\beta_3^{-1}=\frac{1+\alpha}{\alpha}=\frac{\gamma+1}{\gamma-1}.
\end{equation}
The smallness of $\tau_0$ will be used to suppress $M$, $\kappa_0$, $\sigma_\infty$, $\psi_0$, $\alpha$, $\beta_1$, $\beta_2$, $\beta_3$, and any other constants appearing in the proof. We also assume that
\begin{equation}
    \label{absorb r0 sigma with M}
    r_0+r_0^{-1}\le100,\quad\sigma_\infty\le(\ln M)^{\frac{1}{10}}.
\end{equation}

\paragraph{Initial Data of the unknowns} Suppose the initial data of \eqref{rescaled Euler equations} are given by
$$\sigma(0,p)=\sigma_0(p),\quad v(0,p)=v_0(p),\quad\text{for }p\in\mathbb{S}^2.$$
Using the modulation variables, whose initial data have been assigned, we now specify the initial self-similar coordinates. First, the choice of the rotation matrix $O(0)$ ensures that the co-moving coordinates coincide with the original coordinates at $\tilde t=0$. Specifically, we have
$$u_i=\frac{2r_0x_i}{r_0-x_3}=\frac{2r_0\mathring x_i}{r_0-\mathring x_3}.$$
We then introduce the frame $E_i= \partial_{u_i}/|\partial_{u_i}|=\varphi^{-1}\partial_{u_i}$, and expand the vector field $v$ as
\begin{equation}
    v_0(p)=V_{0i}(p)E_i(p).\nonumber
\end{equation}
Next, we define the shock-adapted coordinates by
$$\left\{\begin{aligned}
    &\tilde u_1=u_1-\frac{1}{2}\psi_0u_2^2,\\
    &\tilde u_2=u_2.
\end{aligned}\right.$$
In addition, we introduce the following auxiliary quantities  
$$\lambda_0=\psi_0 u_2,\quad J_0=\varphi^{-1}\langle \lambda_0\rangle,\quad N=\langle \lambda_0\rangle^{-1}(E_1-\lambda_0E_2),\quad T=\langle \lambda_0\rangle^{-1}(\lambda_0E_1+E_2).$$
Now we define the Riemann variables at the initial time as
\begin{equation}
\label{def of w0,z0,a0}
    w_0=v_0\cdot N+\sigma,\quad z_0=v_0\cdot N-\sigma,\quad a_0=v_0\cdot T.
\end{equation}

Subsequently, we set the initial value of the self-similar time as $s_0=-\ln\tau_0$, and define the self-similar coordinates at the initial time by
$y_1=\tilde u_1e^{\frac{3}{2}s_0}$ and $y_2=\tilde u_2e^{\frac{s_0}{2}}$. Now we prescribe the initial data for the self-similar Riemann variables $W,Z$ and $A$ by
\begin{equation}
    \label{def of W0,Z0,A0}
    w_0=e^{-\frac{s_0}{2}}W_0+\kappa_0,\quad z_0=Z_0,\quad a_0=A_0.
\end{equation}

We then state the initial conditions of $W,Z$ and $A$, at $\tilde t=0$, or equivalently, at $s=s_0$. Regarding the spatial support of these unknowns, we assume at the initial time that
\begin{equation}
\label{init of support}
\operatorname{supp}_{y}\left(\nabla_{y} W_0, \nabla_{y} Z_0, \nabla_{y} A_0\right) \subset\left\{\left|y_{1}\right| \le \tau_{0}^{\frac{1}{2}} e^{\frac{3}{2} s},\left|y_{2}\right| \le \tau_{0}^{\frac{1}{6}} e^{\frac{s}{2}}\right\}.
\end{equation}

To state the initial condition for $W_0$, we define the quantity $\mu(\gamma)$ for multi-index $\gamma$ with $|\gamma|\le 3$ as follows:
\begin{equation}
    \label{def of mu}
    \mu(\gamma):=\begin{cases}
        \frac{3\gamma_1+\gamma_2-1}{6},&|\gamma|\le1\\
        \frac{1}{3},&2\le|\gamma|\le3,\gamma_1>0\\
        \frac{1}{6},&2\le|\gamma|\le3,\gamma_1=0
    \end{cases}.
\end{equation}
Recall that $\eta(y):=1+y_1^2+y_2^6$. We choose $W_0(y)$ such that
\begin{equation}
\label{init of W}
|\eta^{\mu(\gamma)}\partial^\gamma W_0(y)| \le \begin{cases}
    \frac{201}{200},&\gamma=(0,0),\\
    16,&|\gamma|=1,\\
    \ln M,&|\gamma|=2,3.
\end{cases}
\end{equation}
For $\widetilde W_0=W_0-\ovl W$, we assume that
\begin{equation}
\label{init of W tilde}
\begin{array}{c}
|\partial^{\gamma} \widetilde{W}_0^0| \le \tau_{0}^{\frac{1}{2}} ,\quad \forall|\gamma|=3; \\
|\partial^{\gamma} \widetilde{W}_0| \mathbbm{1}_{|y| \le l} \le \tau_{0}^{\frac{1}{4}}, \quad \forall|\gamma| \le 4; \\
\eta^{\mu(\gamma)}|\partial^{\gamma} \widetilde{W}_0| \mathbbm{1}_{|y| \le L} \le \tau_{0}^{\frac{1}{\gamma_1+2\gamma_2+15}}, \quad \forall|\gamma| \le 1.
\end{array}
\end{equation}
Here, the constants $l$ and $L$ are defined in \eqref{def of l and L}. The notation $\partial^{\gamma} \widetilde{W}_0^0$ represents $\partial^{\gamma} \widetilde{W}_0(0)$, i.e. the value of $\partial^{\gamma} \widetilde{W}_0$ at $y=0$. We then assign the initial data of $Z$ and $A$ such that
\begin{equation}
\label{init of Z,A}
e^{\left(3\mu(\gamma)+\frac{1}{2}\right)s_0}\left(|\partial^\gamma(Z_0+\sigma_\infty)|+|\partial^\gamma A_0|\right)\le\tau_0.
\end{equation}
Finally, for a fixed integer $k\ge34$, we assume at $\tilde t=0$ (equivalently, at $s=s_0$) that
\begin{equation}
\label{init of Hk norm}
    \tau_0^{\frac{1}{2}}\left\|D^{k} W_0\right\|_{L^{2}}+\left\|D^{k} Z_0\right\|_{L^{2}}+\left\|D^{k} A_0\right\|_{L^{2}} \le M\tau_0^{\frac{1}{2}}.
\end{equation}

\begin{remark}
The imposed initial conditions are admissible in the sense that there exist such data that simultaneously satisfy all the requirements. In particular, one can construct such initial data explicitly by choosing $O(0)=I_3$, $\kappa_0=\sigma_\infty$, $\psi_0=0$, $\tau_0>0$, and
\begin{equation}
    W_0(y) = \chi(\tau_0 y_1,\tau_0^{\frac{1}{3}}y_2)\ovl{W}(y),\quad Z_0(y)+\sigma_\infty = A_0(y) = 0.\nonumber
\end{equation}
Here, the bump function $\chi\in C_c^\infty([-1,1]^2)$ satisfies $\chi=1$ on $[-\frac{1}{2},\frac{1}{2}]^2$, $|\nabla\chi|\le10$, and $|\partial^\gamma\chi|\le100$ for any $|\gamma|=2$. One can verify that
$$|W_0(y)|\le\eta^{\frac{1}{6}}(y),\quad|\partial_1W_0(y)|\le15\eta^{-\frac{1}{3}}(y),\quad|\partial_2W_0(y)|\le15,$$
and
$$|\partial^\gamma W_0(y)|\lesssim_\gamma\eta^{\frac{1}{6}-\frac{\gamma_1}{2}-\frac{\gamma_2}{6}},\quad \text{for any}\,\,\gamma\ge0,$$
provided that the parameter $\tau_0$ is sufficiently small. For $\widetilde W_0$, it holds that
$$|\widetilde W_0(y)|\mathbbm{1}_{|y|\le L}=0.$$
The $\dot H^k$ norm can be estimated using the pointwise bound, yielding that
$$\|D^kW_0\|_{L^2}\lesssim\|\eta^{\frac{1-k}{3}}\|_{L^2}\lesssim_k 1\le M.$$
Next, we verify the conditions for the modulation variables. Using \eqref{eqn of kappa}\eqref{eqn of tau}\eqref{eqn of Q12}\eqref{eqn of Q13,Q23}\eqref{eqn of psi}, we obtain
$$
     \pt\kappa(0)=-\frac{2\beta_3\sigma_\infty\tau_0^2}{12r_0^2-\tau_0^3},\quad \pt\tau(0)=Q_{12}(0)=Q_{23}(0)=0,\quad Q_{13}(0)=\left(1-\frac{\tau_0^3}{12r_0^2}\right)^{-1}\frac{2\beta_3\sigma_\infty}{r_0},$$
$$\pt\psi(0)-\frac{2\beta_3\sigma_\infty}{r_0^2}=\left(1+\frac{\tau_0}{12r_0^3}\right)\left(1-\frac{\tau_0^3}{12r_0^2}\right)^{-1}\frac{\beta_3\sigma_\infty\tau_0^2}{3r_0^2}.$$
By choosing $\tau_0$ to be sufficiently small, the initial conditions for the modulation variables can be fulfilled. 
\end{remark}

\subsection{Statement of the Main Theorem}
Now we state the detailed version of the main theorem. 
\begin{theorem}
    \label{thm: main theorem}
    Consider the system \eqref{rescaled Euler equations} on $\mathbb{S}^2=\{\mathring x\in\mathbb{R}^3:|\mathring{x}|=r_0\}$, where the radius satisfies $r_0+r_0^{-1}\le100$. For any background steady state obeying $\sigma_\infty\ge\beta_3^{-1}=\frac{1+\alpha}{\alpha}=\frac{\gamma+1}{\gamma-1}$, suppose the initial data $(\sigma_0,v_0)$ satisfies the following conditions:
    \begin{itemize}
        \item[(1)] The modulation variables satisfy the initial conditions \eqref{init values of mod variables}\eqref{init conditions for mod variables}, and $\tau_0=\tau_0(\alpha,\sigma_\infty)\ll1$ is sufficiently small.
        \item[(2)] The functions $W_0,Z_0$ and $A_0$, defined by \eqref{def of w0,z0,a0}\eqref{def of W0,Z0,A0} using $(\sigma_0,v_0)$, are subject to the conditions \eqref{init of support}\eqref{init of W}\eqref{init of W tilde} \eqref{init of Z,A}\eqref{init of Hk norm}.
    \end{itemize}
    Then, the following conclusions hold:
    \begin{itemize}
        \item[(1)]\emph{Finite time $\dot C^1$ norm blow-up:} The solution $(\sigma,v)$ blows up at $\tilde t=\widetilde T_*$: 
        \begin{equation}
            \lim_{\tilde t\rightarrow\widetilde T_*}(\|\nabla\sigma\|_{L^\infty}+\|\nabla v\|_{L^\infty})=\infty.\nonumber
        \end{equation}
        Moreover, there exists a time-dependent position $\xi(\tilde t)\in\mathbb{S}^2$ such that the solution becomes unbounded along its trajectory:
        \begin{equation}
            |\nabla\sigma(\tilde t,\xi(\tilde t))|+|\nabla v(\tilde t,\xi(\tilde t))|\gtrsim \frac{1}{\widetilde T_*-\tilde t}.\nonumber
        \end{equation}
        Meanwhile, $\dot C^1$ norm of the solution remains uniformly bounded slightly away from $\xi(\tilde t)$. More precisely, for any $\delta\in(0,\pi r_0/1000)$, it holds that
        \begin{equation}
            \|\nabla\sigma\|_{L^\infty(B_\delta(\xi(\tilde t))^c)}+\|\nabla v\|_{L^\infty(B_\delta(\xi(\tilde t))^c)}\lesssim_{\sigma_\infty,r_0}\delta^{-2}+M.\nonumber
        \end{equation}
        Here, $B_\delta(\xi(\tilde t))=\{p\in\mathbb{S}^2:d(p,\xi(\tilde t))<\delta\}$, and $d$ is the distance induced by the Riemann metric on $\mathbb{S}^2$.
        
        \item[(2)]\emph{Blow-up time:} The blow-up time $\widetilde T_*$ is characterized by $\tau(\widetilde T_*)=\widetilde T_*$, and it satisfies the estimate
        \begin{equation}
            \widetilde T_*=\tau_0+O_{\sigma_\infty,r_0}(\tau_0^2).\nonumber
        \end{equation}
        
        \item[(3)]\emph{Blow-up location:}  Let $(\mathring u_1,\mathring u_2)$ denote the stereographic projection coordinates relative to the original Cartesian frame $(\mathring x_1,\mathring x_2,\mathring x_3)$. Along the trajectory, we have the estimates
        \begin{equation}
            \mathring u_1(\xi(\tilde t))=2\beta_3\kappa_0\tilde t+\mathcal{O}(M\tau_0^{\frac76}),\quad\mathring u_2(\xi(\tilde t))=\mathcal{O}(M\tau_0^{\frac76}).
        \end{equation}
        
        \item[(4)]\emph{Blow-up direction and $C^1$ norm estimates:} The normal vector field $N$ fully characterizes the direction of shock formation. In particular, there holds:
        \begin{equation}
            \|\nabla_N\sigma\|_{L^\infty}+\|\nabla_N\sigma\|_{L^\infty}\approx \frac{1}{\widetilde T_*-\tilde t},\quad \|\nabla_T\sigma\|_{L^\infty}+\|\nabla_T\sigma\|_{L^\infty}\lesssim1.\nonumber
        \end{equation}

        \item[(5)]\emph{No appearance of vacuum:}
        The density remains strictly positive:
        \begin{equation}
            \inf_{(\tilde t,p)\in[0,\widetilde T_*)\times\mathbb{S}^2}\sigma(\tilde t,p)\ge \frac{1}{2}\sigma_\infty.\nonumber
        \end{equation}

        \item[(6)]\emph{Uniform-in-time $\dot C^{1/3}$ regularity:} The Riemann variable $w=v\cdot N+\sigma$ satisfies the following H\"{o}lder estimate:
        \begin{equation}
            \|w\|_{L_{\tilde t}^\infty([0,\widetilde T_*);\dot C^{1/3}(\mathbb{S}^2))}\lesssim1.\nonumber
        \end{equation}
    \end{itemize}
\end{theorem}

\subsection{Framework of Bootstrap Argument}
\label{subsec: Bootstrap Assumptions and Framework of Bootstrap Argument}
Recall from Section \ref{strategy} that the shock formation is attributed to the singularity of self-similar coordinate transformation. And the problem hence reduces to proving global well-posedness for the system governing self-similar variables. To this end, we appeal to a bootstrap argument. We begin by stating the bootstrap assumptions. 

For the modulation variables, we assume that
\begin{equation}
\label{asmp of Q}
\left|Q_{12}\right| \le M \tau_{0}^{\frac{1}{3}} ,\quad\left|Q_{13}-\frac{2\beta_3\kappa}{r_0}\right| \le M\tau_0^{\frac{1}{6}}, \quad\left|Q_{23}\right| \le  M\tau_{0}^{\frac{1}{2}},
\end{equation}
\begin{equation}
\label{asmp of kappa tau}
\left|\pt \kappa\right| \le M \tau_{0}^{\frac{1}{6}} ,\quad \left| \kappa-\kappa_0\right| \le \tau_0,\quad \left|\pt \tau\right| \le M^{\frac12}e^{-s},
\end{equation}
\begin{equation}
\label{asmp of psi}
\left|\pt \psi-2\beta_3\kappa\left(\psi^2+\frac{1}{r_0^2}\right)\right| \le M\tau_0^{\frac{1}{6}}, \quad\left|\psi-\psi_0\right| \le M\tau_0.
\end{equation}
Regarding the spatial support of the unknowns, we impose the following condition:
\begin{equation}
\label{asmp of support}
\operatorname{supp}_{y}\left(\nabla_{y} W, \nabla_{y} Z, \nabla_{y} A\right) \subset\left\{\left|y_{1}\right| \le 2 \tau_{0}^{\frac{1}{2}} e^{\frac{3}{2} s},\left|y_{2}\right| \le 2 \tau_{0}^{\frac{1}{6}} e^{\frac{s}{2}}\right\} =: \mathcal{X}(s).
\end{equation}
We now proceed to state the assumptions on $W$:
\begin{equation}
\label{asmp of W}
|\eta^{\mu(\gamma)}\partial^\gamma W| \le \begin{cases}
    \frac{101}{100},&\gamma=(0,0)\\
    20,&|\gamma|=1\\
    M^{\frac{3|\gamma|+\gamma_2-5}{8}},&|\gamma|=2,3
\end{cases}.
\end{equation}
For $Z$, we assume that
\begin{equation}
\label{asmp of Z}
e^{\left(3\mu(\gamma)+\frac{1}{2}\right)s}|\partial^\gamma(Z+\sigma_\infty)|\le\begin{cases}
    M \tau_{0},&\gamma=(0,0)\\
    \tau_0^{\frac{1}{3}},&|\gamma|=1,2,\gamma_1=0\\
    M^{\frac{|\gamma|+\gamma_2}{4}},&|\gamma|=1,2,\gamma_1>0
\end{cases}.
\end{equation}
For $A$, we require that
\begin{equation}
\label{asmp of A}
e^{\left(3\mu(\gamma)+\frac{1}{2}\right)s}|\partial^\gamma A|\le\begin{cases}
    M\tau_0^{\frac{1}{|\gamma|+1}},&|\gamma|\le 1\text{ and }\gamma_1=0\\
    M,&\gamma\in\{(1,0),(0,2)\}
\end{cases}.
\end{equation}
The assumption on $\widetilde W$, is that
\begin{equation}
\label{asmp of W tilde}
\begin{array}{c}
|\partial^{\gamma} \widetilde{W}^{0}| \le \tau_{0}^{\frac{1}{3}} ,\quad \forall|\gamma|=3; \\
|\partial^{\gamma} \widetilde{W}| \mathbbm{1}_{|y| \le l} \le \tau_{0}^{\frac{1}{4}}(\ln M)^{|\gamma|}, \quad \forall|\gamma| \le 4; \\
\eta^{\mu(\gamma)}|\partial^{\gamma} \widetilde{W}| \mathbbm{1}_{|y| \le L} \le \tau_{0}^{\frac{1}{\gamma_1+2\gamma_2+16}}, \quad \forall|\gamma| \le 1.
\end{array}
\end{equation}
Finally, for a fixed integer $k\ge34$, we impose that
\begin{equation}
\label{asmp of Hk norm}
    e^{-\frac{s}{2}}\left\|D^{k} W\right\|_{L^{2}}+\left\|D^{k} Z\right\|_{L^{2}}+\left\|D^{k} A\right\|_{L^{2}} \le M^{\frac{3k}{2}} e^{-\frac{s}{2}}.
\end{equation}

Provided that all the above bootstrap assumptions \eqref{asmp of Q}--\eqref{asmp of Hk norm} can be improved, we are able to prove the main theorem.
\begin{proof}[Proof of Theorem \ref{thm: main theorem}]
    We proceed by verifying each conclusion of the theorem, leveraging the bootstrap assumptions and the definitions provided in the setup. 

    First, we verify the absence of vacuum:
    $$\sigma=S\overset{\eqref{P expressed by R}}{=}\sigma_\infty+\frac{1}{2}\Big(\underbrace{e^{-\frac{s}{2}}W}_{\eqref{est of W}}+\underbrace{(\kappa-\sigma_\infty)}_{\eqref{init conditions for mod variables}\eqref{asmp of kappa tau}}-\underbrace{(Z+\sigma_\infty)}_{\eqref{asmp of Z}}\Big)\ge\frac{1}{2}\sigma_\infty>0.$$
    
    Next, we estimate the $\dot C^1$ norm of $(v,\sigma)$. We begin by computing the covariant derivatives of the velocity $v$. From \eqref{deriv of N,T}, we have
$$\nabla_Nv=\Big(\nabla_N(v\cdot N)-\theta(N)A\Big)N+\Big(\theta(N)v\cdot N+\nabla_NA\Big)T,$$
    $$\nabla_Tv=\Big(\nabla_T(v\cdot N)-\theta(T)A\Big)N+\Big(\theta(T)v\cdot N+\nabla_TA\Big)T.$$
    If the bootstrap assumptions hold, by \eqref{def of 1-form theta}, the 1-form $\theta$ satisfies that
    $$|\theta(N)|= o_{\tau_0}(1),\quad |\theta(T)|=|\psi|+o_{\tau_0}(1).$$
    Next, we bound the terms involving $v\cdot N$ and $A$ using \eqref{P expressed by R}:
    $$|\nabla_N(v\cdot N)|=\frac{J}{2}|e^s\partial_1W+e^{\frac{3}{2}s}\partial_1Z|+o_{\tau_0}(1)\overset{\eqref{asmp of Z}\eqref{est of D1W}\eqref{est of geo. coef.}}{\le}\Big(\frac{1}{2}+o_{\tau_0}(1)\Big)e^s,$$
    $$|\nabla_T(v\cdot N)|=\frac{J}{2\jpnlambda^2}|\partial_2W+e^{\frac{s}{2}}\partial_2Z|\overset{\eqref{asmp of W}\eqref{asmp of Z}\eqref{est of geo. coef.}}{\le}10+o_{\tau_0}(1),$$
    $$|\nabla_NA|=J|e^{\frac{3}{2}s}A-\lambda\jpnlambda^{-2}e^{\frac{s}{2}}\partial_2A|\overset{\eqref{asmp of A}\eqref{est of geo. coef.}}{\le}M+o_{\tau_0}(1),$$
    $$|\nabla_TA|=\frac{J}{\jpnlambda^2}|\lambda e^{\frac{s}{2}}\partial_2A|\overset{\eqref{asmp of A}\eqref{est of geo. coef.}}{=}o_{\tau_0}(1).$$
    Combining these estimates, for sufficiently small $\tau_0$, we obtain:
    $$|\nabla_Nv|\le e^s,\quad|\nabla_N\sigma|\le e^s,\quad|\nabla_Tv|\le 11,\quad|\nabla_T\sigma|\le11.$$
    The blow-up criterion for classical solutions to the Euler equations (see \cite{Majda-1984}) then guarantees that the solution in the original coordinates remains regular as long as $s$ is finite.
    
    Moreover, we claim that if $s(\tilde t)$ blows up at $\tilde t=\widetilde T_*$, so does the solution. This is shown by establishing a lower bound for the $\dot C^1$ norm of the solution. At the point $\xi(\tilde t)$ where $\tilde u_1 = \tilde u_2 = 0$, we observe that
    $$\begin{aligned}
        |\nabla_Nv(\tilde t,\xi(\tilde t))|=\Big|-\frac{1}{2}(e^s-e^{\frac{3}{2}s}\partial_1Z^0)N^0+e^{\frac{3}{2}}s\partial_1A^0T^0\Big|\overset{\eqref{asmp of Z}\eqref{asmp of A}}{\ge}\Big( \frac{1}{2}-o_{\tau_0}(1)\Big)e^s,
    \end{aligned}$$
    $$\begin{aligned}
        |\nabla_N\sigma(\tilde t,\xi(\tilde t))|=\Big|-\frac{1}{2}(e^s+e^{\frac{3}{2}s}\partial_1Z^0)\Big|\overset{\eqref{asmp of A}}{\ge}\Big( \frac{1}{2}-o_{\tau_0}(1)\Big)e^s.
    \end{aligned}$$
The following estimates then follow:
    $$\frac{1}{3}e^s\le\|\nabla_Nv\|_{L^\infty}\le e^s,\quad\frac{1}{3}e^s\le\|\nabla_N\sigma\|_{L^\infty}\le e^s$$
    Since $\{N,T\}$ are orthonormal, there holds
    $$|\nabla v|^2=|\nabla_Nv|^2+|\nabla_Tv|^2,\quad|\nabla \sigma|^2=|\nabla_N\sigma|^2+|\nabla_T\sigma|^2,$$
which implies
    $$\frac{1}{3}e^s\le\|\nabla v\|_{L^\infty}\le 2e^s,\quad\frac{1}{3}e^s\le\|\nabla \sigma\|_{L^\infty}\le 2e^s.$$
    This proves the claim and concludes the blow-up of the $C^1$ norm and the estimates in (5).
    
    Noting by \eqref{asmp of kappa tau} that $\frac{\mathrm{d}}{\mathrm{d}\tilde{t}} (\tau(\tilde{t}) - \tilde{t}) = \partial_{\tilde{t}} \tau - 1 \leq -\frac{1}{2}$, and $\tau(\tilde{t}) - \tilde{t} = \tau_0 > 0$ at $\tilde{t} = 0$, there exists a finite time $\widetilde{T}_* > 0$ such that $\lim_{\tilde{t} \to \widetilde{T}_*} \tau(\tilde{t}) = \widetilde{T}_*$. By the definition $s=-\ln(\tau(\tilde t)-\tilde t)$, we have $\lim_{\tilde t\rightarrow\widetilde T_*}s(\tilde t)=\infty$. Thus, the solution blows up at a finite time $\widetilde T_*$. Moreover, we obtain the following estimate for the blow-up time:
    $$|\widetilde T_*-\tau_0|=|\tau(\widetilde T_*)-\tau(0)|\le\int_0^{\widetilde T_*}\underbrace{|\pt\tau(\tilde t)|}_{\eqref{asmp of kappa tau}}\dif\tilde t\overset{s=-\ln(\tau(\tilde t)-\tilde t)}{\le}\int_{s_0}^\infty \beta_\tau M^{\frac12}e^{-2s}\dif s\overset{\eqref{est of 1-beta_tau}}{\le}M\tau_0^2.$$
    At the blow-up location $\xi(\tilde t)$, the following estimates result from the above:
    $$|\nabla_N v(\tilde t, \xi(\tilde t))| + |\nabla_N \sigma(\tilde t, \xi(\tilde t))| \gtrsim e^s \gtrsim \frac{1}{\widetilde T_* - \tilde t}.$$
    We now consider the estimates slightly away from $\xi(\tilde t)$. In this case, any point $p\in B_\delta(\xi(\tilde t))^c$ satisfies either $|y_1|\gtrsim \delta e^{\frac{3}{2}s}$ or $|y_2|\gtrsim\delta e^{\frac{s}{2}}$. For such $p$, we derive
    $$|\nabla_Nv|\lesssim e^s\eta^{-\frac{1}{3}}+M\lesssim\delta^{-2}+M.$$
    The proof of $\|\nabla_N\sigma\|_{L^\infty(B_\delta(\xi(\tilde t))^c)}\lesssim\delta^{-2}+M$ follows similarly.

    For the behavior of $\xi(\tilde t)$, we first recall that it corresponds to the origin of the $(\tilde u_1,\tilde u_2)$ coordinates: $\tilde u_1(\xi(\tilde t))=\tilde u_2(\xi(\tilde t))=0$. From the coordinate transformations \eqref{def of tilde u coordinates} and \eqref{def of u coordinates}, we find $u_1(\xi(\tilde t))=u_2(\xi(\tilde t))=0$, and $x_1(\xi(\tilde t))=x_2(\xi(\tilde t))=0$, $x_3(\xi(\tilde t))=-r_0$. Since $x=O(\tilde t)\mathring x$ (see Section \ref{subsect: Rotation and Stereographic Projection}) and $\mathring u_i=\frac{2r_0\mathring x_i}{r_0-\mathring x_3}$, we derive
    $$\mathring u_1(\xi(\tilde t))=-\frac{2r_0 O_{13}(\tilde t)}{1+O_{33}(\tilde t)},\quad\mathring u_2(\xi(\tilde t))=-\frac{2r_0 O_{23}(\tilde t)}{1+O_{33}(\tilde t)}.$$
    From $\dot O^T=QO^T$ and \eqref{asmp of Q}, we have $O=I+\mathcal{O}(M\tau_0)$. Invoking $\dot O^T=QO^T$ one more time, we see that
    $$\dot O^T=\begin{pmatrix}
        0&0&\frac{2\beta_3\kappa}{r_0}\\
        0&0&0\\
        -\frac{2\beta_3\kappa}{r_0}&0&0&
    \end{pmatrix}+\mathcal{O}(M\tau_0^{\frac16}).$$
    Since $O(\tilde t)=I+\int_0^{\tilde t}\dot O(\tilde t')\dif\tilde t'$, we obtain more accurate estimates for $O_{13}$, $O_{23}$, and $O_{33}$ using \eqref{asmp of kappa tau}:
    $$O_{13}(\tilde t)=-\frac{2\beta_3\kappa_0\tilde t}{r_0}+\mathcal{O}(M\tau_0^{\frac76}),\quad O_{23}(\tilde t)=\mathcal{O}(M\tau_0^{\frac76}),\quad O_{33}(\tilde t)=1+\mathcal{O}(M\tau_0^{\frac76}).$$
    Therefore, we derive
    $$\mathring u_1(\xi(\tilde t))=-2r_0\left(-\frac{2\beta_3\kappa_0\tilde t}{r_0}+\mathcal{O}(M\tau_0^{\frac76})\right)\left(\frac12+\mathcal{O}(M\tau_0^{\frac76})\right)=2\beta_3\kappa_0\tilde t+\mathcal{O}(M\tau_0^{\frac76}),$$
    and
    $$\mathring u_1(\xi(\tilde t))=-2r_0\cdot\mathcal{O}(M\tau_0^{\frac76})\cdot\left(\frac12+\mathcal{O}(M\tau_0^{\frac76})\right)=\mathcal{O}(M\tau_0^{\frac76}).$$

    Finally, we estimate the $\dot{C}^{1/3}$ norm of the Riemann variable $w$. For any $p,q \in \operatorname{supp} \nabla w$, note that
    $$d(p,q)\approx |u(p)-u(q)|\approx|\tilde u(p)-\tilde u(q)|.$$
    Employing this relation, we derive the following (uniformly in $\tilde t$) estimate:
    $$\begin{aligned}
        \frac{|w(\tilde t,p)-w(\tilde t,q)|}{d(p,q)^{\frac{1}{3}}}&\lesssim\frac{|w(\tilde t,\tilde u(p))-w(\tilde t,\tilde u(q))|}{|\tilde u(p)-\tilde u(q)|^{\frac{1}{3}}}\\
        &=\frac{e^{-\frac{s}{2}}|W(s,y(p))-W(s,y(q))|}{[e^{-3s}(y_1(p)-y_1(q))^2+e^{-s}(y_2(p)-y_2(q))^2]^{1/6}}\\
        &\le\frac{|W(s,y_1(p),y_2(p))-W(s,y_1(q),y_2(p))|}{|y_1(p)-y_1(q)|^{\frac{1}{3}}}+e^{-\frac{s}{3}}\frac{|W(s,y_1(q),y_2(p))-W(s,y_1(q),y_2(q))|}{|y_2(p)-y_2(q)|^{\frac{1}{3}}}\\
        \overset{\eqref{asmp of W}}&{\le}\frac{\int_{y_1(q)}^{y_1(p)}(1+z^2)^{-\frac{1}{3}}dz}{|y_1(p)-y_1(q)|^{\frac{1}{3}}}+e^{-\frac{s}{3}}|y_2(p)-y_2(q)|^{\frac{2}{3}}\\
        \overset{p,q\in\operatorname{supp}\nabla w}&{\lesssim}1.
    \end{aligned}$$
    Thus, we have
    $$\|w\|_{L_{\tilde t}^\infty([0,\widetilde T_*);\dot C^{1/3}(\mathbb{S}^2))}\lesssim1.$$
    This completes the proof of all conclusions.
\end{proof}

\subsection{Outline for closing the bootstrap argument}
In the following, we sketch the process to close the aforementioned bootstrap argument.

Most of our bootstrap assumptions are $L^\infty$-type bounds for $W$, $Z$, $A$, and $\widetilde W$. Improvements of these assumptions except $\partial_1A$) are achieved by solving the unknowns along the Lagrangian trajectories of the transport equations. In this process, it is crucial to establish estimates for the transport and forcing terms. Typically, the upper bounds take the form of a linear combination of monomials structured as $\eta^{p_1}M^{p_2}(\ln M)^{p_3}\tau_0^{p_4}e^{p_5s}$. For instance, estimate \eqref{est of FW^gamma} yields
$$|\partial_{11}F_W|\lesssim(\ln M)^{7} M \tau_0^{\frac{1}{2}} e^{-s} + (\ln M)^{\frac{1}{10}}(\eta^{-\frac{1}{3}}  M^{\frac{1}{8}} e^{-s} + M^{1+\frac{3k}{2}} e^{-\frac{11}{8}s}).$$
Hence, comparing such ``monomials''—terms of the form $\eta^{p_1}M^{p_2}(\ln M)^{p_3}\tau_0^{p_4}e^{p_5s}$—over the domain $\mathcal{X}(s)$ becomes essential, which motivates the following comparison lemma.
\begin{lemma}
    \label{lem: compare monos}
    On $\mathcal{X}(s)$, the inequality
    $$\eta^{p_1}M^{p_2}(\ln M)^{p_3}\tau_0^{p_4}e^{p_5s}\lesssim \eta^{q_1}M^{q_2}(\ln M)^{q_3}\tau_0^{q_4}e^{q_5s}$$
    holds if and only if the following two conditions are verified:
    \begin{enumerate}
        \item $p_5-q_5+3(p_1-q_1)^+\le0$.
        \item $(2(p_1-q_1)^++p_5-q_5-p_4+q_4,p_2-q_2,p_3-q_3)\le_{\mathrm{lex}}(0,0,0)$
        in lexicographic order.
    \end{enumerate}
\end{lemma}
\begin{proof}
    It suffices to determine the necessary and sufficient conditions under which the following inequality $$\eta^{p_1-q_1}M^{p_2-q_2}(\ln M)^{p_3-q_3}\tau_0^{p_4-q_4}e^{(p_5-q_5)s}\lesssim1$$
    hold on $\mathcal{X}(s)$. Indeed, we have
    \begin{align*}
        &\,\,\,\,\,\,\,\,\,\,\,\,\eta^{p_1-q_1}M^{p_2-q_2}(\ln M)^{p_3-q_3}\tau_0^{p_4-q_4}e^{(p_5-q_5)s}\mathbbm{1}_{\mathcal{X}(s)}\lesssim1\\
        &\Leftrightarrow\sup_{y\in\mathcal{X}(s)}\eta^{p_1-q_1}(y)M^{p_2-q_2}(\ln M)^{p_3-q_3}\tau_0^{p_4-q_4}e^{(p_5-q_5)s}\lesssim1\\
        \overset{\text{Prop }\ref{prop: eta bound}}&{\Leftrightarrow}M^{p_2-q_2}(\ln M)^{p_3-q_3}\tau_0^{p_4-q_4+(p_1-q_1)^+}e^{[p_5-q_5+3(p_1-q_1)^+]s}\lesssim1\\
        &\Leftrightarrow p_5-q_5+3(p_1-q_1)^+\le0\text{ and }M^{p_2-q_2}(\ln M)^{p_3-q_3}\tau_0^{p_4-q_4+(p_1-q_1)^+}e^{[p_5-q_5+3(p_1-q_1)^+]s_0}\lesssim1\\
        \overset{e^{s_0}=\tau_0^{-1}}&{\Leftrightarrow} p_5-q_5+3(p_1-q_1)^+\le0\text{ and }M^{p_2-q_2}(\ln M)^{p_3-q_3}\tau_0^{p_4-q_4-p_5+q_5-2(p_1-q_1)^+}\lesssim1\\
        \overset{\eqref{hierarchy of constants}}&{\Leftrightarrow}p_5-q_5+3(p_1-q_1)^+\le0\text{ and }(2(p_1-q_1)^++p_5-q_5-p_4+q_4,p_2-q_2,p_3-q_3)\le_{\mathrm{lex}}(0,0,0).
    \end{align*}
    This concludes the proof.
\end{proof}

To bound $\partial_1 A$ requires a distinct approach, as it cannot be closed via transport equation estimates. Instead, its control is closely related to the estimate of vorticity $\omega$. In fact, the quantity $\omega/\rho$ is purely transported by the velocity field $v$, and this conservation along particle paths allows us to close the bootstrap assumption for $\partial_1 A$.

However, pure $L^\infty$-type estimates do not close the full bootstrap argument due to a loss of derivatives. To address this issue, we leverage the symmetric hyperbolic structure of the system and utilize complementary energy estimates to recover the lost derivatives.

\subsection{Immediate Corollaries of Bootstrap Assumptions}
The bootstrap assumptions lead to several immediate corollaries. These preliminary results, requiring little beyond the stated assumptions, form the groundwork for the more delicate estimates developed in subsequent sections. We begin with establishing the following proposition:
\begin{proposition}
We have that
\begin{equation}
    \label{est of psi, Dtpsi, kappa, Q13}
    |\psi|\lesssim1,\quad|\pt\psi|+|\kappa|+|Q_{13}|\lesssim(\ln M)^{\frac{1}{3}},
\end{equation}
and
\begin{equation}
\label{est of 1-beta_tau}
    |1-\beta_\tau|\lesssim M^{\frac12}e^{-s}.
\end{equation}
\end{proposition}
\begin{proof}
By \eqref{init conditions for mod variables}\eqref{asmp of psi}, we deduce that $$|\psi|\le|\psi_0|+M\tau_0\le10+ M\tau_0\lesssim1.$$ 
Employing \eqref{asmp of psi}\eqref{asmp of kappa tau}\eqref{asmp of Q}, we obtain
$$|\pt\psi|+|\kappa|+|Q_{13}|\lesssim(1+r_0^{-1}+r_0^{-2})\kappa_0\overset{\eqref{init conditions for mod variables}\eqref{absorb r0 sigma with M}}{\lesssim}(\ln M)^{\frac{1}{3}}.$$
Finally, by \eqref{asmp of kappa tau}, we get 
$$|1-\beta_\tau|=\frac{|\partial_{\tilde t}\tau|}{|1-\partial_{\tilde t}\tau|}\lesssim|\partial_{\tilde t}\tau|\lesssim M^{\frac12}e^{-s}.$$
\end{proof}

\begin{proposition}
\label{prop: eta bound}
    On $\mathcal{X}(s)$, the inequality $\eta \lesssim \tau_0 e^{3s}$ holds.
\end{proposition}
\begin{proof}
    By definition of $\eta$, we observe that $\eta = 1 + y_1^2 + y_2^6\leq 1 + C\tau_0 e^{3s} \lesssim \tau_0 e^{3s}$.
\end{proof}

Using Proposition \ref{prop: eta bound} and the assumptions that $|W|\lesssim\eta^{\frac{1}{6}}$ and $\operatorname{supp}W\subset\mathcal{X}(s)$, we have the following estimate. 
\begin{proposition}
$W$ satisfies the inequality:
\begin{equation}
    \label{est of W}
    |W|\lesssim \tau_0^{\frac{1}{6}}e^{\frac{s}{2}}.
\end{equation}
\end{proposition}
Another important property of $W$ is as follows:  
\begin{proposition}  
The following inequality holds:
\begin{equation}
\label{est of D1W}
    |\partial_1 W| \le 1 + o_{\tau_0}(1),
\end{equation}
where $o_{\tau_0}(1)$ denotes a quantity that vanishes as $\tau_0 \to 0$.  
\end{proposition}
\begin{proof}
If $|y|\le L$, we have that $$|\partial_1 W|\le|\partial_1 \ovl W|+|\partial_1 \widetilde W|\overset{\eqref{est of ovl W and D ovl W}\eqref{asmp of W tilde}}{\le}1+\tau_0^{\frac{1}{17}}.$$
Note that $\eta(y)\ge\frac{1}{2}(1+|y|^2)$. If $|y|\ge L$, we thus deduce $$|\partial_1W|\overset{\eqref{asmp of W}}{\le}20\eta^{-\frac{1}{3}}\lesssim(1+L^2)^{-\frac{1}{3}}\ll1.$$ 
\end{proof}
We also derive an estimate
for the deviation of $Z$ from $-\kappa$:
\begin{proposition}
If the bootstrap assumptions are verified, the following estimate holds:
\begin{equation}
\label{est of Z+kappa}
    |Z+\kappa|\lesssim \tau_0^{\frac15}.
\end{equation}
\end{proposition}
\begin{proof}
This follows from
$$|Z+\kappa|\lesssim|Z+\sigma_\infty|+|\kappa_0-\sigma_\infty|+|\kappa-\kappa_0|\overset{\eqref{init conditions for mod variables}\eqref{asmp of kappa tau}\eqref{asmp of Z}}{\lesssim}\tau_0^{\frac15}.$$
\end{proof}
For higher order derivatives of self-similar variables, we show the proposition below:
\begin{proposition}
Given the bootstrap assumptions, we have
\begin{equation}
\label{higher est of W}
    |D^jW|\overset{3\le j\le5}{\le} M^{2(j-1)}.
\end{equation}
For $\mathcal{R}\in\{Z,A\}$, there holds 
\begin{equation}
\label{higher est of ZA}
|D^j\mathcal{R}|\overset{3\le j\le 5}{\le} e^{-\frac{9}{10}s},\quad
     |D^j\partial_1\mathcal{R}|\overset{2\le j\le 4}{\lesssim}M^{1+\frac{3k}{2}}  e^{-\frac{11}{8}s}.
\end{equation}
\end{proposition}
\begin{proof}
Since $DW$ is compactly supported, by Gagliardo-Nirenberg interpolation inequality, for $2\le j\le k-1$, one can verify that, 
$$\|D^jW\|_{L^\infty}\lesssim\|D^kW\|_{L^2}^{\frac{j-1}{k-2}}\|DW\|_{L^\infty}^{\frac{k-j-1}{k-2}}\overset{\eqref{asmp of W}\eqref{asmp of Hk norm}}{\lesssim}M^{\frac{3k}{2}\cdot\frac{j-1}{k-2}}\overset{k\ge34}{\le}M^{2(j-1)}.$$
Here, we utilize the largeness of $M$ in the last inequality.

For $\mathcal{R}\in\{Z,A\}$ and $1\le j\le 4<k-2$, we obtain
$$
\begin{aligned}
    \|D^j\partial_1\mathcal{R}\|_{L^\infty}&\lesssim\|D^{k-1}\partial_1\mathcal{R}\|_{L^2}^{\frac{j}{k-2}}\|\partial_1\mathcal{R}\|_{L^\infty}^{\frac{k-j-2}{k-2}}\lesssim\left(M^{\frac{3}{2}k}e^{-\frac{s}{2}}\right)^{\frac{j}{k-2}}\left(Me^{-\frac{3}{2}s}\right)^{\frac{k-j-2}{k-2}}\\
    &\lesssim M^{1+\frac{3k}{2}}e^{-\left(\frac{3}{2}-\frac{j}{k-2}\right)s}\lesssim M^{1+\frac{3k}{2}}  e^{-\frac{11}{8}s},
\end{aligned}
$$
which implies $|D^2\mathcal{R}|\lesssim Me^{-s}$. Then, by Gagliardo-Nirenberg inequality, we deduce that
$$
\begin{aligned}
    \|D^j\mathcal{R}\|_{L^\infty}\overset{3\le j\le5}&{\lesssim}\|D^k\mathcal{R}\|_{L^2}^{\frac{j-2}{k-3}}\|D^2\mathcal{R}\|_{L^\infty}^{\frac{k-j-1}{k-3}}\lesssim\left(M^{\frac{3}{2}k}e^{-\frac{s}{2}}\right)^{\frac{j-2}{k-3}}\left(Me^{-s}\right)^{\frac{k-j-1}{k-3}}\\
    &\lesssim M^{3}e^{-\left(1-\frac{j-2}{2(k-3)}\right)s}\lesssim M^{3}e^{-\left(1-\frac{3}{62}\right)s}\le e^{-\frac{9}{10}s}.
\end{aligned}
$$
\end{proof}
Next, we estimate the auxiliary quantities associated with coordinate transformations. These quantities are determined by the modulation variables. To state the resulting inequalities more conveniently, we first define several operations on sets of multi-indices.   
\begin{definition}
For any two sets of multi-indices $E,F\subset\mathbb{Z}_{\geq0}^2$, we define 
$$\operatorname{cone}E = \left\{\sum c_i v_i : c_i \in \mathbb{Z}_{\geq0}, v_i \in E\right\},
$$
$$\operatorname{cone}^+E = \left\{\sum c_i v_i: c_i \in \mathbb{Z}_{\geq0}, \sum c_i > 0, v_i \in E\right\},
$$
$$E+F=\{v_1+v_2:v_1\in E,v_2\in F\}.$$
\end{definition}

We now state the estimates for auxiliary quantities.
\begin{proposition}[Estimates for auxiliary quantities]
\label{prop: est for geo. coeff.}
Suppose that the bootstrap assumptions are valid, the following inequalities hold on the region \(\mathcal{X}(s)\):
\begin{equation}
\label{est of geo. coef.}
\begin{aligned}
    &\big|\partial_{\tilde{u}}^{\gamma}(\dot\psi \tilde{u}_{2}^{2})\big|
    \lesssim \mathbbm{1}_{\gamma\le2e_2} \tau_{0}^{\frac{2-\gamma_2}{6}} \ln M, \quad
    \big|\partial_{\tilde{u}}^{\gamma} u_{1}\big|
    \lesssim \mathbbm{1}_{e_1}(\gamma) + \tau_{0}^{\frac{2-\gamma_2}{6}} \mathbbm{1}_{\gamma\le2e_2}, \\
    &\big|\partial_{\tilde{u}}^{\gamma} \lambda\big|
    \lesssim \tau_{0}^{\frac{1-\gamma_2}{6}} \mathbbm{1}_{\gamma\le e_2}, \quad
    \big|\partial_{\tilde{u}}^{\gamma} (\lambda^2)\big|
    \lesssim \tau_{0}^{\frac{2-\gamma_2}{6}} \mathbbm{1}_{\gamma\le 2e_2}, \\
    &\big|\partial_{\tilde{u}}^{\gamma}(u_{1}^{2})\big|
    \lesssim \mathbbm{1}_{\{2e_1, (1,2), 4e_2\}}(\gamma) + \tau_{0}^{\frac{1}{6}}, \quad
    \big|\partial_{\tilde{u}}^{\gamma}(u_{2}^{2})\big|
    \lesssim \mathbbm{1}_{2e_2}(\gamma) + \tau_{0}^{\frac{1}{6}}, \\
    &\big|\partial_{\tilde{u}}^{\gamma}(u_1 u_{2})\big|
    \lesssim \mathbbm{1}_{\{(1,1),3e_2\}}(\gamma) + \tau_{0}^{\frac{1}{6}}, \quad
    \big|\partial_{\tilde{u}}^{\gamma}(|u|^2)\big|
    \lesssim \mathbbm{1}_{\{2e_1,(1,2),2e_2,4e_2\}}(\gamma) + \tau_{0}^{\frac{1}{6}}, \\
    &\big|\partial_{\tilde{u}}^{\gamma}(\varphi^{-1}-1)\big|
    \lesssim \mathbbm{1}_{\{2e_1,(1,2),2e_2,4e_2\}}(\gamma) + \tau_{0}^{\frac{1}{6}},\\
    &\big|\partial_{\tilde{u}}^{\gamma}(\varphi-1)\big|
    \lesssim_\gamma \mathbbm{1}_{\operatorname{cone}^+\{2e_{1},(1,2),2e_{2},4e_2\}}(\gamma) + \tau_{0}^{\frac{1}{6}}, \\
    &\big|\partial_{\tilde{u}}^{\gamma}(\jpnlambda-1)\big|+\big|\partial_{\tilde{u}}^{\gamma}(\jpnlambda^{-1}-1)\big|
    \lesssim_\gamma \mathbbm{1}_{\gamma=0}\tau_0^{\frac{1}{3}}+\mathbbm{1}_{\operatorname{cone}^+(e_2)}(\gamma)\tau_0^{\frac{1-(-1)^{\gamma_2}}{12}}, \\
    &\big|\partial_{\tilde{u}}^{\gamma}(J-1)\big|
    +\big|\partial_{\tilde{u}}^{\gamma}(\beta_\tau J-1)\big|\lesssim_\gamma \mathbbm{1}_{(\{0,2e_{1},(1,2)\}+\operatorname{cone}(e_2))\setminus\{0\}}(\gamma)\tau_0^{\frac{1-(-1)^{\gamma_2}}{12}} + \tau_{0}^{\frac{1}{6}}.
\end{aligned}
\end{equation}
Additionally, we have:
\begin{equation}
\label{est of hij and gi}
\begin{aligned}
    &\big|\partial_{\tilde{u}}^{\gamma}h_{11}\big|
    \lesssim (\mathbbm{1}_{e_{1}}(\gamma) + \tau_{0}^{\frac{1}{6}}) (\ln M)^{\frac{1}{3}}, \\
    &\big|\partial_{\tilde{u}}^{\gamma}h_{12}\big|
    \lesssim (\mathbbm{1}_{e_{1}}(\gamma) + \tau_{0}^{\frac{1}{6}}) (\ln M)^{\frac{1}{3}} + M^2\tau_0^{\frac{1}{2}}, \\
    &\big|\partial_{\tilde{u}}^{\gamma}(g_1+r_0Q_{13})\big|
    \lesssim \mathbbm{1}_{\{2e_{1}, 2e_{2}\}}(\gamma) (\ln M)^{\frac{1}{3}} + M^{2} \tau_{0}^{\frac{1}{6}}, \\
    &\big|\partial_{\tilde{u}}^{\gamma}g_2\big|
    \lesssim \mathbbm{1}_{e_1+e_2}(\gamma) \ln M + M^{2} \tau_{0}^{\frac{1}{6}}, \\
    &\big|\partial_{\tilde{u}}^{\gamma}(g_1 - \lambda g_2 + r_0Q_{13})\big|
    \lesssim \mathbbm{1}_{\{2e_{1}, 2e_{2}\}}(\gamma) (\ln M)^{\frac{1}{3}} + M^{2} \tau_{0}^{\frac{1}{6}}.
\end{aligned}
\end{equation}
Here, both $\mathbbm{1}_{A}(\gamma)$ and $\mathbbm{1}_{\gamma \in A}$ denote the indicator function, which takes the value $1$ when $\gamma \in A$, and $0$ otherwise.
\begin{proof}
    These bounds follow directly from the definitions of $u_i$, $\tilde u_i$, $\varphi$, $\lambda$, $\jpnlambda$, $J$, $\beta_\tau$, $h_{ij}$, and $g_i$, along with the bootstrap assumptions \eqref{asmp of kappa tau}, \eqref{asmp of psi}, and \eqref{asmp of Q} for the modulation variables.
\end{proof}

\end{proposition}

Next, we consider the lower order estimates on the physical variable $P=(V_1,V_2,S)^T$. We introduce the shifted variables $\widetilde{P} = (V_1, V_2, S - \sigma_\infty)^T = (\widetilde{P}_1, \widetilde{P}_2, \widetilde{P}_3)^T = P - (0, 0, \sigma_\infty)^T$. Based on the support condition \eqref{asmp of support} and the finite speed of propagation property of the system, we have that $\operatorname{supp}_y \widetilde{P} \subset \mathcal{X}(s)$. This compactly supported property allows us to use the Poincar\'e inequality to estimate their $L^2$-norm. The pointwise estimates for the shifted variables $\widetilde{P}$ are given below.
\begin{proposition}
    \label{prop: est of P}
    We have the following estimates for $\widetilde{P}$ on $\mathcal{X}(s)$:
    \begin{equation}
    \label{est of V1}
        |\partial^\gamma V_1|\lesssim\begin{cases}
            \tau_0^{\frac{1}{6}},&\gamma=(0,0),\\
            \eta^{-\frac{1}{3}} e^{-\frac{s}{2}},&\gamma=(1,0)\\
            e^{-\frac{s}{2}},&\gamma=(0,1),\\
            \eta^{-\frac{1}{3}} M^{\frac{1}{8}} e^{-\frac{s}{2}} + M^{1+\frac{3k}{2}} \tau_0^{\frac{1}{6}} e^{-\frac{11}{8}s},&\gamma=(2,0),\\
            \eta^{-\frac{1}{3}} M^{\frac{1}{4}} e^{-\frac{s}{2}} + M^{1+\frac{3k}{2}} \tau_0^{\frac{1}{6}} e^{-\frac{11}{8}s},&\gamma=(1,1),\\
            \eta^{-\frac{1}{6}} M^{\frac{3}{8}} e^{-\frac{s}{2}},&\gamma=(0,2),\\
            M^{1+\frac{3k}{2}} e^{-\frac{11}{8}s} + \eta^{-\frac{1}{3}} M^{\frac{1}{2}} e^{-\frac{s}{2}},&\gamma=(3,0),\\
            \eta^{-\frac{1}{3}} M^{\frac{5}{8}} e^{-\frac{s}{2}} + M^{1+\frac{3k}{2}} e^{-\frac{11}{8}s},&\gamma
            =(2,1),\\
            M^{1+\frac{3k}{2}} e^{-\frac{11}{8}s} + \eta^{-\frac{1}{3}} M^{\frac{3}{4}} e^{-\frac{s}{2}},&\gamma=(1,2),\\
            \eta^{-\frac{1}{6}} M^{\frac{7}{8}} e^{-\frac{s}{2}} + e^{-\frac{9}{10}s},&\gamma=(0,3),\\
            M^{6} e^{-\frac{s}{2}},&|\gamma|=4,\\
            M^{8} e^{-\frac{s}{2}},&|\gamma|=5,
        \end{cases}
    \end{equation}
    \begin{equation}
    \label{est of V2}
        |\partial^\gamma V_2|\lesssim\begin{cases}
            \tau_0^{\frac{1}{3}},&\gamma=(0,0),\\
            \eta^{-\frac{1}{3}} \tau_0^{\frac{1}{6}} e^{-\frac{s}{2}},&\gamma=(1,0),\\
            \tau_0^{\frac{1}{6}} e^{-\frac{s}{2}},&\gamma=(0,1),\\
            \eta^{-\frac{1}{3}} M^{\frac{1}{8}} \tau_0^{\frac{1}{6}} e^{-\frac{s}{2}} + M^{1+\frac{3k}{2}} e^{-\frac{11}{8}s},&\gamma=(2,0),\\
             M^{1+\frac{3k}{2}} e^{-\frac{11}{8}s} + \eta^{-\frac{1}{3}} M^{\frac{1}{4}} \tau_0^{\frac{1}{6}} e^{-\frac{s}{2}},&\gamma=(1,1),\\
            \eta^{-\frac{1}{6}} M^{\frac{3}{8}} \tau_0^{\frac{1}{6}} e^{-\frac{s}{2}} + M e^{-s},&\gamma=(0,2),\\
            M^{1+\frac{3k}{2}} e^{-\frac{11}{8}s} + \eta^{-\frac{1}{3}} M^{\frac{1}{2}} \tau_0^{\frac{1}{6}} e^{-\frac{s}{2}},&\gamma=(3,0),\\
            M^{1+\frac{3k}{2}} e^{-\frac{11}{8}s} + \eta^{-\frac{1}{3}} M^{\frac{5}{8}} \tau_0^{\frac{1}{6}} e^{-\frac{s}{2}},&\gamma
            =(2,1),\\
            M^{1+\frac{3k}{2}} e^{-\frac{11}{8}s} + \eta^{-\frac{1}{3}} M^{\frac{3}{4}} \tau_0^{\frac{1}{6}} e^{-\frac{s}{2}},&\gamma=(1,2),\\
            \eta^{-\frac{1}{6}} M^{\frac{7}{8}} \tau_0^{\frac{1}{6}} e^{-\frac{s}{2}} + e^{-\frac{9}{10}s},&\gamma=(0,3),\\
            M^{6}\tau_0^{\frac{1}{6}} e^{-\frac{s}{2}},&|\gamma|=4,\\
            M^{8}\tau_0^{\frac{1}{6}} e^{-\frac{s}{2}},&|\gamma|=5,
        \end{cases}
    \end{equation}
    \begin{equation}
    \label{est of VN and S}
        |\partial^\gamma(v\cdot N)|+|\partial^\gamma (S-\sigma_\infty))|\lesssim\begin{cases}
            \tau_0^{\frac{1}{6}},&\gamma=(0,0),\\
            \eta^{-\frac{1}{3}} e^{-\frac{s}{2}},&\gamma=(1,0),\\
            e^{-\frac{s}{2}},&\gamma=(0,1),\\
            \eta^{-\frac{1}{3}} M^{\frac{1}{8}} e^{-\frac{s}{2}},&\gamma=(2,0),\\
            \eta^{-\frac{1}{3}} M^{\frac{1}{4}} e^{-\frac{s}{2}},&\gamma=(1,1),\\
            \eta^{-\frac{1}{6}} M^{\frac{3}{8}} e^{-\frac{s}{2}},&\gamma=(0,2),\\
            M^{1+\frac{3k}{2}} e^{-\frac{11}{8}s} + \eta^{-\frac{1}{3}} M^{\frac{1}{2}} e^{-\frac{s}{2}},&\gamma=(3,0),\\
            \eta^{-\frac{1}{3}} M^{\frac{5}{8}} e^{-\frac{s}{2}} + M^{1+\frac{3k}{2}} e^{-\frac{11}{8}s},&\gamma
            =(2,1),\\
            M^{1+\frac{3k}{2}} e^{-\frac{11}{8}s} + \eta^{-\frac{1}{3}} M^{\frac{3}{4}} e^{-\frac{s}{2}},&\gamma=(1,2),\\
            \eta^{-\frac{1}{6}} M^{\frac{7}{8}} e^{-\frac{s}{2}} + e^{-\frac{9}{10}s},&\gamma=(0,3),\\
            M^{6} e^{-\frac{s}{2}},&|\gamma|=4,\\
            M^{8} e^{-\frac{s}{2}},&|\gamma|=5.
        \end{cases}
    \end{equation}
\end{proposition}
\begin{proof}
    These inequalities are derived via using \eqref{P expressed by R} and \eqref{asmp of W}\eqref{asmp of Z}\eqref{asmp of A}\eqref{higher est of W}\eqref{higher est of ZA}\eqref{est of geo. coef.}. Indeed, from \eqref{P expressed by R}, one can observe that
    $$|\partial^\gamma V_1|\overset{\eqref{P expressed by R}}{\lesssim}e^{-\frac{s}{2}}\underbrace{|\partial^\gamma(\jpnlambda^{-1}W)|}_{\eqref{est of geo. coef.}\eqref{asmp of W}\eqref{higher est of W}}+\underbrace{|\partial^\gamma(\jpnlambda^{-1}(Z+\kappa)|}_{\eqref{est of geo. coef.}\eqref{asmp of Z}\eqref{est of Z+kappa}\eqref{higher est of ZA}}+\underbrace{|\partial^\gamma(\lambda\jpnlambda^{-1}A)|}_{\eqref{est of geo. coef.}\eqref{asmp of A}\eqref{higher est of ZA}},$$
    $$|\partial^\gamma V_2|\overset{\eqref{P expressed by R}}{\lesssim}e^{-\frac{s}{2}}\underbrace{|\partial^\gamma(\lambda\jpnlambda^{-1}W)|}_{\eqref{est of geo. coef.}\eqref{asmp of W}\eqref{higher est of W}}+\underbrace{|\partial^\gamma(\lambda\jpnlambda^{-1}(Z+\kappa)|}_{\eqref{est of geo. coef.}\eqref{asmp of Z}\eqref{est of Z+kappa}\eqref{higher est of ZA}}+\underbrace{|\partial^\gamma(\jpnlambda^{-1}A)|}_{\eqref{est of geo. coef.}\eqref{asmp of A}\eqref{higher est of ZA}},$$
    and
    $$|\partial^\gamma (v\cdot N)|+|\partial^\gamma (S-\sigma_\infty)|\overset{\eqref{P expressed by R}}{\lesssim}e^{-\frac{s}{2}}\underbrace{|\partial^\gamma W|}_{\eqref{est of geo. coef.}\eqref{asmp of W}\eqref{higher est of W}}+\underbrace{\mathbbm{1}_{\gamma=(0,0)}|\kappa-\sigma_\infty|}_{\eqref{init conditions for mod variables}\eqref{asmp of kappa tau}}+\underbrace{|\partial^\gamma (Z+\sigma_\infty)|}_{\eqref{asmp of Z}\eqref{higher est of ZA}}.$$
Using the prior estimates indicated under the brackets yields the desired bounds for $V_1$, $V_2$, $v\cdot N$, and $S-\sigma_\infty$.
\end{proof}
Then, we have the following corollary, which is useful in the top-order energy estimates.   
\begin{corollary}
    \label{cor: lower order est of P tilde}
    For $\mu=1,2,3$, there hold the following estimates:
$$\|\widetilde{P}_\mu\|_{L^\infty} \lesssim \tau_0^{\frac{1}{6}}, \quad \|P_\mu\|_{L^\infty} \leq \ln M, \quad \|D^2 P_\mu\|_{L^6_y} \lesssim M e^{-\frac{s}{2}}.$$
\end{corollary}
\begin{proof}
    First, it follows from Proposition \ref{prop: est of P} that $|\widetilde P|\lesssim\tau_0^{\frac{1}{6}}$, and $|D^2P|\lesssim Me^{-\frac{s}{2}}\eta^{-\frac{1}{6}}$. These imply $|P|\le |\widetilde P|+\sigma_\infty\le\ln M$, and $\|D^2 P\|_{L^6}\lesssim Me^{-\frac{s}{2}}\|\eta^{-\frac{1}{6}}\|_{L^6}\lesssim Me^{-\frac{s}{2}}$. 
\end{proof}

In the energy estimates, we also require the following sharper bounds on $|\partial_jP|:=\sqrt{\sum_\mu |\partial_jP_\mu|^2}$. 
\begin{proposition}
    \label{prop: est of DP}
Provided that the bootstrap assumptions are valid, it holds on $\mathcal{X}(s)$ that
    \begin{equation}
    \label{est of DV,DS}
    \max(|\partial_1 V_1|,|\partial_1 S|)\le(\frac{1}{2}    +o_{\tau_0}(1))e^{-\frac{s}{2}},\quad\max(|\partial_2 V_1|,|\partial_2 S|)\le(10+o_{\tau_0}(1))e^{-\frac{s}{2}},\quad |D V_2|\le o_{\tau_0}(1)e^{-\frac{s}{2}}.
    \end{equation}
    In addition, we have the following upper bound estimates for the derivatives of $P$:
    \begin{equation}
        \label{est of DjP}
        |\partial_1 P|\le (\frac{\sqrt{2}}{2}+o_{\tau_0}(1))e^{-\frac{s}{2}},\quad |\partial_2 P|\le (10\sqrt{2}+o_{\tau_0}(1))e^{-\frac{s}{2}}.
    \end{equation} 
\end{proposition}
\begin{proof}
The upper bounds in \eqref{est of DV,DS} can be derived using \eqref{P expressed by R} and the bootstrap assumptions. Then, \eqref{est of DjP} follows from the definition of the vector norm.
\end{proof}

\section{Estimates for Transport Equations}
\label{section: Estimates for Transport Equations}
In this section, we estimate each part of the transport equations satisfied by $\partial^\gamma\mathcal{R}$ with $\mathcal{R}\in\{\widetilde W,W,Z,A\}$. A key observation we make here is that both the transport terms and the forcing terms can be expressed as polynomials of $\widetilde W,W,Z,A$, and certain auxiliary quantities.

\subsection{Estimates for Transport Terms}
In this section, we estimate the transport terms in the equation of $W,Z$ and $A$. We begin with a proposition controlling $e^{-\frac s2}G_{\mathcal{R}}$ on $\mathcal{X}(s)$. 
\begin{proposition}
\label{prop: est of gR,GR, zeroth order}
It holds on $\mathcal{X}(s)$ that
\begin{equation}
    e^{-\frac{s}{2}}G_\mathcal{R}=O(M^2\tau_0^{\frac{1}{6}})-\left\{
    \begin{aligned}
        &0,\quad &\mathcal{R}=W,\\
        &4\beta_3\kappa_0,\quad &\mathcal{R}=Z,\\
        &2\beta_3\kappa_0,\quad &\mathcal{R}=A.
    \end{aligned}
    \right.
\end{equation}
The same conclusion also applies to for $e^{-\frac{s}{2}}g_\mathcal{R}$ with $\mathcal{R}\in\{W,Z,A\}$. 
\end{proposition}
\begin{proof}
From \eqref{asmp of kappa tau} and Lemma \ref{prop: est for geo. coeff.}, we derive that 
$$G=g_1-\lambda g_2-\frac{1}{2}\dot\psi\tilde u_2^2=-2\beta_3\kappa_0+O(M^2\tau_0^{\frac{1}{6}}).$$
Moreover, combining \eqref{asmp of kappa tau}\eqref{est of Z+kappa}\eqref{est of 1-beta_tau} and Proposition \ref{prop: est for geo. coeff.}, we obtain
$$Z=\kappa_0+O(\tau_0^{\frac15}),\quad\kappa=\kappa_0+O(\tau_0),\quad\beta_\tau=1+O(M\tau_0),\quad J=1+O(\tau_0^{\frac{1}{6}}).$$
Then, for any $c_1,c_2\in\mathbb{R}$, it follows from direct calculation that
$$
\begin{aligned}
    \beta_\tau[J(c_1Z+c_2\kappa)+G]&=(1+O(M\tau_0))\left\{(1+O(\tau_0^{\frac{1}{6}}))\left[c_1(-\kappa_0+O(\tau_0^{\frac{1}{6}}))+c_2(\kappa_0+O(\tau_0))\right]-2\beta_3\kappa_0+O(M^2\tau_0^{\frac{1}{6}})\right\}\\
    &=(-c_1+c_2-2\beta_3)\kappa_0+O(M^2\tau_0^{\frac{1}{6}}).
\end{aligned}$$
Therefore, employing the expression \eqref{def of gR,GR,G} of $G_{\mathcal{R}}$, we obtain the desired result. The behavior of $g_{\mathcal{R}}$ then follows directly from \eqref{est of W}. 
\end{proof}
For the derivatives of $G_{\mathcal{R}}$, we have the following pointwise bounds. 
\begin{proposition}
For $\mathcal{R}\in\{W,Z,A\}$, we have the following estimates for higher order derivatives of the transport terms on $\mathcal{X}(s)$: 
\begin{equation}
\label{est of GR}
|\partial^\gamma G_\mathcal{R}|\lesssim
\begin{cases}
    e^{-\frac{|\gamma|-1}{2}s}\tau_0^{\frac{2-|\gamma|}{6}}\ln M,&1\le|\gamma|\le2\text{ and }\gamma_1=0,\\
    M^{\frac{|\gamma|+\gamma_2}{4}}e^{-s},&1\le|\gamma|\le2\text{ and }\gamma_1>0,\\
    e^{-\frac{2}{5}s},&3\le|\gamma|\le5\text{ and }\gamma_1=0,\\
    M^{1+\frac{3k}{2}}e^{-\frac{7}{8}s},&3\le|\gamma|\le5\text{ and }\gamma_1>0,
\end{cases}
\end{equation}
and
\begin{equation}
    \label{est of hR}
    |\partial^\gamma h_{\mathcal{R}}|\lesssim
    \begin{cases}
        (\ln M)^{\frac{1}{10}} \tau_0^{\frac{1}{6}} e^{-\frac{s}{2}},&\gamma=(0,0),\\
        \eta^{-\frac{1}{3}} \tau_0^{\frac{1}{6}} e^{-s},&\gamma=(1,0),\\
        (\ln M)^{\frac{1}{10}} e^{-s}, &\gamma=(0,1),\\
        M e^{-\frac{3}{2}s} + \eta^{-\frac{1}{6}} M^{\frac{3}{8}} \tau_0^{\frac{1}{6}} e^{-s},&|\gamma|=2.
    \end{cases}
\end{equation}
\end{proposition}
\begin{proof}
Via using \eqref{def of gR,GR,G}, we have
$$
\begin{aligned}
    e^{-\frac{s}{2}}|\partial^\gamma G_{\mathcal{R}}|&\lesssim \underbrace{|\partial^\gamma(JZ)|}_{\eqref{asmp of Z}\eqref{higher est of ZA}\eqref{est of geo. coef.}}+\underbrace{|\kappa||\partial^\gamma J|+|\partial^\gamma (g_1-\lambda g_2)|+|\dot\psi||\partial^\gamma (\tilde u_2^2)|}_{\eqref{asmp of kappa tau}\eqref{est of geo. coef.}}.
\end{aligned}
$$
Using the prior estimates indicated under the brackets yields the desired bounds in \eqref{est of GR}. For $h_{\mathcal{R}}$, invoking \eqref{def of hR}, one can see that
$$e^{\frac{s}{2}}|\partial^\gamma h_{\mathcal{R}}|\lesssim\underbrace{|\partial^\gamma(\varphi^{-1}V_2)|+|\partial^\gamma(\varphi^{-2}\lambda\jpnlambda^{-1}S)|}_{\eqref{P expressed by R}\eqref{asmp of W}\eqref{asmp of Z}\eqref{asmp of A}\eqref{higher est of W}\eqref{higher est of ZA}\eqref{est of geo. coef.}}+\underbrace{|\partial^\gamma g_2|}_{\eqref{est of geo. coef.}}.$$
Substituting the estimates indicated under the bracket leads to the estimates for $|\partial^\gamma h_{\mathcal{R}}|$.
\end{proof}

To investigate the transport equations, we introduce the Lagrangian trajectory $\Phi_{\mathcal{R}}(s;s_1,y_0)$ originating from $y_0 \in \mathbb{R}^2$ at initial time $s_1 \ge s_0$, defined by:
\begin{equation}
\label{def of trajectories}
\left\{
\begin{aligned}
    &\partial_s \Phi_{\mathcal{R}}(s;s_1,y_0) = \mathcal{V}_{\mathcal{R}}(s,\Phi_{\mathcal{R}}(s;s_1,y_0)), \\
    &\Phi_{\mathcal{R}}(s_1;s_1,y_0) = y_0.
\end{aligned}
\right.
\end{equation}

We will see that $\Phi_W$ evolves outward exponentially, as long as the starting point is slightly away from the origin. The trajectories $\Phi_Z$ and $\Phi_A$ exhibit different behavior: although they may pass through the origin for certain initial data, they generally move away from the origin at an exponential rate. Although the trajectory intersects the origin, it remains in its vicinity only for a short time, after which the distance increases exponentially. 

The following proposition provides a lower bound on the growth rate of $|\Phi_W|$. 

\begin{proposition}
\label{prop: outgoing property of W trajectories}
For any $s\ge s_1\ge s_0$ and $|y_0|\ge l=(\ln M)^{-5}$, we have $|\Phi_W(s;s_1,y_0)|\ge |y_0|e^{\frac{s-s_1}{5}}$. 
\end{proposition}
\begin{proof}
The proof follows the approach in \cite{buckmaster2023shock}. For the reader's convenience, we provide a brief sketch here. From \eqref{def of gR,GR,G}\eqref{est of hR}\eqref{est of GW0, hW0}, we deduce $|G_W|\lesssim o_{\tau_0}(1)+o_{\tau_0}(1)|y|$ and $|h_W|\lesssim o_{\tau_0}(1)$, which hence imply $y\cdot\mathcal{V}_W(s,y)\ge\left(\frac{1}{2}-o_{\tau_0}(1)\right)|y|^2+y_1^2-\left(1+o_{\tau_0}(1)\right)|y_1W|$. For $l\le|y|\le L$, noting that $\ovl W(0,y_2)=0$, we obtain 
$$\begin{aligned}
|y_1W(s,y)|&\le|y_1|\underbrace{|\widetilde W(s,y)|}_{\eqref{asmp of W tilde}}+|y_1|\underbrace{|\ovl W(y_1,y_2)-\ovl W(0,y_2)|}_{\eqref{evaluation of Wbar and derivs at 0}}\\
    &\le\tau_0^{\frac{1}{5}}|y_1|(1+|y_1|^{\frac{1}{3}}+|y_2|)+y_1^2\le o_{\tau_0}(1)|y|^2+y_1^2.
\end{aligned}$$
Here, we use the condition $|y|\ge l$ in the last inequality. If $|y|\ge L=\tau_0^{-\frac{1}{10}}$, via using \eqref{asmp of W}, we derive that
$$\begin{aligned}
    |y_1W(s,y)|&\le\frac{101}{100}|y_1|(1+|y_1|^{\frac{1}{3}}+|y_2|)\\
    &\le o_M(1)y_1^2+\frac{101}{100}|y_1y_2|\le \left(o_M(1)+\frac{101}{100}\right)y_1^2+\frac{101}{400}y_2^2.
\end{aligned}$$
In both cases, we have $y\cdot\mathcal{V}_W\ge\frac{1}{5}|y|^2$. Since $\partial_s|\Phi_W(s)|^2=(y\cdot\mathcal{V}_W)(s,\Phi_W(s))$, it follows that $|\Phi(s;s_1,y_0)|\ge |y_0|e^{\frac{s-s_1}{5}}$ for $s\ge s_1$.
\end{proof}

The above growth estimate of $|\Phi_W|$ implies the following upper bound for the integral of $\eta^{-p}$ along $\Phi_W$:

\begin{lemma}\label{lemma: estimate of eta to -p along PhiW}
For any $p\in(0,10)$, $s_1\ge s_0$ and $y_0\in\mathbb{R}^2$, there holds
\begin{equation}
    \label{integration along W trajectories}
    \int_{s_1}^{s}\eta^{-p}(\Phi_W(s';s_1,y_0))\dif s'\lesssim_p
    \begin{cases}
        \ln\ln M,&|y_0|\ge l,\\
        \tau_0^{\frac{p}{5}},&|y_0|\ge L.
    \end{cases}
\end{equation}
\end{lemma}
\begin{proof}
Note that $\eta(y)=1+y_1^2+y_2^6\ge\frac{1}{2}(1+|y|^2)$. Setting $\tilde s=2\ln|y_0|+\frac{2}{5}(s'-s_1)$, a direct computation yields
$$\begin{aligned}
    \int_{s_1}^{s}\eta^{-p}(\Phi_W(s';s_1,y_0))\dif s'&\le2^p\int_{s_1}^{s}\left(1+|\Phi_W(s';s_1,y_0)|^2\right)^{-p}\dif s'\\
    \overset{\text{Prop }\ref{prop: outgoing property of W trajectories}}&{\le}2^p\int_{s_1}^{s}\left(1+|y_0|^2e^{\frac{2}{5}(s'-s_1)}\right)^{-p}\dif s'\\
    &\le5\cdot2^{p-1}\int_{2\ln|y_0|}^{\infty}(1+e^{\tilde s})^{-p}\dif\tilde s\\
    &\le\begin{cases}
        5\cdot2^{p-1}\left(\int_{-10\ln\ln M}^{0}+\int_{0}^{\infty}\right)(1+e^{\tilde s})^{-p}\dif\tilde s,&|y_0|\ge l\\
        5\cdot2^{p-1}\int_{-\frac{1}{5}\ln\tau_0}^{\infty}(1+e^{\tilde s})^{-p}\dif\tilde s,&|y_0|\ge L\end{cases}\\
        \overset{M\gg1}&{\le}\begin{cases}
        30\cdot 2^p\ln\ln M,&|y_0|\ge l,\\
        \frac{5\cdot 2^{p-1}}{p}\tau_0^{\frac{p}{5}},&|y_0|\ge L.
        \end{cases}
\end{aligned}$$
\end{proof}    

Next, we control the integral of $\eta^{-p}$ along $\Phi_Z$ and $\Phi_A$. 

\begin{lemma}\label{lemma: integral of eta^-p along PhiZ and A}
Let $\mathcal{R}\in\{Z,A\}$. Suppose the bootstrap assumptions hold on $[s_0,s]$ and $\beta_3\kappa_0\ge1$. Then, for any $p\in(0,10)$, it holds that
\begin{equation}
    \label{integration along ZA trajectories}
    \int_{s_0}^{s}\eta^{-p}(\Phi_{\mathcal{R}}(s';s_0,y_0))\dif s'\lesssim_p 1.
\end{equation}
\end{lemma}
\begin{proof}
According to Proposition \ref{prop: est of gR,GR, zeroth order}, in the region $\Omega(s):=\{y\in\mathbb{R}^2:y_1\le \frac{2}{3}\beta_3\kappa_0e^{\frac{s}{2}}\}$, the trajectories move leftward with speed bounded by $\partial_s\Phi_\mathcal{R}^1(s;s_0,y_0)\le\beta_3\kappa_0e^{\frac{s}{2}}-(2\beta_3\kappa_0-O(M^2\tau_0^{\frac{1}{6}}))e^{\frac{s}{2}}\le-\frac{1}{3}\beta_3\kappa_0 e^{\frac{s}{2}}$. By adapting the argument from \cite[Lemma 8.3]{buckmaster2023shock}, we obtain $|\Phi_\mathcal{R}(s;s_0,y_0)| \ge \frac{2}{3}\beta_3\kappa_0\min(|e^{\frac{s}{2}} - e^{\frac{s_*}{2}}|, e^{\frac{s}{2}})$ for some $s_* = s_*(s_0, y_0)$. Therefore, we deduce that
$$
\begin{aligned}
    \int_{s_0}^{s}\eta^{-p}(\Phi_{\mathcal{R}}(s'))\dif s'&\le\int_{s_0}^{s}(1+(\Phi^1_\mathcal{R}(s'))^2)^{-p}\dif s'\\
    &\lesssim \int_{s_0}^{\infty}\left(1+\frac{2}{3}\beta_3\kappa_0\min(|e^{\frac{s'}{2}} - e^{\frac{s_*}{2}}|, e^{\frac{s'}{2}})\right)^{-2p}\dif s'\\
    &\lesssim \int_{s_0}^{\infty}\max\left[\left(1+\frac{2}{3}\beta_3\kappa_0|e^{\frac{s'}{2}} - e^{\frac{s_*}{2}}|\right)^{-2p}, \left(1+\frac{2}{3}\beta_3\kappa_0e^{\frac{s'}{2}}\right)^{-2p}\right]\dif s'\\
    &\lesssim \underbrace{\int_{s_0}^{\infty}\left(1+\frac{2}{3}\beta_3\kappa_0|e^{\frac{s'}{2}} - e^{\frac{s_*}{2}}|\right)^{-2p}\dif s'}_{r=e^{\frac{s'}{2}}}+\int_{s_0}^{\infty}\left(1+\frac{2}{3}\beta_3\kappa_0e^{\frac{s}{2}}\right)^{-2p}\dif s'\\
    &\lesssim\int_{\tau_0^{-\frac{1}{2}}}^\infty\left(1+\frac{2}{3}\beta_3\kappa_0|r - e^{\frac{s_*}{2}}|\right)^{-2p}r^{-1}\dif r+\tau_0^p(\beta_3\kappa_0)^{-2p}\\
    \overset{\text{Young}}&{\lesssim}\int_{\tau_0^{-\frac{1}{2}}}^\infty r^{-(1+2p)}\dif r+p \int_{-\infty}^\infty\left(1+\frac{2}{3}\beta_3\kappa_0|r - e^{\frac{s_*}{2}}|\right)^{-(1+2p)}\dif r+\tau_0^p(\beta_3\kappa_0)^{-2p}\\
    &\lesssim\tau_0^p\left(\frac{1}{p}+\frac{1}{(\beta_3\kappa_0)^{2p}}\right)+\frac{1}{\beta_3\kappa_0}\lesssim1.
\end{aligned}$$
\end{proof}

\subsection{Estimates for Forcing Terms}
In this section, we establish estimates for the forcing terms $F_W$, $F_Z$, $F_A$, and their derivatives on $\mathcal{X}(s)$. We follow a unified approach: we apply Leibniz's rule to expand derivatives of products, control the terms individually by applying the previously established estimates, and select the maximal terms.

We first give upper bound estimates for the derivatives of $F_W$, $F_Z$, and $F_A$. 
\begin{proposition}
\label{prop: est of FR}
On the domain $\mathcal{X}(s)$, the forcing terms $F_W$, $F_Z$, and $F_A$ obey the following estimates:
\begin{equation}
\label{est of FW,FZ}
    |\partial^\gamma F_W|+e^{\frac{s}{2}}|\partial^\gamma F_Z|\lesssim\left\{\begin{aligned}
   &(\log M)^{\frac{1}{10}} (\tau_0^{\frac{1}{5}} e^{-\frac{s}{2}} + \eta^{\frac{1}{6}}e^{-s}),\quad&\gamma=(0,0),\\
   &(\ln M)^{\frac{1}{10}}(M^{1+\frac{3k}{2}}  e^{-\frac{11}{8}s} + \eta^{-\frac{1}{3}}  e^{-s}),\quad&\gamma=(1,0),\\
   &(\ln M)^{\frac{1}{10}} M e^{-s},\quad&\gamma=(0,1),\\
   &(\ln M)^{\frac{1}{10}} e^{-\frac{9}{10}s},\quad&2\le|\gamma|\le4,\ \gamma_1=0,\\
   &(\ln M)^{7} M \tau_0^{\frac{1}{2}} e^{-s} + (\ln M)^{\frac{1}{10}}(\eta^{-\frac{1}{3}}  M^{\frac{1+\gamma_2}{8}} e^{-s} + M^{1+\frac{3k}{2}} e^{-\frac{11}{8}s}),\quad&|\gamma|=2,\ \gamma_1>0,\\
   &(\ln M)^{\frac{1}{10}}(M^{1+\frac{3k}{2}} e^{-\frac{11}{8}s}+\eta^{-\frac{1}{3}} M^{\frac{4+\gamma_2}{8}} e^{-s}),\quad&|\gamma|=3,\ \gamma_1>0,\\
   &(\ln M)^{\frac{1}{10}} M^{6} e^{-s},\quad&|\gamma|=4,\ \gamma_1>0.
\end{aligned}\right.
\end{equation}
and
\begin{equation}
\label{est of FA}
|\partial^\gamma F_A|\lesssim\left\{\begin{aligned}
   &(\ln M)^{\frac{1}{10}} e^{-s} ,\quad&\gamma=(0,0),\\
   &\eta^{-\frac{1}{6}} (\ln M)^{\frac{1}{10}} M^{\frac{3}{8}} e^{-s},\quad&\gamma=(0,1),\\
   & (\ln M)^{\frac{1}{10}} (\eta^{-\frac{1}{6}} M^{\frac{7}{8}} e^{-s} + e^{-\frac{7}{5}s}),\quad&\gamma=(0,2).
\end{aligned}\right.
\end{equation}
\end{proposition}
\begin{proof}
    Starting from the definitions of the forcing terms (cf. \eqref{def of FR} and \eqref{def of hRR}), we obtain that
    $$\begin{aligned}
        |\partial^\gamma F_W|+e^{\frac{s}{2}}|\partial^\gamma F_Z|\lesssim&e^{-\frac{s}{2}}\underbrace{\left(|\partial^\gamma(\jpnlambda^{-1}|V|^2(u_1-\lambda u_2))|+|\partial^\gamma(J^{-1}(V_1-\lambda V_2)u_j(g_j-2\beta_1\varphi^{-1}V_j))|+|\partial^\gamma(u_jV_jS)|\right)}_{\eqref{est of geo. coef.}\eqref{est of hij and gi}\eqref{est of V1}\eqref{est of V2}\eqref{est of VN and S}}\\
        &+e^{\frac{s}{2}}\underbrace{|\partial^\gamma(\jpnlambda^{-2}h_WA\partial_2\lambda)|}_{\eqref{est of geo. coef.}\eqref{est of hR}\eqref{asmp of A}\eqref{higher est of ZA}}+\underbrace{|\partial^\gamma(\jpnlambda^{-3}\varphi^{-1}v\cdot N\partial_2\lambda)|}_{\eqref{est of geo. coef.}\eqref{est of VN and S}}+e^{-\frac{s}{2}}\underbrace{\left(|\partial^\gamma(h_{11}v\cdot N)|+|\partial^\gamma(h_{12}A)|\right)}_{\eqref{est of hij and gi}\eqref{est of VN and S}\eqref{asmp of A}\eqref{higher est of ZA}}\\
        &+\underbrace{|\partial^\gamma(\varphi^{-1}\jpnlambda^{-1}S\partial_2A)|}_{\eqref{est of geo. coef.}\eqref{est of VN and S}\eqref{asmp of A}\eqref{higher est of ZA}},
    \end{aligned}$$
    and
    $$\begin{aligned}
        |\partial^\gamma F_A|\lesssim&e^{-s}\underbrace{\left(|\partial^\gamma(\jpnlambda^{-1}|V|^2(\lambda u_1+u_2))|+|\partial^\gamma(J^{-1}(\lambda V_1+V_2)u_j(g_j-2\beta_1\varphi^{-1}V_j))|\right)}_{\eqref{est of geo. coef.}\eqref{est of hij and gi}\eqref{est of V1}\eqref{est of V2}}\\
        &+e^{-s}\underbrace{|\partial^\gamma(h_{12}v\cdot N)|}_{\eqref{est of hij and gi}\eqref{est of VN and S}}+e^{-s}\underbrace{|\partial^\gamma(\varphi^{-1}\jpnlambda^{-1}S\partial_2W)|}_{\eqref{est of geo. coef.}\eqref{est of VN and S}\eqref{asmp of W}\eqref{higher est of W}}+e^{-\frac{s}{2}}\underbrace{|\partial^\gamma(\varphi^{-1}\jpnlambda^{-1}S\partial_2Z)|}_{\eqref{est of geo. coef.}\eqref{est of VN and S}\eqref{asmp of Z}\eqref{higher est of ZA}}.
    \end{aligned}$$
    Next, we can expand each term above using Leibniz's rule, and bound it by plugging in the corresponding previous estimates indicated under the braces. As an example, we have
    \begin{align*}
        |\partial^\gamma(\varphi^{-1}\jpnlambda^{-1}S\partial_2Z)|&\lesssim\sum_{\beta_1+\beta_2+\beta_3+\beta_4=\gamma}\underbrace{|\partial^{\beta_1}(\varphi^{-1})|}_{\eqref{est of geo. coef.}}\underbrace{|\partial^{\beta_2}(\jpnlambda^{-1})|}_{\eqref{est of geo. coef.}}\underbrace{|\partial^{\beta_3}S|}_{\eqref{est of VN and S}}\underbrace{|\partial^{\beta_4}\partial_2Z|}_{\eqref{asmp of Z}\eqref{higher est of ZA}}.
    \end{align*}
    This process yields expressions that are ``polynomials'' in $\eta$, $M$, $\ln M$, $\tau_0$, and $e^s$, with each ``monomial'' representing an individual bound. Finally, we compare these ``monomials'' on $\mathcal{X}(s)$ using Lemma \ref{lem: compare monos} and discard all non-dominant terms, thereby arriving at the desired bounds in \eqref{est of FW,FZ} and \eqref{est of FA}.
\end{proof}

Next, we estimate the full forcing terms $F_W^{(\gamma)}$, $F_Z^{(\gamma)}$, and $F_A^{(\gamma)}$.
\begin{proposition}
\label{prop: est of FR^gamma}
We have the following inequalities for $F_W^{(\gamma)}$, $F_Z^{(\gamma)}$, and $F_A^{(\gamma)}$ on the domain $\mathcal{X}(s)$:
\begin{equation}
\label{est of FW^gamma}
|F_W^{(\gamma)}|\lesssim\left\{\begin{aligned}
   &(\log M)^{\frac{1}{10}} (\tau_0^{\frac{1}{5}} e^{-\frac{s}{2}} + \eta^{\frac{1}{6}}e^{-s}),\quad&\gamma=(0,0),\\
   &\eta^{-\frac{1}{3}} M^{\frac{1}{4}} e^{-s} + (\ln M)^{\frac{1}{10}} M^{1+\frac{3k}{2}} e^{-\frac{11}{8}s},\quad&\gamma=(1,0),\\
   &\eta^{-\frac{1}{3}}  \tau_0^{\frac{1}{6}}\ln M,\quad&\gamma=(0,1),\\
   &(\ln M)^{7} M \tau_0^{\frac{1}{2}} e^{-s} + \eta^{-\frac{1}{3}} M^{\frac{1}{2}} e^{-s} + \eta^{-\frac{1}{6}} M^{\frac{3}{8}} \tau_0^{\frac{1}{6}} e^{-s} + (\ln M)^{\frac{1}{10}} M^{1+\frac{3k}{2}} e^{-\frac{11}{8}s},\quad&\gamma=(2,0),\\
   &\eta^{-\frac{1}{3}} M^{\frac{1}{8}},\quad&\gamma=(1,1),\\
   &\eta^{-\frac{1}{3}} M^{\frac{1}{4}} + (\ln M)^{\frac{1}{10}} e^{-\frac{9}{10}s},\quad&\gamma=(0,2),\\
   &\eta^{-\frac{1}{3}} M^{1+\frac{3k}{2}} e^{-\frac{7}{8}s} + \eta^{-\frac{2}{3}} M^{\frac{1}{4}} + (\ln M)^{\frac{1}{10}} M^{1+\frac{3k}{2}} e^{-\frac{11}{8}s},\quad&\gamma=(3,0),\\
   &\eta^{-\frac{1}{3}} M^{\frac{1}{2}},\quad&\gamma=(2,1),\\
   &\eta^{-\frac{1}{3}} M^{\frac{5}{8}},\quad&\gamma=(1,2),\\
   &\eta^{-\frac{1}{3}} M^{\frac{3}{4}} + (\ln M)^{\frac{1}{10}} e^{-\frac{9}{10}s},\quad&\gamma=(0,3),
\end{aligned}\right.
\end{equation}

\begin{equation}
\label{est of FZ^gamma}
|F_Z^{(\gamma)}|\lesssim\left\{\begin{aligned}
   &(\log M)^{\frac{1}{10}} (\tau_0^{\frac{1}{5}} e^{-s} + \eta^{\frac{1}{6}} e^{-\frac{3}{2}s}),\quad&\gamma=(0,0),\\
   &(\ln M)^{\frac{1}{10}} (\eta^{-\frac{1}{3}} e^{-\frac{3}{2}s} +  M^{1+\frac{3k}{2}} e^{-\frac{15}{8}s}),\quad&\gamma=(1,0),\\
   &(\ln M)^{\frac{1}{10}} Me^{-\frac{3}{2}s},\quad&\gamma=(0,1),\\
   &(\ln M)^{7} M\tau_0^{\frac{1}{2}} e^{-\frac{3}{2}s} + \eta^{-\frac{1}{3}} M^{\frac{3}{8}} e^{-\frac{3}{2}s} + (\ln M)^{\frac{1}{10}} M^{1+\frac{3k}{2}} e^{-\frac{15}{8}s},\quad&\gamma=(2,0),\\
   &M^{\frac{1}{2}} e^{-\frac{3}{2}s},\quad&\gamma=(1,1),\\
   &(\ln M)^{\frac{1}{10}} e^{-\frac{7}{5}s},\quad&\gamma=(0,2),
\end{aligned}\right.
\end{equation}

\begin{equation}
\label{est of FA^gamma}
|F_A^{(\gamma)}|\lesssim\left\{\begin{aligned}
   &(\ln M)^{\frac{1}{10}} e^{-s},\quad&\gamma=(0,0),\\
   &\eta^{-\frac{1}{6}} (\ln M)^{\frac{1}{10}} M^{\frac{3}{8}} e^{-s},\quad&\gamma=(0,1),\\
   &\eta^{-\frac{1}{6}} (\ln M)^{\frac{1}{10}} M^{\frac{7}{8}} e^{-s} + M^{1+\frac{3k}{2}} e^{-\frac{11}{8}s},\quad&\gamma=(0,2).
\end{aligned}\right.
\end{equation}
\end{proposition}
\begin{proof}
Using the definitions of the forcing terms, we obtain:
$$\begin{aligned}
    |F_W^{(\gamma)}|\overset{\eqref{def of FW^gamma}}{\lesssim}&\underbrace{|\partial^\gamma F_W|}_{\eqref{est of FW,FZ}}+\underbrace{|\partial_1W[\partial^\gamma,J]W|+\mathbbm{1}_{|\gamma|\ge2}|\partial_2(JW)\partial_1^{\gamma_1+1}\partial_2^{\gamma_2-1}W|}_{\eqref{asmp of W}\eqref{higher est of W}\eqref{est of geo. coef.}}\\
		&+\mathbbm{1}_{|\gamma|\ge3}\sum_{\substack{1\le|\beta|\le|\gamma|-2\\\beta\le\gamma}}\underbrace{|\partial^{\gamma-\beta}(JW)\partial_1\partial^\beta W|}_{\eqref{asmp of W}\eqref{higher est of W}\eqref{est of geo. coef.}}+\sum_{0\le\beta<\gamma}\underbrace{\left(|\partial^{\gamma-\beta}G_W\partial_1\partial^\beta W|+|\partial^{\gamma-\beta}h_W\partial_2\partial^\beta W|\right)}_{\eqref{asmp of W}\eqref{higher est of W}\eqref{est of GR}\eqref{est of hR}},
\end{aligned}$$
$$\begin{aligned}
    |F_Z^{(\gamma)}|\overset{\eqref{def of FZ^gamma}}{\lesssim}&\underbrace{|\partial^\gamma F_Z|}_{\eqref{est of FW,FZ}}+\mathbbm{1}_{|\gamma|\ge2}\underbrace{|\partial_2(JW)\partial_1^{\gamma_1+1}\partial_2^{\gamma_2-1}Z|}_{\eqref{asmp of W}\eqref{asmp of Z}\eqref{higher est of ZA}\eqref{est of geo. coef.}}\\
		&+\mathbbm{1}_{|\gamma|\ge2}\sum_{\substack{0\le|\beta|\le|\gamma|-2\\\beta\le\gamma}}\underbrace{|\partial^{\gamma-\beta}(JW)\partial_1\partial^\beta Z|}_{\eqref{asmp of W}\eqref{asmp of Z}\eqref{higher est of W}\eqref{higher est of ZA}\eqref{est of geo. coef.}}+\sum_{0\le\beta<\gamma}\underbrace{\left(|\partial^{\gamma-\beta}G_Z\partial_1\partial^\beta Z|+|\partial^{\gamma-\beta}h_Z\partial_2\partial^\beta Z|\right)}_{\eqref{asmp of Z}\eqref{higher est of ZA}\eqref{est of GR}\eqref{est of hR}},
\end{aligned}$$
$$\begin{aligned}
    |F_A^{(\gamma)}|\overset{\eqref{def of FA^gamma}}{\lesssim}&\underbrace{|\partial^\gamma F_A|}_{\eqref{est of FA}}+\mathbbm{1}_{|\gamma|\ge2}\underbrace{|\partial_2(JW)\partial_1^{\gamma_1+1}\partial_2^{\gamma_2-1}A|}_{\eqref{asmp of W}\eqref{asmp of A}\eqref{higher est of ZA}\eqref{est of geo. coef.}}\\
		&+\mathbbm{1}_{|\gamma|\ge2}\sum_{\substack{0\le|\beta|\le|\gamma|-2\\\beta\le\gamma}}\underbrace{|\partial^{\gamma-\beta}(JW)\partial_1\partial^\beta A|}_{\eqref{asmp of W}\eqref{asmp of A}\eqref{higher est of W}\eqref{higher est of ZA}\eqref{est of geo. coef.}}+\sum_{0\le\beta<\gamma}\underbrace{\left(|\partial^{\gamma-\beta}G_A\partial_1\partial^\beta A|+|\partial^{\gamma-\beta}h_A\partial_2\partial^\beta A|\right)}_{\eqref{asmp of A}\eqref{higher est of ZA}\eqref{est of GR}\eqref{est of hR}}.
\end{aligned}$$
We bound each term using the prior estimates noted under the braces. Then, we identify the maximal terms on $\mathcal{X}(s)$ using Lemma \ref{lem: compare monos}, producing the desired bounds.
\end{proof}

For convenience, we rewrite the bounds from Proposition \ref{prop: est of FR^gamma} as the following simplified form:

\begin{corollary}
\label{cor: simplified bounds}
The following inequalities hold on $\mathcal{X}(s)$:
\begin{equation}
\label{est of FW^gamma (simplified)}
\eta^{\mu(\gamma)}|F_W^{(\gamma)}|\lesssim\left\{\begin{aligned}
   &o_{\tau_0}(1)\eta^{-\frac{1}{12}},\quad&|\gamma|\le1,\\
   &M^{\frac{3|\gamma|+\gamma_2-6}{8}}\eta^{-\frac{2}{15}\mathbbm{1}_{\gamma_1=0}},\quad&|\gamma|=2,3,
\end{aligned}\right.
\end{equation}
\begin{equation}
\label{est of FZ^gamma (simplified)}
|F_Z^{(\gamma)}|\lesssim\begin{cases}
    e^{-s},&\gamma=(0,0),\\
    M^{\frac{|\gamma|+\gamma_2}{4}-\frac{1}{8}}e^{-\frac{3}{2}s}\eta^{-\frac{1}{8}\mathbbm{1}_{|\gamma|=1}},&|\gamma|=1,2,\ \gamma_1>0,\\
    Me^{-\frac{2}{5}s}e^{-\frac{|\gamma|}{2}s},&|\gamma|=1,2,\ \gamma_1=0,
\end{cases}
\end{equation}
\begin{equation}
\label{est of FA^gamma (simplified)}
|F_A^{(\gamma)}|\lesssim\begin{cases}
    M^{\frac{1}{2}}e^{-\frac{s}{1+|\gamma|}}e^{-\frac{|\gamma|}{2}s},&|\gamma|\le1,\ \gamma_1=0,\\
    \eta^{-\frac{1}{8}}M^{\frac{15}{16}}e^{-s},&\gamma=(0,2).
\end{cases}
\end{equation}
\end{corollary}

\subsection{Estimates Near the Origin}
Building on the global bounds obtained in the previous section, we refine estimates for the transport and forcing terms near the origin, particularly in the regions $\{|y| \le L\}$ and $\{|y| \le l\}$. We begin with improving bounds for $G_W^0$ and $h_W^0$ derived from \eqref{relation of G_W h_W and F_W}.
\begin{proposition}
The transport terms $G_W$ and $h_W$ satisfy the following estimates at the origin:
\begin{equation}
\label{est of GW0, hW0}
    |G_W^0|+|h_W^0|\lesssim M^{\frac{3}{4}}e^{-s},
\end{equation}
\begin{equation}
\label{est of D2GW0}
|\partial_2G_W^0|\lesssim(\ln M)^{\frac{1}{10}}e^{-\frac{9}{10}s}.
\end{equation}
\end{proposition}
\begin{proof}
First, by \eqref{relation of G_W h_W and F_W} and \eqref{evaluation of Wbar and derivs at 0}\eqref{asmp of W tilde}, we have that $\partial_1 \nabla^2 W^0$ is invertible, and $|G_W^0| + |h_W^0| \lesssim |F_W^{(2,0),0}| + |F_W^{(1,1),0}|$.

For $|F_W^{(2,0),0}|$, it follows directly from \eqref{est of FW^gamma} that $|F_W^{(2,0),0}| \lesssim M^{\frac{1}{2}}e^{-s}$. For $|F_W^{(1,1),0}|$, evaluating \eqref{def of FW^gamma} at $y = 0$, there holds $F_W^{(1,1),0} = \partial_{12}^2 F_W^0 + \partial_{12}^2 G_W^0$. Then, using \eqref{est of GR} and \eqref{est of FW,FZ}, we deduce $|F_W^{(1,1),0}| \lesssim M^{\frac{3}{4}}e^{-s}$. Combining these estimates, we obtain $|G_W^0| + |h_W^0| \lesssim M^{\frac{3}{4}}e^{-s}$.

Finally, from \eqref{relation btwn D2FW and D2GW}, we conclude
$$|\partial_2G_W^0|=\underbrace{|\partial_2F_W^0|}_{\eqref{est of FW,FZ}}\lesssim(\ln M)^{\frac{1}{10}}e^{-\frac{9}{10}s}.$$
\end{proof}

By integrating from the origin, we have the following upper bounds for $G_W$ near the origin.  
\begin{corollary}
   Near the origin, the transport term $G_W$ verifies the following estimates:
    \begin{equation}
    \label{est of GW on |y|<L}
        \mathbbm{1}_{|y|\le L}|G_W|\lesssim\tau_0^{\frac{1}{15}}\ln M,
    \end{equation}
    \begin{equation}
    \label{est of GW on |y|<l}
        \mathbbm{1}_{|y|\le l}|\partial^\gamma G_W|\lesssim\begin{cases}
            (\log M)^{-9} e^{-\frac{s}{2}},&\gamma=(0,0),\\
            (\log M)^{-4} e^{-\frac{s}{2}},&\gamma=(0,1).
        \end{cases}
    \end{equation}
\end{corollary}
\begin{proof}
    Employing the fundamental theorem of calculus, for $|y|\le L$, we have that 
$$|G_W|\lesssim\underbrace{|G_W^0|}_{\eqref{est of GW0, hW0}}+\underbrace{(\|\partial_1G_W\|_{L^\infty}+\|\partial_2G_W\|_{L^\infty})}_{\eqref{est of GR}}\cdot L\lesssim\tau_0^{\frac{1}{15}}\ln M.$$
In a similar manner, for $|y|\le l$, we derive that
\begin{equation}
\label{est of D2GW on |y|<l}
|\partial_2G_W|\lesssim\underbrace{|\partial_2G_W^0|}_{\eqref{est of D2GW0}}+\underbrace{(\|\partial_{12}G_W\|_{L^\infty}+\|\partial_{22}G_W\|_{L^\infty})}_{\eqref{est of GR}}\cdot l\lesssim(\ln M)^{-4} e^{-\frac{s}{2}},
\end{equation}
$$|G_W|\lesssim\underbrace{|G_W^0|}_{\eqref{est of GW0, hW0}}+(\underbrace{\|\partial_1G_W\|_{L^\infty}}_{\eqref{est of GR}}+\underbrace{\|\partial_2G_W\|_{L^\infty}}_{\eqref{est of D2GW on |y|<l}})\cdot l\lesssim(\ln M)^{-9} e^{-\frac{s}{2}}.$$
\end{proof}

Next, we consider the forcing term in the equation of $\partial^\gamma\widetilde W$. 

\begin{proposition}
    We have the following bounds for $\partial^\gamma F_{\widetilde W}$ on the region $\{|y|\le L\}$:
    \begin{equation}
    \label{est of FtildeW, |y|<L}
        \mathbbm{1}_{|y|\le L}|\partial^\gamma F_{\widetilde W}|\lesssim\begin{cases}
            \eta^{-\frac{1}{3}} \ln M \tau_0^{\frac{1}{15}} + M \tau_0^{\frac{1}{6}} e^{-\frac{s}{2}} + \eta^{-\frac{1}{6}} \tau_0^{\frac{1}{3}},&\gamma=(0,0),\\
            \eta^{-\frac{1}{2}} (\ln M)^{\frac{1}{10}} \tau_0^{\frac{1}{6}} e^{-\frac{s}{2}} + \eta^{-\frac{2}{3}} \tau_0^{\frac{1}{3}} + \eta^{-\frac{1}{3}} M^{\frac{1}{4}} e^{-s} + \eta^{-\frac{5}{6}} \ln M \tau_0^{\frac{1}{15}} + (\ln M)^{\frac{1}{10}} M^{1+\frac{3k}{2}} e^{-\frac{11}{8}s},&\gamma=(1,0),\\
            \eta^{-\frac{1}{3}} \ln M \tau_0^{\frac{1}{6}} + \eta^{-\frac{1}{2}} \ln M \tau_0^{\frac{1}{15}},&\gamma=(0,1).
        \end{cases}
    \end{equation}
    In addition, on $\{|y|\le l\}$, we have the estimates:
    \begin{equation}
    \label{est of FtildeW, |y|<l}
        \mathbbm{1}_{|y|\le l}|\partial^\gamma F_{\widetilde W}|\lesssim\begin{cases}
            \tau_0^{\frac{1}{3}},&|\gamma|\le4,\,\gamma_1\text{ is odd and }\gamma_2\text{ is even},\\
            (\ln M)^{-5} \tau_0^{\frac{1}{3}},&|\gamma|\le4,\gamma_1\text{ is even or }\gamma_2\text{ is odd}.
        \end{cases}
    \end{equation}
\end{proposition}
\begin{proof}
    For $|y| \le L$, we have
    $$\begin{aligned}
        |\partial^\gamma F_{\widetilde W}|\lesssim&\underbrace{|\partial^\gamma F_W|}_{\eqref{est of FW,FZ}}+\underbrace{\Big|\partial^\gamma\left[(1-\beta_\tau J)\ovl{W}\partial_1\ovl{W}\right]\Big|}_{\eqref{est of geo. coef.}\eqref{est of Dk ovl W}}+\underbrace{|\partial^\gamma (G_W\partial_1\ovl W)|}_{\eqref{est of GW on |y|<L}\eqref{est of Dk ovl W}}+\underbrace{|\partial^\gamma (h_W\partial_2\ovl{W})|}_{\eqref{est of hR}\eqref{est of Dk ovl W}}+e^{-\frac{s}{2}}\underbrace{|\partial_{\tilde t}\kappa|}_{\eqref{asmp of kappa tau}}\mathbbm{1}_{\gamma=(0,0)}.
    \end{aligned}$$
    Substituting the prior estimates indicated under the braces yields \eqref{est of FtildeW, |y|<L}. For the case $|y|\le l$, notice that for any multi-index $\gamma\ge0$, we have
    \begin{equation}
    \label{est of Dk ovl W, |y|<l}
        |\partial^\gamma\ovl W|\lesssim\begin{cases}
            (\ln M)^{-5},&{\gamma_1\text{ is even or }\gamma_2\text{ is odd}},\\
            1,&{\gamma_1\text{ is odd and }\gamma_2\text{ is even}}.
        \end{cases}
    \end{equation}
    Subsequently, employing \eqref{def of F_Wtilde}, we derive that
    $$\begin{aligned}
        \mathbbm{1}_{|y|\le l}|\partial^\gamma F_{\widetilde W}|\overset{|\gamma|\le4}{\lesssim}&\underbrace{|\partial^\gamma F_W|}_{\eqref{est of FW,FZ}}+\underbrace{\Big|\partial^\gamma\left[(1-\beta_\tau J)\ovl{W}\partial_1\ovl{W}\right]\Big|}_{\eqref{est of geo. coef.}\eqref{est of Dk ovl W, |y|<l}}+\underbrace{|\partial^\gamma (G_W\partial_1\ovl W)|}_{\eqref{est of GR}\eqref{est of GW on |y|<l}\eqref{est of Dk ovl W, |y|<l}}+\underbrace{|\partial^\gamma h_W\partial_2\ovl{W}|}_{\eqref{est of hR}\eqref{est of Dk ovl W, |y|<l}}.
    \end{aligned}$$
     The stated bounds then follow by substituting the brace-marked estimate and applying Lemma~\ref{lem: compare monos} to compare monomials.
\end{proof}

Building upon the above proposition, we have the following estimates for the full forcing term $F_{\widetilde W}^{(\gamma)}$.
\begin{proposition}
On the region $\{|y|\le L\}$, it holds for $F_{\widetilde W}^{(\gamma)}$ that
\begin{equation}
\label{est of F_tildeW^gamma, |y|<L}
    |F_{\widetilde W}^{(\gamma)}|\lesssim
    \begin{cases}
        \eta^{-\frac{1}{3}}  \tau_0^{\frac{1}{15}} \ln M+ M \tau_0^{\frac{1}{6}} e^{-\frac{s}{2}} + \eta^{-\frac{1}{6}} \tau_0^{\frac{1}{3}},&\gamma=(0,0)\\
        \eta^{-\frac{2}{3}} \tau_0^{\frac{1}{16}} + (\ln M)^{\frac{1}{10}} M^{1+\frac{3k}{2}} e^{-\frac{11}{8}s},&\gamma=(1,0)\\
        \eta^{-\frac{1}{3}} \tau_0^{\frac{1}{17}},&\gamma=(0,1)
    \end{cases}.
\end{equation}
While on $\{|y|\le l\}$, we have that
\begin{equation}
\label{est of F_tildeW^gamma, |y|<l}
    |F_{\widetilde W}^{(\gamma)}|\overset{|\gamma|=4}{\lesssim} (\ln M)^{2} \tau_0^{\frac{1}{4}}.
\end{equation}
Moreover, at the origin $y=0$, there holds
\begin{equation}
    \label{est of F_W^gamma at 0}
    |F_{W}^{(\gamma),0}|\overset{|\gamma|=3}{\lesssim}e^{-\frac{2}{5}s}.
\end{equation}
\end{proposition}
\begin{proof}
    For $|y|\le l$ and $|\gamma|\le4$, we note that
    \begin{equation}
    \label{est of DkW, |y|<l}
        |\partial^\gamma W|\le\underbrace{|\partial^\gamma\ovl W|}_{\eqref{est of Dk ovl W, |y|<l}}+\underbrace{|\partial^\gamma\ovl W|}_{\eqref{asmp of W tilde}}\lesssim\begin{cases}
            (\ln M)^{-5},&{\gamma_1\text{ is even or }\gamma_2\text{ is odd}},\\
            1,&{\gamma_1\text{ is odd and }\gamma_2\text{ is even}}.
        \end{cases}
    \end{equation}
    Then, by \eqref{def of F_tildeW^gamma} and \eqref{def of FW^gamma}, we derive the following expansions, adapted to the cases $\{|y| \le L\}$, $\{|y| \le l\}$, and $y = 0$, respectively:

$$\begin{aligned}
    |F_{\widetilde W}^{(\gamma)}|\overset{|y|\le L,\ |\gamma|\le1}{\lesssim}&\underbrace{|\partial^\gamma F_{\widetilde{W}}|}_{\eqref{est of FtildeW, |y|<L}}+\sum_{0\le\beta<\gamma}\Big(\underbrace{|\partial^{\gamma-\beta}G_W\partial_1\partial^\beta \widetilde{W}|}_{\eqref{est of GR}\eqref{est of GW on |y|<L}\eqref{asmp of W tilde}}+\underbrace{|\partial^{\gamma-\beta}h_W\partial_2\partial^\beta \widetilde{W}|}_{\eqref{est of hR}\eqref{asmp of W tilde}}+\underbrace{|\partial^{\gamma-\beta}(J\partial_1\ovl{W})\partial^\beta \widetilde{W}|}_{\eqref{est of Dk ovl W}\eqref{est of geo. coef.}\eqref{asmp of W tilde}}\Big)\\
		&+\mathbbm{1}_{\gamma_2>0}\underbrace{|\partial_2(JW)\partial_1^{\gamma_1+1}\partial_2^{\gamma_2-1}\widetilde{W}|}_{\eqref{asmp of W}\eqref{asmp of W tilde}\eqref{est of geo. coef.}}+\mathbbm{1}_{\gamma_1>0}\underbrace{|\partial_1JW\partial^\gamma\widetilde{W}|}_{\eqref{asmp of W}\eqref{asmp of support}\eqref{est of geo. coef.}}+\mathbbm{1}_{|\gamma|\ge2}\sum_{\substack{0\le|\beta|\le|\gamma|-2\\\beta\le\gamma}}\underbrace{|\partial^{\gamma-\beta}(JW)\partial_1\partial^\beta\widetilde{W}|}_{\eqref{asmp of W}\eqref{asmp of W tilde}\eqref{est of geo. coef.}},
\end{aligned}$$

$$\begin{aligned}
    |F_{\widetilde W}^{(\gamma)}|\overset{|y|\le l,\ |\gamma|=4}{\lesssim}&\underbrace{|\partial^\gamma F_{\widetilde{W}}|}_{\eqref{est of FtildeW, |y|<l}}+\sum_{0\le\beta<\gamma}\Big(\underbrace{|\partial^{\gamma-\beta}G_W\partial_1\partial^\beta \widetilde{W}|}_{\eqref{est of GR}\eqref{est of GW on |y|<l}\eqref{asmp of W tilde}}+\underbrace{|\partial^{\gamma-\beta}h_W\partial_2\partial^\beta \widetilde{W}|}_{\eqref{est of hR}\eqref{asmp of W tilde}}+\underbrace{|\partial^{\gamma-\beta}(J\partial_1\ovl{W})\partial^\beta \widetilde{W}|}_{\eqref{est of Dk ovl W}\eqref{est of geo. coef.}\eqref{asmp of W tilde}}\Big)\\
		&+\mathbbm{1}_{\gamma_2>0}\underbrace{|\partial_2(JW)\partial_1^{\gamma_1+1}\partial_2^{\gamma_2-1}\widetilde{W}|}_{\eqref{est of DkW, |y|<l}\eqref{asmp of W tilde}\eqref{est of geo. coef.}}+\mathbbm{1}_{\gamma_1>0}\underbrace{|\partial_1JW\partial^\gamma\widetilde{W}|}_{\eqref{est of DkW, |y|<l}\eqref{asmp of support}\eqref{est of geo. coef.}}+\mathbbm{1}_{|\gamma|\ge2}\sum_{\substack{0\le|\beta|\le|\gamma|-2\\\beta\le\gamma}}\underbrace{|\partial^{\gamma-\beta}(JW)\partial_1\partial^\beta\widetilde{W}|}_{\eqref{est of DkW, |y|<l}\eqref{asmp of W tilde}\eqref{est of geo. coef.}},
\end{aligned}$$

$$\begin{aligned}
    |F_{W}^{(\gamma),0}|\overset{|\gamma|=3}{\lesssim}&\underbrace{|\partial^\gamma F_W^0|}_{\eqref{est of FW,FZ}}+\underbrace{|([\partial^\gamma,J]W)^0|}_{\eqref{match W and ovl W}\eqref{est of geo. coef.}}+\mathbbm{1}_{\gamma_2>0}\underbrace{|\partial_2(JW)^0\partial_1^{\gamma_1+1}\partial_2^{\gamma_2-1}W^0|}_{\eqref{match W and ovl W}\eqref{est of geo. coef.}\eqref{est of DkW, |y|<l}}\\
		&+\sum_{\substack{1\le|\beta|\le|\gamma|-2\\\beta\le\gamma}}\underbrace{|\partial^{\gamma-\beta}(JW)^0\partial_1\partial^\beta W^0|}_{\eqref{match W and ovl W}\eqref{est of geo. coef.}}+\sum_{0\le\beta<\gamma}\Big(\underbrace{|\partial^{\gamma-\beta}G_W^0\partial_1\partial^\beta W^0|+|\partial^{\gamma-\beta}h_W^0\partial_2\partial^\beta W^0|}_{\eqref{match W and ovl W}\eqref{est of GR}\eqref{est of hR}\eqref{est of D2GW0}\eqref{est of DkW, |y|<l}}\Big),
\end{aligned}$$
We then bound the corresponding terms with the indicated estimates, and apply Lemma~\ref{lem: compare monos} to get the stated results.
\end{proof}

For convenience, we rewrite the bounds of $|F_{\widetilde W}^{(\gamma)}|$ in the following simplified form.
\begin{corollary}
    \label{cor: est of F_tildeW^gamma (simplified)}
    The following estimates hold  for $|F_{\widetilde W}^{(\gamma)}|$:
    \begin{equation}
\label{est of F_tildeW^gamma (simplified)}
    \eta^{\mu(\gamma)}|F_{\widetilde W}^{(\gamma)}|\mathbbm{1}_{|y|\le L}\overset{|\gamma|\le1}{\lesssim}\eta^{-\frac{1}{8}}\tau_0^{\frac{1}{\gamma_1+2\gamma_2+15}}\ln M.
\end{equation}
\end{corollary}

\section{Estimates for Modulation Variables}
\label{section: Estimates for Modulation Variables}
In this section, we bound the modulation variables, and improve the corresponding bootstrap assumptions \eqref{asmp of Q}\eqref{asmp of kappa tau}\eqref{asmp of psi}. 

\subsection{Estimates on \texorpdfstring{$\kappa$}{kappa}}
By \eqref{eqn of kappa}, we have that
$$
\begin{aligned}
    |\partial_{\tilde t}\kappa|&\lesssim e^{\frac{s}{2}}(\underbrace{|G_W^0|}_{\eqref{est of GW0, hW0}}+\underbrace{|F_W^0|}_{\eqref{est of FW,FZ}})\lesssim (\ln M)^{\frac{1}{10}} \tau_0^{\frac{1}{6}}\le\frac{1}{2}M\tau_0^{\frac{1}{6}},
\end{aligned}$$
and
$$|\kappa-\kappa_0|\le\int_{s_0}^s|\partial_s\kappa|\dif s'\le\int_{s_0}^s\beta_\tau e^{-s'}\underbrace{|\partial_{\tilde t}\kappa|}_{\eqref{asmp of kappa tau}}\dif s'\lesssim\int_{s_0}^sM\tau_0^{\frac{1}{6}}e^{-s'}\dif s'\lesssim M\tau_0^{\frac{7}{6}}\le\frac{1}{2}\tau_0.$$

\subsection{Estimate on \texorpdfstring{$\tau$}{tau}}
Using \eqref{eqn of tau}, we deduce that
$$\begin{aligned}
    |\partial_{\tilde t}\tau|&\lesssim\underbrace{|\partial_1F_W^0|}_{\eqref{est of FW,FZ}}+\underbrace{|\partial_1G_W^0|}_{\eqref{est of GR}}\lesssim M^{\frac{1}{4}} e^{-s}\le\frac{1}{2}M^{\frac12}e^{-s}.
\end{aligned}$$

\subsection{Estimate on \texorpdfstring{$Q$}{Q}}
Via using \eqref{eqn of Q12}, we obtain that
\begin{equation}
\begin{aligned}
    |Q_{12}|\lesssim\underbrace{|\partial_2F_W^0|}_{\eqref{est of FW,FZ}}+\underbrace{e^{\frac{s}{2}}|\partial_2Z^0|}_{\eqref{asmp of Z}}+\underbrace{|\psi r_0Q_{23}|}_{\eqref{absorb r0 sigma with M}\eqref{asmp of Q}\eqref{asmp of psi}}\lesssim \tau_0^{\frac{1}{3}}\le\frac{1}{2}M\tau_0^{\frac{1}{3}}.
\end{aligned}
\end{equation}
Then, it follows from \eqref{eqn of Q13,Q23} that
\begin{equation}
\begin{aligned}
    |Q_{23}|\lesssim\underbrace{e^{\frac{s}{2}}|h_W^0|}_{\eqref{est of GW0, hW0}}+\underbrace{|A^0|}_{\eqref{asmp of A}}\lesssim M\tau_0+ M^{\frac{3}{4}} e^{-\frac{s}{2}}\le\frac{1}{2}M\tau_0^{\frac{1}{2}},
\end{aligned}
\end{equation}
and
\begin{equation}
    \label{est of Q13}
    \begin{aligned}
    \left|Q_{13}-\frac{2\beta_3\kappa_0}{r_0}\right|&\lesssim e^{-\frac{s}{2}}\underbrace{|G
    _W^0|}_{\eqref{est of GW0, hW0}}+
    \underbrace{|Z^{0}+\kappa|}_{\eqref{asmp of Z}\eqref{asmp of kappa tau}\eqref{init conditions for mod variables}}+\underbrace{|\kappa-\kappa_0|}_{\eqref{asmp of kappa tau}}\lesssim\tau_0^{\frac{1}{6}}\le\frac{1}{2}M\tau_0^{\frac{1}{6}}.
\end{aligned}
\end{equation}

\subsection{Estimate on \texorpdfstring{$\psi$}{psi}}
Employing \eqref{eqn of psi}, we derive that
$$\begin{aligned}
    \left|\pt \psi-\left(2\psi^2+\frac{3}{r_0^2}\right)\beta_3\kappa_0\right|&\lesssim e^{\frac{s}{2}}\left(\underbrace{|\partial_{22}F_W^0|}_{\eqref{est of FW,FZ}}+\underbrace{|G_W^0||\partial_{122}W^0|+|h_W^0||\partial_{222}W^0|}_{\eqref{evaluation of Wbar and derivs at 0}\eqref{asmp of W tilde}\eqref{est of GW0, hW0}}\right)+e^s\underbrace{|\partial_{22}Z^0|}_{\eqref{asmp of Z}}\\
    &\quad+\underbrace{\left(\frac{1}{2 r_{0}^{2}}+\psi^{2}\right)}_{\eqref{asmp of psi}}\left[\underbrace{|Z^0+\kappa|}_{\eqref{asmp of Z}\eqref{est of Z+kappa}}+\underbrace{|\kappa-\kappa_0|}_{\eqref{asmp of kappa tau}}\right]+\underbrace{\left|Q_{13}-\frac{2\beta_3\kappa_0}{r_0}\right|}_{\eqref{est of Q13}}\\
    &\lesssim \tau_0^{\frac{1}{6}}\ln M\le\frac{1}{2}M\tau_0^{\frac{1}{6}}.
\end{aligned}$$
Here, we use the facts that $|\psi_0|\le10$ and $\left(\frac{1}{2 r_{0}^{2}}+\psi^{2}\right)\le\ln M$. The above estimate for $\pt\psi$ also implies $|\pt\psi|\le\ln M$. Therefore, for $\psi$ itself, we have that
$$|\psi-\psi_0|\le\int_{s_0}^s\ln Me^{-s'}\dif s'\lesssim\tau_0\ln M\le\frac{1}{2}M\tau_0.$$

\section{Estimates for \texorpdfstring{$Z$}{Z} and \texorpdfstring{$A$}{A}}
\label{section: Estimates for Z and A}
In this section, we improve the bootstrap assumptions \eqref{asmp of Z} and \eqref{asmp of A} for \texorpdfstring{$Z$}{Z} and \texorpdfstring{$A$}{A}. Let $\mathcal{R}\in\{Z,A\}$. We define 
\begin{equation}
    \widetilde{\mathcal{R}}=\begin{cases}
        Z+\sigma_\infty,&\mathcal{R}=Z,\\
        A,&\mathcal{R}=A. 
    \end{cases}
\end{equation}
From \eqref{eqns of deriv of W,Z,A}, $\partial^\gamma\widetilde{\mathcal{R}}$ with $\gamma\ge0$ solves the transport-type equation
$$\partial_s\partial^\gamma\widetilde{\mathcal{R}}+D_\mathcal{R}^{(\gamma)}\partial^\gamma\widetilde{\mathcal{R}}+\mathcal{V}_\mathcal{R}\cdot\nabla \partial^\gamma\widetilde{\mathcal{R}}=F_{\mathcal{R}}^{(\gamma)}.$$
By solving this equation along the characteristic, we get
$$
\begin{aligned}
    \partial^\gamma\widetilde{\mathcal{R}}(s,\Phi_{\mathcal R}(s;s_1,y_0))=&\partial^\gamma\widetilde{\mathcal{R}}(s_1,y_0)\exp\left(-\int_{s_1}^sD_\mathcal{R}^{(\gamma)}(s',\Phi_{\mathcal{R}}(s';s_1,y_0))\dif s'\right)\\
    &+\int_{s_1}^sF_{\mathcal{R}}^{(\gamma)}(s',\Phi_{\mathcal R}(s';s_1,y_0))\exp\left(-\int_{s'}^sD_\mathcal{R}^{(\gamma)}(s'',\Phi_{\mathcal{R}}(s'';s_1,y_0))\dif s''\right)\dif s'.
\end{aligned}$$
Note that by \eqref{def of DRgamma}, \eqref{asmp of W} and \eqref{est of geo. coef.}, we have $D_\mathcal{R}^{(\gamma)}=\frac{3\gamma_1+\gamma_2}{2}+O_\gamma(\eta^{-\frac{1}{3}})$, where the notation $f=O_\gamma(\eta^{-\frac{1}{3}})$ means that $|f|\le C_\gamma\eta^{-\frac{1}{3}}$, with constant $C_\gamma$ depending on $\gamma$. Hence, by invoking Lemma \ref{integration along ZA trajectories}, and setting $s_1=s_0$, we obtain
$$\begin{aligned}
    |\partial^\gamma\widetilde{\mathcal{R}}(s,\Phi_\mathcal{R}(s;s_0,y_0))|\lesssim&|\partial^\gamma\widetilde{\mathcal{R}}(s_0,y_0)|e^{-\frac{3\gamma_1+\gamma_2}{2}(s-s_0)}+\int_{s_0}^s|F_{\mathcal{R}}^{(\gamma)}(s',\Phi_\mathcal{R}(s';s_0,y_0))|e^{-\frac{3\gamma_1+\gamma_2}{2}(s-s')}\dif s'.
\end{aligned}$$
Since the map $y_0\mapsto\Phi_\mathcal{R}(s;s_0,y_0)$ is a bijection on $\mathbb{R}^2$, the above identity implies that
$$
    \|\partial^\gamma\widetilde{\mathcal{R}}(s,\cdot)\|_{L^\infty}\lesssim\|\partial^\gamma\widetilde{\mathcal{R}}(s_0,\cdot)\|_{L^\infty}e^{-\frac{3\gamma_1+\gamma_2}{2}(s-s_0)}+\sup_{y_0}\int_{s_0}^s|F_{\mathcal{R}}^{(\gamma)}(s',\Phi_\mathcal{R}(s';s_0,y_0))|e^{-\frac{3\gamma_1+\gamma_2}{2}(s-s')}\dif s',
$$
and
\begin{equation}
    \label{est of partial^gamma R}
    \begin{aligned}
    e^{\left(3\mu(\gamma)+\frac{1}{2}\right)s}\|\partial^\gamma\widetilde{\mathcal{R}}(s,\cdot)\|_{L^\infty}\lesssim&\underbrace{\|\partial^\gamma\widetilde{\mathcal{R}}(s_0,\cdot)\|_{L^\infty}}_{\eqref{init of Z,A}}\tau_0^{-\left(3\mu(\gamma)+\frac{1}{2}\right)}\underbrace{e^{\left(3\mu(\gamma)-\frac{3\gamma_1+\gamma_2-1}{2}\right)(s-s_0)}}_{\le1,\text{ see }\eqref{def of mu}}\\
    &+e^{\left(3\mu(\gamma)-\frac{3\gamma_1+\gamma_2-1}{2}\right)s}\sup_{y_0}\int_{s_0}^s|F_{\mathcal{R}}^{(\gamma)}(s',\Phi_\mathcal{R}(s';s_0,y_0))|e^{\frac{3\gamma_1+\gamma_2}{2}s'}\dif s'\\
    &\lesssim\tau_0+e^{\left(3\mu(\gamma)-\frac{3\gamma_1+\gamma_2-1}{2}\right)s}\sup_{y_0}\int_{s_0}^s|F_{\mathcal{R}}^{(\gamma)}(s',\Phi_\mathcal{R}(s';s_0,y_0))|e^{\frac{3\gamma_1+\gamma_2}{2}s'}\dif s'.
\end{aligned}
\end{equation}
The second term, representing the contribution of the forcing term, is now the primary concern.

\begin{proposition}

    We have the following estimates for $Z$ and $A$:
    \begin{equation}
    e^{\left(3\mu(\gamma)+\frac{1}{2}\right)s}|\partial^\gamma(Z+\sigma_\infty)|\le\frac{1}{2}\times\begin{cases}
        M \tau_{0},&\gamma=(0,0),\\
        \tau_0^{\frac{1}{3}},&|\gamma|=1,2,\gamma_1=0,\\
        M^{\frac{|\gamma|+\gamma_2}{4}},&|\gamma|=1,2,\gamma_1>0,
\end{cases}
\end{equation}
\begin{equation}
    e^{\left(3\mu(\gamma)+\frac{1}{2}\right)s}|\partial^\gamma A|\le\frac{1}{2}\times\begin{cases}
        M\tau_0^{\frac{1}{|\gamma|+1}},&|\gamma|\le 1\text{ and }\gamma_1=0,\\
        M,&\gamma=(0,2).
\end{cases}.
\end{equation}
\end{proposition}
\begin{proof}
    We deal with $Z$ first. From \eqref{est of partial^gamma R}, we deduce that
    $$\begin{aligned}
        |Z+\sigma_\infty|&\lesssim\tau_0+\sup_{y_0}\int_{s_0}^s\underbrace{|F_{Z}^{(0,0)}(s',\Phi_Z(s';s_0,y_0))|}_{\eqref{est of FZ^gamma (simplified)}}\dif s'\lesssim\tau_0+\int_{s_0}^se^{-s'}\dif s'\lesssim\tau_0\le\frac{1}{2}M\tau_0,
    \end{aligned}$$
    $$\begin{aligned}
        e^{\left(3\mu(\gamma)+\frac{1}{2}\right)s}|\partial^\gamma Z|\overset{|\gamma|=1,2\ \gamma_2=0}&{\lesssim}\tau_0+\underbrace{e^{\left(3\mu(\gamma)-\frac{3\gamma_1+\gamma_2-1}{2}\right)s}}_{=1,\text{ see }\eqref{def of mu}}\sup_{y_0}\int_{s_0}^s\underbrace{|F_{\mathcal{R}}^{(\gamma)}(s',\Phi_Z(s';s_0,y_0))|}_{\eqref{est of FZ^gamma (simplified)}}e^{\frac{\gamma_2}{2}s'}\dif s'\\
        &\lesssim\tau_0+\int_{s_0}^sMe^{-\frac{2}{5}s'}\dif s'\lesssim M\tau_0^{\frac{2}{5}}\le\frac{1}{2}\tau_0^{\frac{1}{3}},
    \end{aligned}$$
    $$\begin{aligned}
        e^{\left(3\mu(\gamma)+\frac{1}{2}\right)s}|\partial^\gamma Z|\overset{|\gamma|=1,2\ \gamma_2>0}&{\lesssim}\tau_0+\underbrace{e^{\left(3\mu(\gamma)-\frac{3\gamma_1+\gamma_2-1}{2}\right)s}}_{\eqref{def of mu}}\sup_{y_0}\int_{s_0}^s\underbrace{|F_{\mathcal{R}}^{(\gamma)}(s',\Phi_Z(s';s_0,y_0))|}_{\eqref{est of FZ^gamma (simplified)}}e^{\frac{3\gamma_1+\gamma_2}{2}s'}\dif s'\\
        &\lesssim\tau_0+e^{-\frac{3\gamma_1+\gamma_2-3}{2}s}\int_{s_0}^sM^{\frac{|\gamma|+\gamma_2}{4}-\frac{1}{8}}e^{\frac{3\gamma_1+\gamma_2-3}{2}s'}\eta^{-\frac{1}{8}\mathbbm{1}_{|\gamma|=1}}(\Phi_Z(s';s_0,y_0))\dif s'\\
        \overset{\text{Lemma }\ref{lemma: integral of eta^-p along PhiZ and A}}&{\lesssim}\tau_0+M^{\frac{|\gamma|+\gamma_2}{4}-\frac{1}{8}}\le\frac{1}{2}M^{\frac{|\gamma|+\gamma_2}{4}},
    \end{aligned}$$
    Next, we estimate derivatives of $A$. Using \eqref{est of partial^gamma R} again, we obtain
    $$\begin{aligned}
        e^{\left(3\mu(\gamma)+\frac{1}{2}\right)s}|\partial^\gamma A|\overset{|\gamma
        |\le1,\gamma_1=0}&{\lesssim}\tau_0+\underbrace{e^{\left(3\mu(\gamma)-\frac{3\gamma_1+\gamma_2-1}{2}\right)s}}_{=1,\text{ see }\eqref{def of mu}}\sup_{y_0}\int_{s_0}^s\underbrace{|F_{A}^{(\gamma)}(s',\Phi_A(s';s_0,y_0))|}_{\eqref{est of FA^gamma (simplified)}}e^{\frac{3\gamma_1+\gamma_2}{2}s'}\dif s'\\
        &\lesssim\tau_0+\int_{s_0}^sM^{\frac{1}{2}}e^{-\frac{s'}{1+|\gamma|}}e^{-\frac{|\gamma|}{2}s'}e^{\frac{|\gamma|}{2}s'}\dif s'\lesssim M^{\frac{1}{2}}\tau_0^{\frac{s'}{1+|\gamma|}}\le\frac{1}{2}M\tau_0^{\frac{s'}{1+|\gamma|}},
    \end{aligned}$$
    $$\begin{aligned}
        e^{\left(3\mu(\gamma)+\frac{1}{2}\right)s}|\partial^\gamma A|\overset{\gamma=(0,2)}&{\lesssim}\tau_0+\underbrace{e^{\left(3\mu(\gamma)-\frac{3\gamma_1+\gamma_2-1}{2}\right)s}}_{=1,\text{ see }\eqref{def of mu}}\sup_{y_0}\int_{s_0}^s\underbrace{|F_{A}^{(\gamma)}(s',\Phi_A(s';s_0,y_0))|}_{\eqref{est of FA^gamma (simplified)}}e^{\frac{3\gamma_1+\gamma_2}{2}s'}\dif s'\\
        &\lesssim\tau_0+\int_{s_0}^sM^{\frac{15}{16}}\eta^{-\frac{1}{8}}(\Phi_A(s';s_0,y_0)))e^{-s'}e^{s'}\dif s'\overset{\text{Lemma }\ref{lemma: integral of eta^-p along PhiZ and A}}{\lesssim}M^{\frac{15}{16}}\le\frac{1}{2}M.
    \end{aligned}$$
\end{proof}

We have now improved the bootstrap assumptions \eqref{asmp of Z} and \eqref{asmp of A}, except for the bound on $\partial_1 A$. This remaining estimate builds on the analysis of the vorticity $\omega$.

\begin{proposition}
The vorticity $\omega$ satisfies the identity:
\begin{equation}
    \omega=\varphi^{-2}\left[\varphi\jpnlambda\partial_{\tilde u_1}a-\partial_{\tilde u_2}(\varphi V_1)\right],\nonumber
\end{equation}
Furthermore, it admits the estimate:
\begin{equation}
    |\omega|\approx e^{\frac{3}{2}s}|\partial_1A|+O(\tau_0^{\frac{1}{2}}),\nonumber
\end{equation}
where $\approx$ denotes equivalence up to multiplicative constants, i.e. 
$$|\omega|\lesssim e^{\frac{3}{2}s}|\partial_1A|+\tau_0^{\frac{1}{2}}\quad\text{and}\quad e^{\frac{3}{2}s}|\partial_1A|\lesssim |\omega|+\tau_0^{\frac{1}{2}}. $$
\end{proposition}
\begin{proof}
In terms of the definition of $\omega:=\epsilon^{12}(\nabla_{u_1}v_2^\flat-\nabla_{u_2}v_1^\flat)$, we derive that
$$
\begin{aligned}
    \omega&=\epsilon^{12}(\nabla_{u_1}v_2^\flat-\nabla_{u_2}v_1^\flat)=\epsilon^{12}(\partial_{u_1}v_2^\flat-\partial_{u_2}v_1^\flat)=\varphi^{-2}(\partial_{u_1}(\varphi^2v_2)-\partial_{u_2}(\varphi^2v_1))=\varphi^{-2}(\partial_{u_1}(\varphi V_2)-\partial_{u_2}(\varphi V_1))\\
    &=\varphi^2\left[\partial_{\tilde u_1}\left(-\frac{\lambda}{2J}(w+z)+\frac{a}{J}\right)-(\partial_{\tilde u_2}-\lambda\partial_{\tilde u_1})\left(\frac{w+z}{2J}+\frac{\lambda}{J}a\right)\right]\\
    &=\varphi^{-2}\left[\varphi\jpnlambda\partial_{\tilde u_1}a-\partial_{\tilde u_2}(\varphi V_1)\right].
\end{aligned}$$
Recalling the relations $\partial_{\tilde u_1}=e^{\frac{3}{2}s}\partial_1$ and $\partial_{\tilde u_2}=e^{\frac{s}{2}}\partial_2$, we obtain
$$
\begin{aligned}
    \omega&=\underbrace{\varphi^{-2}}_{\eqref{est of geo. coef.}}[\underbrace{\varphi\jpnlambda}_{\eqref{est of geo. coef.}}e^{\frac{3}{2}s}\partial_1A-\underbrace{e^\frac{s}{2}\partial_2(\varphi V_1)}_{\eqref{est of geo. coef.},\text{ Corollary }\ref{cor: lower order est of P tilde}}]\\
    &=\left(1+O(\tau_0^{\frac{1}{3}})\right)\left[\left(1+O(\tau_0^{\frac{1}{6}})\right)e^{\frac{3}{2}s}\partial_1A+O(\tau_0^{\frac{1}{2}})\right],
\end{aligned}$$
which yields $|\omega|\lesssim e^{\frac{3}{2}s}|\partial_1A|+\tau_0^{\frac{1}{2}}$ and $e^{\frac{3}{2}s}|\partial_1A|\lesssim |\omega|+\tau_0^{\frac{1}{2}}$.
\end{proof}
Now, we are ready to estimate $\partial_1A$.
\begin{proposition}
The quantity $e^{\frac{3}{2}s}\partial_1 A$ is almost conserved. Specifically, we have
$$
e^{\frac{3}{2}s}\|\partial_1A(s,\cdot)\|_{L^\infty} \lesssim e^{\frac{3}{2}s_0}\|\partial_1A(s_0,\cdot)\|_{L^\infty} + \tau_0^{\frac{1}{2}}.
$$
\end{proposition}
\begin{proof}
From \eqref{eqn of specific vorticity}, we note that $\frac{\omega}{\rho} = \frac{\omega}{(\alpha S)^{\frac{1}{\alpha}}}$ is transported by the velocity field, and hence $\|\frac{\omega}{\rho}(s,\cdot)\|_{L^\infty} = \|\frac{\omega}{\rho}(s_0,\cdot)\|_{L^\infty}$. 

By taking $\tau_0$ small enough and employing Corollary \ref{cor: lower order est of P tilde}, we have $|S| \approx \sigma_\infty$. Therefore, we derive
$$
\begin{aligned}
    e^{\frac{3}{2}s}\|\partial_1A(s,\cdot)\|_{L^\infty} &\lesssim \|\omega(s,\cdot)\|_{L^\infty} + \tau_0^{\frac{1}{2}} \\
    &\lesssim \|\rho(s,\cdot)\|_{L^\infty} \left\|\frac{\omega}{\rho}(s,\cdot)\right\|_{L^\infty} + \tau_0^{\frac{1}{2}} \\
    &\lesssim (\alpha\sigma_\infty)^{\frac{1}{\alpha}} \|\omega(s_0,\cdot)\|_{L^\infty} \left\|\frac{1}{\rho(s_0,\cdot)}\right\|_{L^\infty} + \tau_0^{\frac{1}{2}} \\
    &\lesssim (\alpha\sigma_\infty)^{\frac{1}{\alpha}} \left(e^{\frac{3}{2}s_0} \|\partial_1A(s_0,\cdot)\|_{L^\infty} + C\tau_0^{\frac{1}{2}}\right)(\alpha\sigma_\infty)^{-\frac{1}{\alpha}} + \tau_0^{\frac{1}{2}} \\
    &\lesssim e^{\frac{3}{2}s_0} \|\partial_1A(s_0,\cdot)\|_{L^\infty} + \tau_0^{\frac{1}{2}}.
\end{aligned}
$$
\end{proof}

Using the above lemma and the initial condition on $\partial_1A$, we deduce that 
\[
\|\partial_1A(s,\cdot)\|_{L^\infty} \lesssim \|\partial_1A(s_0,\cdot)\|_{L^\infty}e^{\frac{3}{2}(s-s_0)} + \tau_0^{\frac{1}{2}}e^{-\frac{3}{2}s} \le \frac{1}{2}M e^{-\frac{3}{2}s}.
\]
This improves the bootstrap assumption for $\partial_1A$.

\section{Estimates for \texorpdfstring{$W$}{W} and \texorpdfstring{$\widetilde W$}{W tilde}}
\label{section: Estimates for W}
This section is devoted to improving the bootstrap assumptions for the Riemann variable $W$ and its deviation $\widetilde{W}$, i.e., \eqref{asmp of W} and \eqref{asmp of W tilde}.

Let $\mathcal{R}\in\{W,\widetilde W\}$. Via using \eqref{eqn of W,Z,A}\eqref{eqns of deriv of W,Z,A}\eqref{eqn of tilde W}\eqref{eqn of derivatives of tilde W}, the equation for $\partial^\gamma\mathcal{R}$ can be written in the form
$$\partial_s\partial^\gamma\mathcal{R}+D_\mathcal{R}^{(\gamma)}\partial^\gamma\mathcal{R}+\mathcal{V}_W\cdot\nabla \partial^\gamma\mathcal{R}=F_\mathcal{R}^{(\gamma)},$$
with $D_\mathcal{R}^{(\gamma)}=\frac{3\gamma_1+\gamma_2-1}{2}+O_{\mathcal{R},\gamma}(\eta^{-\frac{1}{3}})$. 

\subsection{Estimates Near the Origin}
We begin with the estimates for $W$ at $y=0$. By letting $|\gamma|=3$ and $y=0$ in the equation \eqref{eqns of deriv of W,Z,A} of $\partial^\gamma W$, we obtain
$$\partial_s\partial^\gamma W^0+(1+\gamma_1)(1-\beta_\tau)\partial^\gamma W^0+G_W^0\partial_1\partial^\gamma W+h_W^0\partial_2\partial^\gamma W=F_W^{(\gamma),0}.$$
It follows from \eqref{asmp of W tilde} that $|\partial^\beta W^0|\le|\partial^\beta W^0|+|\partial^\beta \ovl W^0|\lesssim1$ for $|\beta|\le4$. Thus, $|\partial_s\partial^\gamma W^0|$ is bounded as below:
$$\begin{aligned}
    |\partial_s\partial^\gamma W^0|&\lesssim \underbrace{|1-\beta_\tau|}_{\eqref{est of 1-beta_tau}}+\underbrace{|G_W^0|+|h_W^0|}_{\eqref{est of GW0, hW0}}+\underbrace{|F_W^{(\gamma),0}|}_{\eqref{est of F_W^gamma at 0}}\lesssim e^{-\frac{2}{5}s},
\end{aligned}$$
and $|\partial^\gamma \widetilde W^0|$ is controlled in the following manner:
$$|\partial^\gamma \widetilde W^0(s)|\le\underbrace{|\partial^\gamma \widetilde W^0(s_0)|}_{\eqref{init of W tilde}}+\int_{s_0}^s|\partial_s\partial^\gamma W^0(s')|\dif s'\lesssim\tau_0^{\frac{2}{5}}\le\frac{1}{2}\tau_0^{\frac{1}{3}}.$$
Here, we use the fact $\partial_s\partial^\gamma W^0=\partial_s\partial^\gamma \widetilde W^0$. This closes the bootstrap argument for $|D^3\widetilde W^0|$.

Next, we deal with $\widetilde W$ in the region $\{|y|\le l\}$. For $|\gamma|=4$, we derive a lower bound estimate for the damping term $D_{\widetilde W}^{(\gamma)}$ in $\{|y|\le l\}$:
$$\begin{aligned}
    D_{\widetilde W}^{(\gamma)}\overset{\eqref{eqn of derivatives of tilde W}}&{=}\frac{3\gamma_1+\gamma_2-1}{2}+\beta_\tau J(\partial_1\ovl{W}+\gamma_1\partial_1W)\\
    &=\frac{3\gamma_1+\gamma_2-1}{2}+\underbrace{\beta_\tau J}_{\eqref{est of geo. coef.}}\underbrace{((1+\gamma_1)\partial_1\ovl{W}+\gamma_1\partial_1\widetilde W)}_{\eqref{est of ovl W and D ovl W}\eqref{asmp of W tilde}}\\
    &\ge\frac{1}{2}+(1+\gamma_1)-\gamma_1 o_{\tau_0}(1)\ge\frac{1}{3}.
\end{aligned}$$
Using Proposition \ref{prop: outgoing property of W trajectories} and a contradiction argument, we conclude that $|\Phi(s';s,y)| \le l$ for $s' \in [s_0, s]$ and $|y| \le l$. Hence, by solving along trajectories, we obtain for $|\gamma|=4$ that
$$\begin{aligned}
    |\partial^\gamma \widetilde W(s,y)|\overset{|y|\le l}&{\le}|\partial^\gamma \widetilde W(s_0,\Phi(s_0;s,y))|e^{-\frac{s-s_0}{2}}+\int_{s_0}^s\underbrace{|F_{\widetilde W}^{(\gamma)}(s',\Phi(s';s,y))|}_{\eqref{est of F_tildeW^gamma, |y|<l}}e^{-\frac{s-s'}{2}}\dif s'\\
    &\lesssim\underbrace{\|\partial^\gamma\widetilde W(s_0,\cdot)\|_{L^\infty(|y|\le l)}}_{\eqref{init of W tilde}}+\int_{s_0}^se^{-\frac{s-s'}{2}}(\ln M)^{2} \tau_0^{\frac{1}{4}}\dif s'\\
    &\lesssim(\ln M)^{2} \tau_0^{\frac{1}{4}}\le\frac12(\ln M)^{4} \tau_0^{\frac{1}{4}}.
\end{aligned}$$
For $|\gamma| = 3$ and $|y| \le l$, we have  
$$|\partial^\gamma \widetilde W(s,y)| \lesssim |\partial^\gamma \widetilde W^0| + l \cdot \|D^4 \widetilde W\|_{L^\infty(|y| \le l)} \lesssim \tau_0^{\frac{1}{3}} + (\ln M)^{-5}\cdot (\ln M)^{4} \tau_0^{\frac{1}{4}} \le \frac{1}{2}(\ln M)^{3} \tau_0^{\frac{1}{4}}.$$
For $|\gamma| \le 2$, since $\partial^\gamma \widetilde W^0 = 0$, we derive consecutively that  
$$|\partial^\gamma \widetilde W(s,y)| \overset{|y|\le l}{\lesssim}l \cdot \|D^{|\gamma|+1} \widetilde W\|_{L^\infty(|y| \le l)}\lesssim l\cdot (\ln M)^{|\gamma|+1}\tau_0^{\frac{1}{4}} \le \frac{1}{2} (\ln M)^{|\gamma|}\tau_0^{\frac{1}{4}}.
$$

\subsection{Estimates Away from the Origin}
In this section, we continue to improve \eqref{asmp of W} and \eqref{asmp of W tilde} outside $\{|y|\le l\}$. Recall that $\mu=\mu(\gamma)$ was defined in \eqref{def of mu}. Let $q=\eta^{\mu(\gamma)}\partial^\gamma\mathcal{R}$, where $\mathcal{R}\in\{W,\widetilde W\}$. Then, $q$ solves the equation
$$\partial_sq+D_qq+\mathcal{V}_W\cdot\nabla q=F_q,$$
with $D_q=D_\mathcal{R}^{(\gamma)}-\mu\eta^{-1}\mathcal{V}_W\cdot\nabla \eta$ and $F_q=\eta^\mu F_\mathcal{R}^{(\gamma)}$. 

\begin{lemma}
For the damping term, it holds that
\begin{equation}
    D_q=\frac{3\gamma_1+\gamma_2-1}{2}-3\mu+O(e^{-\frac{s}{2}})+O_{\gamma}(\eta^{-\frac{1}{6}}).\nonumber
\end{equation}
\end{lemma}
\begin{proof}
Since $D_\mathcal{R}^{(\gamma)}=\frac{3\gamma_1+\gamma_2-1}{2}+O_{\gamma}(\eta^{-\frac{1}{3}})$, it suffices to show that $\mu\eta^{-1}\mathcal{V}_W\cdot\nabla \eta=3\mu+O(e^{-\frac{s}{2}})+O(\eta^{-\frac{1}{6}})$. Note that
$$\begin{aligned}
    |G_W|&\lesssim\underbrace{|G_W^0|}_{\eqref{est of GW0, hW0}}+\underbrace{\|\partial_1G_W\|_{L^\infty}|y_1|+\|\partial_2G_W\|_{L^\infty}|y_2|}_{\eqref{est of GR}}\\
    &\lesssim Me^{-s}+e^{-s}\eta^{\frac{1}{2}}+M^2\tau_0^{\frac{1}{6}}\eta^{\frac{1}{6}}\lesssim Me^{-s}\eta^{\frac{1}{2}}+M^2\tau_0^{\frac{1}{6}}\eta^{\frac{1}{6}},
\end{aligned}$$
and that $|h_W|\lesssim e^{-\frac{s}{2}}$ by \eqref{est of hR}. Therefore, using the assumption $|W|\lesssim\eta^{\frac{1}{6}}$, we have
$$\begin{aligned}
    \mu\eta^{-1}\mathcal{V}_W\cdot\nabla \eta&=\mu\eta^{-1}\left[\left(\frac{3}{2}y_1+\beta_\tau JW+G_W\right)\cdot 2y_1+\left(\frac{1}{2}y_2+h_W\right)\cdot6y_2^5\right]\\
    &=\mu\eta^{-1}\left[3\eta-3+O(\eta^{\frac{2}{3}})+O(Me^{-s}\eta)+o_{\tau_0}(1)\eta^{\frac{5}{6}}\right]\\
    &=3\mu+O(e^{-\frac{s}{2}})+O(\eta^{-\frac{1}{6}}).
\end{aligned}$$
\end{proof}

Combining the above lemma and \eqref{integration along W trajectories}, one can immediately obtain that
\begin{corollary}
There exists a universal constant $C$, such that for any $s'\ge s_1\ge s_0$ and any $y_0\in\mathbb{R}^2$, the following inequality holds:
\begin{equation}
    \exp\left(-\int_{s'}^sD_q(s'',\Phi_{W}(s'';s_1,y_0))\dif s''\right)\le e^{\left(-\frac{3\gamma_1+\gamma_2-1}{2}+3\mu\right)(s-s')}\times\begin{cases}
        (\ln M)^C,&|y_0|\ge l,\\
        (1+C\tau_0^{\frac{1}{30}}),&|y_0|\ge L.
    \end{cases}\nonumber
\end{equation}
\end{corollary}

Next, we consider the cases $|y_0|\ge l$ and $|y_0|\ge L$ separately. 

\subsubsection{Estimates in $\{|y|\ge l\}$} 
We aim to estimate $|q(s,y)|$ for $|y|\ge l$. By tracing the flow backward in time, we assume $y=\Phi_W(s;s_1,y_0)$. In view of Proposition \ref{prop: outgoing property of W trajectories}, we further assume that either $s_1=s_0$ or $|y_0|=l$. Solving along the trajectories then yields the estimate
\begin{equation}
\label{est in |y| > l}
\begin{aligned}
    |q(s,y)|&\le(\ln M)^C\max\Big(\|q(s_0,\cdot)\|_{L^\infty(|y|\ge l)},\sup_{|y'|=l}|q(s_1,y')|\Big)\\
    &\quad+(\ln M)^C\int_{s_1}^s|F_q(s',\Phi_W(s',s_1,y_0))|e^{\left(-\frac{3\gamma_1+\gamma_2-1}{2}+3\mu\right)(s-s')}\dif s'.
\end{aligned}
\end{equation}
In the following proposition, we control the lower order derivatives of $W$.
\begin{proposition}
For $|\gamma|\le1$ and $l\le|y|\le L$, there holds
\begin{equation}
\label{integration of F_tW^gamma}
\eta^{\mu(\gamma)}|\partial^\gamma \widetilde W(s,y)|\le\frac{1}{2}\tau_0^{\frac{1}{\gamma_1+2\gamma_2+16}}.
\end{equation}
\end{proposition}
\begin{proof}
By \eqref{est in |y| > l}, we have
\begin{equation}
\begin{aligned}
    \eta^{\mu(\gamma)}(y)|\partial^\gamma\widetilde W(s,y)|&\le(\ln M)^C\max\Big(\underbrace{\|\eta^{\mu(\gamma)}\partial^\gamma\widetilde W(s_0,\cdot)\|_{L^\infty(|y|\ge l)}}_{\eqref{init of W tilde}},\sup_{|y'|=l}\underbrace{\eta^{\mu(\gamma)}|\partial^\gamma\widetilde W(s_1,y')|}_{\eqref{asmp of W tilde}}\Big)\\
    &\quad+(\ln M)^C\int_{s_1}^s\underbrace{(\eta^{\mu(\gamma)}|F_{\widetilde W}^{(\gamma)}|)(s',\Phi(s';s_1,y_0))|}_{\eqref{est of F_tildeW^gamma (simplified)}}\underbrace{e^{\left(-\frac{3\gamma_1+\gamma_2-1}{2}+3\mu\right)(s-s')}}_{=1,\ \text{see }\eqref{def of mu}}\dif s'\\
    &\le(\ln M)^C\tau_0^{\frac{1}{\gamma_1+2\gamma_2+15}}+(\ln M)^C\int_{s_1}^s\tau_0^{\frac{1}{\gamma_1+2\gamma_2+15}}\eta^{-\frac{1}{8}}(s',\Phi(s';s_1,y_0))\dif s\\
    \overset{\eqref{integration along W trajectories}}&{\le}(\ln M)^C\tau_0^{\frac{1}{\gamma_1+2\gamma_2+15}}\le\frac{1}{2}\tau_0^{\frac{1}{\gamma_1+2\gamma_2+16}}.
\end{aligned}\nonumber
\end{equation}
\end{proof}
Next, we prove the following estimates for higher order derivatives:
\begin{proposition}
For $|y|\ge l$ and $|\gamma|=2,3$, we have the estimates
\begin{equation}
    \eta^{\mu(\gamma)}|\partial^\gamma W(s,y)|\le\frac{1}{2}M^{\frac{3|\gamma|+\gamma_2-5}{8}}.
\end{equation}
\end{proposition}
\begin{proof}
    Direct computation using \eqref{est in |y| > l} leads us to
    $$\begin{aligned}
    \eta^{\mu(\gamma)}|\partial^\gamma W(s,y)|&\le(\ln M)^C\max\Big(\underbrace{\|\eta^{\mu(\gamma)}\partial^\gamma W(s_0,\cdot)\|_{L^\infty(|y|\ge l)}}_{\eqref{init of W}},\underbrace{\sup_{|y'|=l}|\eta^{\mu(\gamma)}\partial^\gamma W(s_1,y')|}_{W=\ovl W+\widetilde W,\ \eqref{est of Dk ovl W}\eqref{asmp of W tilde}}\Big)\\
    &\quad+(\ln M)^C\int_{s_1}^s\underbrace{(\eta^{\mu(\gamma)}|F_W^{(\gamma)}|)(s',\Phi(s',s_1,y_0))}_{\eqref{est of FW^gamma (simplified)}}e^{\left(-\frac{3\gamma_1+\gamma_2-1}{2}+3\mu\right)(s-s')}\dif s'\\
    &\le(\ln M)^C+(\ln M)^C\int_{s_1}^sM^{\frac{3|\gamma|+\gamma_2-6}{8}}\eta^{-\frac{2}{15}\mathbbm{1}_{\gamma_1=0}}e^{-\left(\frac{3\gamma_1+\gamma_2-1}{2}-1+\frac{1}{2}\mathbbm{1}_{\gamma_1=0}\right)(s-s')}\dif s'\\
    \overset{|\gamma|=2,3}&{\le}(\ln M)^C+(\ln M)^CM^{\frac{3|\gamma|+\gamma_2-6}{8}}\int_{s_1}^s\eta^{-\frac{2}{15}\mathbbm{1}_{\gamma_1=0}}e^{-\frac{1}{2}(s-s')\mathbbm{1}_{\gamma_1>0}}\dif s'\\
    \overset{\eqref{integration along W trajectories}}&{\le}(\ln M)^C+(\ln M)^CM^{\frac{3|\gamma|+\gamma_2-6}{8}}\le\frac{1}{2}M^{\frac{3|\gamma|+\gamma_2-5}{8}}.
\end{aligned}$$
\end{proof}

\subsubsection{Estimates in $\{|y|\ge L\}$}
We proceed to estimates in the region $\{|y|\ge L\}$. According to Proposition \ref{prop: outgoing property of W trajectories}, we can require either $s_1=s_0$ or $|y_0|=L$ in this case. Consequently, we gain the bound on $q(s,y)$ for $|y|\ge L$:
\begin{equation}
\label{est in |y| > L}
\begin{aligned}
    |q(s,y)|&\le(1+C\tau_0^{\frac{1}{30}})\max\Big(\|q(s_0,\cdot)\|_{L^\infty(|y|\ge L)},\sup_{|y'|=L}|q(s_1,y')|\Big)\\
    &\quad+(1+C\tau_0^{\frac{1}{30}})\int_{s_1}^s|F_q(s',\Phi(s',s_1,y_0))|e^{\left(-\frac{3\gamma_1+\gamma_2-1}{2}+3\mu\right)(s-s')}\dif s'.
\end{aligned}
\end{equation}

\begin{proposition}
    For $|\gamma|\le1$ and $|y|\ge L$, we derive the following estimates:
    \begin{equation}
        \eta^{\mu(\gamma)}|\partial^\gamma W|\le\begin{cases}
            \frac{151}{150},&\gamma=(0,0),\\
            18,&|\gamma|=1.
        \end{cases}
    \end{equation}
\end{proposition}
\begin{proof}
    Using \eqref{est in |y| > L}, we compute:
    $$\begin{aligned}
    \eta^{\mu(\gamma)}|\partial^\gamma W(s,y)|&\le(1+C\tau_0^{\frac{1}{30}})\max\Big(\underbrace{\|\eta^{\mu(\gamma)}\partial^\gamma W(s_0,\cdot)\|_{L^\infty(|y|\ge l)}}_{\eqref{init of W}},\underbrace{\sup_{|y'|=L}|\eta^{\mu(\gamma)}\partial^\gamma W(s_1,y')|}_{W=\ovl W+\widetilde W,\ \eqref{est of ovl W and D ovl W}\eqref{asmp of W tilde}}\Big)\\
    &\quad+(1+C\tau_0^{\frac{1}{30}})\int_{s_1}^s\underbrace{(\eta^{\mu(\gamma)}|F_W^{(\gamma)}|)(s',\Phi(s',s_1,y_0))}_{\eqref{est of FW^gamma (simplified)}}\underbrace{e^{\left(-\frac{3\gamma_1+\gamma_2-1}{2}+3\mu\right)(s-s')}}_{=1,\ \text{see }\eqref{def of mu}}\dif s'\\
    &\le(1+C\tau_0^{\frac{1}{30}})\max\Big(\frac{201}{200}\mathbbm{1}_{\gamma=0}+16\mathbbm{1}_{|\gamma|=1},1+\tau_0^{\frac{1}{18}}\Big)\\
    &\quad+(1+C\tau_0^{\frac{1}{30}})\int_{s_1}^so_{\tau_0}(1)\eta^{-\frac{1}{12}}(s',\Phi(s',s_1,y_0))\dif s'\\
    \overset{\eqref{integration along W trajectories}}&{\le}(1+C\tau_0^{\frac{1}{30}})\Big(\frac{201}{200}\mathbbm{1}_{\gamma=0}+16\mathbbm{1}_{|\gamma|=1}\Big)+o_{\tau_0}(1)\\
    &\le\begin{cases}
            \frac{201}{200}+o_{\tau_0}(1),&\gamma=(0,0),\\
            16+o_{\tau_0}(1),&|\gamma|=1.
        \end{cases}
\end{aligned}$$
This yields the desired bounds for $\eta^{\mu(\gamma)}|\partial^\gamma W(s,y)|$.
\end{proof}

\section{Top Order Energy Estimates}
\label{section: energy estimates}
In this section, we conclude our proof by improving the bootstrap assumption \eqref{asmp of Hk norm} on the top order $\dot H^k$ norm.

For a fixed integer $k\ge34$, define the energy $E_k(s)$ for the physical variable $P=(V_1,V_2,S)^T$ as
\begin{equation}
    E_k^2(s):=\sum_{|\gamma|=k}(\gamma_2!)^{-2}\left(\|\partial^\gamma V_1(s,\cdot)\|_{L^2}^2+\|\partial^\gamma V_2(s,\cdot)\|_{L^2}^2+\|\partial^\gamma S(s,\cdot)\|_{L^2}^2\right).\nonumber
\end{equation}
This energy is equivalent to the $\dot{H}^k$-norm of the unknowns as follows:
\begin{equation}
    \label{norm equivalence}
    E_k^2\le\|P\|_{\dot H^k}^2=\|\partial^\gamma V_1(s,\cdot)\|_{L^2}^2+\|\partial^\gamma V_2(s,\cdot)\|_{L^2}^2+\|\partial^\gamma S(s,\cdot)\|_{L^2}^2\le (k!)^2E_k^2.
\end{equation}
We will establish the following energy estimate in this section:
\begin{lemma}[$\dot H^k$ bound]
    \label{lem: Hk est of P}
    Let $k\ge 34$ be a fixed integer. Then, we have the following energy inequality:
    \begin{equation}
        E_k^2(s)\lesssim M^{3k-5}e^{-s}.\nonumber
    \end{equation}
    Here, the implicit constant is universal, and the large constant $M\gg1$ depends on $k$.
\end{lemma}

\subsection{Equation for $\widetilde P$}
First, we derive the governing equations for the physical variables $\widetilde P$. Recall that $\widetilde P=(V_1,V_2,S-\sigma_\infty)^T=P-(0,0,\sigma_\infty)^T$. To derive the equations of $\widetilde P$, we note that the third row of $D_P$ vanishes and $\partial P = \partial \widetilde{P}$. Together with \eqref{eqn of P, self-similar}, this leads us to the governing equation of $\widetilde{P}$ as below: 
\begin{equation}
\label{eqn of P tilde, self-similar}
\partial_{s}\widetilde P+\beta_{\tau} e^{-s} D_{P} \widetilde P+\left(\frac{3}{2} y_{1}+\beta_{\tau} e^{\frac{s}{2}} A_{P,\tilde u_1}\right) \partial_{1} \widetilde P+\left(\frac{1}{2} y_{2}+\beta_{\tau} e^{-\frac{s}{2}} A_{P, \tilde u_2}\right) \partial_{2} \widetilde 
P=\beta_{\tau} e^{-s} F_{P}.
\end{equation}

To facilitate the proof of energy estimates, we analyze in detail the structures of the damping, transport, and forcing terms. Here, Greek indices $\mu,\nu$ range over $1,2,3$, and Latin indices $i,j$ range over $1,2$. 
\begin{enumerate}
    \item From the definitions \eqref{def of DP} of $D_P$ and \eqref{def of hij} of $h_{ij}$, the damping term $D_P\widetilde P$ contains solely linear terms with the form $\Lambda \widetilde P_\mu$, where $\Lambda\in\{h_{11},h_{12}\}$.

    \item From the expression \eqref{def of AP tilde u1}\eqref{def of AP tilde u2} of $A_{P,\tilde u_1}$ and $A_{P,\tilde u_2}$, the transport terms $A_{P,\tilde u_j}\partial_j\widetilde P$ consist of 
    \begin{itemize}
        \item Linear terms $\Lambda D\widetilde P_\mu$ with $\Lambda\in\{g_1-\lambda g_2,g_2,\frac{1}{2}\dot\psi\tilde u_2^2,2\beta_3\sigma_\infty\varphi^{-1},2\beta_3\sigma_\infty\lambda\varphi^{-1}\}$;
        \item Quadratic terms $\Lambda \widetilde P_\mu D\widetilde P_\nu$ with $\Lambda\in\{2\beta_1\varphi^{-1},2\beta_1\lambda\varphi^{-1},2\beta_3\varphi^{-1},2\beta_3\lambda\varphi^{-1}\}$.
    \end{itemize}

    \item By \eqref{def of FP}, $F_P$ is composed of
    \begin{itemize}
        \item Linear terms $\Lambda \widetilde P_\mu$ with $\Lambda\in\{\frac{\varphi}{2r_0^2}u_jg_j,\frac{\beta_3}{r_0^2}\sigma_\infty u_i\}$;
        \item Quadratic terms $\Lambda \widetilde P_\mu \widetilde P_\nu$ with $\Lambda\in\{\frac{\beta_1}{r_0^2}u_i,\frac{\beta_3}{r_0^2}u_i\}$.
    \end{itemize}
\end{enumerate}
The auxiliary quantities mentioned above obey the bounds: $|\Lambda|\le\ln M$ and $\|\Lambda\|_{C_{\tilde u}^k}\le M^2$.

\subsection{Useful Lemmas}
Before establishing the energy estimates, we present several auxiliary lemmas.
\begin{lemma}
\label{lem: est of linear and quad. terms}
Suppose the auxiliary quantities $\Lambda$ satisfy the estimates $\|\Lambda\|_{L^{\infty}\left(\mathcal{X}(s)\right)} \le A_{0}$, $\|\Lambda\|_{C_{\tilde{u}}^{k}\left(\mathcal{X}(s)\right)}\le M^{2}$ (which yields $\left\|D_{y}^{j} \Lambda\right\|_{\left.L^{\infty}\mathcal{X}(s)\right)} \le M^2 e^{-\frac{j}{2} s}$ for $0 \le j \le k$). Then, for $1\le\mu,\nu\le3$, we have:
\begin{enumerate}
    \item $k$-th order estimates for linear terms:
    \begin{equation}
    \label{k th order linear est, single term}
        \|\Lambda \partial^{\gamma} \widetilde{P}_\mu\|_{L^{2}} \le A_{0}\|\partial^{\gamma} \widetilde{P}_\mu\|_{L^{2}},\quad \|D^{j} \Lambda D^{k-j} \widetilde{P}_\mu\|_{L^{2}}\overset{1\le j\le k}{\le} o_{\tau_0}(1)E_k,
    \end{equation}
    which further implies that
    \begin{equation}
    \label{k th order linear est, whole term}
        \|\partial^{\gamma}(\Lambda \widetilde{P}_\mu)\|_{L^{2}} \overset{|\gamma|=k}{\le}A_{0}\|\partial^\gamma\widetilde{P}_\mu\|_{L^2}+o_{\tau_{0}}(1)E_k.
    \end{equation}
    \item $(k+1)$-th order estimates for linear terms:
    \begin{equation}
    \label{k+1 th order linear est}
        \|D^{j} \Lambda D^{k+1-j} \widetilde{P}_\mu\|_{L^{2}} \overset{2 \le j \le k+1}{\lesssim} o_{\tau_{0}}(1)e^{-\frac{s}{2}}E_k.
    \end{equation}
    \item $k$-th order estimates for quadratic terms:
    \begin{equation}
    \label{k th order quadratic est, single term}
        \|D^{j}(\Lambda \widetilde{P}_\mu) D^{k-j} \widetilde{P}_\nu\|_{L^{2}} \overset{0 \leq j \leq k}{\lesssim_k} (A_0+o_{\tau_0}(1)) \ln M \cdot E_k.
    \end{equation}
    This implies that
    \begin{equation}
    \label{k th order quadratic est, whole term}
        \|D^{k}\left(\Lambda \widetilde{P}_\mu \widetilde{P}_\nu\right)\|_{L^{2}} \lesssim_k (A_0+o_{\tau_0}(1)) \ln M \cdot E_k.
    \end{equation}
    \item $(k+1)$-th order estimates for quadratic terms:
    \begin{equation}
    \label{k+1 th order quadratic est}
        \|D^{j}(\Lambda \widetilde{P}_\mu) D^{k+1-j} \widetilde{P}_\nu\|_{L^{2}}\overset{2 \leq j \leq k-1}{\lesssim_k}(1+A_0)^{2} M^{1+\frac{2}{3 k-8}} e^{-\frac{1}{2}\left(1+\frac{2}{3 k-8}\right) s}E_k^{1-\frac{2}{3 k-8}}.
    \end{equation}
\end{enumerate}
\end{lemma}
\begin{proof}
(1) The $j=0$ case is trivial. 
For $j \ge 1$, noticing the fact that $\operatorname{supp} \widetilde P\subset\mathcal{X}(s)$, we obtain that
$$
\|D^j \Lambda D^{k-j} \widetilde P_\mu\|_{L^2} \leq M^2 e^{-\frac{j}{2} s}\left(\tau_0^{\frac{1}{6}} e^{\frac{s}{2}}\right)^jE_k
\le M^2 \tau_0^{\frac{1}{6}} E_k=o_{\tau_0}(1)\|D^k \widetilde P_\mu\|_{L^2}.
$$
Thus, \eqref{k th order linear est, whole term} can be obtained by applying Leibniz rule. 

(2)Via direct computation, we arrive at:
$$
\begin{aligned}
    \|D^{j} \Lambda D^{k+1-j} \widetilde{P}_\mu\|_{L^{2}}&\lesssim M^2e^{-\frac{j}{2}s}\left(\tau_0^{\frac{1}{6}}e^{\frac{j}{2}s}\right)^{k-(k+1-j)}E_k\overset{2\le j\le k+1}{\lesssim} o_{\tau_{0}}(1)e^{-\frac{s}{2}}E_k.
\end{aligned}$$

(3)Employing H\"older inequality and Gagliardo-Nirenberg interpolation inequality, we deduce that
$$\begin{aligned}
    \|D^{j}(\Lambda \widetilde{P}_\mu) D^{k-j} \widetilde{P}_\nu\|_{L^{2}} \overset{0 \leq j \leq k}&{\lesssim} \|D^{j}(\Lambda \widetilde{P}_\mu)\|_{L^{{\frac{2k}{j}}}}\|D^{k-j} \widetilde{P}_\nu\|_{L^{{\frac{2k}{k-j}}}}\\
    &\lesssim_k\|\Lambda\widetilde P_\mu\|_{L^\infty}^{1-\frac{j}{k}}\|D^k(\Lambda\widetilde P_\mu)\|_{L^2}^{\frac{j}{k}}\|\widetilde P_\nu\|_{L^\infty}^{\frac{j}{k}}\|D^k\widetilde P_\nu\|_{L^2}^{1-\frac{j}{k}}\\
    &\lesssim(A_0\ln M)^{1-\frac{j}{k}}\left((A_0+o_{\tau_0}(1))E_k\right)^{\frac{j}{k}}(\ln M)^{\frac{j}{k}}E_k^{1-\frac{j}{k}}\\
    &\lesssim (A_0+o_{\tau_0}(1))\ln M\cdot E_k.
\end{aligned}$$
Inequality \eqref{k th order quadratic est, whole term} hence follows via using Leibniz rule.

(4)Leveraging Corollary \ref{cor: lower order est of P tilde} and Leibniz rule, one can deduce that $\|D^2(\Lambda \widetilde P_\mu)\|_{L^6}\lesssim(A_0+o_{\tau_0}(1))Me^{-\frac{s}{2}}$. For $2\le j\le k-1$, we define 
$$p_k(j):=\left(\frac{j-2}{3(k-2)}+\frac{1}{12}\right)^{-1},\quad a_k(j):=\frac{\frac{j-2}{2}+\frac{1}{6}-\frac{1}{p_k(j)}}{\frac{k-3}{2}+\frac{1}{6}}.$$
It follows directly that $\frac{1}{p_k(j)}+\frac{1}{p_k(k+1-j)}=\frac{1}{2}$ and $a_k(j)+a_k(k+1-j)=1-\frac{2}{3k-8}$. By Gagliardo-Nirenberg interpolation inequality, we have that
$$\begin{aligned}
    \|D^{j}(\Lambda \widetilde{P}_\mu) D^{k+1-j} \widetilde{P}_\nu\|_{L^{2}}&\le\|D^{j}(\Lambda \widetilde{P}_\mu)\|_{L^{p_k(j)}}\|D^{k+1-j} \widetilde{P}_\nu\|_{L^{p_k(k+1-j)}}\\
    &\lesssim_k\|D^2(\Lambda \widetilde{P}_\mu)\|_{L^6}^{1-a_k(j)}\|D^k(\Lambda \widetilde{P}_\mu)\|_{L^2}^{a_k(j)}\|D^2\widetilde P_\nu\|_{L^6}^{1-a_k(k+1-j)}\|D^k\widetilde P_\nu\|_{L^2}^{a_k(k+1-j)}\\
    \overset{\eqref{k th order linear est, whole term}}&{\lesssim} (A_0+o_{\tau_0}(1))^{1-a_k(j)}\left(Me^{-\frac{s}{2}}\right)^{1+\frac{2}{3k-8}}\left((1+A_0)E_k\right)^{1-\frac{2}{3k-8}}\\
    &\lesssim (1+A_0)^2M^{1+\frac{2}{3 k-8}} e^{-\frac{1}{2}\left(1+\frac{2}{3 k-8}\right) s}E_k^{1-\frac{2}{3 k-8}}.
\end{aligned}$$
\end{proof}

\begin{corollary}
Let $\Lambda,\Lambda_1,\Lambda_2,\Lambda_3$ be auxiliary quantities satisfying $|\Lambda_\mu|\le1+o_{\tau_0}(1)$ for some $1\le\mu\le3$, $|\Lambda_\nu|\le o_{\tau_0}(1)$ for $\nu\ne\mu$, and $\|\Lambda\|_{C_{\tilde u}^k}+\sum_\nu\|\Lambda_\nu\|_{C_{\tilde u}^k}\le M^2$. Then, for any multi-index $\gamma$ with order $k$, the following estimate holds:
\begin{align}
\label{est of linear comb.}
    \|\sum_\nu\partial^\gamma(\Lambda_\nu\widetilde P_\nu)+\partial^\gamma\Lambda\|_{L^2(\mathcal{X}(s))}\le\|\partial^\gamma\widetilde P_\mu\|_{L^2}+o_{\tau_0}(1)E_k+o_{\tau_0}(1)e^{-\frac{k-2}{2}s}.
\end{align}
\end{corollary}
\begin{proof}
Using Lemma \ref{lem: est of linear and quad. terms}, we derive that
\begin{align*}
    \|\sum_\nu\partial^\gamma(\Lambda_\nu\widetilde P_\nu)+\partial^\gamma\Lambda\|_{L^2(\mathcal{X}(s))}&\le\underbrace{\|\partial^\gamma(\Lambda_\nu\widetilde P_\mu)\|_{L^2}}_{\eqref{k th order linear est, whole term}}+\sum_{\nu\ne\mu}\underbrace{\|\partial^\gamma(\Lambda_\nu\widetilde P_\nu)\|_{L^2}}_{\eqref{k th order linear est, whole term}}+\|\partial^\gamma\Lambda\|_{L^\infty(\mathcal{X}(s))}|\mathcal{X}(s)|^{\frac12}\\
    &\le\|\partial^\gamma\widetilde P_\mu\|_{L^2}+o_{\tau_0}(1)E_k+CM^2\tau_0^{\frac13}e^{-\frac{k-2}{2}s}.
\end{align*}
This yields the desired inequality. 
\end{proof}

Referring to Lemma \ref{lem: est of linear and quad. terms}, along with the previous discussion on the structures of the damping, transport, and forcing terms, we conclude the following inequality:

\begin{corollary}
\label{cor: k-th order est of forcing and damping}For any multi-index $\gamma$ with $|\gamma|=k$, we have the following estimates for the forcing term and the damping term: 
    $$\|\partial^\gamma F_P\|_{L^2}+\|\partial^\gamma (D_P\widetilde P)\|_{L^2}\lesssim\ln M\cdot E_k.$$
\end{corollary}
\begin{proof}
    By Proposition \ref{prop: est for geo. coeff.}, for $\Lambda\in\{h_{11},h_{12},\frac{\varphi}{2r_0^2}u_jg_j,\frac{\beta_3}{r_0^2}\sigma_\infty u_i,\frac{\beta_1}{r_0^2}u_i,\frac{\beta_3}{r_0^2}u_i\}$, it holds that $|\Lambda|\le\ln M$ and $\|\Lambda\|_{C_{\tilde u}}^k\le M^2$. Hence, the desired result follows from Lemma \ref{lem: est of linear and quad. terms}. 
\end{proof}

For the transport terms, we prove the following estimates:

\begin{lemma}
    For $i=1,2$ and $|\gamma|=k$, there holds that
    \begin{equation}
    \label{k th order est of AP tilde ui}
        \|\partial^\gamma A_{P,\tilde u_i}\|_{L^2(\mathcal{X}(s))}:=\||\partial^\gamma A_{P,\tilde u_i}|_F\|_{L^2}\le 4\|\partial^\gamma\widetilde P\|_{L^2}+o_{\tau_0}(1)E_k+o_{\tau_0}(1)e^{-\frac{k-2}{2}s}.
    \end{equation}
    Here, $|(a_{\mu\nu})|_F:=(\sum_{\mu,\nu}a_{\mu\nu}^2)^{1/2}$ is the Frobenius norm of a matrix. Moreover, for $1\le j\le k-2$, the following estimate holds:
    \begin{equation}
    \label{k+1 th order est of AP tilde ui}
        \|D^{k-j}A_{P,\tilde u_i}D^{j+1}\widetilde P\|_{L^2}\lesssim M^{1+\frac{2}{3 k-8}} e^{-\frac{1}{2}\left(1+\frac{2}{3 k-8}\right) s}E_k^{1-\frac{2}{3 k-8}}+o_{\tau_0}(1)e^{-\frac{s}{2}}E_k.
    \end{equation}
\end{lemma}
\begin{proof}
    Employing \eqref{def of AP tilde u1}\eqref{def of G}\eqref{k th order linear est, whole term} and Proposition \ref{prop: est for geo. coeff.}, we deduce that:
    $$\begin{aligned}
        \|\partial^\gamma A_{P,\tilde u_1}\|_{L^2(\mathcal{X}(s))}^2&\le3\|\left[2\beta_1(\partial^\gamma(\varphi^{-1}V_1)-\partial^\gamma(\lambda\varphi^{-1}V_2))+\partial^\gamma G\right\|_{L^2}^2+8\beta_3^2(\|\partial^\gamma(\varphi^{-1}S)\|_{L^2}^2+\|\partial^\gamma(\lambda\varphi^{-1}S)\|_{L^2}^2)\\
        \overset{\eqref{est of linear comb.}}&{\le}12\beta_1^2\|\partial^\gamma V_1\|_{L^2}^2+8\beta_3^2\|\partial^\gamma S\|_{L^2}^2+o_{\tau_0}(1)E_k^2+o_{\tau_0}(1)e^{-(k-2)s}\\
        &\le 16\|\partial^\gamma \widetilde P\|_{L^2}^2+o_{\tau_0}(1)E_k^2+o_{\tau_0}(1)e^{-(k-2)s}.
    \end{aligned}$$
    Similarly, by \eqref{def of AP tilde u2}, we obtain that
    $$\begin{aligned}
        \|\partial^\gamma A_{P,\tilde u_2}\|_{L^2}^2&\le12\beta_1^2\|\partial^\gamma V_2\|_{L^2}^2+8\beta_3^2\|\partial^\gamma S\|_{L^2}^2+o_{\tau_0}(1)E_k^2+o_{\tau_0}(1)e^{-(k-2)s}\\
        &\le 16\|\partial^\gamma \widetilde P\|_{L^2}^2+o_{\tau_0}(1)E_k^2+o_{\tau_0}(1)e^{-(k-2)s}.
    \end{aligned}$$
    Then, \eqref{k th order est of AP tilde ui} follows from the elementary inequality $(\sum_i a_i)^{1/2}\le\sum_i a_i^{1/2}$ for finite sums. To prove \eqref{k+1 th order est of AP tilde ui}, note that $D^{k-j}A_{P,\tilde u_i}D^{j+1}\widetilde P$ consists of linear terms $D^{k-j}\Lambda D^{j+1}\widetilde P_\mu$ with $\Lambda\in\{g_1-\lambda g_2,g_2,\frac{1}{2}\dot\psi\tilde u_2^2,2\beta_3\sigma_\infty\varphi^{-1},2\beta_3\sigma_\infty\lambda\varphi^{-1}\}$, and quadratic terms $D^{k-j}(\Lambda \widetilde P_\mu)D^{j+1}\widetilde P_\nu$ with $\Lambda\in\{2\beta_1\varphi^{-1},2\beta_1\lambda\varphi^{-1},2\beta_3\varphi^{-1},2\beta_3\lambda\varphi^{-1}\}$. Via applying Proposition \ref{prop: est for geo. coeff.} and \eqref{k+1 th order linear est}\eqref{k+1 th order quadratic est}, these terms can be bounded by $M^{1+\frac{2}{3 k-8}} e^{-\frac{1}{2}\left(1+\frac{2}{3 k-8}\right) s}E_k^{1-\frac{2}{3 k-8}}+o_{\tau_0}(1)e^{-\frac{s}{2}}E_k$.
\end{proof}
\begin{lemma}
\label{lem: est of quadratic form}
    For any $a,b\in\mathbb{R}^3$ and Latin indices $i=1,2$, we have
    \begin{equation}
        |a^T\partial_1A_{P,\tilde u_i}b|\le(1+o_{\tau_0}(1))e^{-\frac{s}{2}}|a||b|, \quad |a^T\partial_2A_{P,\tilde u_i}b|\le(20+o_{\tau_0}(1))e^{-\frac{s}{2}}|a||b|.\nonumber
    \end{equation}
\end{lemma}
\begin{proof}
    According to the definition \eqref{def of AP tilde u1} of $A_{P,\tilde u_1}$, by using  bootstrap assumptions on $Z$ and $A$, and Proposition \ref{prop: est of DP}, we get
    $$\begin{aligned}
        |a^T\partial_j A_{P,\tilde u_1}b|&\le 2\beta_1\varphi^{-1}|\partial_j V_1||a^Tb|+2\beta_3\varphi^{-1}|\partial_jS||a_1b_3+a_3b_1|+o_{\tau_0}(1)e^{-\frac{s}{2}}|a||b|\\
        \overset{\eqref{est of geo. coef.}\eqref{est of DV,DS}}&{\le}\begin{cases}
            \beta_1e^{-\frac{s}{2}}|a^Tb|+\beta_3e^{-\frac{s}{2}}|a||b|+o_{\tau_0}(1)e^{-\frac{s}{2}}|a||b|,&j=1\\
            20\beta_1e^{-\frac{s}{2}}|a^Tb|+20\beta_3e^{-\frac{s}{2}}|a||b|+o_{\tau_0}(1)e^{-\frac{s}{2}}|a||b|,&j=2
        \end{cases}\\
        &\le\begin{cases}
            (1+o_{\tau_0}(1))e^{-\frac{s}{2}}|a||b|,&j=1\\
            (20+o_{\tau_0}(1))e^{-\frac{s}{2}}|a||b|,&j=2
        \end{cases}.
    \end{aligned}$$
    Similarly, by definition \eqref{def of AP tilde u2}, we derive that
    $$\begin{aligned}
        |a^T\partial_j A_{P,\tilde u_2}b|&\le 2\beta_3\varphi^{-1}|\partial_jS||a_2b_3+a_3b_2|+o_{\tau_0}(1)e^{-\frac{s}{2}}|a||b|\\
        \overset{\eqref{est of geo. coef.}\eqref{est of DV,DS}}&{\le}\begin{cases}
            \beta_3 e^{-\frac{s}{2}}|a||b|+o_{\tau_0}(1)e^{-\frac{s}{2}}|a||b|,&j=1\\
            20\beta_3 e^{-\frac{s}{2}}|a||b|+o_{\tau_0}(1)e^{-\frac{s}{2}}|a||b|,&j=2
        \end{cases}\\
        &\le\begin{cases}
            (1+o_{\tau_0}(1))e^{-\frac{s}{2}}|a||b|,&j=1\\
            (20+o_{\tau_0}(1))e^{-\frac{s}{2}}|a||b|,&j=2
        \end{cases}.
    \end{aligned}$$
    In both cases, we use the fact that $0\le\beta_1,\beta_3\le1$.
\end{proof}

\subsection{Proof of the Energy Estimate}
Now we are ready to prove Lemma \ref{lem: Hk est of P}. 

\begin{proof}[Proof of Lemma \ref{lem: Hk est of P}]
    By applying $\partial^\gamma$ to \eqref{eqn of P tilde, self-similar}, left-multiplying by $\partial^\gamma \widetilde P^T$, summing over $|\gamma|=k$ with the weight $(\gamma_2!)^{-2}$, and then integrating over $\mathbb{R}^2$, we arrive at
    $$
    \begin{aligned}
        \frac{1}{2}\frac{\dif}{\dif s}E_k^2&+\sum_{|\gamma|=k}(\gamma_2!)^{-2}\int\partial^\gamma\widetilde P^T\partial^\gamma\left[\left(\frac{3}{2} y_{1}+\beta_{\tau} e^{\frac{s}{2}} A_{P,\tilde u_1}\right) \partial_{1} \widetilde P+\left(\frac{1}{2} y_{2}+\beta_{\tau} e^{-\frac{s}{2}} A_{P,\tilde u_2}\right) \partial_{2} \widetilde P\right]\dif y\\
        &=\beta_\tau e^{-s}\sum_{|\gamma|=k}(\gamma_2!)^{-2}\int\partial^\gamma\widetilde P^T\partial^\gamma(F_P-D_P\widetilde P)\dif y.
    \end{aligned}$$
    \underline{\emph{Damping and forcing terms.}} From Corollary \ref{cor: k-th order est of forcing and damping} and H\"older inequality, the right-hand side is bounded by
    $$\begin{aligned}
        \quad\beta_\tau e^{-s}\sum_{|\gamma|=k}(\gamma_2!)^{-2}\int|\partial^\gamma\widetilde P^T\partial^\gamma(F_P-D_P\widetilde P)|\dif y&\lesssim e^{-s}E_k\sum_{|\gamma|=k}\underbrace{(\|\partial^\gamma F_P\|_{L_2}+\|\partial^\gamma(D_P\widetilde P)\|_{L^2})}_{\text{Corollary }\ref{cor: k-th order est of forcing and damping}}\\
        \overset{\eqref{hierarchy of constants}}&{\lesssim}o_{\tau_0}(1)E_k^2
    \end{aligned}$$
    
    \underline{\emph{Transport terms.}} Integration by part yields that
    $$
    \begin{aligned}
        &\sum_{|\gamma|=k}(\gamma_2!)^{-2}\int\partial^\gamma\widetilde P^T\partial^\gamma\left(\frac{3}{2} y_{1} \partial_{1} \widetilde P+\frac{1}{2} y_{2}\partial_{2} \widetilde P\right)\dif y
        =\left(\frac{k}{2}-1\right)E_k^2+\sum_{|\gamma|=k}\gamma_1(\gamma_2!)^{-2}\|\partial^\gamma\widetilde P\|_{L^2}^2,
    \end{aligned}$$
    and
    $$\begin{aligned}
        \int\partial^\gamma\widetilde P^T\partial^\gamma \left(A_{P,\tilde u_i} \partial_{i} \widetilde P\right)\dif y=&\left(\gamma_i-\frac{1}{2}\right)\int\partial^\gamma\widetilde P^T\partial_i A_{P,\tilde u_i}\partial^\gamma\widetilde P\dif y+\sum_{j\ne i}\gamma_j \int\partial^\gamma\widetilde P^T\partial_j A_{P,\tilde u_i}\partial^{\gamma-e_j+e_i}\widetilde P\dif y\\
        &+\sum_{\substack{1\le|\beta|\le k-2\\ \beta\le\gamma}}\binom{\gamma}{\beta}\int\partial^\gamma\widetilde P^T\partial^{\gamma-\beta} A_{P,\tilde u_i}\partial_i\partial^\beta\widetilde P\dif y+\int\partial^\gamma\widetilde P^T\partial^\gamma A_{P,\tilde u_i}\partial_i\widetilde P\dif y\\
        =&I_1+I_2+I_3+I_4.
    \end{aligned}$$
    By Lemma \ref{lem: est of quadratic form}, we have
    $$|I_1|\le
    \begin{cases}
        (\gamma_i+\frac{1}{2})e^{-\frac{s}{2}}\|\partial^\gamma\widetilde P\|_{L^2}^2+o_{\tau_0}(1)e^{-\frac{s}{2}}E_k^2,&i=1\\
        o_{\tau_0}(1)E_k^2,&i=2
    \end{cases},$$
    and
    $$|I_2|\le
    \begin{cases}
        e^{-\frac{s}{2}}\sum_{j\ne i}\gamma_j\|\partial^\gamma\widetilde P\|_{L^2}\|\partial^{\gamma-e_j+e_i}\widetilde P\|_{L^2}+o_{\tau_0}(1)e^{-\frac{s}{2}}E_k^2, &i=1\\
        o_{\tau_0}(1)E_k^2, &i=2
    \end{cases}.$$
    By using \eqref{k+1 th order est of AP tilde ui} and the notation $\delta:=\frac{2}{3k-8}$, we derive that 
    $$\begin{aligned}
        |I_3|&\lesssim M^{1+\delta} e^{-\frac{1}{2}\left(1+\delta\right) s}E_k^{2-\delta}+o_{\tau_0}(1)e^{-\frac{s}{2}}E_k^2\\
        &\lesssim\left(\frac{1}{\ln M}e^{-\frac{2-\delta}{4}s}E^{2-\delta}\right)^{\frac{2}{2-\delta}}+\left(e^{-\frac{3}{4}\delta s}M^{1+\delta}\ln M\right)^{\frac{2}{\delta}}+o_{\tau_0}(1)e^{-\frac{s}{2}}E_k^2\\
        &\lesssim o_{\tau_0,M}(1)e^{-\frac{s}{2}}E_k^2+(\ln M)^{3k-8}M^{3k-6}e^{-\frac{3}{2}s}.
    \end{aligned}$$
    Here, $o_{\tau_0,M}(1)$ denotes a positive quantity that tends to zero as $\tau_0\rightarrow0$ and $M\rightarrow\infty$. Moreover, it follows from \eqref{est of DjP} and \eqref{k th order est of AP tilde ui} that
    $$\begin{aligned}
        |I_4|&\le\||\partial^\gamma\widetilde P^T||\partial^\gamma A_{P,\tilde u_i}|_F|\partial_i\widetilde P|\|_{L^1}\\
        &\le\|\partial^\gamma\widetilde P\|_{L^2}\|\partial^\gamma A_{P,\tilde u_i}\|_{L^2}\|\partial_i\widetilde P\|_{L^\infty}\\
        &\le \|\partial^\gamma\widetilde P\|_{L^2}(4\|\partial^\gamma\widetilde P\|_{L^2}+o_{\tau_0}(1)E_k+o_{\tau_0}(1)e^{-\frac{k-2}{2}s})\|\partial_i\widetilde P\|_{L^\infty}\\
        &\le\begin{cases}
        4e^{-\frac{s}{2}}\|\partial^\gamma\widetilde P\|_{L^2}^2+o_{\tau_0}(1)e^{-\frac{s}{2}}E_k^2+o_{\tau_0}(1)e^{-(k-2)s},&i=1\\
        o_{\tau_0}(1)E_k^2+o_{\tau_0}(1)e^{-(k-2)s},&i=2
    \end{cases}.
    \end{aligned}$$ 
Recall that $\beta_\tau=1+o_{\tau_0}(1)$ (see \eqref{est of 1-beta_tau}), we obtain
    $$\begin{aligned}
        &\beta_{\tau}e^{\frac{s}{2}}\sum_{|\gamma|=k}(\gamma_2!)^{-2}\int\partial^\gamma\widetilde P^T\partial^\gamma\left( A_{P,\tilde u_1}\partial_{1} \widetilde P\right)\dif y\\
        \ge&-\sum_{|\gamma|=k}(\gamma_2!)^{-2}\gamma_1\|\partial^\gamma\widetilde P\|_{L^2}^2-\sum_{|\gamma|=k}(\gamma_2!)^{-2}\gamma_2\|\partial^\gamma\widetilde P\|_{L^2}\|\partial^{\gamma-e_2+e_1}\widetilde P\|_{L^2}\\
        &-\frac{9}{2}E_k^2-o_{\tau_0,M}(1)E_k^2-C(\ln M)^{3k-8}M^{3k-6}e^{-s}\\
        \overset{\text{AM-GM}}&{\ge}-\left(\frac{13}{2}+o_{\tau_0,M}(1)\right)E_k^2-\sum_{|\gamma|=k}(\gamma_2!)^{-2}\gamma_1\|\partial^\gamma\widetilde P\|_{L^2}^2-M^{3k-5}e^{-s},
    \end{aligned}$$
    and
    $$\beta_{\tau}e^{-\frac{s}{2}}\sum_{|\gamma|=k}(\gamma_2!)^{-2}\int\partial^\gamma\widetilde P^T\partial^\gamma\left( A_{P,\tilde u_2}\partial_{1} \widetilde P\right)\dif y\ge-o_{\tau_0,M}(1)E_k^2-M^{3k-5}e^{-s}.$$
    Invoking all the bounds on transport terms, we arrive at
    $$\begin{aligned}
        &\sum_{|\gamma|=k}(\gamma_2!)^{-2}\int\partial^\gamma\widetilde P^T\partial^\gamma\left[\left(\frac{3}{2} y_{1}+\beta_{\tau} e^{\frac{s}{2}} A_{P,\tilde u_1}\right) \partial_{1} \widetilde P+\left(\frac{1}{2} y_{2}+\beta_{\tau} e^{-\frac{s}{2}} A_{P,\tilde u_2}\right) \partial_{2} \widetilde P\right]\dif y\\
        \ge&\left(\frac{k-15}{2}-o_{\tau_0,M}(1)\right)E_k^2-CM^{3k-5}e^{-s}.
    \end{aligned}$$

    Finally, gathering all terms in the energy estimate, we obtain
    $$\frac{\dif}{\dif s}E_k^2 + \left(k - 15 - o_{\tau_0,M}(1)\right)E_k^2 \leq M^{3k-4}e^{-s}.
    $$
    Since $k - 15 - o_{\tau_0,M}(1) \geq 1$, via using Gr\"onwall's inequality, we have
    $$E_k^2(s) \leq e^{-(s-s_0)}E_k^2(s_0) + M^{3k-4}e^{-s}.
    $$
    Plugging in the initial condition \eqref{init of Hk norm}, we conclude that $E_k^2(s_0) \lesssim M\tau_0$ and $E_k^2(s) \lesssim M^{3k-4}e^{-s}$.
\end{proof}

\begin{proof}[Closure of bootstrap assumption \eqref{asmp of Hk norm}]
    Notice that $\|D^k\widetilde P\|_{L^2}\lesssim E_k(s)\lesssim M^{\frac{3k-4}{2}}e^{-\frac{s}{2}}$, it suffices to show
    $$e^{-\frac{s}{2}}\left\|D^{k} W\right\|_{L^{2}}+\left\|D^{k} Z\right\|_{L^{2}}+\left\|D^{k} A\right\|_{L^{2}}\lesssim\|D^k\widetilde P\|_{L^2}.$$
    By \eqref{R expressed by P}, this is an immediately follows from \eqref{k th order linear est, whole term}. Therefore, we obtain
    $$e^{-\frac{s}{2}}\left\|D^{k} W\right\|_{L^{2}}+\left\|D^{k} Z\right\|_{L^{2}}+\left\|D^{k} A\right\|_{L^{2}}\lesssim M^{\frac{3k}{2}-2}e^{-\frac{s}{2}}\le\frac{1}{2}M^{\frac{3k}{2}}e^{-\frac{s}{2}}.$$
This completes the proof of the energy estimate.
\end{proof}

\appendix

\section{Shock Formation for Equivariant Euler Equations} \label{append}

In the appendix\footnote{The proofs of this appendix are based on Fulin Qi's UROPS college research project.}, we present a proof of shock formation for the compressible Euler equations on $\mathbb{S}^2$ under equivariant symmetry, that is, under the assumptions $v_\phi = 0$ and $\partial_\phi \mathcal W = 0$ for $\mathcal W \in \{v_\theta, v_\phi, \rho\}$. From a physical perspective, the equations considered in the appendix describe axisymmetric flows without swirl, a regime of significant physical relevance. While the methods parallel to those in the main text, the coordinate choice and the modulation variables differ herehere. The appendix also provides a streamlined illustration of the core ideas developed in the paper. 

\subsection{Preliminaries}\label{section: preliminaries}
In this section, we derive the Euler equations on $\mathbb{S}^2$ under equivariant symmetry, introduce the co-moving coordinates, define the Riemann variables, and perform the self-similar transformations.

\subsubsection{The Equivariant Euler's Equations}
We begin with the derivation of equivariant Euler's equations on $\mathbb{S}^2$. To capture the equivariant motion of fluids, we adopt the spherical coordinates on $\mathbb S^2$, instead of the stereographic coordinates used in the previous argument. Setting $\mathcal R:=(0,2\pi)\times\left(-\frac\pi2,\frac\pi2\right)$, the standard spherical coordinates are given by
\begin{align}
    (\phi,\theta)\in\mathcal R\mapsto (\cos\phi\cos\theta,\sin\phi\cos\theta, \sin\theta)\in\mathbb S^2.\nonumber
\end{align}
The isentropic compressible Euler equations can then be written in spherical coordinates as
\begin{align}\label{eq: EV_Euler equations under spherical coordinates}
	\begin{cases}
        \displaystyle\partial_t\rho+\frac1{\cos\theta}(\partial_\phi\rho)(v_\phi)+(\partial_\theta\rho)(v_\theta)+\frac\rho{\cos\theta}\partial_\phi v_\phi+\rho \partial_\theta v_\theta-\rho\tan\theta v_\theta=0,\\
		\partial_tv_\phi+\left(\dfrac{v_\phi\partial_\phi v_\phi}{\cos\theta}-v_\phi v_\theta\tan\theta+ v_\theta\partial_\theta v_\phi\right)+\dfrac{\rho^{\gamma-2}}{\cos\theta}\partial_\phi\rho=0,\\
		\partial_tv_\theta+\left(\dfrac{v_\phi\partial_\phi v_\theta}{\cos\theta}+v_\phi^2\tan\theta+v_\theta\partial_\theta v_\theta\right)+\rho^{\gamma-2}\partial_\theta\rho=0,
	\end{cases}
\end{align}
where $v=v_\theta\textbf{e}_\theta+v_\phi\textbf{e}_{\phi}$, with $\textbf{e}_\theta$ and $\textbf{e}_\phi$ being unit vectors defined by $\textbf{e}_\theta:=\partial_\theta/\abs{\partial_\theta}$ and $\textbf{e}_\phi:=\partial_\phi/\abs{\partial_\phi}$. 

Assuming that the flow is equivariant (i.e., $v_\phi = 0$ and $\partial_\phi \mathcal W = 0$ for $\mathcal W \in \{v_\theta, v_\phi, \rho\}$), we abuse the notation and write $v_\theta$ as $v$ for simplicity. Denote the rescaled sound speed as $\sigma = \frac{1}{\alpha} \rho^\alpha$ with $\alpha = \frac{\gamma - 1}{2} > 0$. The system \eqref{eq: EV_Euler equations under spherical coordinates} is reduced to
\begin{align}\label{eq: EV_Euler equations under EV}
    \begin{cases}
        \partial_t \sigma + v \partial_\theta \sigma + \alpha \sigma \partial_\theta v = \alpha \sigma  v\tan\theta, \\
        \partial_t v + v \partial_\theta v + \alpha \sigma \partial_\theta \sigma = 0.
    \end{cases}
\end{align}
Although the spherical coordinates do not cover the line $\{\phi=0\}$ on $\mb S^2$, this poses no issue under the equivariant assumptions since all unknowns are independent of $\phi$.

\subsubsection{Co-moving Coordinates}
In the following, we introduce the co-moving coordinates. To keep track of the shock, we first define three time-dependent modulation variables. Namely, we use $\xi\in\mb R$ to trace the location of the shock formation, $\tau\in\mathbb R$ to track the slope of the Riemann invariant $w$ (defined in \eqref{eq: EV_Riemann invariant}) at $x=\xi(t)$, and $\kappa\in\mathbb R$ to record the value of $w$ at $x=\xi(\tilde t)$. Since the slope of $w$ is approximately $\frac{c}{T^*-t}$ for some constant $c\in\mathbb R$, $\tau$ can also be understood as a variable tracking the time of shock formation. With these modulation variables, we can define the co-moving coordinates as
\begin{align}
	\tilde\theta:=\theta-\xi(\tilde t), \quad\tilde t:=\dfrac{1+\alpha}2t.\nonumber
\end{align}
Under the co-moving coordinates, the governing equations of $(v, \sigma)$ become
\begin{align}\label{eq: EV_Euler equations under coordinates adapted to the shock}
	\begin{cases}
		\pt v+\left( 2\beta_1v-\pt\xi\right) \partial_{\tilde\theta}v+2\alpha\beta_1\sigma\partial_{\tilde\theta}\sigma=0,\\
		\pt \sigma+\left( 2\beta_1v-\pt\xi\right) \partial_{\tilde\theta}\sigma+2\alpha\beta_1\sigma\partial_{\tilde\theta}v=2\alpha\beta_1\sigma v\tan\left(\tilde\theta+\xi\right),\\
	\end{cases}	
\end{align}
where the constants $\beta_1,\beta_2$ and $\beta_3$ are defined in \eqref{beta param}.

\subsubsection{Riemann variables}
We define Riemann variables by
\begin{align}\label{eq: EV_Riemann invariant}
	w:=v_\theta+\sigma,\quad z:=v_\theta-\sigma.
\end{align}
The system \eqref{eq: EV_Euler equations under coordinates adapted to the shock} can then be rewritten in $(w, z)$ as
\begin{align}\label{eq: EV_Euler in w and z}
	\begin{cases}
		\pt w+\left(w+\beta_2 z-\pt\xi\right)\partial_{\tilde\theta} w=\frac12\beta_3(w^2-z^2)\tan\left(\tilde\theta+\xi\right),\\
		\pt z+\left(\beta_2 w+ z-\pt\xi\right)\partial_{\tilde\theta}z=\frac12\beta_3(z^2-w^2)\tan\left(\tilde\theta+\xi\right).
	\end{cases}
\end{align}

\subsubsection{Self-similar transformations}
We now introduce the self-similar coordinates
\begin{align}
	y:=\tilde\theta e^{\frac32s},\quad s:=-\ln\left(\tau(\tilde t)-\tilde t\right)\nonumber,
\end{align}
together with the corresponding self-similar variables
\begin{align}
	W(s,y):=e^{\frac s2}\left(w(\tilde t,\tilde\theta)-\kappa(\tilde t)\right), \quad Z(s,y):=z(\tilde t,\tilde\theta).\nonumber
\end{align}
Letting $\beta_\tau=(1-\partial_{\tilde t}\tau)^{-1}$, we can then reformulate the system \eqref{eq: EV_Euler in w and z} as
\begin{align}\label{eq: EV_Euler in W and Z}
    \begin{cases}
	    \left(\partial_s-\frac12\right)W+\left(\frac32y+g_W\right)\partial_yW=F_W-\beta_\tau e^{-\frac s2}\partial_{\tilde t}\kappa,\\
        \partial_sZ+\left(\frac32y+g_Z\right)\partial_yZ=F_Z.
	\end{cases}
\end{align}
The transport and forcing terms in the above system \eqref{eq: EV_Euler in W and Z} are 
\begin{align}
    &\begin{cases}\label{eq: EV_g_W and g_Z}
        g_W=\beta_\tau W+G_W=\beta_\tau W+\beta_\tau e^{\frac s2}\left(\kappa+\beta_2Z-\partial_{\tilde t}\xi\right),\\
	g_Z=\beta_2\beta_\tau W+G_Z=\beta_2\beta_\tau W+\beta_\tau e^{\frac s2}\left(\beta_2\kappa+Z-\partial_{\tilde t}\xi\right),
    \end{cases}\\
    &\begin{cases}\label{eq: EV_F_W and F_Z}
        F_W=\frac{\beta_\tau\beta_3}2e^{-\frac s2}\tan(ye^{-\frac{3}2s}+\xi)\left(\left(e^{-\frac s2}W+\kappa\right)^2-Z^2\right),\\
	F_Z=\frac{\beta_\tau\beta_3}2 e^{-s}\tan(ye^{-\frac32s}+\xi)(Z^2-(e^{-\frac s2}W+\kappa)^2).
    \end{cases}	
\end{align}

\subsubsection{1D self-similar Burgers profile}
The profile we use is the one-dimensional self-similar Burgers profile introduced in Section \ref{subsection: blow-up profile}. For convenient reference, we restate some of its key properties here.

First, $\overline W$ satisfies the following estimates:
\begin{align}
	\begin{cases}\label{eq: EV_W_bar decay}
		|\overline W|\leq |y|^\frac13\leq \langle y\rangle^\frac13, &|\partial_y\overline W|\leq\langle y\rangle^{-\frac23},\quad
		|\partial_y^2\overline W|\leq\langle y\rangle^{-\frac53},\\
        |\partial_y^3\overline W|\leq6\langle y\rangle^{-\frac83},&|\partial_y^4\overline W|\leq35\langle y\rangle^{-\frac{11}3},
	\end{cases}
\end{align}
where $\langle y\rangle =\sqrt{1+y^2}$. At the origin $y=0$, we also have
\begin{align}\label{eq: EV_W_bar initial}
	\begin{cases}
		\overline W(0)=0, &\partial_y\overline W(0)=-1,\\
		\partial_y^2 \overline W(0)=0,&\partial_y^3\overline W(0)=6.\\
	\end{cases}
\end{align}

\subsubsection{Evolution of unknowns}
For $k\geq 1$, the evolution equations for $(\partial_y^k W, \partial_y^kZ)$ are given by
\begin{align}
    &\begin{cases}\label{eq: EV_dkW and dkZ}
        \left(\partial_s+\frac{3k-1}2+(1+\mathbbm 1_{k\geq 2}k)\beta_\tau \partial_y W\right)\partial_y^kW+(\beta_\tau W+\frac32 y+G_W)\partial_y^{k+1}W=F^{(k)}_W,\\
		\left(\partial_s+\frac{3k}2+\beta_2\beta_\tau k\partial_y Z\right)\partial_y^kZ+\left(\frac32y+\beta_2\beta_\tau W+G_Z\right)\partial_y^{k+1}Z=F^{(k)}_Z,
    \end{cases}
    \end{align}
where $\mathbbm 1$ is the indicator function, and the corresponding forcing terms take the following form:
    \begin{align}
    &\begin{cases}\label{eq: EV_F^k_W and F^k_Z}
        F^{(k)}_W=\partial_y^kF_W-\mathbbm 1_{k\geq 3}\beta_\tau\left(\sum_{j=1}^{k-2}\binom kj\partial_y^{k-j}W\partial_y^{j+1}W\right)-\sum_{j=0}^{k-1}\binom kj\partial_y^{k-j}G_W\partial_y^{j+1}W,\\
		F^{(k)}_Z=\partial_y^kF_Z-\mathbbm 1_{k\geq 2}\beta_2\beta_\tau\left(\sum_{j=0}^{k-2}\binom kj\partial_y^{k-j}W\partial_y^{j+1}Z\right)-\sum_{j=0}^{k-1}\binom kj\partial_y^{k-j}G_Z\partial_y^{j+1}Z,
    \end{cases}\\
    &\begin{cases}\label{eq: EV_dFkW and dFkZ}
        \partial_y^k F_W=\frac{\beta_\tau\beta_3}2e^{-\frac s2}\sum_{j=0}^{k}\binom kj\partial_y^{k-j}\tan(ye^{-\frac{3}2s}+\xi)\left(\partial_y^j\left(e^{-\frac s2}W+\kappa\right)^2-\partial_y^jZ^2\right),\\
		\partial_y^k F_Z=\frac{\beta_\tau\beta_3}2e^{-s}\sum_{j=0}^{k}\binom kj\partial_y^{k-j}\tan(ye^{-\frac{3}2s}+\xi)\left(\partial_y^jZ^2-\partial_y^j\left(e^{-\frac s2}W+\kappa\right)^2\right).
    \end{cases}
\end{align}

\subsubsection{Evolution of $\widetilde W$}

We define $\widetilde W:=W-\overline W$, and one can check that $\widetilde W$ satisfies
\begin{align}
    \label{eq: EV_tilde W}\begin{cases}
        \left(\partial_s-\frac12+\beta_\tau\partial_y\overline W\right)\widetilde W+\left(\frac32 y+g_W\right)\partial_y\widetilde W=F_{\widetilde W},\\
	   F_{\widetilde W}=F_W-G_W\partial_y\overline W+(1-\beta_\tau)\overline W\partial_y\overline W-\beta_\tau e^{-\frac s2}\partial_{\tilde t}\kappa.
    \end{cases}
\end{align}
For $k\geq 1$, the evolution equation of $\partial_y^k\widetilde W$ can be deduced from \eqref{eq: EV_dkW and dkZ}, \eqref{eq: EV_F^k_W and F^k_Z}, and \eqref{eq: EV_dFkW and dFkZ}:
\begin{align}\label{eq: EV_dyW_tilde_k}
    \left(\partial_s+\frac{3k-1}2+\beta_\tau\left(\partial_y\overline W+k\partial_yW\right)\right)\partial_y^k\widetilde W+\left(\frac32y+\beta_\tau W+G_W\right)\partial_y^{k+1}\widetilde W=F^{(k)}_{\widetilde W},
\end{align}
where the forcing term is defined by
\begin{align}\label{eq: EV_F^k tilde W}
 	\begin{cases}
        F^{(k)}_{\widetilde W}=\partial_y^kF_{\widetilde W}-\beta_\tau\Sigma_1-\beta_\tau \mathbbm 1_{k\geq2}\Sigma_2,\\
	   \Sigma_1=\sum_{j=0}^{k-1}\binom kj\left(\partial_y^{k-j+1}\overline W\partial_y^j\widetilde W+\partial^{k-j}_yG_W\partial_y^{j+1}\widetilde W\right),\\
	   \Sigma_2=\sum_{j=0}^{k-2}\binom{k}{j}\partial^{k-j}_y W\partial_y^{j+1}\widetilde W,\\
	   \partial_y^k F_{\widetilde W}=\partial_y^kF_W+\sum_{j=0}^k\binom kj\left((1-\beta_\tau)\partial_y^{k-j}\overline W-\partial_y^{k-j}G_W\right)\partial_y^{j+1}\overline W.
    \end{cases}	
\end{align}

\subsection{Main results}\label{section: main result}

In this section, we present the the main shock formation result for the 2D equivariant compressible Euler equations on $\mathbb S^2$. 

\subsubsection{Initial conditions}

We begin by specifying the initial state. We take the initial time $t_0$ to be $0$, and we set the initial values of modulation variables as
\begin{align}
    \label{eq: EV_initial conditions for modulation variables}\kappa(0)=\kappa_0=\sigma_\infty>0,\quad \xi(0)=\xi_0\in\left[\frac{\pi}{16},\frac{\pi}8\right],\quad \tau(0)=\tau_0>0.
\end{align}
Furthermore, we require $\sigma_\infty$ to satisfy
\begin{align}\label{eq: EV_requirement of kappa_0}
	\kappa_0=\sigma_\infty>\frac{\xi_0^{\frac13}}{2\beta_3}.	
\end{align}
\begin{remark}
Note that the choice of $\xi_0$ here is made for technical convenience. With modified bootstrap arguments, this interval can be further enlarged.
\end{remark}

We now prescribe the initial data of $(v,\sigma)$ in the $\tilde\theta$-coordinate as below
\begin{align}
	v(0,\tilde\theta):=v_0( \tilde\theta),\quad \sigma(0,\tilde\theta):=\sigma_0(\tilde\theta),\nonumber
\end{align}
where $v_0,\sigma_0$ are to be carefully chosen such that the Riemann variables satisfy the conditions stated in this section.

The initial data of the Riemann variables are denoted as
\begin{align}
	w_0(\tilde \theta):= v_0(\tilde \theta)+ \sigma_0(\tilde \theta),\quad z_0(\tilde \theta):= v_0(\tilde \theta)-\sigma_0(\tilde \theta).\nonumber
\end{align}
Furthermore, we assume at the initial time that 
\begin{align}
	\Supp_{\tilde \theta}(w_0-\sigma_\infty,z_0+\sigma_\infty)\subseteq \left(-\frac1{10}\xi_0,\frac{1}{10}\xi_0\right),\nonumber
\end{align}
which corresponds to 
\begin{align}
	\Supp_{\theta}(\tilde w_0-\sigma_\infty,z+\sigma_\infty)\subseteq \left(\frac9{10}\xi_0,\frac{11}{10}\xi_0\right).\nonumber
\end{align}

The self-similar Riemann variables would initially verify that:
\begin{align}
	&\begin{cases}\tag{IB0-$W$}\label{eq: EV_IB0 on W}
		|W(s_0,y)|\leq (1+\frac1{200}\tau_0^\frac1{23})\langle y\rangle^{\frac13}, &
		|\partial_y W(s_0,y)|\leq 12\langle y\rangle^{-\frac23},\\
		|\partial_y^2 W(s_0,y)|\leq \frac1{10}M^\frac18\langle y\rangle^{-\frac23},&
		|\partial^3_y W(s_0,y)|\leq \frac14M^\frac13,\\
		|\partial_y^4 W(s_0,y)|\leq \frac14M^{\frac9{10}},
	\end{cases}\\
	&\begin{cases}\tag{IB0-$\widetilde W$}\label{eq: EV_IB0 on W_tilde}
		|\widetilde W(s_0,y)|\mathbbm{1}_{|y|\leq L}\leq \frac14\tau_0^\frac{1}{2}\langle y\rangle^{\frac13},\\
		|\partial_y \widetilde W(s_0,y)|\mathbbm{1}_{|y|\leq L}\leq \frac14\tau_0^\frac1{3}\langle y\rangle^{-\frac23},\\
		 |\partial_y^2\widetilde W(s_0,y)|\mathbbm{1}_{|y|\leq L}\leq \frac14\tau_0^\frac1{4}\langle y\rangle^{-\frac23},\\
         \abs{\partial_y^3\widetilde W^0(s_0)}\leq\frac1{10}\tau_0^\frac45,\\
         \abs{\partial_y^k\widetilde W(s_0,y)}\mathbbm{1}_{|y|\leq l}\leq\frac{1}{10}M^2\tau_0^{\frac12}|y|^{4-k}+\frac{1}{10}\tau_0^{\frac35}\abs{y}^{3-k}\mathbbm{1}_{k\leq3},
	\end{cases}\\
	\tag{IB0-$Z$}\label{eq: EV_IB0 on Z}
	&\begin{cases}
		|Z(s_0,y)+\sigma_\infty|\leq \frac14M\tau_0,&
		|\partial_y Z(s_0,y)|\leq \frac14M\tau_0^\frac32,\\
		|\partial_y^2 Z(s_0,y)|\leq \frac14M^{\frac43}\tau_0^\frac32,&
		|\partial_y^3 Z(s_0,y)|\leq \frac14M^6\tau_0^\frac32,\\
		|\partial_y^4 Z(s_0,y)|\leq \frac14M^7\tau_0^\frac32.
	\end{cases}
\end{align}	
The choice of constants $M$ and $\tau_0$ follows Section \ref{hierarchy of constants}. In principle, $M$ is a large constant that dominates all universal constants and initial conditions on physical variables, while $\tau_0$ is set to be small enough so that its reciprocal can suppress the largeness resulting from $M$.

\subsubsection{Main theorem}
We now state the main theorem.
\begin{theorem}\label{theorem: EV_main theorem}
	Consider the equivariant compressible Euler equations \eqref{eq: EV_Euler equations under coordinates adapted to the shock}. Assume that the following conditions hold:
	\begin{enumerate}
		\item For any $\kappa_0$ satisfying \eqref{eq: EV_requirement of kappa_0}, the initial data of $(\widetilde v_\theta,\widetilde \sigma)$ is smooth, and the corresponding self-similar variables satisfy \eqref{eq: EV_IB0 on Z}, \eqref{eq: EV_IB0 on W}, and \eqref{eq: EV_IB0 on W_tilde};
		\item The parameter $\tau_0$ satisfies $0<\tau_0\leq \epsilon$ for a sufficiently small $\epsilon\in(0,1)$.
	\end{enumerate}
	Then, the corresponding solution $(\widetilde v_\theta,\tilde\sigma)$ to \eqref{eq: EV_Euler equations under coordinates adapted to the shock} blows up at a finite time $\widetilde T^*<+\infty$. Furthermore, the solution exhibits the following properties:
	\begin{enumerate}
		\item \emph{Blow-up speed:} The following estimates for $(\partial_{\tilde \theta} \widetilde v_{\theta}, \partial_{\tilde \theta} \widetilde \sigma)$ holds:
		\begin{align}\label{eq: EV_blow up speed}
			\frac{c}{\widetilde T^*-\tilde t}\leq\max\left(\norm{\partial_{\tilde\theta}\widetilde v_{\theta}}_\infty, \norm{\partial_{\tilde\theta}\widetilde\sigma}_\infty\right)\leq \frac{C}{\widetilde T^*-\tilde t}.
		\end{align}
		\item \emph{Blow-up time:} The maximal lifespan $\widetilde T^*$ satisfies
			\begin{align}\label{eq: EV_blow up time}
			|\tau_0-\widetilde T^*|\leq 2M\tau_0^2.
			\end{align} 
		\item \emph{Blow-up location:} For $\delta\in(0,1)$, there exists a constant $C(\delta)$ such that
		\begin{align}\label{eq: EV_blow-up location}
				\left\|\partial_{\tilde \theta} \widetilde v_\theta(\tilde t,\cdot)\right\|_{L^{\infty}\left(B_\delta(\xi(\tilde t))^c\right)}+\left\|\partial_{\tilde \theta} \widetilde \sigma(\tilde t,\cdot)\right\|_{L^{\infty}\left(B_\delta(\xi(\tilde t))^c\right)} \leq M+2\delta^{-\frac23},
		\end{align}
		while the gradient is unbounded along $\xi(\tilde t)$ as shown in (1). Here, $B_\delta(\xi(\tilde t))$ is the ball centered at $\xi(\tilde t)$ with radius $\delta$, and $B_\delta(\xi(\tilde t))^c$ denotes its set complement. Moreover, the limit $\lim_{\tilde t\rightarrow \widetilde T^*}\xi(\tilde t)=\xi_*$ exists.
		\item \emph{1/3-H\"older continuity:} The solution admits a uniform-in-time $C^\frac13$ bound. More precisely, we have
		\begin{align}\label{eq: EV_blow-up solution holder continuity}
			(\tilde v_\theta,\tilde \sigma)\in L^\infty_t([0,T^*))\cap C_{\tilde\theta}^\frac13.
		\end{align}
	\end{enumerate}
\end{theorem}

\subsection{Bootstrap argument}\label{section: bootstrap argument}

In this section, we set up the bootstrap argument, by which the global well-posedness in self-similar ansatz can be established.

\subsubsection{Bootstrap assumptions}

We first state the bootstrap assumptions. For the modulation variables, we assume that 
\begin{align}\tag{BA-M}\label{eq: EV_BA on modulation variables}
    \begin{cases}
		|\kappa-\kappa_0|\leq M\tau_0,&\left|\partial_{\tilde t}\kappa\right|\leq M,\\
		\left|\xi-\xi_0-2\beta_3\kappa_0\tilde t\right|\leq M^2\tau_0^2, &\left|\partial_{\tilde t}\xi-2\beta_3\kappa_0\right|\leq M^2\tau_0,\\
		\left|\tau-\tau_0\right|\leq 2M\tau_0^2,&\left|\partial_{\tilde t}\tau\right|\leq 2Me^{-s}.
	\end{cases}	
\end{align}

For spatial support, we define $\mathcal X(s):=\left\{e^{\frac{3s}2}(\frac{\xi_0}{2}-\xi(\tilde t))\leq y\leq e^{\frac{3s}2}(\frac{3\xi_0}{2}-\xi(\tilde t))\right\}$ and assume that 
\begin{align}\tag{BA-$\operatorname{Supp}$}\label{eq: EV_BA on spatial support}
	\operatorname{supp}\left(\partial_yW, \partial_yZ\right)\subseteq\mathcal X(s).
\end{align}

\def\BAWtildeSmall{10M^2\tau_0^{\frac12}|y|^{4-k}+\tau_0^{\frac35}\abs{y}^{3-k}\mathbbm{1}_{k\leq3}}
For $W$ and $\widetilde W$, we impose that
\begin{align}
	&\begin{cases}\tag{BA-$W$}\label{eq: EV_BA on W}
		|W|\leq (1+\tau_0^\frac1{23})\langle y\rangle^{\frac13}, &
		|\partial_y W|\leq 15\langle y\rangle^{-\frac23},\\
		|\partial_y^2 W|\leq M^\frac16\langle y\rangle^{-\frac23},&
		|\partial^3_y W|\leq M^\frac12,\\
		|\partial_y^4 W|\leq M,
	\end{cases}\\
	&\begin{cases}\tag{BA-$\widetilde W$}\label{eq: EV_BA on W_tilde}
		|\widetilde W|\mathbbm{1}_{|y|\leq L}\leq \tau_0^\frac{1}{3}\langle y\rangle^{\frac13},\\
		|\partial_y \widetilde W|\mathbbm{1}_{|y|\leq L}\leq \tau_0^\frac1{4}\langle y\rangle^{-\frac23},\\
		 |\partial_y^2\widetilde W|\mathbbm{1}_{|y|\leq L}\leq \tau_0^\frac1{5}\langle y\rangle^{-\frac23},\\
		 |\partial^3_y\widetilde{W}^0|\leq \tau_0^{\frac45},\\
		 |\partial_y^k\widetilde W|\mathbbm{1}_{|y|\leq l}\leq \BAWtildeSmall,\quad (k\leq 4),
	\end{cases}
\end{align}
where quantities with a superscript $0$ should be understood as their values at $y=0$. 

For $Z$, we require that
\begin{align}\tag{BA-$Z$}\label{eq: EV_BA on Z}
	\begin{cases}
		|Z+\sigma_\infty|\leq M\tau_0,&
		|\partial_y Z|\leq Me^{-\frac32s},\\
		|\partial_y^2 Z|\leq M^{\frac43}e^{-\frac32s},&
		|\partial_y^3 Z|\leq M^6e^{-\frac32s},\\
		|\partial_y^4 Z|\leq M^7e^{-\frac32s}.
	\end{cases}
\end{align}	

\subsubsection{Bootstrap procedure}
We now state the improved bootstrap inequalities (IB), which will be later deduced from bootstrap assumptions and initial conditions:
\begin{align}\tag{IB-M}\label{eq: EV_IB on modulation variables}
    &\begin{cases}
		 |\kappa-\kappa_0|\leq\frac12M\tau_0,&\left|\partial_{\tilde t}\kappa\right|\leq M^\frac13,\\
		 \left|\xi-\xi_0-2\beta_3\kappa_0\tilde t\right|\leq M^{\frac74}\tau_0^2, &\left|\partial_{\tilde t}\xi-2\beta_3\kappa_0\right|\leq M^\frac32\tau_0,\\
		\left|\tau-\tau_0\right|\leq \frac32M\tau_0^2,&\left|\partial_{\tilde t}\tau\right|\leq Me^{-s}.
	\end{cases}\\
	\tag{IB-$\operatorname{Supp}$}\label{eq: EV_IB on spatial support}
	&\operatorname{Supp}\left(\partial_yW, \partial_yZ\right)\subseteq\frac78\mathcal X(s),\\
	&\begin{cases}\tag{IB-$W$}\label{eq: EV_IB on W}
		|W|\leq (1+\frac12\tau_0^\frac1{23})\langle y\rangle^{\frac13}, &
		|\partial_y W|\leq 13\langle y\rangle^{-\frac23},\\
		|\partial_y^2 W|\leq M^\frac18\langle y\rangle^{-\frac23},&
		|\partial^3_y W|\leq \frac12M^\frac12,\\
		|\partial_y^4 W|\leq \frac12M,
	\end{cases}\\
	&\begin{cases}\tag{IB-$\widetilde W$}\label{eq: EV_IB on W_tilde}
		|\widetilde W|\mathbbm{1}_{|y|\leq L}\leq \frac12\tau_0^\frac{1}{3}\langle y\rangle^{\frac13},\\
		|\partial_y \widetilde W|\mathbbm{1}_{|y|\leq L}\leq \frac12\tau_0^\frac1{4}\langle y\rangle^{-\frac23},\\
		 |\partial_y^2\widetilde W|\mathbbm{1}_{|y|\leq L}\leq \frac12\tau_0^\frac1{5}\langle y\rangle^{-\frac23},\\
		 |\partial^3_y\widetilde{W}^0|\leq \frac12\tau_0^{\frac45},\\
		 |\partial_y^k\widetilde W|\mathbbm{1}_{|y|\leq l}\leq M^2\tau_0^\frac12|y|^{4-k}+\tau_0^{\frac45}\abs{y}^{3-k}\mathbbm{1}_{k\leq3},\quad (k\leq 4).
	\end{cases}\\
	\tag{IB-$Z$}\label{eq: EV_IB on Z}
	&\begin{cases}
		|Z+\sigma_\infty|\leq \frac12M\tau_0,&
		|\partial_y Z|\leq \frac12Me^{-\frac32s},\\
		|\partial_y^2 Z|\leq \frac12M^{\frac43}e^{-\frac32s},&
		|\partial_y^3 Z|\leq \frac12M^6e^{-\frac32s},\\
		|\partial_y^4 Z|\leq \frac12M^7e^{-\frac32s}.
	\end{cases}
\end{align}

\subsection{Immediate corollaries of bootstrap assumptions}\label{section: Immediate corollaries of bootstrap assumptions}

In this section, we will present some immediate results derived from the bootstrap assumptions.

\subsubsection{Relation between $\langle y\rangle$, $\xi_0$, and $\mathcal X(s)$}\

We define $\mathcal X(s):=\left\{e^{\frac{3s}2}(\frac{\xi_0}{2}-\xi(\tilde t)), e^{\frac{3s}2}(\frac{3\xi_0}{2}-\xi(\tilde t))\right\}$, on which we assume the supports of $\partial_y W$ and $\partial_y Z$ lie. When $y\in\mathcal X(s)$, by our choice of $\tau_0$ and $\xi_0$, and the bootstrap assumption \eqref{eq: EV_BA on modulation variables} on $\xi$ , we derive the following relationship between $\langle y \rangle$ and $\xi_0$:
\begin{align}\label{eq: EV_relation_<y>_xi0}
	\langle y\rangle = \sqrt{1+y^2}
	&\leq \sqrt{1+\xi_0^2e^{3s}}\leq (\tau_0^{\frac32}+\xi_0)e^{\frac{3s}2},
\end{align}
which immediately implies
\begin{align}
	\langle y\rangle\leq (\tau_0^{\frac32}+\xi_0)e^{\frac{3s}2}\lesssim\xi_0e^{\frac{3}{2}s},\quad e^{-s}\lesssim\xi_0^\frac23\langle y\rangle^{-\frac23}.\nonumber
\end{align}

\subsubsection{Absence of Vacuum}
Using \eqref{eq: EV_Riemann invariant} and bootstrap assumptions, we estimate the deviation of $\sigma$ from the background state $\sigma_\infty$:
\begin{align*}
    |\sigma-\sigma_\infty|&\le\frac{1}{2}\Big(\underbrace{e^{-\frac s2}|W|}_{\eqref{eq: EV_BA on W}\eqref{eq: EV_relation_<y>_xi0}}+\underbrace{|\kappa-\sigma_\infty|}_{\eqref{eq: EV_initial conditions for modulation variables}\eqref{eq: EV_BA on modulation variables}}+\underbrace{|Z+\sigma_\infty|}_{\eqref{eq: EV_BA on Z}}\Big)\\
    &\le\frac12(\xi_0^{\frac13}+O(\tau_0^{\frac12})).
\end{align*}
Therefore, it follows from \eqref{eq: EV_requirement of kappa_0} that
\begin{align}
    \sigma\ge\sigma_\infty-|\sigma-\sigma_\infty|\ge\frac{1}{\alpha}\xi^{\frac13}-O(\tau_0^{\frac12})\ge\frac{1}{2\alpha}\xi^{\frac13}>0.\nonumber
\end{align}

\subsubsection{Estimates of $\partial_y^k F_W$ at $y=0$}\

Appealing to bootstrap assumptions, \eqref{eq: EV_F_W and F_Z} and \eqref{eq: EV_F^k_W and F^k_Z}, we obtain
\begin{align}
	|F_W^0|&\leq M^{\frac14}e^{-\frac s2},\nonumber\\
	|\partial_y F_W^0|
	&\lesssim M^\frac14 e^{-s},\nonumber\\
	\abs{\partial_y^2 F_W^0}
	&\lesssim \tan(\xi) e^{-\frac32s}\lesssim M^\frac 54e^{-\frac32s}.\nonumber
\end{align}

\subsubsection{Evolution equations of modulation variables}\

We first derive the equations governing the modulation variables. By evaluating \eqref{eq: EV_Euler in W and Z} and \eqref{eq: EV_dkW and dkZ} at $y=0$ for $k=1,2$, and invoking the boundary conditions of $W$, we have
\begin{align}
	\eqref{eq: EV_Euler in W and Z}&\implies \partial_{\tilde t}\kappa=\frac{e^\frac s2}{\beta_\tau}\left(F_W^0+G_W^0\right),
	\label{eq: EV_dt_kappa}\\
	(\eqref{eq: EV_dkW and dkZ},k=1)
	&\implies \partial_{\tilde t}\tau=\frac1{\beta_\tau}\left(\partial_yF_W^0+\partial_y G_W^0\right),
	\label{eq: EV_dt_tau}\\
	(\eqref{eq: EV_dkW and dkZ},k=2)
	&\implies G_W^0=\left(\partial_y^3W^0\right)^{-1}\left(\partial_y^2F_W^0+\partial_y^2G_W^0\right)\label{eq: EV_GW0_alternative}.
\end{align}

The last equation, along with \eqref{eq: EV_g_W and g_Z}, gives
\begin{align}\label{eq: EV_dt_xi}
	\partial_{\tilde t}\xi=\kappa+\beta_2 Z^0-\frac{1}{\beta_\tau}e^{-\frac s2}\left(\partial_y^3W^0\right)^{-1}\left(\partial_y^2F_W^0+\partial_y^2G_W^0\right).
\end{align}

We remark that, since $\partial_y^3 \widetilde W^0 = \mathcal{O}(\tau_0^{1/4})$ and $\partial_y^3 \overline W(0)=6$, it follows from \eqref{eq: EV_GW0_alternative} that
\begin{equation}
    \abs{G_W^0} \leq M^{\frac43} e^{-s}.
\end{equation}

\subsubsection{Improvement of bootstrap assumptions on $\kappa$}\ 

Using the estimates obtained and bootstrap assumptions, for all $T\in[0,T^*]$, we deduce
\begin{align*}
	\abs{\partial_{\tilde t}\kappa}&\leq\frac{e^{\frac s2}}{\beta_\tau}\left(|F_W^0|+|G_W^0|\right)
	\leq 2M^\frac14\leq \frac12 M,\\
	|\kappa-\kappa_0|&\leq\int_0^T\abs{\partial_{\tilde t}\kappa}\;dt'\leq\frac12MT
	<\frac12M\tau_0.
\end{align*}

\subsubsection{Improvement of bootstrap assumptions on $\tau$}\

Similarly, we can show that
\begin{align*}
	|\partial_{\tilde t}\tau|&\leq \frac1{\beta_\tau}\left(\abs{\partial_yF_W^0}+\abs{\partial_y G_W^0}\right)\lesssim(1+\epsilon)\left(M^\frac14e^{-s}+e^{\frac s2}\abs{\partial_y Z^0}\right)\lesssim Me^{-s},\\
	|\tau-\tau_0|&
	\leq \int_0^{T_*}Me^{-s}\;dt'
	\leq \frac32MT_*\tau_0<2M\tau_0^2,
\end{align*} 
where $\epsilon, \tau_0\ll1$.

\subsubsection{Improvement of bootstrap assumptions on $\xi$}

By direct computation and the previous estimates, we can prove
\begin{align*}
	\abs{\partial_{\tilde t}\xi-2\beta_3\kappa_0    }&=\abs{(\kappa-\kappa_0)+\beta_2 (Z^0+\kappa_0)-\frac{1}{\beta_\tau}e^{-\frac s2}\left(\partial_y^3W^0\right)^{-1}\left(\partial_y^2F_W^0+\partial_y^2G_W^0\right)}\\
	&\lesssim M\tau_0+M\tau_0+e^{-\frac s2}\abs{\left(\partial_y^3W^0\right)^{-1}}\left(M^\frac 54e^{-\frac32s}+M^\frac43e^{-s}\right)\lesssim M^\frac32\tau_0<M^2\tau_0,\\
	|\xi-\xi_0-\beta_3\kappa_0\tilde t|&\leq\int_0^{T^*}\abs{\xi-\xi_0-\beta_3\kappa_0}\;\dif t'\leq M^{\frac32}T^*\tau_0\leq M^{\frac 74}\tau_0^2<M^2\tau_0^2.
\end{align*}

Therefore, we have improved all bootstrap assumptions on modulation variables. 

\subsection{A priori estimates}\label{section: prior estimates}

In this section, we derive the estimates of transport terms and forcing terms step by step for the equations of $W$, $Z$, and $\widetilde W$. These estimates will be used to close the bootstrap arguments on $W$, $Z$, and $\widetilde W$. 

\subsubsection{Estimates on the transport term}

We first provide estimates on transport terms and their derivatives.

\begin{lemma}\label{lemma: EV_estimate of G_W and G_Z} 
There hold the following estimates for $|G_W|$ and $|G_Z|$:
	\begin{align}
	|G_W^0|&\lesssim M^\frac43e^{-s},\label{eq: EV_estimate_GW0}\\
	|G_W|&\lesssim \begin{cases}
		\tau_0^{\frac45}, & |y|\leq L=\tau_0^{-\frac1{10}};\\
		M\xi_0e^{\frac s2}, & y\in\mathcal X(s);
	\end{cases}\label{eq: EV_estimate_GW}\\
	e^{-\frac s2}G_Z&=-4\beta_3\kappa_0+\mathcal{O}(M^2\tau_0),\quad y\in\mathcal X(s);
	\label{eq: EV_estimate_GZ}
\end{align}
\end{lemma}

\begin{proof}
	Regarding $G_W$ as a function of $y$, we obtain
\begin{align*}
	|G_W|&\leq|G_W^0|+|y||\partial_yG_W|\lesssim Me^{-s}+|y|Me^{-s}.
\end{align*}

When $|y|\leq L=\tau_0^{-\frac1{10}}$, we have the following estimates:
\begin{align*}
	|G_W|&\leq Me^{-s}+\tau_0^{-\frac1{10}}Me^{-s}\lesssim Me^{-\frac9{10}s}\lesssim\tau_0^{\frac45}.
\end{align*}

When $y\in\mathcal X(s)$, we derive 
\begin{align*}
	|G_W|\leq Me^{-s}+C\left(\frac32\xi_0-\xi\right)e^{\frac32 s}\left(Me^{-s}\right)
	&\leq Me^{-s}+C\left(\frac{\xi_0}2+|\xi_0-\xi|\right)Me^{\frac s2}\\
	&\leq Me^{-s}+C\left(\frac{51\xi_0}{100}\right)Me^{\frac s2}\lesssim M\xi_0e^{\frac s2}.
\end{align*}
The proof of \eqref{eq: EV_estimate_GZ} then proceeds in a similar manner as that for Proposition \ref{prop: est of gR,GR, zeroth order}, and we omit the details here. 
\end{proof}

\begin{lemma}\label{lemma: EV_estimate of dkG_W and dkG_Z}
	For $1\leq k\leq 4$, the following estimates on $\abs{\partial_y^kG_W}$ and $\abs{\partial_y^kG_Z}$ hold: 
	\begin{align}\label{eq: EV_estimate_dGW}
	\begin{cases}
		\abs{\partial_y G_W}\lesssim Me^{-s},&
		\abs{\partial_y^2 G_W}\lesssim M^{\frac43}e^{-s},\\
		\abs{\partial_y^3 G_W}\lesssim M^6e^{-s},&
		\abs{\partial_y^4 G_W}\lesssim M^7e^{-s},
	\end{cases}	
\end{align}
and
\begin{align}\label{eq: EV_estimate_dGZ}
	\begin{cases}
		\abs{\partial_y G_Z}\lesssim Me^{-s},&
		\abs{\partial_y^2 G_Z}\lesssim M^{\frac43}e^{-s},\\
		\abs{\partial_y^3 G_Z}\lesssim M^6e^{-s},&
		\abs{\partial_y^4 G_Z}\lesssim M^7e^{-s}.
	\end{cases}	
\end{align}
\end{lemma}

\begin{proof}
	The estimates follow from the relation $|\partial_y^kG_W|\approx|\partial_y^kG_Z|\approx e^{\frac s2}\abs{\partial_y^kZ}$ and \eqref{eq: EV_BA on Z}.
\end{proof}

\subsubsection{Estimates on the derivative of forcing terms}

In this section, we will provide estimates on derivatives of forcing terms, namely $\partial_y^kF_W$,  $\partial_y^kF_Z$, and $\partial_y^k F_{\widetilde W}$. We start by introducing the estimates of several recurring terms in these three derivatives.

\begin{lemma}\label{lemma: EV_miscellaneous terms in derivatives of forcing}
    The following estimates hold:
    \begin{align}
	&\begin{cases}\label{eq: EV_estimate_eWk}
	\abs{e^{-\frac s2}W+\kappa}\leq M^\frac14\\
		\abs{\partial_y\left(e^{-\frac s2}W+\kappa\right)}\leq e^{-\frac s2}\langle y\rangle^{-\frac23},&
	\abs{\partial_y^2\left(e^{-\frac s2}W+\kappa\right)}\leq e^{-\frac s2}\langle y\rangle^{-\frac23},\\
	\abs{\partial_y^3\left(e^{-\frac s2}W+\kappa\right)}\leq M^\frac12e^{-\frac s2},&
	\abs{\partial_y^4\left(e^{-\frac s2}W+\kappa\right)}\leq M^5e^{-\frac s2},
	\end{cases}\\
	&\begin{cases}\label{eq: EV_estimate_eWk2}
		\left(e^{-\frac s2}W+\kappa\right)^2\leq M^\frac12\\
	\abs{\partial_y\left(e^{-\frac s2}W+\kappa\right)^2}\leq 2M^\frac14 e^{-\frac s2}\langle y\rangle^{-\frac 23},&
	\abs{\partial_y^2\left(e^{-\frac s2}W+\kappa\right)^2}\leq M^\frac13 e^{-\frac s2}\langle y\rangle^{-\frac 23},\\
	\abs{\partial_y^3\left(e^{-\frac s2}W+\kappa\right)^2}\leq M^\frac{11}{16}e^{-\frac s2},&
	\abs{\partial_y^4\left(e^{-\frac s2}W+\kappa\right)^2}\leq M^\frac{19}{16}e^{-\frac s2}.
	\end{cases}\\
	&\begin{cases}\label{eq: EV_estimate_dyZ2}
		|Z|^2\lesssim M^\frac98,\\
		\abs{\partial_y Z^2}\lesssim M^\frac{13}8e^{-\frac32s}, & \abs{\partial_y^2 Z^2}\lesssim M^2e^{-\frac32s},\\
		\abs{\partial_y^3 Z^2}\lesssim M^7e^{-\frac32s}, & \abs{\partial_y^4 Z^2}\lesssim M^{7+\frac16}e^{-\frac32s}.
	\end{cases}
\end{align}
\end{lemma}

\begin{proof}
    These estimates follow directly from \eqref{eq: EV_BA on W} and \eqref{eq: EV_BA on Z}.
\end{proof}

Then, we prove the estimates for $\partial_y^kF_W$ and $\partial_y^kF_Z$:
\begin{lemma}
    The derivatives of forcing terms, namely $\partial_y^kF_W$ and $\partial_y^kF_Z$, obey the bounds below:
    \begin{align}\label{eq: EV_estimate_dkFW}
	&\begin{cases}
		|F_W|\lesssim M^\frac{3}4e^{-\frac s2},\\
	\abs{\partial_yF_W}\lesssim M^2e^{-2s}+M^\frac13 e^{-s}\langle y\rangle^{-\frac23},&
	\abs{\partial_y^2F_W}\lesssim M^{\frac{17}{8}} e^{-2s}+M^\frac23 e^{-s}\langle y\rangle^{-\frac 23},\\
	\abs{\partial_y^3F_W}\lesssim M^{3}e^{-s},&
	\abs{\partial_y^4F_W}\lesssim M^\frac{19}{16}e^{-s}.
	\end{cases}	\\
    \label{eq: EV_estimate_dkFZ}
	&\begin{cases}
		|F_Z|\lesssim M^\frac{3}4e^{-s},\\
	\abs{\partial_yF_Z}\lesssim M^2e^{-\frac{5s}2}+M^\frac13 e^{-\frac{3s}2}\langle y\rangle^{-\frac23},&
	\abs{\partial_y^2F_Z}\lesssim M^{\frac{17}{8}} e^{-\frac{5s}2}+M^\frac23 e^{-\frac{3s}2}\langle y\rangle^{-\frac 23},\\
	\abs{\partial_y^3F_Z}\lesssim M^{3}e^{-\frac{3s}{2}},&
	\abs{\partial_y^4F_Z}\lesssim M^\frac{19}{16}e^{-\frac{3s}{2}}.
	\end{cases}	
\end{align}
\end{lemma}

\begin{proof}
    Estimates in \eqref{eq: EV_estimate_dkFW} can be deduced from \eqref{eq: EV_dFkW and dFkZ}, \eqref{eq: EV_BA on W}, and lemma \ref{lemma: EV_miscellaneous terms in derivatives of forcing}, while estimates in \eqref{eq: EV_estimate_dkFZ} follows directly from \eqref{eq: EV_estimate_dkFW} and the relation $F_Z=-e^{-\frac s2}F_W$.
\end{proof}

We now proceed to estimates of $\abs{F_{\widetilde W}}$, which will be carried out in two regions $\{|y|\leq l\}$ and $\{|y|\leq L\}$ separately. We also present an estimate for $\abs{\partial_y^3F_{\widetilde W}^0}$, which is crucial for closing the bootstrap argument on $\widetilde W$ later.
\begin{lemma}
    For $1\leq k\leq 4$, we have for $\partial_y^kF_{\widetilde W}$ that
\begin{align}\label{eq: EV_estimates of dy F_tilde W within L}
	&\begin{cases}
		\abs{F_{\widetilde W}}\mathbbm{1}_{\abs{y}\leq L}\lesssim Me^{-\frac s2}+10M\tau_0^{\frac12}\langle y\rangle^{-\frac23},\\
		\abs{\partial_y F_{\widetilde W}}\mathbbm{1}_{\abs{y}\leq L}\lesssim 2M^2 e^{-s}\langle y\rangle^{-\frac23}+\tau_0^\frac45\langle y\rangle^{-\frac53},\\
		\abs{\partial_y^2 F_{\widetilde W}}\mathbbm{1}_{\abs{y}\leq L}\lesssim M^\frac94e^{-s}\langle y\rangle^{-\frac23}+\tau_0^\frac45\langle y\rangle^{-\frac83},\\
		\abs{\partial_y^3 F_{\widetilde W}}\mathbbm{1}_{\abs{y}\leq L}\lesssim M^{6+\frac18}e^{-s}+\tau_0^\frac45\langle y\rangle^{-\frac{11}3},
	\end{cases}	\\
    &\abs{\partial_y^k F_{\widetilde W}}\mathbbm{1}_{\abs{y}\leq l}\lesssim M\tau_0^{\frac12}\label{eq: EV_estimates of dy F_tilde W within l},\\
    &\abs{\partial_y^3F_{\widetilde W}^0}
	\lesssim M^{\frac{13}2}e^{-s}.\label{eq: EV_estimates of dy3FtW0}
\end{align}
\end{lemma}

\begin{proof}
 Notice that by \eqref{eq: EV_F^k tilde W}, the equation of $\partial_y^kF_{\widetilde W}$ contains $\partial_y^kF_W$, $\partial_y^k\overline W$, and $\partial_y G_W$, which are controlled in \eqref{eq: EV_W_bar decay}, \eqref{eq: EV_estimate_GW}, \eqref{eq: EV_estimate_dGW}, and \eqref{eq: EV_estimate_dkFW} respectively. With these estimates, \eqref{eq: EV_estimates of dy F_tilde W within L} can be verified. For $1\leq k\leq 4$ and $|y|<l$, we employ \eqref{eq: EV_W_bar decay}, \eqref{eq: EV_estimate_GW} and \eqref{eq: EV_estimate_dGW}, and \eqref{eq: EV_estimate_dkFW} to get
\begin{align}
	\abs{\partial_y^k F_{\widetilde W}}\mathbbm{1}_{\abs{y}\leq l}&\leq \left(\abs{\partial_y^kF_W}+\sum_{j=0}^k\abs{\binom kj\left((1-\beta_\tau)\partial_y^{k-j}\overline W-\partial_y^{k-j}G_W\right)\partial_y^{j+1}\overline W}\right)\mathbbm{1}_{\abs{y}\leq l}\lesssim 2M\tau_0^{\frac12},\nonumber
\end{align}
where the last inequality holds because all terms decay with respect to $s$ for $y$ being small. For the last inequality \eqref{eq: EV_estimates of dy3FtW0}, it follows from  \eqref{eq: EV_estimate_GW0}, \eqref{eq: EV_estimate_dGZ}, and \eqref{eq: EV_BA on modulation variables} that
\begin{align}
	\abs{\partial_y^3F_{\widetilde W}^0}&\lesssim \abs{\partial_y^3F_W^0}+\abs{\partial_y^3G_W^0}+9\abs{\partial_y^3G_W^0}+18\abs{G_W^0}+12\beta_\tau|\partial_{\tilde t}\tau|
	\lesssim M^{\frac{13}2}e^{-s}.\nonumber
\end{align}
\end{proof}

\subsubsection{Estimates on forcing terms}

In this section, we estimate the forcing terms in the evolution equations of self-similar Riemann invariants.

\begin{lemma}\label{lemma: EV_estimate of F_W^(k)}
	For $0\leq k\leq 4$, we have estimates for $\abs{F_W^{(k)}}$ in the following sense:
	\begin{align}\label{eq: EV_estimate_FWk}
	\begin{cases}
		|F_W|\leq M^\frac34e^{-\frac s2},\\
	\abs{F_W^{(1)}}\lesssim M^{2+\frac18}e^{-s}\langle y\rangle^{-\frac23},&
	\abs{F_W^{(2)}}\lesssim M^2 e^{-s}\langle y\rangle^{-\frac23},\\
	\abs{F_W^{(3)}}\lesssim M^{6+\frac14}e^{-s}+M^\frac5{12}\langle y\rangle^{-\frac{4}3},&
	\abs{F_W^{(4)}}\lesssim M^{7+\frac14}e^{-s}+M^\frac56\langle y\rangle^{-\frac23}.
	\end{cases}
	\end{align}
\end{lemma}

\begin{proof}
	From the expression of $F_W^{(k)}$ in \eqref{eq: EV_F^k_W and F^k_Z}, we observe that $F_W^{(k)}$ only depends on $\partial_y^kF_W$, $\partial_y^kG_W$, and $\partial_y^k W$, which are controlled in \eqref{eq: EV_estimate_dkFW}, \eqref{eq: EV_estimate_GW0} and \eqref{eq: EV_estimate_dGW}, and \eqref{eq: EV_BA on W}. The estimates in \eqref{eq: EV_estimate_FWk} follow then.
\end{proof}

\begin{lemma}\label{lemma; EV_estimate_FZ^k}
	For $0\leq k\leq 4$, the following estimates for $\abs{F_Z^{(4)}}$ hold:
	\begin{align}\label{eq: EV_estimate_FZk}
	\begin{cases}
		|F_Z|\leq M^\frac{3}4e^{-s},\\
	\abs{F_Z^{(1)}}\lesssim M^2e^{-\frac52s}+M^\frac13e^{-\frac32s}\langle y\rangle^{-\frac23},&
	\abs{F_Z^{(2)}}\lesssim M^\frac{11}4e^{-\frac52s}+M^\frac54e^{-\frac32s}\langle y\rangle^{-\frac23},\\
	\abs{F_Z^{(3)}}\lesssim M^{\frac{7}2}e^{-\frac32s},&
	\abs{F_Z^{(4)}}\lesssim M^{2}e^{-\frac32s}+M^{6+\frac{1}5}\langle y\rangle^{-\frac23}e^{-\frac{3s}2}.
	\end{cases}
\end{align}
\end{lemma}

\begin{proof}
	From \eqref{eq: EV_F^k_W and F^k_Z}, the expression of $F_Z^{(k)}$ only involves $\partial_y^kF_Z,\partial_y^kG_Z$, $\partial_y^kW$, and $\partial_y^kZ$, whose estimates have been obtained in \eqref{eq: EV_estimate_dkFZ}, \eqref{eq: EV_estimate_GZ} and \eqref{eq: EV_estimate_dGZ}, \eqref{eq: EV_BA on W}, and \eqref{eq: EV_BA on Z}. This then leads to the desired estimates.
\end{proof}

\begin{lemma}\label{lemma: EV_estimate_FtWk}
	There holds for $\abs{F_{\widetilde W}^{(k)}}$ that:
	\begin{align}\label{eq: EV_estimate_FtWk}
	\begin{cases}
	   \abs{F_{\widetilde W}}\mathbbm{1}_{\abs{y}\leq L}\lesssim Me^{-\frac s2}+10M\tau_0^{\frac12}\langle y\rangle^{-\frac13},\\
	   \abs{F_{\widetilde W}^{(1)}}\mathbbm{1}_{\abs{y}\leq L}\lesssim M^\frac94e^{-s}\langle y\rangle^{-\frac23}+\langle y\rangle^{-\frac43}\tau_0^\frac{7}{24}+Me^{-s}\tau_0^\frac14\langle y\rangle^{-\frac23},\\
	   \abs{F_{\widetilde W}^{(2)}}\mathbbm{1}_{\abs{y}\leq L}\lesssim \tau_0^\frac{9}{40}\langle y\rangle^{-\frac73}+M^\frac52e^{-s}\langle y\rangle^{-\frac23}+M^{\frac{1}{6}}\tau_0^{\frac{1}{4}}\langle y\rangle^{-\frac43},\\
	   {\abs{F_{\widetilde W}^{(4)}}\mathbbm{1}_{\abs{y}\leq l}\lesssim M^2l\tau_0^{\frac12}}.
	\end{cases}
\end{align}
\end{lemma}

\begin{proof}
	Invoking \eqref{eq: EV_BA on W_tilde}, \eqref{eq: EV_BA on W}, \eqref{eq: EV_W_bar decay}, \eqref{eq: EV_estimate_GW0}, \eqref{eq: EV_estimate_dGW}, and \eqref{eq: EV_estimates of dy F_tilde W within L} into the expression of $F_{\widetilde W}^{(k)}$ in \eqref{eq: EV_F^k tilde W}, we obtain the above estimates of $\abs{F_{\widetilde W}^{(k)}}$.  
\end{proof}

\subsection{Estimates on Lagrangian Trajectory }\label{sec_weighted estimate}

In this section, we define the Lagrangian trajectories $\Phi_W$ and $\Phi_Z$, and establish their lower bound estimates. We also compute the boundary values of higher-order derivatives of $\widetilde W$.

\subsubsection{Lagrangian trajectories}

We first define Lagrangian trajectories.  For a general transport equation
\begin{align}\label{eq: EV_general equation in R}
	\partial_s R+D_RR+\mathcal V_R\partial_yR=F_R,
\end{align}
with $R$ depending on $s$, $y$, and $D_R$, $V_R$ depending on $s$, $y$, and lower order terms, given an initial point $y_0$ and an initial time $s_1\geq s_0$, we can define the Lagrangian trajectory $\Phi_R(s;s_1,y_0)$ by
\begin{align}
	\begin{cases}
		\frac{\dif}{\dif s}\Phi_R(s;s_1,y_0)=\mathcal V_R(s, \Phi_R(s;s_1,y_0)),\\
		\Phi_R(s_1;s_1,y_0)=y_0.
	\end{cases}\nonumber
\end{align}
Then, the solution to \eqref{eq: EV_general equation in R} is given by
\begin{align}\label{eq: EV_solution of R}
	R\circ \Phi_R(s;s_1,y_0)=&R(s_1,y_0)\exp\left(-\int_{s_1}^s D_R(s', \Phi_R(s';s_1,y_0))\;\dif s'\right)\\
	\notag &+\int_{s_1}^s F_R (s', \Phi_R(s';s_1,y_0))\exp\left(-\int^{s}_{s'} D_R(s'', \Phi_R(s'';s_1,y_0))\;\dif s''\right)\dif s'.	
\end{align}

\subsubsection{Lower bound of $\Phi_W$}

Recall that $\mathcal V_W = \tfrac{3}{2}y + g_W$ in \eqref{eq: EV_general equation in R}. 
We now state a lower bound estimate for $\Phi_W(s;s_1,y_0)$ (Lemma \ref{lemma: EV_lower bound of Phi_W}), in the region $\abs{y_0}\geq l$, together with its corollary (Lemma \ref{lemma: estimate of y-weighted R}). These propositions will be used to improve the bootstrap assumptions on $W$ and $\widetilde W$ later. Since the arguments are analogous to those in Proposition \ref{prop: outgoing property of W trajectories} and Lemma \ref{lemma: estimate of eta to -p along PhiW}, we omit the details and refer the reader to the proofs given there.

\begin{lemma}\label{lemma: EV_lower bound of Phi_W}
	If $|y_0|\geq l$ and $s_0\geq-\ln\epsilon$, then we have for all $s\geq s_0$ that
\begin{align}\label{eq: EV_lower bound of Phi_W}
		\abs{\Phi_W(s;s_1,y_0)}\geq|y_0|e^{\frac{s-s_1}3}.
	\end{align}
\end{lemma}

\begin{lemma}\label{lemma: estimate of y-weighted R}
	For any $p\in(\frac{1}{10},10)$, $s_1>s_0$, and $y_0\in\mathbb R$, the following estimate holds:
	\begin{align}
		\int_{s_1}^s\langle y\rangle^{-p}\circ \Phi_W(s';s_1,y_0)\; \dif s^{\prime}\leq 
		\begin{cases}
			-4\ln l, & \text{$|y_0|\geq l$},\\
			\tau_0^\frac{p}{11}, & \text{$|y_0|\geq L$}.
		\end{cases}\nonumber
	\end{align} 
\end{lemma}

\subsubsection{Lower bound of $\Phi_{Z}$}

In this section, we derive a lower bound for $\Phi_Z$, where $\mathcal V_R$ is $\mathcal V_Z=\frac32y+g_Z$. The derivation is based on Lemma 8.3 in \cite{Buckmaster2023}, and we reproduce it here for completeness.

\begin{proposition}[\cite{Buckmaster2023}, Lemma 8.3]\label{lemma: EV_estimate of dPhi_Z}
	If $\Phi_Z(s;s_1,y_0)\leq \beta_3\kappa_0e^\frac s2$, then it holds that $\frac{d}{ds}\Phi_Z(s;s_1,y_0)\leq-\frac12\beta_3\kappa_0 e^\frac s2$.
\end{proposition}

\begin{proof}
	By \eqref{eq: EV_requirement of kappa_0} and \eqref{eq: EV_estimate_GZ}, we compute
	\begin{align*}
		\frac{d}{ds}\Phi_Z(s;s_1,y_0)
		&=\frac32\Phi_Z(s;s_1,y_0)+\beta_2\beta_\tau W\circ\Phi_Z(s;s_1,y_0)+G_Z\circ\Phi_Z(s;s_1,y_0)\\
		&\le \left(\frac 32\beta_3\kappa_0+\beta_2\xi_0^{\frac13}-4\beta_3\kappa_0+O(\tau_0^{\frac{1}{23}})\right)e^{\frac s2}\le-\frac12\beta_3\kappa_0 e^\frac s2.
	\end{align*}
	This leads us to the desired conclusion.
\end{proof}

\begin{proposition}[\cite{Buckmaster2023}, Lemma 8.3]
	For all $y_0\in\mathcal X_0$ and $s_1\ge s_0$, there exists $s_*(s_1,y_1)\geq s_0$ such that 
	\begin{align}\label{eq: EV_estimate of Phi_Z}
		\abs{\Phi_Z(s;s_1,y_0)}\geq\beta_3\kappa_0\min\left(\abs{e^{\frac{s}{2}}-e^{\frac{s_*}2}}, e^{\frac s2}\right).
	\end{align}
\end{proposition}

\begin{proof}
    We first notice that if $\Phi_Z(s;s_1,y_0)\leq 0$ for all $s\in[s_0,\infty)$, then Proposition \ref{lemma: EV_estimate of dPhi_Z} suggests that
	\begin{align}
		\Phi_Z(s;s_1,y_0)\leq \Phi_Z(s_1;s_1,y_0)-\frac12\beta_3\kappa_0\int_{s_1}^se^{\frac{s'}2}\;\dif s'
		\leq -\beta_3\kappa_0\left(e^{\frac s2}-e^{\frac{s_1}2}\right),
	\end{align}
	which yields the desired result by taking $s_*=s_1$. Therefore, we can assume $\Phi_Z^{y_0}(s) > 0$ for certain $s \in [s_0,\infty)$. 
    
    Consider the set
    \begin{align*}
        S_{Z}^{y_0}:=\left\{s\in[s_0,\infty):\Phi_Z(s;s_1,y_0)<\beta_3\kappa_0e^{\frac{ s}2}\right\}.
    \end{align*}
    If $S_{Z}^{y_0}$ is empty, then \eqref{eq: EV_estimate of Phi_Z} holds trivially. Thus, we can assume $S_{Z}^{y_0}$ is nonempty. 

    Since $S_{Z}^{y_0}$ is nonempty and bounded below, the quantity $\bar s := \inf S_{Z}^{y_0}$ is well defined. We then have:
    \begin{enumerate}
        \item $0<\Phi_Z(\bar s;s_1,y_0)\leq \beta_3 \kappa_0 e^{\tfrac{\bar s}{2}}$ by continuity of $\Phi_{Z}$ and the assumption that $\Phi_{Z}$ is not always non-positive.
        \item $\Phi_Z(s;s_1,y_0) > \beta_3 \kappa_0 e^{\tfrac{s}{2}}$ for all $s \in [s_0,\bar s)$. 
        \item $\Phi_Z(s;s_1,y_0) < \beta_3 \kappa_0 e^{\tfrac{s}{2}}$ and $\tfrac{d}{ds}\Phi_Z(s;s_1,y_0)\leq 0$ for all $s \in (\bar s,\infty)$ by Proposition \ref{lemma: EV_estimate of dPhi_Z}.
    \end{enumerate}
    We note that the interval $[s_0,\bar s)$ might be empty. From 1 and 3, by continuity of $\Phi_Z^{y_0}(s)$, there exists a unique $s_* > \bar s$ such that $\Phi_Z^{y_0}(s_*) = 0$. Consequently, in view of Proposition \ref{lemma: EV_estimate of dPhi_Z}, we obtain for $s > s_*$ that
    \begin{align}
		\Phi_Z(s;s_1,y_0)=\Phi_Z(s_*;s_1,y_0)+\int_{s_*}^s\Phi_Z(s';s_1,y_0)\;\dif s'\leq -\beta_3\kappa_0\left(e^\frac s2-e^{\frac{s_*}{2}}\right),\nonumber
	\end{align}
    and for $\bar s \leq s < s_*$,
    \begin{align}
		\Phi_Z(s;s_1,y_0)=\Phi_Z(s_*;s_1,y_0)-\int_s^{s^*}\frac{d}{ds}\Phi_Z(s';s_1,y_0)\;\dif s'\leq -\beta_3\kappa_0\left(e^{\frac{s_*}2}-e^\frac s2\right).\nonumber
	\end{align}
    We therefore conclude the estimate \eqref{eq: EV_estimate of Phi_Z}.
\end{proof}

The above proposition results in the following estimate for the integral of $\langle y\rangle^{-p}$ along $\Phi_Z$, which plays a crucial role in improving the bootstrap assumptions on $Z$. As the proof closely follows that of Lemma \ref{lemma: integral of eta^-p along PhiZ and A}, we omit it here for brevity.
\begin{lemma}
    For any $p>\frac{1}{10}$, $s_1\geq s_0$, and $y_0\in\mathbb R$, the following estimate holds:
	\begin{align}
		\label{eq: EV_weighted estimate of Phi_Z}
        \int_{s_0}^s\langle y\rangle^{-p}\circ \Phi_Z(s';s_1,y_0)\; \dif s'\lesssim_p1.
	\end{align}
\end{lemma}

\subsection{Improvement of bootstrap assumptions on $W$ and $\widetilde W$}\label{section: close bootstrap assumptions on W}

In this section, we provide improved bootstrap estimates for $W$, $\widetilde W$, and their derivatives. 
Since the transport terms in the equation for $W$ coincide with those in the equation for $\widetilde W$, one can verify that $\Phi_W = \Phi_{\widetilde W}$. For notational convenience, we will therefore use $\Phi_W$ only.

\subsubsection{Estimates of $\abs{\partial_y^3\widetilde W^0}$}

We first improve the bootstrap assumption on $\abs{\partial_y^3\widetilde W^0}$.

\begin{lemma}\label{lemma: EV_estimate of d_y^3W_tilde at y=0}
	There holds the following inequality:
	\begin{align}
		\abs{\partial_y^3\widetilde W^0}\leq \frac12\tau_0^{\frac35}.\nonumber
	\end{align}
\end{lemma}
\begin{proof}
	We first derive a ``good" equation satisfied by $\partial_y^3\widetilde W$. By directly computing $(\ref{eq: EV_dyW_tilde_k}, k=3)$ at $y=0$, we get
	\begin{align}\label{eq: EV_d3tW0}
	\partial_s\partial_y^3W^0=\partial_s\partial_y^3\widetilde W^0=\partial_y^3F_W^0+\partial_y^3G_W^0-G_W^0\partial_y^4W^0-3\partial_yG_W^0\partial_y^3W^0-4(1-\beta_\tau)\partial_y^3W^0,
	\end{align}
	which, by \eqref{eq: EV_BA on W}, \eqref{eq: EV_estimate_GW0}, \eqref{eq: EV_estimate_dGW}, and \eqref{eq: EV_estimate_dkFW}, implies
	\begin{align}
	\left|\partial_s \partial_y^3 W^0\right| & \leq \frac{1}{100} M^{\frac{3}{2}} e^{-s} \tau_0^{\frac{1}{5}}+M^3 e^{-s}+M^6 e^{-s}+3M^2M^5 e^{-s}+4 \beta_\tau\left|\partial_{\tilde t} \tau\right|\partial_y^3 W^0 
		\lesssim M^{7+\frac{1}{4}} e^{-s}.\nonumber
	\end{align}
	Hence, we can conclude that 
	\begin{align*}
		\left|\partial_y^3 \widetilde{W}^0(s)\right|&\lesssim\left|\partial_y^3\widetilde{W}^0(s_0)\right|+C \int_{s_0}^{+\infty}\left|\partial_s \partial_y^3 W^0\left(s'\right)\right| \;\dif s' 
		\leq \frac1{10} \tau_0^{\frac{4}{5}}+C M^{\frac{29}{4}} \tau_0
		\leq \frac{1}{2} \tau_0^{\frac{4}{5}},
	\end{align*}
	where the last inequality holds for a sufficiently small $\tau_0$.
\end{proof}

\subsubsection{Estimates of $\partial_y^k\widetilde W\mathbbm{1}_{|y|\leq l}$}
We proceed to establish the improved bootstrap bounds for $\partial_y^k\widetilde W$ in the region $\{|y|\leq l\}$.
\begin{lemma}
	For $0\leq k\leq 4$, the following estimate holds:
	\begin{align}
		\left|\partial_y^k \widetilde{W}(s, y)\right| \mathbbm1_{|y| \leq l}\leq M^2\tau_0^\frac12|y|^{4-k}+\tau_0^{\frac45}\abs{y}^{3-k}\mathbbm{1}_{k\leq3}.\nonumber
	\end{align}
\end{lemma}

\begin{proof} When $k=4$, by setting $D_{\widetilde W}^{(4)}:=\frac{11}{2}+\beta_\tau(\partial_y\overline W+4\partial_y W)$, we obtain 
\begin{align}
	D^{(4)}_{\widetilde W}\geq\frac{11}2-\beta_\tau\left(5+4\tau_0^{\frac1{23}}\right)\geq\frac13.\nonumber
\end{align}
Hence, from \eqref{eq: EV_solution of R}, by setting $R=\partial_y^4\widetilde W$, we derive
\begin{align*}
	\left|\partial_y^4 \widetilde{W}(s, y)\right| \mathbbm1_{|y| \leq l} &\leq\left|\partial_y^4 \widetilde{W}\left(y_0, s_0\right)\right| e^{-\int_{s_0}^s \frac{1}{3} \;\dif s'}+e^{-\int_{s_0}^{s} \frac{1}{3} \;\dif s'} \int_{s_0}^{s} e^{\int_{s_0}^{s'} \frac{1}{3} \;\dif s''}\abs{F_{\widetilde W}^{(4)}}\mathbbm{1}_{|y| \leq l} \circ \Phi_{W}^{y_0}(s';s_0,y_0)\;\dif s' \\
	& \lesssim\abs{\partial_y^4 \widetilde{W}(s_0, y_0)}e^{-\frac{s-s_0}{3}}+e^{-\frac{s-s_0}{3}} \int_{s_0}^s e^{\frac{s^{\prime}-s_0}{3}}\left(M^2l\tau_0^\frac12\right)\;\dif s' \\
	& \leq \frac{1}{10}M^2\tau_0^{\frac{1}{2}}+100M^2l\tau_0^{\frac12}
	\leq M^2\tau_0^\frac12.
\end{align*}

	For the case of $k=3$, it follows that 
	\begin{align}
			\left|\partial_y^3 \widetilde{W}(s, y)\right| \mathbbm{1}_{|y| \leq l} &\leq\left|\partial_y^3 \widetilde{W}^0(s)\right|+\left\|\partial_y^4 \widetilde{W}\right\|_\infty\mathbbm{1}_{|y|\leq l}|y|\leq \tau_0^{\frac45} + M^2\tau_0^{\frac12}|y|<10M^2\tau_0^\frac12|y|+\tau_0^\frac35.\nonumber
	\end{align}

	The case of $0\leq k\leq 2$ follows recursively by employing
	\begin{align}
		\left|\partial_y^k \widetilde{W}(s, y)\right|\mathbbm{1}_{\abs{y} \leq l}&\leq\left|\partial_y^k \widetilde{W}(0, s)\right|+\left\|\abs{\partial_y^{k+1}\widetilde{W}}\mathbbm{1}_{|y| \leq l}\right\|_\infty|y|, \nonumber
	\end{align}
	together with the fact that $\partial_y^kW(0,s)=\partial_y^k\overline W(0,s)$, i.e., $\left|\partial_y^k \widetilde{W}(0, s)\right|=0$, for $0\leq k\leq 2$.
\end{proof}

\subsubsection{$\langle y\rangle^\mu$-weighted estimates}
In this section, we derive weighted estimates for $W$ and $\widetilde W$. Notice that the equations of $R\in\{\partial_y^kW,\partial_y^k\widetilde W\}$ take the form of 
\begin{align}
    \partial_sR+D_RR+\left(\frac32y+g_W\right)\partial_yR=F_R.	\nonumber
\end{align}
Defining $q:=\langle y\rangle^\mu R$ with $\mu$ being a constant, we have
\begin{align}
	\partial_s q+D_R q+\left(\frac{3}{2} y+g_W\right) \partial_y q-\left(\frac{3}{2} y+g_W\right)\left(\partial_y\langle y\rangle^\mu\right) R=\langle y\rangle^\mu F_R.\nonumber
\end{align}
Namely, $q$ satisfies the system:
\begin{align}\label{eq: EV_equation of q}
	\begin{cases}
		\partial_s q+D_q q+\left(\frac{3}{2} y+g_W\right) \partial_y q=F_q, \\
		D_q = D_R-\left(\frac{3}{2} y+g_W\right)\left(\partial_y\langle y\rangle^\mu\right)\langle y\rangle^{-\mu}, \\
		F_q =\langle y\rangle^\mu F_R.
	\end{cases}
\end{align}

We start by controlling the transport term $D_R$.
\begin{lemma}\label{lemma: EV_estimate on D_R}
	The transport term $D_q$ satisfies the following estimate:
	\begin{align}
	D_q
	\geq D_R-\frac32\mu-6|\mu|\langle y\rangle^{-\frac12}.\nonumber
	\end{align}
\end{lemma}

\begin{proof}
	Via direct computation, one can rewrite $D_q$ as
	\begin{align*}
	D_q
	=D_R-\left(\frac{3}{2} y+g_W\right)(\partial_y \ln (\langle y\rangle^ \mu))
	& =D_R-\mu\left(\frac{3}{2} y+g_W\right)\left(\frac{y}{1+y^2}\right) \\
	& =D_R-\frac{3}{2} \mu-\mu g_W y\langle y\rangle^{-2}+\frac{3}{2} \mu\langle y\rangle^{-2}.
\end{align*}
It remains to estimate $D_\eta:=y\langle y\rangle^{-2}g_W$. Using \eqref{eq: EV_BA on W},\eqref{eq: EV_estimate_GW0}, and \eqref{eq: EV_estimate_dGW}, we deduce
\begin{align*}
	\left|D_\eta\right| &\leq|y|\langle y\rangle^{-2}\left(\beta_\tau|W|+\left|G_W^0\right|+\abs{\partial_y G_W}\abs{y}\right)\\
& \leq|y|\langle y\rangle^{-2}\left(\beta_\tau\left(1+\tau_0^{\frac{1}{23}}\right)\langle y\rangle^{\frac{1}{3}}+M^{\frac{3}{2}}e^{-s}+M e^{-s}|y|\right)
\leq 4\langle y\rangle^{-\frac{1}{2}},
\end{align*}
which yields the desired inequality as $\langle y\rangle^{-2}$ decays faster than $\langle y\rangle^{-\frac12}$.
\end{proof}
This immediately implies the following estimate for $q$:
\begin{corollary}\label{corollary: EV_estimate of q}
	The quantity $q$ is bounded from above as follows:
	\begin{align}\label{eq: EV_estimate of general q}
		\abs{q(s,y)}\leq&\abs{q(s_1,y_0)}\exp\left(-\int_{s_1}^s \left(D_R-\frac32\mu-6|\mu|\langle y\rangle^{-\frac12}\right) \circ \Phi_W(s';s_1,y_0)\;\dif s'\right)\\
	\notag &+\int_{s_1}^s \left[\langle y\rangle^{\mu}F_R \circ \Phi_W(s';s_1,y_0)\right]\exp\left(-\int^{s}_{s'} \left(D_R-\frac32\mu-6|\mu|\langle y\rangle^{-\frac12}\right) \circ \Phi_W(s'';s_1,y_0)\;\dif s''\right)\;\dif s'.
	\end{align}
\end{corollary}
\begin{proof}
	One can verify that $\Phi_q=\Phi_W$. The corollary then follows from Lemma \ref{lemma: EV_estimate on D_R} and \eqref{eq: EV_solution of R} for $R=q$. 
\end{proof}

In the following sections, we will make use of this corollary to derive weighted estimates for $R\in\{\partial_y^kW,\partial_y^k\widetilde W\}$ in the regions $\{l \leq |y| \leq L\}$ and $\{|y| \geq L\}$. For the region $\{l \leq |y| \leq L\}$, in view of the trajectory estimates in Lemma \ref{lemma: estimate of y-weighted R}, if we start from a point $(s,y)$ and trace the flow backwards in time, there are two possible scenarios. First, the trajectory may reach $(s_0,y_0)$ with $\abs{y_0}>l$. Second, before the time reaches $s_0$, the trajectory may intersect the boundary $\abs{y_0}=l$. Consequently, when performing estimates in the region $\{l \leq |y| \leq L\}$, we assume either $s_1=s_0$ or $\abs{y_0}=l$. Similarly, when estimating in the region $\{|y| \geq L\}$, we assume either $s_1=s_0$ or $\abs{y_0}=L$.

\subsubsection{Improvement of bootstrap assumptions on $\abs{\partial_y^k\widetilde W}\mathbbm{1}_{|y|\leq L}$}

Since we have improved the bootstrap assumptions on $\abs{\partial_y^k \widetilde W}\mathbbm{1}_{|y|\leq l}$, we only consider the region $\{l\leq|y|\leq L\}$.

\begin{lemma}
	For $|y_0|\geq l$, it holds that:
	\begin{align}
		\begin{cases}
			\langle y\rangle^{-\frac{1}{3}}|\widetilde{W}(s, y)| \mathbbm{1}_{\abs{y}\leq L}<\frac12\tau_0^\frac13,\\
			\langle y\rangle^{\frac{2}{3}}|\partial_y\widetilde{W}(s, y)| \mathbbm{1}_{\abs{y}\leq L}<\frac12\tau_0^\frac14,\\
			\langle y\rangle^{\frac{2}{3}}|\partial_y^2\widetilde{W}(s, y)| \mathbbm{1}_{\abs{y}\leq L}<\frac12\tau_0^\frac15.\\
		\end{cases}	\nonumber
	\end{align}
\end{lemma}

\begin{proof}
When $k=0$, $D_R=-\frac12+\beta_\tau\partial_y\overline W$. Let $q=\langle y\rangle^{-\frac{1}{3}}|\widetilde{W}(s, y)| \mathbbm{1}_{\abs{y}\leq L}$. Employing \eqref{eq: EV_W_bar decay}, \eqref{eq: EV_estimate_FtWk}, Lemma \ref{lemma: estimate of y-weighted R}, and Corollary \ref{corollary: EV_estimate of q} leads to
\begin{align*}
    &\phantom{\lesssim\quad}\langle y\rangle^{-\frac{1}{3}}|\widetilde{W}(s, y)| \mathbbm{1}_{\abs{y}\leq L}\\
    &\leq\langle y_0\rangle^{-\frac13}\abs{\widetilde W(s_1,y_0)}\exp\left(-\beta_\tau\int_{s_1}^s \partial_y\overline{W}\circ \Phi_W(s';s_1,y_0)\;\dif s'+2\int_{s_1}^s \langle y\rangle^{-\frac12} \circ \Phi_W(s';s_1,y_0)\;\dif s'\right)\\
	&\quad+\int_{s_1}^s \left[\langle y\rangle^{-\frac13}F_{\widetilde W} \circ \Phi_W(s';s_1,y_0)\right]\exp\left(-\beta_\tau\int_{s'}^s \partial_y\overline{W}\circ \Phi_W(s'';s_1,y_0)\;\dif s''+2\int_{s'}^s \langle y\rangle^{-\frac12} \circ \Phi_W(s'';s_1,y_0)\;\dif s''\right)\;\dif s'\\
    &\overset{\abs{y_0}\geq l}{\leq}l^{-16}M^2\tau_0^\frac12+l^{-100}M\tau_0^\frac12\leq\frac12\tau_0^{\frac13}.
\end{align*}

When $k=1$, $D_R$ is given by $D_{\widetilde W}^{(1)}=1+\beta_\tau(\partial_y\overline W+\partial_y W)$. Then, for $q=\langle y\rangle^{\frac{2}{3}}|\partial_y \widetilde{W}(s, y)|$, using \eqref{eq: EV_W_bar decay}, \eqref{eq: EV_estimate_FtWk}, \eqref{eq: EV_BA on W}, Corollary \ref{corollary: EV_estimate of q}, and Lemma \ref{lemma: estimate of y-weighted R}, we obtain
\begin{align*}
	\langle y\rangle^{\frac{2}{3}}|\partial_y \widetilde{W}(s, y)| 
    & \leq l^{-4}\left\langle y_0\right\rangle^{\frac{2}{3}}\left|\partial_y \widetilde{W}(s_1,y_0)\right| \exp \left(3 \beta \tau \int_{s_1}^s\langle y\rangle^{-\frac{2}{3}} \circ \Phi_W(s';s_1,y_0) \;\dif s'\right) \\
    &\quad +l^{-8} \int_{s_1}^s\langle y\rangle^{\frac{2}{3}}\left(M^{\frac{9}{4}} e^{-s}\langle y\rangle^{-\frac{2}{3}}+\tau_0^{\frac{7}{24}}\langle y\rangle^{-\frac{4}{3}}+M \tau_0^{\frac{1}{4}} e^{-s}\langle y\rangle^{-\frac{2}{3}}\right) \circ \Phi_W(s';s_1,y_0)) \;\dif s' \\
    & \leq 2l^{-8} \tau_0^{\frac{1}{3}}+l^{-8}\left(M^{\frac{9}{4}} \tau_0^{\frac{1}{3}} \xi_0^{\frac{2}{9}} l^{-4}+\tau_0^{\frac{7}{24}} l^{-4}+M \tau_0^{\frac{5}{4}}\right) <\frac{1}{2} \tau_0^{\frac{1}{4}}.
\end{align*}

In a similar fashion, we derive for $k=2$ and $q=\langle y\rangle^{\frac{2}{3}}\abs{\partial_y^2\widetilde{W}}$ that
\begin{align*}
	\langle y\rangle^{\frac{2}{3}}\abs{\partial_y^2\widetilde{W}}&\leq l^{-8}\tau_0^{\frac{1}{4}}+l^{-8} \int_{s_0}^s\langle y\rangle^{\frac{2}{3}}\left(M^{\frac{5}{2}} e^{-s}\langle y\rangle^{-\frac{2}{3}}+\tau^{\frac{9}{40}}\langle y\rangle^{-\frac{7}{3}}+M^{\frac{1}{6}}\tau_0^{\frac{1}{4}}\langle y\rangle^{-\frac43}\right) \circ \Phi_W(s';s_1,y_0) \;\dif s' \\
& \leq l^{-8} \tau_0^{\frac{1}{4}}+l^{-8}\left(M^{\frac{5}{2}} \tau_0^{\frac{1}{3}} \xi_0^{\frac{2}{9}} l^{-4}+\tau_0^{\frac{9}{40}} l^{-4}+M^{\frac16}\tau_0^{\frac14}l^{-4}\right) \le \frac{1}{2} \tau_0^{\frac{1}{5}}.
\end{align*}

We hence conclude the estimates stated in the lemma.
\end{proof}

\subsubsection{Improvement of bootstrap assumptions on $W$}

To close the bootstrap argument on $W$, we investigate the estimates of $\partial_y^k W$, $0\leq k\leq 4$, on three regions: $\{|y|<l\}$, $\{l\leq |y|<L\}$, and $\{|y|\geq L\}$:
\begin{enumerate}
	\item  In the region $\{|y|<l\}$, since $W=\widetilde W+\overline W$, we can directly close the bootstrap argument by combining the estimates of $\abs{\partial_y^k\widetilde W}\mathbbm{1}_{|y|<l}$ and the estimates of $\abs{\partial_y^k\overline W}$.
	\item In the region $\{l\leq |y|<L\}$, the bootstrap assumptions of $\partial_y^k W$, $k\in\{0,1,2\}$ can be directly improved by using the bootstrap assumptions on $\partial_y^k\widetilde W$, $k\in\{0,1,2\}$ and \eqref{eq: EV_W_bar decay}.
\end{enumerate}
Hence, it remains to address the bootstrap assumptions of $\partial_y^k W$, $k\in\{0,1,2\}$ on $\{|y|\geq L\}$ and the bootstrap assumptions of $\partial_y^k W$, $k\in\{3,4\}$ on $\{|y|\geq l\}$.

\begin{lemma}
	For $\partial_y^k W$, $k\in\{0,1,2\}$ on $\{|y|\geq L\}$, we have the following estimates:
	\begin{align}
		\begin{cases}
			|W|\leq (1+\frac12\tau_0^\frac1{23})\langle y\rangle^{\frac13},\\
			|\partial_y W|\leq 13\langle y\rangle^{-\frac23},\\
			|\partial_y^2 W|\leq M^\frac18\langle y\rangle^{-\frac23}.
		\end{cases}	\nonumber
	\end{align}
\end{lemma}

\begin{proof}
	For the case of $k=0$, $D_W$ is given by $D_W=-\frac12$. Hence, by setting $q=\langle y\rangle^{-\frac{1}{3}}|W(s, y)|$, and combining \eqref{eq: EV_estimate_FWk}, Lemma \ref{lemma: estimate of y-weighted R}, and Corollary \ref{corollary: EV_estimate of q}, we get
\begin{align*}
	\langle y\rangle^{-\frac{1}{3}}|W(s, y)| &\stackrel{\abs{y_0}\geq L}{\leq} e^{2 \tau_0^{\frac{1}{22}}}\left\langle y_0\right\rangle^{-\frac{1}{3}}\left|W(s_1, y_0)\right|+e^{2 \tau_0^{\frac{1}{22}}} \int_{s_1}^s\left|\langle y\rangle^{-\frac{1}{3}} F_W\right|\circ \Phi_W(s';s_1,y_0) \;\dif s' \\
	& \leq e^{2 \tau_0^{\frac{1}{22}}}\left\langle y_0\right\rangle^{-\frac{1}{3}}\left|W(s_1, y_0)\right|+e^{2 \tau_0^{\frac{1}{22}}} M^{\frac{3}{4}} \tau_0^{\frac{1}{4}} \int_{s_1}^s \xi_0^{\frac{1}{6}}\langle y\rangle^{-\frac{1}{6}} \circ \Phi_W(s';s_1,y_0) \;\dif s'\\
	& \leq e^{2 \tau_0^{\frac{1}{22}}}\left(1+\tau_0^{\frac{1}{22}}\right)+e^{2 \tau_0^{\frac{1}{22}}} M^{\frac{3}{4}} \tau_0^{\frac{1}{4}} \xi_0^{\frac{1}{6}} e^{\tau_0^\frac{1}{66}} \le 1+\frac{1}{2} \tau_0^{\frac{1}{23}}.
\end{align*}

For $k=1$, we take $q=\langle y\rangle^\frac23\abs{\partial_y W}$, and  $D_W^{(1)}=1+\beta_\tau\partial_y W$. Hence, with a similar approach, we have
\begin{align*}
	|q(s, y)| &\leq e^{4 \tau_0^{\frac{1}{22}}}\left\langle y_0\right\rangle^{\frac{2}{3}}\left|\partial_y W\left(s_1, y_0\right)\right| \exp \left(\beta_\tau \int_{s_1}^s\left|\partial_y W\right| \circ \Phi_W(s';s_1,y_0) \;\dif s'\right) \\
    & \quad +e^{4\tau_0^{\frac{1}{22}}} \int_{s_1}^s \exp \left(\beta_\tau \int_{s^{\prime}}^s|\partial_y W| \circ \Phi_W(s'';s_1,y_0)\;\dif s''\right)\abs{\langle y\rangle^{\frac{2}{3}} F_{W}^{(1)}} \circ \Phi_W(s';s_1,y_0) \;\dif s' \\
    & \leq 12e^{4 \tau_0^\frac{1}{22}} e^{\tau_0 ^\frac{2}{34}}+e^{4 \tau_0^{\frac{1}{22}}} e^{\tau_0^\frac{1}{34}}M^{\frac{17}{8}} \tau_0^{\frac{1}{4}} \int_{s_1}^s \xi_0^{\frac{1}{6}}\langle y\rangle^{-\frac{1}{6}} \circ \Phi_W(s';s_1,y_0) \;\dif s' <13.
\end{align*}

For $k=2$ and $q=\langle y\rangle^\frac23\abs{\partial_y^2W}$, direct computation yields $D_W-\frac32\mu=3\beta_\tau\partial_yW+\frac32$ and 
\begin{align}
	D_q\geq \frac32+3\beta_\tau\partial_y W-10\langle y\rangle ^{-\frac12}.\nonumber
\end{align}
Utilizing \eqref{eq: EV_BA on W}, \eqref{eq: EV_estimate_FWk}, Lemma \ref{lemma: estimate of y-weighted R}, and Corollary \ref{corollary: EV_estimate of q} then leads to
\begin{align*}
	|q(s, y)| &\leq e^{4 \tau_0^{\frac{1}{22}}}\left\langle y_0\right\rangle^{\frac{2}{3}}\left|\partial_y^2 W\left(s_1, y_0\right)\right| \exp \left(3 \beta_\tau \int_{s_1}^s\left|\partial_y W\right| \circ \Phi_W(s';s_1,y_0) \;\dif s'\right) e^{-\frac{3}{2}\left(s-s_0\right)}\\
    &\quad +e^{4 \tau_0^{\frac{1}{22}}} \int_{s_1}^s \exp \left(3 \beta_\tau \int_{s^{\prime}}^s|\partial_y W| \circ \Phi_W(s'';s_1,y_0) \;\dif s'\right)\left|\langle y\rangle^{\frac{2}{3}} F_W^{(2)}\right| \circ \Phi_W(s';s_1,y_0) \;\dif s' \\
    & \leq e^{4 \tau_0^{\frac{1}{22}}}\left\langle y_0\right\rangle^{\frac{2}{3}} M^{\frac{1}{10}}\left\langle y_0\right\rangle^{-\frac{2}{3}} e^{\tau_0^\frac{2}{34}}+e^{4\tau_0^\frac{1}{22}} e^{\tau_0^{\frac{2}{34}}} \int_{s_1}^s\left|\langle y\rangle^{\frac{2}{3}} M^2\langle y\rangle^{-\frac{2}{3}} e^{-s^{\prime}}\right| \circ\Phi_W(s';s_1,y_0) \;\dif s' \\
    & \lesssim M^{\frac{1}{9}}+M^2 e^{4 \tau_0^\frac{1}{22}} e^{\tau_0^\frac{2}{34}} \tau_0^{\frac{1}{4}} \int_{s_1}^s \xi_0^{\frac{1}{6}}\langle y\rangle^{-\frac{1}{6}} \circ \Phi_W(s';s_1,y_0) \;\dif s' \leq M^{\frac{1}{8}}.
\end{align*}
We therefore conclude our proof of the lemma.
\end{proof}

\begin{lemma}
	For $\partial_y^k W$, $k\in\{3,4\}$ on $\{|y|\geq l\}$, we have the following estimates:
	\begin{align}
		|\partial^3_y W|\leq \frac12M^\frac12,\quad 
			|\partial_y^4 W|\leq \frac12M.\nonumber
	\end{align}
\end{lemma}

\begin{proof}
For the case of $k=3$, we define $q=\left|\partial_y^3 W(s, y)\right|$. By \eqref{eq: EV_BA on W}, \eqref{eq: EV_estimate_FWk}, Lemma \ref{lemma: estimate of y-weighted R}, and Corollary \ref{corollary: EV_estimate of q}, we derive 
\begin{align*}
	\left|\partial_y^3 W(s, y)\right|&\stackrel{\left|y_0\right|\geq l}{\leq}\left|\partial_y^3 W\left(s_1, y_0\right)\right| \exp ^{-4\left(s-s_1\right)} \exp \left(4 \beta_\tau \int_{s_1}^s|\partial_y W| \circ \Phi_W(s';s_1,y_0) \;\dif s'\right)\\
    &\quad +\int_{s_1}^s e^{-4\left(s -s^{\prime}\right)}\left|F_W^{(3)}\right| \circ \Phi_W(s';s_1,y_0) \exp \left(4 \beta_\tau \int_{s^{\prime}}^s|\partial_y W| \circ \Phi_W(s'';s_1,y_0)\; \dif s''\right) \;\dif s' \\
    & \leq \frac14M^{\frac{1}{3}}\left(1+\tau_0^{\frac{1}{2}}\right)^4+2^4 \int_{s_1}^s\left|M^{\frac{25}{4}} e^{-s'}+M^{\frac{5}{12}}\langle y\rangle^{-\frac{4}{3}}\right|\circ \Phi_W(s';s_1,y_0) \;\dif s' \\
    & \leq \frac14M^{\frac{1}{3}}\left(1+\tau_0^{\frac{1}{2}}\right)^4+16 M^{\frac{25}{4}} \tau_0^{\frac{1}{4}} \int_{s_1}^s \xi_0^{\frac{1}{6}}\langle y\rangle^{-\frac{1}{6}} \circ \Phi_W(s';s_1,y_0)\;\dif s'+16 M^{\frac{5}{12}}\int_{s_1}^s\langle y\rangle^{-\frac{4}{3}} \circ \Phi_W(s';s_1,y_0)\;\dif s' \\
    & \leq \frac3{8}M^{\frac{1}{3}}+16 M^{\frac{25}{4}} \tau_0^{\frac{1}{4}} l^{-4}+16 M^{\frac{5}{12}} l^{-4} \leq \frac{1}{2} M^{\frac{1}{2}}.
\end{align*}
Similarly, for the case of $k=4$, we define $q=\left|\partial_y^4 W(s, y)\right|$. Direct computation gives
\begin{align}
	D_W^{(4)}\geq\frac{11}2-5(1+4M\tau_0)\left(1+\tau_0^\frac1{23}\right)>\frac13.\nonumber
\end{align}
Then, using \eqref{eq: EV_BA on W}, \eqref{eq: EV_estimate_FWk}, Corollary \ref{corollary: EV_estimate of q}, and Lemma \ref{lemma: estimate of y-weighted R}, we get
\begin{align*}
	\left|\partial_y^4 W(s, y)\right| &\stackrel{\left|y_0\right|\geq l}{\leq}\left|\partial_y^4 W(s_1, y_0)\right| e^{-\frac{s-s_1}{3}}+\int_{s_1}^s e^{-\frac{s-s^{\prime}}{3}}\left|F_W^{(4)}\right| \circ \Phi_W(s';s_1,y_0) \;\dif s' \\
& \leq\left|\partial_y^4 W(s_1, y_0)\right| e^{-\frac{s-s_1}{3}}+e^{-\frac{s}{3}} M^{\frac{29}{4}} \int_{s_1}^s e^{-\frac{2 s^{\prime}}{3}}\;\dif s' +e^{-\frac{s}{3}} M^{\frac{5}{6}}\int_{s_1}^s e^{\frac{s^{\prime}}{3}}\langle y\rangle^{-\frac{2}{3}} \circ \Phi_W(s';s_1,y_0) \;\dif s' \\
& \leq\left|\partial_y^4 W\left(s_1, y_0\right)\right| e^{-\frac{s-s_1}{3}}+\tau_0 M^{\frac{29}{4}} \xi_0^{\frac{4}{9}} l^{-4}+e^{-\frac{s}{3}} M^{\frac{5}{6}} 2^{-\frac{1}{3}} \int_{s_1}^s l^{-\frac{1}{3}} e^{-\frac{2}{9}\left(s'-s_1\right)+\frac{s'}{3}} \;\dif s'
\leq \frac{1}{2} M.
\end{align*}
This completes our proof of the lemma.
\end{proof}

\subsection{Improvement of bootstrap assumptions on $Z$}\label{section: closure of bootstrap assumptions on Z}

In this section, we improve the bootstrap assumptions posed on $Z$ and its higher-order derivatives.

\subsubsection{Decay estimates for $\partial_y^kZ$}

\begin{lemma}\label{lemma: EV_decay estimate}
For $1\leq k\leq 4$, the following decay estimates on $\partial_y^kZ$ hold:
	\begin{align}\label{eq: EV_decay estimates}
		\left|e^{\frac{3}{2}s} \partial_y^k Z(s, y) \right|\lesssim e^{-\frac{3}{2}(k-1)(s-s_0)} e^{\frac{3}{2} s_0}\abs{\partial_y^kZ(s_0, y_0)}+e^{-\frac{3}{2}(k-1) s} \int_{s_0}^s e^{\frac{3}{2}k s^{\prime}} \abs{F_Z^{(k)}(s')}\circ \Phi_Z(s';s_0,y_0)\;\dif s'.
	\end{align}
\end{lemma}

\begin{proof}
	We define $q=e^{\frac32s}\partial_y^kZ$ for $0\leq k\leq 4$. It can be verified that $q$ satisfies the system
\begin{align}
	\begin{cases}
		\partial_s q+D_q q+\left(\frac{3}{2}+g_Z\right) \partial_y q=F_q ,\\ 
		F_q=e^{\frac{3}{2} s} F_Z, \\ 
		D_q=D_Z-\frac{3}{2}=\frac{3}{2}(k-1)+k \beta_2 \beta_\tau\partial_y W.
	\end{cases}\nonumber
\end{align}

In terms of the equation of $D_q$, via direct computation and using \eqref{eq: EV_BA on W}, we derive
\begin{align*}
	\exp \left(-\int_{s^{\prime}}^s D_q \circ \Phi_Z(s'';s_0,y_0)\;\dif s''\right) &= e^{-\frac{3}{2}(k-1)\left(s-s^{\prime}\right)} \exp \left(-k \beta_2 \beta_\tau \int_{s^{\prime}}^s \partial_y W \circ \Phi_Z(s'';s_0,y_0)\; \dif s''\right)\\
	& \leq e^{-\frac{3}{2}(k-1)\left(s-s^{\prime}\right)}\exp \left(2(k+1) \int_{s_0}^{+\infty}\langle y\rangle^{-\frac{2}{3}}\circ\Phi_Z(s';s_0,y_0) \;\dif s'\right) \\
	& =e^{-\frac{3}{2}(k-1)\left(s-s^{\prime}\right)}\exp \left(2(k+1) \int_{s_0}^{+\infty}\left(1+\left|\Phi_Z(s';s_0,y_0)\right|^2\right)^{-\frac{1}{3}} \;\dif s'\right) \\
	& \leq C_ke^{-\frac{3}{2}(k-1)\left(s-s^{\prime}\right)},
\end{align*}
where $C_k$ is a constant only depending on $k$. The last inequality holds because of \eqref{eq: EV_weighted estimate of Phi_Z}. Combining this with the general solution \eqref{eq: EV_solution of R}, we conclude the desired inequality. 
\end{proof}

\subsubsection{Improvement of bootstrap assumptions on $Z$}
We first close the bootstrap argument on $Z$. When $k=0$, we put $\mu=0$ to get
\begin{align*}
	|(Z+\sigma_\infty)(s, y)| & \le\left|(Z+\sigma_\infty)\left(y_0, s_0\right)\right|+ C\int_{s_0}^s M^{\frac{3}{4}} e^{-s^{\prime}} \;\dif s'\\ \overset{\eqref{eq: EV_IB0 on Z}}&{\le}\frac{1}{4}M\tau_0+ M^{\frac{3}{4}} \tau_0 \leq \frac{1}{2} M\tau_0.
\end{align*}
The improved bootstrap estimates for the remaining cases are stated in the following lemma.
\begin{lemma}
	For $1\leq k\leq 4$, we have the following estimates on $\abs{e^{\frac32 s}\partial_y^kZ}$:
	\begin{align}
		\begin{array}{cc}
			|e^{\frac{3}{2}s} \partial_yZ(s, y)|\leq \frac12 M, & |e^{\frac{3}{2}s} \partial_y^2Z(s, y)|\leq \frac12 M^\frac43,\\
			|e^{\frac{3}{2}s} \partial_y^3Z(s, y)|\leq \frac12 M^4, & |e^{\frac{3}{2}s} \partial_y^4Z(s, y)|\leq \frac12 M^4.\\
		\end{array}\nonumber
	\end{align}
\end{lemma}
\begin{proof} 
	All estimates follow directly from Lemma \ref{lemma: EV_decay estimate} and the forcing estimates \eqref{eq: EV_estimate_FZk}.
	In particular, when $k=1$, we derive
\begin{align*}
	\left|e^{\frac{3}{2}s} \partial_yZ(s, y)\right|& \lesssim e^{\frac{3}{2} s_0}\left|\partial_y Z(s_0, y_0)\right|+ \int_{s_0}^{s_0} e^{\frac{3}{2} s^{\prime}}\left|F_Z^{(1)}(s^{\prime})\right| \circ \Phi_Z^{y_0}\left(s^{\prime}\right) \dif s' \\
& \lesssim M^{\frac{1}{2}}+ M^{\frac{1}{3}} \int_{s_0}^s\langle y\rangle^{-\frac{2}{3}} \circ \Phi_Z^{y_0}\left(s^{\prime}\right) \;\dif s'+ M^2 \int_{s_0}^s e^{-s} \circ \Phi_Z^{y_0}\left(s^{\prime}\right) \;\dif s' \\
& \lesssim M^{\frac{1}{2}}+M^{\frac{1}{3}}+M^2 \xi_0^{\frac{1}{3}} \tau_0^{\frac{1}{2}} \int_{s_0}^s\langle y\rangle^{-\frac{1}{3}} \circ \Phi_Z^{y_0}(s^{\prime})\;\dif s'\\
&\lesssim M^\frac12+M^\frac13+M^2\xi_0^\frac13\tau_0^\frac12 \le\frac{1}{2} M.
\end{align*}
When $k=2$, we obtain
\begin{align*}
	\left|e^{\frac{3}{2}} \partial_y^2 Z(s, y)\right|
 	&\lesssim e^{-\frac{3}{2}(s-s_0)} e^{\frac{3}{2} s_0}\left|\partial_y^2 Z(s_0, y_0)\right|+ e^{-\frac{3}{2}s} \int_{s_0}^s e^{3 s^{\prime}}\left|F_Z^{(2)}(s^{\prime})\right| \circ \Phi_Z(s';s_0,y_0) \dif s' \\
	& \lesssim M^{\frac{2}{3}}+M^{\frac{11}{4}} e^{-\frac{3}{2} s} e^{3 s-2 s} \xi_0^{\frac{1}{3}} \int_{s_0}^s \langle y\rangle^{-\frac{1}{3}} \circ \Phi_Z(s';s_0,y_0) \;\dif s' +M^{\frac{5}{4}} e^{-\frac{3}{2} s+3 s-\frac{3}{2} s} \int_{s_0}^s\langle y\rangle^{-\frac{2}{3}} \circ \Phi_Z(s';s_0,y_0)\;\dif s' \\
	& \lesssim M^{\frac{1}{3}}+M^{\frac{11}{4}} \tau_0^{\frac{1}{2}} \xi_0^{\frac{1}{3}}+M^{\frac{5}{4}} \leq \frac{1}{2} M^{\frac{4}{3}}.
\end{align*}
For the case of $k=3$, one can verify that
\begin{align*}
	\left|e^{\frac{3}{2}s} \partial_y^3 Z(s, y)\right| &\lesssim e^{-3\left(s-s_0\right)} e^{\frac{3}{2} s_0}\left|\partial_y^3 Z(s_0, y_0)\right|+e^{-3 s} \int_{s_0}^s e^{\frac{9}{2} s^{\prime}}\left|F_Z^{(3)}(s^{\prime})\right| \circ \Phi_Z(s';s_0,y_0) \;\dif s' \\
	& \lesssim M^3+e^{-3 s} \int_{s_0}^s e^{\frac{9}{2} s^{\prime}}\left|M^{\frac{7}{2}} e^{-\frac{3}{2} s^{\prime}}\right| \circ \Phi_Z(s';s_0,y_0)\;\dif s' \\
	& \lesssim M^3+ e^{-3s} M^{\frac{7}{2}} \left(e^{3s}-e^{-3s_0}\right)<M^3+M^{\frac{7}{2}} \leq \frac{1}{2} M^4.
\end{align*}
Finally, for $k=4$, the following estimate can be derived:
\begin{align*}
	\abs{e^{\frac32 s}\partial_y^4Z(s,y)}&\lesssim e^{-\frac{9}{2}(s-s_0)} e^{\frac{3}{2} s_0}\abs{\partial_y^4 Z(s_0, y_0)}+ e^{-\frac{9}{2} s} \int_{s_0}^s e^{12 s}\abs{F_Z^4(s')}\circ\Phi_Z(s';s_0,y_0)\; \dif s' \\
&\lesssim M^3+e^{-\frac{9}{2} s}\left(e^{\frac{9}{2} s}-e^{\frac{9}{2} s_0}\right) M^2\lesssim  M^3+ M^2\leq \frac{1}{2} M^4.
\end{align*}
We therefore finish the proof of the lemma.
\end{proof}
\subsection{Proof of the main theorem}\label{section: proof of the main theorem}
	
In this section, we prove Theorem \ref{theorem: EV_main theorem}. To show that the blow-up must happen in finite time, we use bootstrap assumptions on $\tau$. To demonstrate that $\widetilde v_{\theta}$ and $\widetilde\sigma$ remain bounded, it suffices to verify that $W$ stays bounded near the origin $(\tilde\theta=0)$.

\subsubsection{Blow-up speed estimate}
Building upon the bootstrap assumptions on modulation variables \eqref{eq: EV_BA on modulation variables}, we know that $|1-\pt\tau|$ is bounded. Hence, from the equation $\tau(\tilde t)-\tilde t=\int_{\tilde t}^{\widetilde T^*}(1-\pt\tau(t'))\;dt'$, we can conclude that 
\begin{align*}
	c(\widetilde T^*-\tilde t)\leq \tau(\tilde t)-\tilde t= e^{-s}\leq C(\widetilde T^*-\tilde t).
\end{align*}
\begin{align*}
	W(s,y):=e^{\frac s2}\left(\widetilde w(\tilde\theta,\tilde t)-\kappa(\tilde t)\right), \quad Z(s,y):=\widetilde z(\tilde\theta,\tilde t)
\end{align*}
Using \eqref{eq: EV_BA on W}, \eqref{eq: EV_BA on Z}, and the definition of $y$, we get
\begin{align*}
\abs{\partial_{\tilde\theta}\widetilde w}= e^{s}\abs{\partial_y W}\leq (1+\tau_0^\frac1{23})\langle y\rangle^{-\frac23}e^s\leq C_1e^s, \quad
	\abs{\partial_{\tilde\theta}\widetilde z}\leq e^{\frac{3s}2}\abs{\partial_yZ}\leq M,
\end{align*}
which then implies that 
\begin{align*}
	\frac{C_2}{\widetilde T^*-\tilde t}\leq \abs{\partial_{\tilde\theta}\widetilde w}-\abs{\partial_{\tilde\theta}\widetilde z}\leq \abs{\partial_{\tilde\theta}R}\leq \abs{\partial_{\tilde\theta}\widetilde w}+\abs{\partial_{\tilde\theta}\widetilde z}\leq\frac{C_3}{\widetilde T^*-\tilde t},
\end{align*}
where $R\in\{\widetilde v_\theta, \widetilde \sigma\}$. Here, we use the fact that $|\partial_y W|\le2$.

\subsubsection{Blow-up location estimate}

The limit $\xi_*$ exists since by \eqref{eq: EV_BA on modulation variables}, $\abs{\partial_{\tilde t}\xi}$ is bounded. To prove \eqref{eq: EV_blow-up location}, we first notice that close to $\theta=\xi(\tilde t)$, we have for $\abs{\partial_{\tilde\theta} \widetilde v_\theta}$ that
\begin{align*}
	\abs{\partial_{\tilde\theta} \widetilde v_\theta}&\lesssim\frac12\norm{e^s\partial_y W}_{L^\infty\left(B_\delta^c(\xi(t))\right)}+\frac12\norm{e^\frac32\partial_yZ}_{L^\infty\left(B_\delta^c(\xi(t))\right)}\\
	&\leq e^s\norm{\langle y\rangle^{-\frac23}}_{L^\infty\left(\left\{|y|\geq e^{\frac{3s}{2}}\delta\right\}\right)}+\frac14M\leq \frac14M+\delta^{-\frac23},
\end{align*}
provided that $\delta\in(0,1)$. The bound for $\abs{\partial_{\tilde\theta} \widetilde \sigma}$ can be derived similarly.

\subsubsection{Blow-up time estimate}

To establish the blow-up time estimate, since $\widetilde T^*=\tau(\widetilde T^*)$, we have
\begin{align*}
	|\tau_0-\widetilde T^*|=|\tau_0-\tau(\widetilde T^*)|=\int_0^{\widetilde T^*}\abs{\partial_{\tilde t}\tau}\;\dif s\leq \int_0^{\widetilde T^*}2M\tau_0\;\dif s\leq 2M\tau_0^2.
\end{align*}

\subsubsection{Proof of 1/3-H\"older continuity}

We first show that $\widetilde w\in C^{1/3}_{\tilde\theta}$. Via direct calculation, for all $x_1,x_2\in\mathbb R$ (which makes sense as we can assume $\widetilde w$ and $\widetilde z$ to be $2\pi$-periodic in $\theta$), we obtain
\begin{align*}
	\frac{\left|\widetilde w\left(\tilde{t}, x_1\right)-\widetilde w\left(\tilde{t}, x_2\right)\right|}{\left|x_1-x_2\right|^{\frac{1}{3}}} & =e^{-\frac{s}{2}} \frac{\left|W\left(s,\left(x_1-\xi\right) e^{\frac{3}{2} s}\right)-W\left(s,\left(x_2-\xi\right) e^{\frac{3}{2} s}\right)\right|}{\left|x_1-x_2\right|^{\frac{1}{3}}} \\
& \leq e^{-\frac{s}{2}}\left|x_1-x_2\right|^{-\frac{1}{3}} \int_{\left(x_1-\xi\right) e^{\frac{3}{2} s}}^{\left(x_2-\xi\right) e^{\frac{3}{2} s}}\left|\partial_y W(s, y)\right| d y \\
& \lesssim e^{-\frac{s}{2}}\left|x_1-x_2\right|^{-\frac{1}{3}}\int_0^{\left|x_1-x_2\right| e^{\frac{3}{2} s}}\abs{y}^{-\frac{2}{3}} d y \\
& \lesssim e^{-\frac{s}{2}}\left|x_1-x_2\right|^{-\frac{1}{3}}\left(\left|x_1-x_2\right|e^{\frac{3s}2}\right)^{-\frac{1}{3}} = 1 .
\end{align*}

Similarly, we can conclude that $\widetilde z\in C^{1/3}_{\tilde\theta}$ as below:
\begin{align*}
	\frac{\left|\widetilde z\left(\tilde{t}, x_1\right)-\widetilde z\left(\tilde{t}, x_2\right)\right|}{\left|x_1-x_2\right|^{\frac{1}{3}}} & =\frac{\left|Z\left(s,\left(x_1-\xi\right) e^{\frac{3}{2} s}\right)-Z\left(s,\left(x_2-\xi\right) e^{\frac{3}{2} s}\right)\right|}{\left|x_1-x_2\right|^{\frac{1}{3}}} \\
& \leq \left|x_1-x_2\right|^{-\frac{1}{3}} \int_{\left(x_1-\xi\right) e^{\frac{3}{2} s}}^{\left(x_2-\xi\right) e^{\frac{3}{2} s}}\left|\partial_y Z(s, y)\right| d y \\
& \lesssim \left|x_1-x_2\right|^{-\frac{1}{3}}\int_0^{\left|x_1-x_2\right| e^{\frac{3}{2} s}}Me^{-\frac{3s}2} d y \\
& \lesssim \left|x_1-x_2\right|^{-\frac{1}{3}}\int_0^{\left|x_1-x_2\right| e^{\frac{3}{2} s}}Me^{-\frac{s}{2}}\xi_0^\frac32\langle y\rangle^{-\frac23} d y \lesssim 1.
\end{align*}

Since $\widetilde v_\theta$ and $\widetilde \sigma$ are linear combinations of two $C^{1/3}_{\tilde\theta}$ functions, it follows that they also belong to $C^{1/3}_{\tilde\theta}$.

\bibliographystyle{plain} 

\bibliography{reference_Euler_shock_formation_on_S2}

\end{document}